\documentclass[draft,a5paper,twoside,russian,10pt]{article}

\usepackage[intlimits]{amsmath}
\usepackage{amsthm,amsfonts}
\usepackage{amssymb}
\usepackage{mathrsfs}
\usepackage[final]{graphicx,epsfig} 
\usepackage{longtable}
\usepackage{indentfirst}
\usepackage[utf8]{inputenc}
\usepackage[T2A]{fontenc}
\usepackage[english]{babel}
\usepackage[usenames]{color}
\usepackage{esint}
\usepackage{bbm}
\usepackage{styles/mystyle}

\usepackage{caption}
\makeatletter
\def\@seccntformat#1{\csname the#1\endcsname.\quad}
\makeatother

\def\maketitle{
  \begin{center}

\begin{titlepage}
\begin{center}
{\footnotesize МИНИСТЕРСТВО\,ОБРАЗОВАНИЯ\,И\,НАУКИ\,РОССИЙСКОЙ\,ФЕДЕРАЦИИ

\vspace{0.5em}

 МОСКОВСКИЙ ФИЗИКО-ТЕХНИЧЕСКИЙ ИНСТИТУТ

(ГОСУДАРСТВЕННЫЙ УНИВЕРСИТЕТ)

}

\vspace{24mm}

\vspace{7mm}

{\LARGE {\bfseries{СТОХАСТИЧЕСКИЙ АНАЛИЗ  \\[5pt] В ЗАДАЧАХ\\
[25pt]}}  Часть I 
  }

\vspace{9.5mm}

\begin{center}
{\largeПод редакцией\ \ А.\,В.~Гасникова}
\end{center}

\vspace{14mm}

{\footnotesize \itshape

Допущено  \\ Учебно-методическим объединением высших учебных заведений РФ \\
по образованию в области прикладных математики и физики в~качестве 
учебного\,пособия\,для\,студентов\,вузов,\,обучающихся\:по\,направлению\,подготовки <<Прикладные математика и физика>>,
а~также по~другим математическим и~естественнонаучным направлениям и~специальностям

}

\vspace{35mm} МОСКВА \\ МФТИ \\ 2016
\end{center}
\end{titlepage}

\thispagestyle{empty}
{\small
\begin{flushleft}
{\bf УДК 519.21(075)\\
ББК 22.317я73} \\
\hspace{27pt}{\bf С81}
\end{flushleft}}

\vspace*{-0,25cm}

\begin{center}
{\normalsize А\,в\,т\,о\,р\,ы: \\
Н.\,О.~Бузун, А.\,В.~Гасников, Ф.\,О.~Гончаров, О.\,Г.~Горбачев,\\ С.\,А.~Гуз, Е.\,А.~Крымова, А.\,А.~Натан, Е.\,О.~Черноусова} 
\end{center}

\begin{center}
Р\,е\,ц\,е\,н\,з\,е\,н\,т\,ы: \\
{\small Факультет математики НИУ ВШЭ \\
доктор физико-математических наук, доцент  \emph{А.\,В.~Колесников} }\\[3pt]
{\small Независимый Московский универститет\\
директор ИППИ РАН,
доктор физико-математических наук, \\
профессор РАН \emph{А.\,Н.~Соболевский}}
\end{center}

\noindentС81\hspace{22pt}{\bfСтохастический анализ в задачах.} Ч.\,I~: учеб. пособие\,/\\ 
\hspace*{27,6pt}Н.\,О.~Бузун, А.\,В.~Гасников, Ф.\,О.~Гончаров и др.\,; под ред.\\
\hspace*{27,6pt}А.\,В.~Гасникова. --~М.~: МФТИ, 2016. --~212~с.

\noindent\hspace{40pt}ISBN 978-5-7417-0593-3 (Ч.\,I)

\vspace{3mm}
    
{\small

Содержит программу, список литературы и задачи одноименного курса, читаемого сотрудниками кафедры математических основ управления студентам факультета управления и прикладной математики Московского физико-технического института (государственного университета). Задачи могут быть использованы в качестве упражнений на семинарских занятиях, сдачах заданий, экзаменах, а также при самостоятельном освоении курса. 

Включены задачи, отражающие  опыт преподавания вероятностных дисциплин в Независимом Московском университете и опыт наших коллег из ПреМоЛаб МФТИ, преподающих вероятностные дисциплины в ведущих научных центрах Германии и Франции. 

Учебное пособие содержит задачи, разработанные с учетом направлений подготовки студентов МФТИ, студентов НМУ, магистров НИУ~ВШЭ, СколТеха, Иннополиса, БФУ им. И.~Канта.

}

{\small
\hspace{220pt}{\bf УДК 519.21(075)}

\hspace{220pt}{\bf ББК 22.317я73}

\vspace{7mm}}

{\footnotesize \parindent=0mm  

{\bf ISBN 978-5-7417-0593-3 (Ч.\,I)}  \hspace{0.5em} \copyright\  Бузун~Н.\,О., Гасников~А.\,В.,

{\bf ISBN 978-5-7417-0595-7} \hspace{5.3em}Гончаров~Ф.\,О., Горбачев~О.\,Г.,

\hspace{17.8em}Гуз~С.\,А., Крымова~Е.\,А.,

\hspace{17.8em}Натан~А.\,А., Черноусова~Е.\,О., 2016

\hspace{16.1em} \copyright\  Федеральное государственное автономное

\hspace{17.8em}образовательное учреждение

\hspace{17.8em}высшего профессионального образования

\hspace{17.8em}<<Московский\,физико-технический\,институт

\hspace{17.8em}(государственный университет)>>, 2016

}\newpage

  \end{center}
}

\renewcommand{\refname}{Литература}

\graphicspath{{images/}}
\newcommand{\imgh}[3]{\begin{figure}[!h]\center{\includegraphics[width=#1]{#2}}\caption{#3}\label{Fig:#2}\end{figure}}

\addto{\captionsrussian}{
	
}

\usepackage{wasysym}
\usepackage{verbatim}

\makeatletter
\g@addto@macro\th@definition{\thm@headpunct{.}}
\makeatother
\theoremstyle{definition}

\newtheorem{problem}{\noindent\normalsize\bf{}}[section]

\newtheorem*{example}{\noindent\normalsize\bf{}Пример}

\newtheorem*{remark}{\noindent\normalsize\bf{}Замечание}

\newtheorem*{lemma}{\noindent\normalsize\bf{}Лемма}

\newtheorem*{ordre}{\noindent\normalsize\bf{}Указание}

\newcommand{\atoc}[1]{\addtocontents{toc}{\hspace{0.4\linewidth}#1\par}}

\newcommand{\fixme}[1]{}

\usepackage{enumitem}

\RequirePackage{enumitem}
\renewcommand{\alph}[1]{\asbuk{#1}} 

\setenumerate[1]{label=\alph*), fullwidth, itemindent=\parindent, 
  listparindent=\parindent} 
\setenumerate[2]{label=\arabic*), fullwidth, itemindent=\parindent, 
  listparindent=\parindent, leftmargin=\parindent}

\newcommand{\diag}{\ensuremath{\mathrm{diag}}}
\newcommand{\const}{\ensuremath{\mathop{\mathrm{const}}}\nolimits}
\newcommand{\Var}{\ensuremath{\mathrm{{\mathbb D}}}}
\newcommand{\Exp}{\ensuremath{\mathrm{{\mathbb E}}}}
\newcommand{\PR}{\ensuremath{\mathrm{{\mathbb P}}}}
\newcommand{\Be}{\ensuremath{\mathrm{Be}}}
\newcommand{\Po}{\ensuremath{\mathrm{Po}}}
\newcommand{\Beta}{\ensuremath{\mathrm{Beta}}}
\newcommand{\Logn}{\ensuremath{\mathrm{Log}\mathcal{N}}}
\newcommand{\Dir}{\ensuremath{\mathrm{Dir}}}

\newcommand{\mes}{\ensuremath{\mathrm{mes}}}
\newcommand{\cov}{\ensuremath{\mathrm{cov}}}
\newcommand{\I}{\ensuremath{\mathrm{I}}}
\newcommand{\N}{\ensuremath{\mathcal{N}}}
\newcommand{\KL}{\ensuremath{\mathcal{KL}}}
\newcommand{\Star}{\hspace{-6pt}*\hspace{3pt}}
\newcommand{\DStar}{\hspace{-6pt}**\hspace{3pt}}
\newcommand{\Nbb}{\mathbb{N}}

\usepackage{watermark}

\usepackage{url}
\makeatletter
\g@addto@macro{\UrlBreaks}{\UrlOrds}
\makeatother

\begin{document}

\selectlanguage{russian}
\def\contentsname{Оглавление}
\def\figurename{\small Рис.}



\setcounter{page}{3}

\tableofcontents

\section{Введение}

Данное учебное пособие, состоящее из двух частей, было задумано в 2005 г. в ходе бесед с одним из основателей кафедры МОУ ФУПМ МФТИ А.\,А.\,Натаном (06.02.1918--09.01.2009). Андрей Александрович вместе со своими учениками (С.\,А.\,Гузом и~О.\,Г.\,Горбачевым) сформировал цикл вероятностных  дисциплин, по которому сейчас обучаются студенты факультета управления и прикладной математики МФТИ. Задачи, накопленные авторским коллективом Натан--Гуз--Горбачев более чем за 30 лет преподавания на Физтехе, легли в основу настоящего сборника. Однако большая часть собранных в пособии задач была предложена (в основном отобрана из различных известных источников) за последние 10 лет в ходе методической работы, проводимой сотрудниками кафедры МОУ в рамках семинара ``Стохастический анализ в задачах''. Задачи сборника  предлагались студентам МФТИ и Независимого Московского университета в качестве домашних заданий, разбирались на семинарах. С 2015 года часть собранных в пособии задач вошла в состав методических материалов для магистратур факультета компьютерных наук НИУ ВШЭ и факультета прикладной математики и информационных технологий БФУ им.~И.~Канта.

По замыслу авторов, предлагаемый сборник задач должен обеспечить читателей необходимым ``вероятностным фундаментом'' для последующего изучения более специализированных курсов: прикладной статистики, машинного обучения, стохастической оптимизации, стохастических дифференциальных уравнений, эффективных алгоритмов, экспериментальной экономики (финансовой математики), математического моделирования и др. Несмотря на широкий спектр представленных тем, основной акцент делается на формирование у читателей геометрической интуиции, восходящей к Пуанкаре, которая позволяет с единых позиций понять многообразие асимптотических результатов стохастической теории как проявление одного общего принципа концентрации меры (см. Часть 2). 
Большое внимание в сборнике уделяется различным приемам (идеям) доказательства предельных теорем -- асимптотических результатов теории вероятностей. Для этого прежде всего используется аппарат производящих функций и теория функций комплексного переменного (см. Часть 1).

Более половины представленных в сборнике задач не являются стандартными. Для таких задач даны указания, комментарии (замечания), ссылки на публикации. 
Желающие более глубоко изучить представленные темы могут ознакомиться со списком литературы, приведенным в конце книги. 

При работе с пособием может возникнуть впечатление сильной перегруженности дополнительным материалом. Однако авторы и не ставили перед собой цели  точно подстроиться под какие-либо существующие вероятностные курсы. Напротив, намерение авторов склонялось к  составлению достаточно полного и замкнутого набора задач (в смысле агрегации основных идей/приемов, нашедших применение в различных современных приложениях). 

Одной из основных мотиваций к широкому охвату тематик  послужило следующее наблюдение. За последние десять лет практически повсеместно заметно возросла востребованность в выпускниках ведущих технических вузов с хорошим знанием вероятностных дисциплин, что в свою очередь оказало влияние на учебные планы вузов. В частности, это привело к увеличению числа часов, отводимых вероятностным дисциплинам, открытию новых (специализированных) магистерских программ, появлению новых оригинальных курсов с ``вероятностной составляющей''.

Прежде всего это вызвано широким распространением задач анализа больших массивов данных. Веяния последних нескольких лет (биоинформатика, Интернет) в основном связаны с изучением сверхбольших массивов данных (Big Data).  
Другая, не менее важная причина, -- разработка эффективных (приближенных, рандомизированных) алгоритмов решения трудных задач. 
Сложно, например, представить себе современного специалиста по моделированию, который бы не использовал методы Монте-Карло. Также сложно представить себе специалиста в области Computer Science, которому не приходилось применять рандомизированные алгоритмы и проводить вероятностный анализ алгоритмов (например, для оценки сложности в среднем). 
Еще одна причина связана с возросшей ролью вероятностных методов в анализе и разработке экономических моделей.
Наконец, можно заметить, что задачи анализа больших компьютерных, социальных, транспортных сетей в последнее время выходят на передний план во многих приложениях. Огромную роль в изучении таких сетей играют вероятностные модели, которые отчасти будут описаны в предлагаемом сборнике задач. 

Важную роль в подготовке настоящего сборника сыграли сотрудники Лаборатории структурных методов анализа данных в предсказательном моделировании (ПреМоЛаб) и особенно ее заведующий проф. В.\,Г.\,Спо\-кой\-ный. Лаборатория была открыта на базе ФУПМ МФТИ в~2011~году. Благодаря этой лаборатории у студентов, например, появилась возможность посмотреть на сайте www.mathnet.ru и канале Premolab на YouTube видеозаписи выступлений ведущих ученых, посвященные ряду нестандартных задач из этого сборника. 

Мы благодарны нашим коллегам Д.~В. Беломестному, Я.\,И.\,Бело\-поль\-ской, М.\,Л.\,Бланку, А.\,И.\,Буфетову, Н.\,Д. Введенской, В.\,В. Веденяпину, А.\,М. Вер\-ши\-ку, К.\,В. Во\-рон\-цову, В.\,В. Высоцкому, В.\,В. Вьюгину, М.\,Н. Вялому, М.\,С. Гельфанду, Э.\,Х. Гимади, Г.\,К. Голубеву, А.\,Б. Дайняку, Н.\,Х. Иб\-ра\-ги\-мо\-ву, М.\,И. Исаеву, В.\,А. Зоричу, Г.\,А. Кабатянскому, А.\,В.Ка\-лин\-ки\-ну, А.\,В. Ко\-лес\-ни\-ко\-ву, Т.\,Л. Коробковой, А.\,В. Леонидову, Г. Лугоши, В.\,А. Малышеву, В.\,Д. Мильману, В.\,В. Моттлю, Т.\,А. Нагапетяну, А.\,В. Назину, А.\,А. Наумову, Ю.\,Е.\,Нес\-те\-ро\-ву, А.\,С. Немировскому, В.\,И. Опойцеву, Ф.\,В. Пет\-ро\-ву, С.\,А. Пи\-ро\-го\-ву, Б.\,Т. Поляку, И.\,Г. Поспелову, В.\,Н. Разжевайкину, А.\,М. Райгородскому, М.\,А. Раскину, В.\,Г. Редько, A.\,Е. Ромащенко, А.\,Н. Рыб\-ко, А.\,В. Сав\-ва\-тееву, А.\,Н.\,Соболевскому, А. Содину, У. Сэнд\-холь\-му, С.\,П. Тарасову, М.\,А. Тихонов, И.\,О. Толстихину, О.\,С. Федь\-ко, Ю.\,А. Флё\-ро\-ву, М.\,Ю. Хачаю, А.\,X. Шеню, оказавшим заметное влияние на формирование различных разделов этого учебного пособия, и А.\,А.\,Шананину, во многом способствовавшему развитию на ФУПМ базового цикла вероятностных дисциплин. Также отметим большую помощь старшекурсников, аспирантов ФУПМ и~участников нашего стохастического семинара в НМУ и Физтехе в вычитке этого сборника задач (в~особенности Д.\,Бабичева, А.\,Балицкого, П.\,Двуреченского, С.\,Дов\-га\-ля, Ю.\,Дорна, Н.\,Животовского, Д.\,Камзолова, Е.\,Клочкова, А.\,Крошнина, А.\,Макарова, Ю.\,Максимова, К.\,Мищенко, П.\,Мищенко, Е.\,Молчанова, М.\,Панова, Д.\,Пет\-раш\-ко, Л.\,Прохоренковой\,(Ос\-тро\-ум\-ов\-ой), А.\,Ро\-ди\-на, А.\,Сувориковой, И.\,Усмановой, М.\,Ши\-ро\-бо\-ко\-ва, Д. Щед\-ри\-ной).

Особую благодарность хотелось бы выразить Йордану Стоянову, который скрупулезно вычитывал первыее версии этой книги, снабжая авторов ценными замечаниями. Авторы книги несут полную ответственность за любые оставшиеся неаккуратности в тексте, о которых, надеемся, читатели нам сообщат (mou@mail.mipt.ru).

В конце 2016 года было решено готовить новое (второе) издание данного пособия (Часть 1), не меняя его структуры и набора задач. Тем не менее, поскольку в 2016--2017 учебном году авторы активно использовали данное пособие (первое издание, вышедшее в издательстве МФТИ) в занятиях, проводимых со студентами ФУПМ МФТИ, Факультета Компьютерных наук НИУ ВШЭ, СколТеха, Иннополиса, то было решено учесть получнную в ходе этих занятий обратную связь. Это заставило немного скорректировать несколько десятков задач, сделав их (надеемся) более понятными. В целом, полученная обратная связь указала на необходимость более подробно во введении пояснить замысел данного пособия. 

Настоящее новое издание первой части пособия (состоящего из двух частей), предназначено для углубленного изучения цикла вероятностных дисциплин (теории вероятностей, случайных процессов, математической статистики). При составлении задач авторы исходили из того, что читатели предварительно ознакомились с лекционным материалом в объеме стандартных университетских курсов. Задачи раделов 2 и 3 представляют собой своеобразную ``солянку'', набранную в основном из стандартных задавальников, по которым на протяжении многих лет обучались студенты Физтеха. Все последующие разделы пособия (в обеих частях) являются специализированными, и большая часть задач этих разделов редко встречается в стандартных сборниках задач по вероятностным дисциплинам. Отбор тем разделов осуществлялся, прежде всего, исходя из восстребованности соответствующих направлений в современных приложениях. Отбор задач осуществлялся исходя из наличия достаточно красивой математической составляющей, полезной для решения других задач. Отличительной особенностью данного пособия (помимо наличия большего числа нестандартных задач) является междисципдинарность и требования к читателю. Нам представляется, что данное пособие стоит использовать в первую очередь ``ангажированным'' читателям, которые готовы потратить достаточно времени на работу с пособием. Написанное выше можно резюмировать следующим образом: пособие является своеобразным (задачным) путеводителем по математическим жемчужинам стохастического анализа, при отборе которых авторы старались максимально учесть возможность последующего использования полученных знаний на практике.

\begin{center}
\textbf{
Список обозначений
}
\end{center}

\noindent $( \Omega, \mathcal{F}, \PR )$ -- вероятностное пространство ($\Omega$ -- множество исходов,\linebreak $\mathcal{F}$ -- $\sigma$-алгебра, $\PR$ -- вероятностная мера). \\
$p(x),  f(x),  f_X(x)$ -- плотность распределения случайной величины $X$.\\
$\Exp X$ -- математическое ожидание случайной величины $X$. \\
$\Exp_p X$ -- математическое ожидание $X$ с плотностью распределения $p$. \\
$\Exp_\PR X$ -- математическое ожидание $X$ по мере $\PR$. \\
$\Var X$ -- дисперсия случайной величины $X$. \\
$\N(\mu,\sigma^2)$ -- нормальное распределение. \\
$\Logn (\mu, \sigma^2)$ -- логнормальное распределение.\\
$\Po(\lambda)$ -- распределение Пуассона. \\
$\Dir(\alpha_1,...,\alpha_n)$ -- распределение Дирихле.\\
$\Beta(\alpha,\beta)$ -- бета-распределение.\\
$\Be(p)$ -- распределение Бернулли. \\
$R[a,b]$ -- равномерное распределение на отрезке $[a,b]$. \\
$\Phi(x)$ -- функция стандартного нормального распределения $\N(0,1)$. \\
$\mathrm{Exp}(\lambda)$ -- показательное распределение с параметром $\lambda$. \\
$\overset{d}{\longrightarrow}$  -- сходимость по распределению, в ряде случаев $n\to\infty$ опущено во избежание громоздких обозначений.\\
$\overset{p}{\longrightarrow}$ -- сходимость по вероятности. \\
$\overset{\text{п.н.}}{\longrightarrow}$ -- сходимость с вероятностью 1.\\
$[x^{n}]\varphi(x)$ -- коэффициент при $x^n$ в разложении в степенной ряд функции $\varphi(x)$.\\
с.в. -- случайная величина. \\
ЗБЧ -- закон больших чисел. \\
х.ф. -- характеристическая функция. \\
ЦПТ -- центральная предельная теорема. \\
п.ф. -- производящая функция. \\
э.п.ф -- экспоненциальная производящая функция.\\
$\langle \cdot, \cdot\rangle$ -- скалярное произведение. \\
Индикаторная функция:
\[
\I(\text{true}) = [\text{true}] = 1, \quad \I(\text{false}) = [\text{false}] = 0.
\]
Простая выборка с плотностью распределения $p_X$:
\[
X_1,\ldots,X_n \sim p_X.
\]
В схожих ситуациях используется символ $\in$ вместо $\sim$ (например, $X \in \N(m, \sigma^2)$), если подразумевается принадлежность случайной величины  семейству распределений. Знак $\sim$ также обозначает пропорциональность. \\ 
Ненормированная плотность распределения:
\[
p(x) \propto g(x), \quad p(x) = \frac{g(x)}{\int g(x) dx}.
\]
Математическое ожидание по заданной переменной:
\[
\Exp_X \left[h(X,Y)\right] = \int_{-\infty}^{+\infty} h(x,Y) dF(x). 
\]
Дивергенция Кульбака--Лейблера для распределений $\PR_1$ и $\PR_2$ с общим носителем~$\Omega$ (соответствующие плотности распределений обозначены как $p_1$ и $p_2$):
\[
KL(\PR_1 \Vert \PR_2) = \KL(\PR_1 \Vert \PR_2) = \int_{\Omega} \log \left( \frac{p_1(x)}{p_2(x)} \right) p_1(x) dx
\]
(в тексте могут встречаться два различных обозначения дивергенции $KL(\PR_1 \Vert \PR_2)$ и $\mathcal{KL}(\PR_1 \Vert \PR_2)$).  
Дивергенция Кульбака--Лейблера для распределения $\PR_\theta$, зависящего от параметра (соответствующая плотность распределения обозначена как $p(x| \theta)$):
\[
\KL(\theta_1 \Vert \theta_2) = \KL(\PR_{\theta_1} \Vert \PR_{\theta_2}).
\]
$L_p(X,\mu)$ -- пространство {\it p}-интегрируемых функций по мере $\mu$ (если $\mu$  не указана, то рассматривается мера Лебега на множестве $X$). Норма в $L_p(X,\mu)$ определяется следующим образом:
\[
    \|f\|_{L_p} = \left(\int_X |f(x)|^pd\mu(x) \right)^{1/p}.
\]
$A \triangle B = A \setminus B + B \setminus A$ -- симметрическая разность. \\
$\log = \ln$. \\
$\mathcal{B}$ -- борелевская сигма-алгебра.

\newpage

\section{Стандартные задачи}
\label{standart}


\begin{problem}
Некто имеет $N$ ключей, из которых только один от его двери. Какова вероятность того, что, используя ключи в случайном порядке, 
он откроет дверь 
\begin{enumerate}
\item[а)] первым ключом; 
\item[б)] последним ключом? 
\end{enumerate}
Требуется найти вероятность того, что потребуется не менее $k$ попыток, чтобы открыть дверь, если ключи, которые не подошли 
\begin{enumerate}
\item[в)] откладываются; 
\item[г)] не откладываются. 
\end{enumerate}

\begin{ordre}
В пункте в) воспользоваться заменой вероятности бинарной величины на математическое ожидание.
\end{ordre}

\end{problem}

\begin{problem}
Ребенок играет с десятью буквами разрезной азбуки: А, А, А, Е, И, К, М, М, Т, Т. 
Какова вероятность того, что при случайном расположении букв в ряд он получит слово << МАТЕМАТИКА >>? 
\end{problem}

\begin{problem}
Казалось бы, при бросании двух игральных костей как 9, так и 10 можно получить двумя разными способами: $9 = 3+6 = 4+5$, $10 = 4+6 = 5+5$. В случае бросания трех костей 9 и 10 получаются 6 различными способами. Почему тогда при бросании двух костей 9 появляется чаще 10, а при бросании трех костей -- 10 чаще 9? 
\end{problem}

\begin{problem}(Парадокс раздела ставки \cite{2013}, \cite{book12}.)
Два игрока играют в честную игру (с вероятностью победы в каждой из партий 1/2). Игроки договорились, что тот, кто первым выиграет 6 партий, получит  весь приз. Предположим, что игроки не успели доиграть: первый выиграл 5 партий, а второй -- 3. Как следует разделить приз?   
\end{problem}

\begin{problem}(Вероятностный парадокс лжеца \cite{book2012}.)
Если Вы выберете ответ на этот вопрос случайно, какова вероятность правильного ответа:

а) 25{\%}; б) 50{\%}; в) 60{\%}; г) 25{\%}?
\end{problem}

\begin{problem}
Из урны, содержащей $a$ белых, $b$ черных и $c$ красных шаров, последовательно извлекаются три шара. Найдите 
вероятность следующих событий: 
\begin{enumerate}
\item[а)] все три шара разного цвета; 
\item[б)] шары извлечены в последовательности ``белый, черный, красный''; 
\item[в)] шары извлечены в обратной последовательности. 
\end{enumerate}
\end{problem}

\begin{problem}
Из урны, содержащей $a$ белых и $b$ черных шаров, извлекается наугад один шар и откладывается в сторону. Какова вероятность 
того, что извлеченный наугад второй шар окажется белым, если: 
\begin{enumerate}
\item[а)] первый извлеченный шар белый; 
\item[б)] цвет  первого извлеченного шара остается неизвестным? 
\end{enumerate}
\end{problem}

\begin{problem}(Задача Стефана Банаха \cite{2013}, \cite{29}.)
В двух спичечных коробках имеется по $n$ спичек. На каждом шаге наугад выбирается коробок, и из него удаляется (используется) 
одна спичка. Найдите вероятность того, что в момент, когда один из коробков опустеет, в другом останется $k$ спичек. 
\end{problem}

\begin{ordre}
Событие, удовлетворяющее условию задачи: из выбранного коробка взяли последнюю спичку, а в другом коробке имеется $k$ спичек. 
\end{ordre}

\begin{problem}
Партия продукции состоит из десяти изделий, среди которых два изделия дефектные. Какова вероятность того, что из пяти отобранных 
наугад и проверенных изделий: 
\begin{enumerate}
\item[а)] ровно одно изделие дефектное; 
\item[б)] ровно два изделия дефектные; 
\item[в)] хотя бы одно изделие дефектное? 
\end{enumerate} 

\begin{ordre}
В пункте в) удобнее искать вероятность противоположного события.
\end{ordre}

\end{problem}

\begin{problem}
Известно, что в результате бросания десяти игральных костей выпала  по крайней мере одна «шестерка». Какова вероятность того, что число выпавших «шестерок» больше единицы?
\end{problem}

\begin{problem}
Найдите вероятность того, что из $50$ студентов, присутствующих на лекции по теории вероятностей, хотя бы двое имеют одну и ту же дату рождения. 

\begin{remark}
См. Мостеллер~Ф. Пятьдесят занимательных вероятностных задач с решениями. М.: Едиториал УРСС, 2009.
\end{remark}

\end{problem}

\begin{problem}
В урне находится $m$ шаров, из которых $m_1$ белых и $m_2$ черных $(m_1 + m_2 = m)$. 
Производится $n$ извлечений одного шара с возвращением его (после определения его цвета) обратно в урну. Найдите вероятность того, 
что ровно $r$ раз из $n$ будет извлечен белый шар. 
\end{problem}

\begin{problem}
Найдите вероятность того, что при размещении $n$ различных шаров по $N$ различимым ящикам заданный ящик будет содержать ровно 
$k$: $0\leqslant k\leqslant n$ шаров (все различимые размещения равновероятны). 
\end{problem}

\begin{problem}
В урне находится $m$ шаров, из которых $m_1$ -- первого цвета, $m_2$ -- второго цвета, $\ldots$, $m_s$ -- $s$-го цвета 
$(m_1+m_2+\ldots +m_s=m)$. 
Производится $n$ извлечений одного шара с возвращением его (после определения его цвета) обратно в урну. Найдите вероятность того, 
что $r_1$ раз будет извлечен шар первого цвета, $r_2$ раз -- шар второго цвета, $\ldots$, $r_s$ раз -- шар $s$-го цвета 
$(r_1+r_2+\ldots +r_s=n)$. 
\end{problem}

\begin{problem}
\label{sec:clubok}
Из клубка с $n$ разными веревочками выходит $2n$ концов этих веревочек. Свяжем их попарно в случайном порядке (все возможные варианты связок равновероятны); получится несколько (зацепленных друг за друга) веревочных колец разной длины. Найдите математическое ожидание числа колец.
\end{problem}

\begin{ordre}
Эта задача  легко сводится (\cite{book2012}) к подсчету математического ожидания числа циклов в случайной перестановке \linebreak(см.~задачу~\ref{cycle} из раздела~\ref{genF}). 
\end{ordre}

\begin{problem}
\label{vkl_iskl}
Пусть $A_1,\ldots,A_n$ -- последовательность событий и $J\subset \{1,\ldots,n\}$. Обозначим
\begin{align*}
&\PR_1(J) = \sum\limits_{i\in J} \PR(A_i),\\
&\PR_2(J) = \underset{i < j, (i,j)\in J}{\sum} \PR(A_i \cap A_j),\\
&\PR_k(J) = \underset{i_1 < \ldots < i_k, i_{s}\in J}{\sum} \PR(A_{i_1} \cap \ldots \cap A_{i_k}).
\end{align*}

 Покажите, что для всех нечетных $k \le n$ справедливо неравенство
 \[
 \PR(A_{i_1} \cup \ldots \cup A_{i_k}) \leq \underset{j = 1}{\overset{k}{\sum}} (-1)^{j+1}\PR_j\left(\left\{i_1,i_2,..,i_k\right\}\right),    
 \]
  а для всех четных $k \le n$
 
 \[
 \PR(A_{i_1} \cup \ldots \cup A_{i_k}) \geq \underset{j = 1}{\overset{k}{\sum}} (-1)^{j+1}\PR_j\left(\left\{i_1,i_2,..,i_k\right\}\right).	
 \]

\end{problem}

\begin{remark}
При $k=n$ получим равенство, известное как \textit{формула включений-исключений}.
\end{remark}

\begin{problem}
В гардеробе случайным образом перепутались $N$ одинаковых шляп посетителей. Какова вероятность того, что хотя бы один посетитель получит свою шляпу? Рассмотреть случаи $N=4, \, N=10\,000$.  
\end{problem}

\begin{ordre}
Воспользоваться формулой включений-исключений \linebreak(см.~задачу~\ref{vkl_iskl}). См.~также задачу~\ref{permloop} из раздела~\ref{genF}.
\end{ordre}

\begin{problem}
Несколько раз бросается игральная кость. Какое событие более вероятно: 
\begin{enumerate}
\item[а)] сумма выпавших очков четна; 
\item[б)] сумма выпавших очков нечетна? 
\end{enumerate}
\end{problem}

\begin{problem}
Для уменьшения общего количества игр $2n$ команд спортсменов разбиваются на две подгруппы. Определите вероятности того, что 
две наиболее сильные команды окажутся: 
\begin{enumerate}
\item[а)] в одной подгруппе; 
\item[б)] в разных подгруппах. 
\end{enumerate}
\end{problem}

\begin{problem}
В урне находятся $a$ белых и $b$ черных шаров, которые наугад по одному без возвращения извлекаются из урны до тех пор, пока урна не опустеет. 
Какое событие более вероятно: 
\begin{enumerate}
\item[а)] первый извлеченный шар белый; 
\item[б)] последний извлеченный шар белый?
\end{enumerate}
\end{problem}

\begin{problem}
В урне находятся $a$ белых и $b$ черных шаров. Шары наугад по одному извлекаются из урны без возвращения. Найдите вероятность того, 
что $k$-й вынутый шар оказался белым. 
\end{problem}

\begin{problem}
$30$ шаров размещаются по $8$ ящикам так, что для каждого шара вероятности попадания в каждый ящик одинаковы. Найдите вероятность 
размещения, при котором будет $3$ пустых ящика, $2$ ящика -- с~тремя, $2$ ящика -- с~шестью и $1$ ящик -- с~двенадцатью шарами. 
\end{problem}

\begin{problem}(Вероятностная формулировка Большой теоремы Ферма.)
Две урны содержат одинаковое количество шаров. Шары окрашены в белый и черный цвета. Из каждой урны вынимают по $n$ шаров с возвращением, $n\geq 3$. Найдите $n$ и ``состав'' каждой урны, если вероятность того, что все шары, взятые из 1-й урны, только белые, равна вероятности того, что все шары, взятые из 2-й урны, либо только белые, либо только черные.
\end{problem}

\begin{problem} (Мультиномиальное распределение.)
В $k$ ячейках случайно и независимо друг от друга размещаются $N$ частиц так, что каждая из них попадает 
в $i$-ю ячейку с вероятностью $p_i$ $(i=1,\ldots,k, \sum\limits_{i=1}^{k} p_i=1)$. Найдите вероятность того, что число частиц в ячейках 
примет заданные значения $n_1$, $\ldots$, $n_i$, $\ldots$, $n_k$. 
\end{problem}

\begin{problem}
Из $n$ лотерейных билетов $k$ -- выигрышные $(n\geqslant 2k)$. Какова вероятность, что среди $k$ купленных билетов по крайней мере 
один будет выигрышным? 
\end{problem}

\begin{problem}(Какие автобусы переполнены? \cite{book2012}.) Двое заканчивают рабочий день 
одновременно и идут к автобусной остановке вместе. Каждый уезжает на первом 
подъехавшем автобусе: одному нужен автобус №~24, другому -- №~25. При этом 
первый считает, что в среднем автоубы №~24 более полные, а второй -- что автобусы №~25. Почему так может быть?
\end{problem}

\begin{problem} 
Что более вероятно: при одновременном бросании четырех игральных костей получить хотя бы одну единицу или при $24$ бросаниях 
 по две игральные кости одновременно получить хотя бы один раз две единицы? Найдите вероятности указанных событий. 
\end{problem}

\begin{problem}
Из $2^N$ множеств -- совокупности всех подмножеств множества $\{1,2,\ldots,N\}$ -- случайно и независимо выбираются два множества $A$ и $B$. Найти вероятность, что $A$ и $B$ не пересекаются.
\end{problem}

\begin{problem}(Н.Н. Константинов \cite{book2012}.)
В самолете $n$ мест. Есть $n$ пассажиров, выстроившихся друг за другом в очередь. Во главе очереди -- ``заяц'' (пассажир без билета). У всех, 
кроме ``зайца'', есть билет, на котором указан номер посадочного места. Так как ``заяц'' входит первым, он случайным образом занимает 
некоторое место. Каждый следующий пассажир, входящий в салон самолета, действует по такому принципу: если его место свободно, то 
садится на него, если занято, то занимает с равной вероятностью любое свободное. Найдите вероятность того, что последний пассажир 
сядет на свое место. 
\end{problem}


\begin{problem}
Опыт состоит в подбрасывании симметричной монеты до тех пор, пока два раза подряд она не выпадет одной и той же стороной. Подбрасывания независимы в совокупности. 
Постройте пространство элементарных событий и найдите вероятности следующих событий: 
\begin{enumerate}
\item опыт окончится до шестого бросания; 
\item для завершения опыта потребуется четное число бросаний. 
\end{enumerate}
\end{problem}

\begin{problem}(Pаспределение Паскаля или отрицательно биномиальное распределение.)
Найдите распределение дискретной случайной величины, равной количеству произошедших неудач в последовательности испытаний Бернулли с вероятностью успеха $p$, проводимой до $r$-го успеха.

\end{problem}

\begin{problem}
При каждом подбрасывании монеты орел выпадает с вероятностью $p>0$. Пусть $\pi _{n} $ -- вероятность того, что за $n$ независимых подбрасываний орел выпадет четное число раз. Покажите, что $\pi _{n+1} =\left(1-p\right)\cdot \pi _{n} +p\cdot \left(1-\pi _{n} \right)$, $n\in {\mathbb N}$. Найдите~$\pi _{n}$. 
\end{problem}

\begin{problem}

В орлянку играют два игрока. Начальный капитал первого игрока равен $z$ руб., второго равен $\infty$. В каждой партии первый игрок выигрывает или проигрывает 1 руб. с вероятностью $0,5$, независимо от предшествующих партий. Обозначим за $z + s_k$ капитал первого игрока после $k$-й итерации, $\eta(z)$ -- число шагов до разорения:
\[
\eta(z) = \min \{k: z + s_k = 0 \}.
\]     
Покажите, что с вероятностью 1 первый игрок разорится, т.е. \linebreak что $\PR(\eta(z)<\infty) = 1,$ но при этом $\Exp \eta(z) = \infty$.  

\end{problem}

\begin{ordre}
Получите рекуррентные соотношения для вероятности разорения и  $\Exp \eta(z)$. Допустив  $\Exp \eta(z) < \infty$, придите к противоречию, доказав $\Exp \eta(z) < 0$ при больших $z$.
\end{ordre}

\begin{problem}
Сколько раз в среднем нужно подбросить монету, чтобы решка выпала два раза подряд?
\end{problem}

\begin{problem} (Честная игра \cite{book2012}.)
Два человека играют в {\it орлянку} (кидают симметричную монету и ставят по рублю: один на орла, другой на решку). У каждого есть своя монетка, и каждый подозревает партнера в несимметричности его монетки. Предложите правила (как, имея эти две монеты, ``сгенерировать'' симметричную монету), по которым игра будет гарантированно честной?
\end{problem}

\begin{problem} (Как играть в проигрышную игру? \cite{book12})
Пусть игрок А и Б играют в {\it орлянку} с четным числом розыгрышей $n$. Игрок А выигрывает в каждом розыгрыше с вероятностью $p = 0,45$ (игрок Б -- с вероятностью $0,55$). Чтобы в итоге выиграть в этой игре, надо выиграть больше чем в половине розыгрышей. Кажется, что если у игрока А есть возможность выбирать число розыгрышей, то он выберет $n = 2$, желая максимизировать вероятность своего выигрыша. Однако это не лучший выбор. Покажите, что оптимальный выбор для игрока А -- $n = 10$ розыгрышей. 
\end{problem}

\begin{problem}
Приведите пример вероятностного пространства и трёх событий на этом пространстве, которые попарно независимы, но зависимы в совокупности. Предложите обобщение полученной конструкции, в котором любые $n$ из $n+1$ событий независимы в совокупности, а все $n+1$ -- зависимы в совокупности.
\end{problem}

\begin{ordre}
См. задачу \ref{bernshtein}, а также источники \cite{book12}, \cite{stoianov}.
\end{ordre}

\begin{problem}
Покажите, что из независимости событий $A$ и $B$ следует независимость событий $A$ и $\overline B$, $\overline A$ и $B$, 
$\overline A$ и $\overline B$. 
\end{problem}

\begin{problem}
Покажите, что из равенства ${\mathbb P}(A\, |\, B)={\mathbb P}(A\, |\, \overline B)$ для ненулевых событий $A$ и $B$ следует 
равенство ${\mathbb P}(AB)={\mathbb P}(A){\mathbb P}(B)$, т.е. их независимость. 
\end{problem}

\begin{problem}(Пример Бернштейна \cite{stoianov}.)
\label{bernshtein}
Подбрасываются три игральные кости. События $A$, $B$ и $C$ означают выпадение одинакового числа очков (соответственно) на первой и 
второй, на второй и третьей, на первой и третьей костях. Являются ли эти события независимыми 
\begin{enumerate}
\item[а)] попарно, 
\item[б)] в совокупности? 
\end{enumerate}
\end{problem}

\begin{problem} (См. \cite{book2012}, \cite{2013}.)
Имеется три картонки. На одной с обеих сторон нарисована буква $A$, на другой -- $B$, а на третьей с одной стороны -- $A$,\linebreak с другой -- $B$. Одна из картонок выбирается наугад и кладется на стол. Предположим, что на видимой стороне картонки оказывается буква $A$. Какова вероятность, что на второй стороне тоже $A$?
\end{problem}

\begin{problem}
Пусть $A$, $B$, $C$ -- заданные события. Докажите справедливость неравенств: 
\begin{enumerate}
\item ${\mathbb P}(AB)+{\mathbb P}(AC)+{\mathbb P}(BC)\geqslant {\mathbb P}(A)+{\mathbb P}(B)+{\mathbb P}(C)-1$; 
\item ${\mathbb P}(AB)+{\mathbb P}(AC)-{\mathbb P}(BC)\leqslant {\mathbb P}(A)$; 
\item ${\mathbb P}(A\bigtriangleup B)\leqslant {\mathbb P}(A\bigtriangleup C)+{\mathbb P}(C\bigtriangleup B).$ 
\end{enumerate}
\end{problem}

\begin{problem}(Задача A.A. Натана \cite{5}.)
Вероятности конъюнкции событий $A$, $B$ и $C$ приведены в таблице истинности: 
\vspace{0.3cm}

\begin{tabular}{|c|c|c|c|}
\hline
$\quad A\quad$ & $\quad B\quad$ & $\quad C\quad$ & $\quad {\mathbb P}(ABC)\quad$ \\
\hline
$0$ & $0$ & $0$ & $20/36$ \\
\hline
$0$ & $0$ & $1$ & $5/36$ \\
\hline
$0$ & $1$ & $0$ & $5/36$ \\
\hline
$0$ & $1$ & $1$ & $0$ \\
\hline
$1$ & $0$ & $0$ & $5/36$ \\
\hline
$1$ & $0$ & $1$ & $0$ \\
\hline
$1$ & $1$ & $0$ & $0$ \\
\hline
$1$ & $1$ & $1$ & $1/36$ \\
\hline
\end{tabular}

\vspace{0.3cm}

Можно ли наблюдаемые события $A$ и $B$ использовать как признаки, используемые для обнаружения события $C$? 
\end{problem}

\begin{problem}
Юноша собирается сыграть три теннисных матча со своими родителями, и для выигрыша он должен победить два раза подряд. 
Порядок матчей может быть следующим: отец--мать--отец, мать--отец--мать. Известно,  что отец играет лучше матери. Юноше нужно решить, какой порядок для него предпочтительней? 

\end{problem}

\begin{problem}
Имеются две урны. В одной из них находится один белый шар, в другой -- один черный шар (других шаров урны не содержат). Выбирается 
наугад одна урна. В нее добавляется один белый шар и после перемешивания один из шаров извлекается. Извлеченный шар оказался белым. 
Определите апостериорную вероятность того, что выбранной оказалась урна, которая первоначально содержала белый шар. 
\end{problem}

\begin{problem}
В первой урне содержится $a$ белых и $b$ черных шаров (и только они), во второй -- $c$ белых и $d$ черных шаров 
(и только они). Из выбранной наугад урны извлекается один шар, который обратно не возвращается. Извлеченный шар оказался белым. 
Найдите вероятность того, что и второй шар, извлеченный из той же урны, окажется белым. 
\end{problem}

\begin{problem}(Парадокс Монти Холла \cite{book12}.)
Представьте, что вы стали участником шоу. Перед Вами три закрытых двери. Ведущий поместил за одной из трех 
пронумерованных дверей автомобиль, а за двумя другими дверями -- по козе (козы тоже пронумерованы) случайным образом -- это значит, 
что все $3! = 6$ вариантов расположения автомобиля и коз за пронумерованными дверями равновероятны. У~вас нет никакой информации 
о том, что за какой дверью находится. Ведущий говорит: <<Сначала вы должны выбрать одну из дверей. После этого я открою одну из 
оставшихся дверей (при этом если вы выберете дверь, за которой находится автомобиль, то я с вероятностью $1/2$ выберу дверь, 
за которой находится коза номер $1$, и с вероятностью $1-1/2=1/2$ дверь, за которой находится коза номер $2$). Затем я предложу 
вам изменить свой первоначальный выбор и выбрать оставшуюся закрытую дверь вместо той, которую вы выбрали сначала. Вы можете 
последовать моему совету и выбрать другую дверь либо подтвердить свой первоначальный выбор. После этого я открою дверь, 
которую вы выбрали, и вы выиграете то, что находится за этой дверью.>> Вы выбираете дверь номер $3$. Ведущий открывает дверь номер $1$ 
и показывает, что за ней находится коза. Затем ведущий предлагает вам выбрать дверь номер $2$. Увеличатся ли ваши шансы 
выиграть автомобиль, если вы последуете его совету? 
\end{problem}

\begin{problem}
Известно, что $96\%$ выпускаемой продукции соответствует стандарту. Упрощенная схема контроля признает годным с вероятностью 
$0,98$ каждый стандартный экземпляр аппаратуры и с вероятностью $0,05$ -- каждый нестандартный экземпляр аппаратуры. Найдите вероятность того, 
что изделие, прошедшее контроль, соответствует стандарту. 
\end{problem}

\begin{problem}
Пусть отличник правильно решает задачу с вероятностью $0,95$, а двоечник с вероятностью $0,15$. Сколько задач достаточно дать на зачете и сколько требовать решить, чтобы отличник не сдал зачет с вероятностью, не большей $0,01$, а при этом двоечник сдал зачет с вероятностью, не большей $0,1$?
\end{problem}

\begin{problem}(Задача про мужика.)
Пусть некоторый мужик говорит правду с вероятностью 75\% и лжет с 
вероятностью 25\%. Он подбрасывает симметричную кость и говорит, 
что ``выпала 6''. C какой вероятностью выпала 6?
\end{problem}

\begin{problem}
\label{exp_eps}
Пусть с.в. $X$ имеет показательное распределение $\mathrm{Exp}(\lambda)$. Покажите, что имеет место ``отсутствие последействия'': 
$$
{\mathbb P}(X \ge x+y\,|\, X \ge x)={\mathbb P}(X \ge y) . 
$$
Верно ли обратное утверждение? Приведите пример дискретной с.в. со свойством ``отсутствия последействия''.
\end{problem}
\begin{ordre}
В качестве дискретной с.в. рассмотрите с.в., имеющую \textit{геометрическое распределение}. А именно, с.в. $X$ -- момент первого выпадения орла, в схеме испытаний Бернулли $Be(p)$ (независимо подкидывается монетка с вероятностью выпадения орла равной $p$). Для такой с.в.
$$\PR\left(X\ge n\right) = \left(1 - p\right)^n, \quad n\in\mathbb{N}.$$
Следует сравнить со случаем, когда $X \in \mathrm{Exp}(\lambda)$. В этом случае
$$\PR\left(X\ge t\right) = \exp\left(-\lambda t \right), \quad t\ge 0.$$
Оказывается, что этими примерами исчерпывается список распределений, удовлетворяющих свойству отсутствия последействия.

Подобно тому как сумма $n$ независимых с.в., с расределением $Be(p)$ порождает биномиальную с.в. $Bi(p,n)$ -- число успехов в схеме испытаний Бернулли (независимо $n$ раз подкидывается монетка с вероятностью успеха $p$), сумма $n$ независимых с.в., распределенных по геометрическому закону с параметром $p$, порождает \textit{отрицательное биномиальное распределение} с параметрами $p$, $n$, также называемое \textit{распределенинем Паскаля}.
\end{ordre}
\begin{remark}
Следует также обратить внимание на задачу~\ref{sec:poisson} раздела~\ref{zb4}. Описанные в задаче распределения естественным образом возникают при построении марковских процессов в непрерывном и дискретном времени (см. раздел 6 и Часть 2). 
\end{remark}

\begin{problem}
В $m+1$ урнах содержится по $m$ шаров (в каждой), причем урна с номером $n$ содержит $n$ белых и $m-n$ черных шаров $(n = 0,1,\ldots,m)$. 
Случайным образом выбирается урна и из нее $k$ раз с возвращением извлекаются шары. Найдите: 
\begin{enumerate}
\item[а)] вероятность того, что следующим также будет извлечен белый шар при условии, что все $k$ шаров оказались белыми, 
\item[б)] ее предел при $m\to\infty$. 
\end{enumerate}
\end{problem}

\begin{ordre}
Применить \textit{формулу полной вероятности} в следующем виде: 
$$
{\mathbb P}(B\, |\, A)=\sum\limits_{n=1} {\mathbb P}(B\, |\, H_n A){\mathbb P}(H_n\, |\, A).
$$
\end{ordre}

\begin{problem}
    Береговая охрана замечает вражеский корабль на расстоянии 1~км от побережья и начинает его обстреливать, делая 1 выстрел в минуту. После первого же выстрела корабль разворачивается и движется в сторону от побережья со скоростью 60 км/ч. Вероятность попадания обратно пропорциональна квадрату расстояния от побережья до корабля (если цель находится на расстоянии $x$, тогда вероятность попадания $0,75x^{-2}$). Вероятность того, что корабль выдержит $n$ попаданий и не утонет, равна $4^{-n}$. Посчитайте вероятность того, что кораблю удастся уйти.
\end{problem}

\begin{problem}
Параметр $p$ (вероятность выпадения <<орла>> при бросании монеты) в схеме испытаний Бернулли является равномерно распределенной случайной величиной на отрезке $[0.1, 0.9]$ и разыгрывается до начала испытаний, причем $p$ не изменяется от опыта к опыту. В серии из $n=1000$ бросаний было подсчитано число успехов $r=777$.  Найдите условную плотность распределения $p(x|r=777)$. Оцените, как изменится ответ, если известно, что точное значение числа успехов $r\in[750, 790]$.
\begin{remark}
Эта задача на формулу Байеса и формулу полной вероятности. С применением этих формул связано и решение следующей задачи.

Отметим, что во многих современных приложениях вероятностых методов (машинное обучнение, теория информации и др.) байесовские подход играет важную роль MacKay D. Information Theory, Inference, and Learning Algorithms. Cambridge University Press, 2003. 640 p.;  Jaynes E.T. Probability Theory: The Logic of Science. Cambridge University Press, 2003. 753 p.; Мамфорд Д. На заре эры стохастичности. В сборнике ``Математика: границы и перспективы''. М.: ФАЗИС, 2005. С. 327--354; Bishop C. Pattern Recognition and Machine Learning. Springer, 2006. 738 p. Более подробно о ряде таких приложений будет написано в Части 2. В России в этом направлении активно работают Дмитрий Ветров и Дмитрий Кропотов с ФКН НИУ ВШЭ (ГУ). 

Отметим также, что формула Байеса и формула полной вероятности совершенно естественным образом возникает в генетики, например, при решении задач на применение \textit{законов Менделя, Харди--Вайнберга} (см., например, Тейлор Д., Грин Н., Стаут У. Биология. Под ред. Р. Сопера. Т. 3. М.: Лаборатория Базовых знаний, 2016. 451 с.). Популярно об этом написано, например, в книге Болтянский В.Г., Савин А.П. Беседы о математике. Кинига 1. Дискретные объекты. М.: ФИМА, МЦНМО, 2002. 368 с.
\end{remark}
\end{problem}

\begin{problem}(Задача Лапласа.)
\label{laplace}
Пусть имеется монетка, о которой в начальный момент известно, что вероятность выпадения орла $p$ может быть произвольным числом из отрезка $[0,1]$ с равной вероятностью. Таким образом априорное распределение вероятности выпадения орла $p$ -- равномерное распределение на отрезеке $[0,1]$, т.е. $R[0,1]$. Следовательно, априорно ожидаемое значение $p = 0,5$. Монетку независимо подбросили $n$ раз и посчитали общее число успехов $r$. Определите апостериорное распределение $p|\left(n,r\right)$ и ожидаемое значение $p$ с учетом результатов опыта $\left(n,r\right)$. 
\begin{ordre} (Бета-распределение, как сопряженное распределение к биномиальному распределению)
Для решения задачи полезно ввести бета-распределение $\Beta$ с параметрами $\alpha_1, \alpha_2$, задаваемое плотностью ($\Gamma(\alpha)$ -- гамма-функция Эйлера, в частности $\Gamma(m+1) = m!$)
$$p_{\Beta}(x) = \frac{\Gamma(\alpha_1+\alpha_2)}{\Gamma(\alpha_1)\Gamma(\alpha_2)}x^{\alpha_1 - 1}(1-x)^{\alpha_2 - 1}, \quad x\in[0,1].$$
И воспользоваться следующим общим результатом: если априорное распределение $p \in \Beta(\alpha_1,\alpha_2)$ (в условиях задачи $\alpha_1 = \alpha_2 = 1$), то апостериорное распределение $p|\left(n,r\right) \in \Beta(\alpha_1 + r,\alpha_2 + n - r)$. Также отметим, что для подсчета ожидаемого (среднего) значения с.в. $\xi$, распределенной по закону $\Beta(\alpha_1,\alpha_2)$ можно воспользоваться формулой $\Exp \xi = \alpha_1/(\alpha_1 + \alpha_2)$.
\end{ordre}
\begin{remark} (Распределение Дирихле, гамма-распределение \cite{202}.)
В более общем случае рассмотрим не монетку, а многогранник, при бросании которого он выпадает на одну из своих $m$ граней (в такой интерпретации $m > 3$, однако приводимые далее формулы справедливы при $m\ge 2$, где случай $m=2$ сводится к разобранному в указании бета-распределению). До начала проведения опыта можно считать, что вероятности упасть на каждую грань равны между собой. Требуется определить апостериорные вероятности упасть на каждую из граней после проведения опыта $\left(n_1,r_1; ...; n_m,r_m\right)$. Для решения этой задачи полезно ввести распределение Дирихле $\Dir$ с параметрами $\alpha_1, ..., \alpha_m$, задаваемое плотностью (приводимая ниже формула внешне похожа на формулу, определяющую мультиномиальное распределение, см., например, задачу 19 раздела 6; это совпадение совсем не случайно -- чуть ниже будет показано, что распределение Дирихле является сопряженным распределением для мультиномиального распределения в том же смысле, что бета-распределение является сопряженным распределением для биномиального)
$$p_{\Dir}(x_1,...,x_m) = \frac{\Gamma(\alpha_1+...+\alpha_m)}{\Gamma(\alpha_1)\cdot ...\cdot\Gamma(\alpha_m)}x_1^{\alpha_1 - 1}\cdot ...\cdot x_m^{\alpha_m - 1},$$ где $\sum_{k=1}^m{x_k}=1, x_k\in[0,1], k=1,...,m$, и воспользоваться следующим общим результатом: если априорное распределение $p \in \Dir(\alpha_1,...,\alpha_m)$ (в условиях новой задачи $\alpha_1 =...= \alpha_m = 1$), то апостериорное распределение $p|\left(n_1,r_1; ...; n_m,r_m\right)\in \Dir(\alpha_1 + r_1,...,\alpha_m + r_m)$. 

Введем новый класс с.в. (гамма-распределение) $X\in \Gamma(\lambda,\alpha)$ -- неотрицательные случайные величины с плотностью распределения
$$
p_{\Gamma}(x) = \dfrac{\left(x/\lambda\right)^{\alpha - 1}{\rm e}^{-x/\lambda}}{\Gamma(\alpha)}, \quad  x\geq 0. 
$$
В частности, $\Gamma(\lambda^{-1},1) = \mbox{Exp}(\lambda)$ (определение показательного распределения см. в задаче \ref{exp_eps}). Более того, гамма-распределение $\Gamma(\lambda^{-1},\alpha)$ в случае натурального $\alpha = n$ можно представить в виде суммы $n$ независимых одинаково распределенных показательных с.в. $\mbox{Exp}\left( \lambda \right)$:
 $$\Gamma \left( {\lambda^{-1} ,n} 
\right) \mathop =\limits^d \underbrace {\mbox{Exp}\left( \lambda \right)+...+\mbox{Exp}\left( 
\lambda \right)}_n.$$
В математической статистике также бывает полезно следующее представление распределение $\chi^2(n) = \Gamma \left(2 ,n/2\right)$ в виде суммы $n$ независимых с.в., распределенных по закону $\N(0,1)^2$ (где $\N(0,1)$ -- стандартная нормальная с.в.):
$$\Gamma \left(2 ,n/2\right) \mathop =\limits^d \underbrace {\N(0,1)^2 +...+\N(0,1)^2}_n.$$
Если $Y_k \in\Gamma(1,\alpha_k), k = 1,...,m$ -- независимые с.в. и $Y = Y_1 + ... + Y_m$, то (см. также задачу 13 раздела 7) $$\left(\frac{Y_1}{Y},...,\frac{Y_m}{Y} \right)^T \in \Dir(\alpha_1,...,\alpha_m).$$
Гамма-распределение внешне похоже на распределение Пуассона $Z\in\Po(\lambda)$ (см. также задачу \ref{sec:poisson} раздела \ref{zb4}):
$$
Z = \max\Bigl\{ n:\; \sum\limits_{k=1}^{n} X_k<1 \Bigr\}, 
$$
где $X_1, X_2, X_3,\ldots$ -- независимые одинаково распределенные по закону ${\rm Exp}(\lambda)$ с.в. Действительно,
$$\PR\left(Z = n \right) = \dfrac{\lambda^n{\rm e}^{-\lambda}}{\Gamma(n+1)}, \quad n = 0,1,2,...$$
Оказывается, эта связь не только ``внешняя''. Если априорно было известно, что $\lambda \in \Gamma(\lambda,\alpha)$, то после наблюдения $n$ независимых реализаций $Z_k \in \Po(\lambda), k = 1,...,n$ имеем $$\lambda|\left(Z_1,...,Z_n\right) \in \Gamma\left(\frac{\lambda}{\lambda\sum_{k=1}^n{Z_k} + 1},\alpha + n\right).$$
\end{remark}
\end{problem}

\begin{problem}
Пусть $X\in\Po(\lambda)$ и $Y\in\Po(\mu)$ -- независимые случайные величины,  имеющие распределение Пуассона. 
Доказать, что случайная величина $Z = X + Y$ имеет распределение Пуассона $\Po(\lambda + \mu)$. 
Найдите вид условного распределения случайной величины $X$ при фиксированном значении случайной величины $Z$.
\end{problem}

\begin{problem}
Покажите, что борелевская $\sigma$-алгебра в ${\mathbb R}^1$, содержащая все числовые промежутки вида $[a,b)$, 
содержит все промежутки вида $(a,b)$, $(a,b]$, $[a,b]$ и отдельные точки прямой. 

\begin{ordre}
Учесть свойство замкнутости $\sigma$-алгебры относительно операций 
объединения, пересечения и вычитания. 
\end{ordre}

\end{problem}

\begin{problem}
Пусть $\Omega = [a,b]\times...\times[a,b]$, $\mathcal{F}$ -- $\sigma$–алгебра, содержащая всевозможные множества вида 
$[\alpha_1,\beta_1]\times...\times[\alpha_n,\beta_n]$ $(a \leqslant \alpha_1 < \beta_1 \leqslant b, \ldots, a \leqslant \alpha_n < \beta_n \leqslant b)$. На этих множествах задана вероятностная мера 
\[
{\mathbb P}\left(\omega_1  \in [\alpha_1,\beta_1], ..., \omega_n  \in [\alpha_n,\beta_n]  \right)=\dfrac{\mes[\alpha_1, \beta_1]}{\mes[a,b]}  \ldots  \dfrac{\mes[\alpha_n, \beta_n]}{\mes[a,b]}.
\]
Покажите, что 
\begin{enumerate}
\item ${\mathbb P}\left( \omega_1=c=\const\right)=0$; 
\item ${\mathbb P}\left( \omega_1=\omega_2\right)=0$.
\end{enumerate}
Найдите вероятность ${\mathbb P}\left( \omega_1 \leq \omega_3 \leq \omega_2\right)$. 

\end{problem}

\begin{problem}
Число элементарных событий некоторого вероятностного $\left<\Omega,\mathcal{F},\PR\right>$ пространства равно $n$. Укажите минимальное и максимальное возможные значения для числа событий. Другими словами, известно, что $|\Omega| = n$, чему может быть равно в этом случае $|\mathcal{F}|$?
\end{problem}

\begin{problem}
Может ли число всех событий какого-либо вероятностного пространства быть равным $128$; $129$; $130$? 
\end{problem}

\begin{problem}(Задача о встрече \cite{2}.)
\label{L_extension}
Два лыжника условились о встрече между $10$ и $11$ часами утра у подножия горы (склона), причем договорились ждать друг друга не более $10$ минут, чтобы не замерзнуть. Считая, что момент прихода на встречу каждым выбирается наудачу в пределах указанного часа, найдите вероятность того, что встреча состоится. 
\end{problem}

\begin{remark}
В подобного рода задачах, где явно не определена вероятностная тройка $\left< \Omega, \mathcal{F}, \PR \right>$, в общем случае следует придерживаться такой последовательности действий: 
\begin{enumerate}
\item задать вероятностную меру на множествах простой структуры;
\item продолжить вероятностную меру на множества более сложной структуры, рассматриваемые в задаче;  
\item проверить условия, гарантирующие $\sigma$-аддитивность продолженной меры ($\PR(\bigcup_{i=1}^{n} A_i) \to \PR(A)$ при $n\to\infty$, где $\bigcup_{i=1}^{\infty} A_i = A$).   
\end{enumerate}
Детали описаны во многих учебниках (см., например, \cite{3}, \cite{21}), см. также
http://www.mccme.ru/ium/postscript/s12/gasnikov-teoria-mery.pdf. 

}
\end{remark}

\begin{problem}(Парадокс Бертрана \cite{2}, \cite{book12}.)
Рассмотрим окружность, описанную вокруг равностороннего треугольника. Какова вероятность, что случайным образом проведенная хорда будет иметь длину большую, чем сторона этого треугольника?
\begin{remark}
Парадокс Бертрана -- классическая (и, вероятно, ``самая главная'') задача на ``геометрические вероятности'' (см., например, Кендалл М., Моран П. Геометрические вероятности. М.: Наука, 1972. 192 с.). Это направление представляется достаточно важным не только с точки зрения демонстрации аксиоматики теории вероятностей, построенной А.Н. Колмогоровым \cite{3}, но и для решения большого количества практических задач. Следующая задача также на ``геометрические вероятности'' и на первый взгляд кажется достаточно практичной, в этой связи заметим, что так как описано в этой задаче число $\pi$ на практике, на самом деле, никто не оценивает. В Части 2 ``геометрическим вероятностям'' планируется посвятить отдельный раздел.
\end{remark}
\end{problem}

\begin{problem}(Задача Бюффона.)
Как с помощью листа тетрадки в линейку (с расстояними между линиями 1 см.) и иголки длиной 1,5 см., согнутой в кочергу (буквой Г), определить число $\pi$? Сколько раз достаточно бросить иголку, чтобы с вероятностью не менее 0,99 можно было оценить число $\pi$ (по общему числу пересечений иголки с линиями тертадки) с точностью до 3 знаков после запятой?
\begin{remark}
См. Кендалл М., Моран П. Геометрические вероятности.  М.: Наука, 1972.
\end{remark}
\end{problem}
\begin{problem}
Имеется монетка с вероятностью выпадения орла $p=1/2$, которую независимо подкидывают. Покажите, что событие $A$, заключающееся в том, что отношение числа успехов к общему числу бросаний $n$ стремится к $1/2$ при $n \to \infty$ измеримо. Найдите вероятность этого события $\PR\left(A\right)$.
\end{problem}
}

\begin{problem} (Случайное направление \cite{book2012}.) 
вы приходите на станцию метро в случайный момент времени $T$ (например, $T$ имеет равномерное распределение на интервале между 10 и 11 часами) и садитесь в первый пришедший поезд (в ту или другую сторону). Будем считать, что поезда в обе стороны ходят одинаково часто (например, каждые 10 минут) согласно расписанию. Одинаковы ли вероятности событий, что вы поедете в одну или в другую сторону? Рассмотрите случай, когда в одну сторону поезда идут в 10:00, 10:10, 10:20,\dots , а в другую -- в 10:02, 10:12, 10:22,\dots 
\end{problem}

\begin{remark} Казалось бы, интуиция подсказывает, что оба направления ``равноправны'' -- ведь мы приходим в ``случайный'' момент времени. ``Парадокс'' легко разрешается введением вероятностного пространства (конкретизацией ``случайности'' времени прихода на станцию).
\end{remark}

\begin{problem}
В урне находится $3$ белых и $2$ черных шара. 
Эксперимент состоит в последовательном извлечении из урны всех шаров по одному наугад без возвращения. Постройте вероятностное пространство. 
Описать $\sigma$-алгебру, порожденную случайной величиной $X$, если: 
\begin{enumerate}
\item $X$ -- число белых шаров, предшествующих первому черному шару; 
\item $X$ -- число черных шаров среди извлеченных; 
\item $X=X_1+X_2$, где $X_1$ -- число белых шаров, предшествующих первому черному шару, 
$X_2$ -- число черных шаров, предшествующих белому шару. 
\end{enumerate}
\end{problem}

\begin{problem}\Star
\label{SigmaAlgebra}
Пусть $(\Omega,\mathcal{F},{\mathbb P})$ -- некоторое вероятностное пространство и $A$ -- алгебра подмножеств $\Omega$ такая, что 
$\sigma(A)=\mathcal{F}$\quad $(\sigma(A)$ -- наименьшая \mbox{$\sigma$-алгебра}, содержащая алгебру $A)$. Докажите, что 
$$
\forall\,\varepsilon>0,\, B\in\mathcal{F}\quad \exists A_{\varepsilon}\in A:\quad {\mathbb P}(A_{\varepsilon}\bigtriangleup B)
\leqslant\varepsilon . 
$$
\end{problem}

\begin{ordre}
Рассмотрите совокупность множеств (см. детали в \cite{21} Т.~2): 
$$
{\mathcal B}=\bigl\{ B\in\Sigma\, | \, \forall\varepsilon>0 \; \exists A_B\in A:\; {\mathbb P}(A_B\bigtriangleup B)
\leqslant\varepsilon \bigr\} . 
$$

\noindent Покажите, что ${\mathcal B}$ является минимальной $\sigma$-алгеброй (см. \cite{22}, \cite{220}).

\end{ordre}

 \begin{problem}
 \begin{enumerate}
 
 \item Может ли функция
\[
F\left( {x_1 ,x_2 } \right)=\left\{ {\begin{array}{l}
 1,\quad \min \left( {x_1 ,x_2 } \right)>1, \\ 
 0,\quad \min \left( {x_1 ,x_2 } \right)\le 1 \\ 
 \end{array}} \right.
\]
быть функцией распределения некоторого двумерного случайного вектора? 
Приведите необходимое и достаточное условие того, чтобы $F\left( {x_1 
,...,x_n } \right)$ была функцией распределения некоторого случайного 
вектора.

\item На ${\mathbb R}^n$ задана функция $F\left( {x_1,\ldots,x_n } \right)$, 
неубывающая по каждому аргументу, такая, что $F\left( {x_1 ,\ldots,-\infty 
,\ldots, x_n } \right)=0$, $F\left( {\infty,\ldots,\infty } \right)=1$. Что еще 
надо потребовать от этой функции, чтобы с ее помощью  можно было бы 
естественным образом задать счетно-аддитивную вероятностную меру на 
минимальной $\sigma $-алгебре, содержащей всевозможные прямоугольные 
параллелепипеды в ${\mathbb R}^n$?

\item Может ли функция распределения $F_X \left( x \right)$ с.в. $X$ 
иметь более чем счетное множество точек разрыва?

\end{enumerate}
\end{problem}

\begin{problem}(Обобщённое неравенство Чебышёва.)
\label{mark_cheb}
Пусть $g$ -- неотрицательная и неубывающая функция такая, что $\Exp g(|\xi|) < \infty$. Покажите, что для всех $x$, удовлетворяющих условию $g(x) > 0$,
\[
\PR(|\xi| \geq x) \leq \frac{ \Exp g(|\xi|) }{ g(x)}.
\]
Докажите следствия этого неравенства (см. \cite{21}):
\begin{enumerate}
\item (Неравенство Маркова) Пусть $\Exp |\xi|^p < \infty$, $p > 0$, $x > 0$. Тогда
\[
\PR(|\xi| \geq x) \leq \frac{ \Exp |\xi|^p }{ x^p}.
\]
\item (Неравенство Чебышёва) Пусть $\Var \xi < \infty$, $x > 0$. Тогда
\[
\PR(|\xi - \Exp \xi| \geq x) \leq \frac{ \Var \xi }{ x^2}.
\]
\item (Экспоненциальное неравенство Чебышёва, неравенство Чернова) Пусть $\Exp \left[e^{t \xi}\right] < \infty$, $t > 0$. Тогда
\[
\PR(\xi \geq x) \leq \frac{ \Exp \left[e^{t \xi}\right] }{ e^{t x} }.
\]
\item (Неравенство Кантелли) Пусть $\Var \xi < \infty$, $x > 0$. Тогда
\[
\PR( |\xi - \Exp \xi | \geq x) \leq \frac{ 2\Var \xi }{ x^2 + \Var \xi}.
\]
\end{enumerate}

\end{problem}

\begin{problem}
Докажите \textit{неравенство Ляпунова} для $0 < p < q$ 
\[
(\Exp |\xi|^p)^{1/p} \leq (\Exp |\xi|^q)^{1/q}. 
\]
\end{problem}

\begin{ordre}
Воспользуйтесь неравенством Йенсена: $\Exp |\xi| < \infty$, $g$ -- выпуклая функция:
\[
g(\Exp \xi) \leq \Exp g(\xi). 
\]

\end{ordre}

\begin{problem}
Докажите \textit{неравенство Гёльдера}. Пусть $1 < p < \infty$, $1/p + 1/q = 1$, $\Exp |\xi|^p < \infty$, $\Exp |\eta|^q < \infty$. Тогда  
\[
\Exp |\xi \eta | \leq (\Exp |\xi|^p)^{1/p} (\Exp |\eta|^q)^{1/q}. 
\]
\end{problem}

\begin{ordre}
Воспользуйтесь неравенством Юнга: $x > 0$, $y > 0$: 
\[
x y \leq \frac{x^p}{p} + \frac{y^q}{q}.  
\]

\end{ordre}

\begin{problem}
Докажите \textit{неравенство Минковского}: $1 \geq p < \infty$, $\Exp |\xi|^p < \infty$, $\Exp |\eta|^p < \infty$:  
\[
(\Exp |\xi + \eta |^p)^{1/p} \leq (\Exp |\xi|^p)^{1/p} + (\Exp |\eta|^p)^{1/p}. 
\]
\end{problem}

\begin{ordre}
Воспользуйтесь неравенством $p > 0$: 
\[
\left| \sum \limits_{k=1}^n x_i \right|^p \leq \max\{1, n^{p-1}\} \sum \limits_{k=1}^n |x_i|^p  
\]
и неравенством Гёльдера.
\end{ordre}

\begin{problem}
Покажите, что неравенство Чебышёва:
\[P\left(\left|X-EX\right|>\varepsilon \right)\le \frac{DX}{\varepsilon ^{2} } \] 
принципиально не улучшаемо.

\begin{ordre} 
Рассмотрите с.в. 
\[X=\left\{\begin{array}{cc} {a,} & {{\raise0.7ex\hbox{$ p $}\!\mathord{\left/ {\vphantom {p 2}} \right. \kern-\nulldelimiterspace}\!\lower0.7ex\hbox{$ 2 $}} } \\ {0,} & {1-p} \\ {-a,} & {{\raise0.7ex\hbox{$ p $}\!\mathord{\left/ {\vphantom {p 2}} \right. \kern-\nulldelimiterspace}\!\lower0.7ex\hbox{$ 2 $}} } \end{array}.\right.\] 
Положите $\varepsilon =a-\delta $, где $\delta \to 0+$. См. также Червоненкис~А.Я. Компьютерный анализ данных. М.: Яндекс, 2009.
\end{ordre}

\begin{remark}
В частности, данное неравенство не улучшаемо в классе $X \in L_2$. 
\end{remark}

\end{problem}

\begin{problem}
Существует ли случайная величина с конечным вторым моментом и бесконечным первым моментом?
\end{problem}

\begin{problem}
В начале карточной игры принято с помощью жребия определять первого сдающего. Для этого колода хорошо тасуется, и игрокам сдается по одной карте до появления первого туза. Кому выпал туз, тот и сдает в игре. На каком месте в среднем появляется первый туз, если в колоде 32 карты (то есть нужно найти математическое ожидание случайной величины «Число карт, сданных до первого туза»)?
\end{problem}

\begin{problem}(А. Шень \cite{book2012}.)
 В лотерее на выигрыш уходит 40\% от стоимости проданных билетов. Каждый билет стоит 100 рублей. Докажите, что вероятность выиграть 5000 рублей (или больше) меньше 1\%.

Искомая вероятность зависит, конечно, от правил лотереи, но ни при каких условиях она не превосходит $1\%$.
Приведите пример правил лотереи, где искомая вероятность минимальна и максимальна.

\begin{ordre} 
Воспользуйтесь неравенством Маркова (см. задачу \ref{mark_cheb}).
\end{ordre}

\end{problem}

\begin{problem}
Восемь мальчиков и семь девочек купили билеты в кинотеатр на $15$ подряд идущих мест. Все $15!$ возможных способов рассадки равновероятны. Вычислите среднее число пар рядом сидящих мальчика и девочки. Например, (\mars, \female, \mars, \female) содержит три такие пары. 
\end{problem}

\begin{problem}
На первом этаже семнадцатиэтажного общежития в лифт вошли десять человек. Предполагая, что каждый из вошедших (независимо от остальных) может с равной вероятностью жить на любом из шестнадцати этажей (со 2-го по 17-й), найдите математическое ожидание числа остановок лифта.
\end{problem}

\begin{problem}
\label{sec:latters}
Имеется $n$ пронумерованных писем и $n$ пронумерованных конвертов. Письма случайным образом раскладываются по конвертам (все $n!$ 
способов равновероятны). 
Найдите математическое ожидание числа совпадений номеров письма и конверта (письмо лежит в конверте с тем же номером). 
\end{problem}

\begin{problem}
\label{mom_ineq}
С.в. $\xi$ имеет ограничение $\PR(\xi \geq h(x)) \leq e^{-x}$ при $x \geq x_0$, $h(x)$ -- абсолютно непрерывная функция. Докажите следующие неравенства для первого и второго моментов $\xi$:
\[
\Exp \xi \leq h(x_0) + \int_{x_0}^{\infty}h'(x)e^{-x}dx,
\]
\[
\Exp \xi^2 \leq h^2(x_0) +  2\int_{x_0}^{\infty}h(x) h'(x)e^{-x}dx.
\]
\end{problem}

\begin{ordre}
Удобно воспользоваться следующей формулой для подсчета математического ожидания: 
\[
\Exp \xi = \int\limits_{0}^{+\infty} \PR(\xi \geq x)dx - 
    \int\limits_{-\infty}^{0}\PR(\xi < x)dx.
\]
\end{ordre}

\begin{problem}
Случайным образом выстраиваются в шеренгу  $n$ человек разного роста. Найдите вероятность того, что
\begin{enumerate}
\item[а)] самый низкий окажется $i$-м слева; 
\item[б)] самый высокий окажется первым слева, а самый низкий --- последним слева; 
\item[в)] самый высокий и самый низкий окажутся рядом; 
\item[г)] между самым высоким и самым низким расположатся более $k$ человек. 
\end{enumerate}
\end{problem}

\begin{problem}
Урна содержит $N$ шаров с номерами от $1$ до $N$. Пусть $K$ -- наибольший номер, полученный при $n$ их поштучных извлечениях 
с возвращением. Найдите: 
\begin{enumerate}
\item распределение $K$;
\item асимптотику математического ожидания ${\mathbb E}K$ при $N\to\infty$. 
\end{enumerate}
\end{problem}

\begin{ordre}
\[
{\mathbb E}K=\sum\limits_{j=0}^{N} {\mathbb P}(K>j).
\]
\end{ordre}

\begin{problem}
Найдите математическое ожидание ${\mathbb E}Z$, где $Z$ -- $k$-я по величине из координат $n$ точек, взятых наудачу на отрезке 
$[0, 1]$ $(k \leqslant n)$. 
\end{problem}

\begin{problem}
Покажите, что если независимые с.в. $X_1,\dots,X_n$ имеют   
\textit{показательное распределение}, т.е. 
$$
f_{X_i}(x)=\begin{cases}
\lambda_i\exp(-\lambda_i x), \; x\geqslant 0, \\
0,\; x<0
\end{cases}
$$
(часто пишут $X_i\in \mathrm{Exp}(\lambda_i)$), то 
$$
\min\{ X_1,\ldots, X_n\}\in \mathrm{Exp}\Bigl( \sum\limits_{i=1}^{n}\lambda_i\Bigr) . 
$$
\end{problem}

\begin{problem}
В некотором вузе проходит экзамен. Количество экзаменационных билетов $N$. Перед экзаменационной аудиторией выстроилась очередь из 
студентов, которые не знают, чему равно $N$. Согласно этой очереди студенты вызываются на экзамен (второй студент заходит в аудиторию 
после того, как из нее выйдет первый, и т.д.). Каждый студент с равной вероятностью может выбрать любой из $N$ билетов (в независимости 
от других студентов). Проэкзаменованные студенты, выходя из аудитории, сообщают оставшейся очереди номера своих билетов. 
Оцените (сверху), сколько студентов должно быть проэкзаменовано, чтобы оставшаяся к этому моменту очередь смогла оценить число 
экзаменационных билетов с точностью $10\%$ с вероятностью, не меньшей $0.95$. 
\end{problem}


\begin{problem}
Допустим, что вероятность столкновения молекулы с другими молекулами в промежутке времени $[t,t + \Delta t)$ 
равна $p = \lambda\Delta t + o(\Delta t)$ и не зависит от времени, прошедшего после предыдущего столкновения $(\lambda = \const)$. 
Найдите распределение времени свободного пробега молекулы и вероятность того, что это время превысит заданную \linebreak величину $t^*$. 
\end{problem}

\begin{ordre}
Разобьем интервал $\Delta=[0,t)$ на $n$ отрезков равной длины.
Пусть $A_i$ -- событие, означающее, что на отрезке $\Delta_i$  молекула претерпит столкновение с другими молекулами. Можно представить вероятность отсутствия столкновения в виде произведения вероятностей событий $\overline{A_i}$.
\end{ordre}

\begin{problem}(Задача об оптимальном моменте замены оборудования.)
Пусть некоторая система состоит из $n$ элементов, выходящих из строя независимо друг от друга через случайное время $\tau _{i} $, имеющее показательное распределение с параметром $\lambda _{i} $ ($i$ -- номер элемента, $1\le i\le n$). Вышедший из строя элемент немедленно заменяется новым. Пусть $c_{i} $ -- убытки, связанные с выходом из строя и заменой элемента $i$-го типа. Через промежуток времени $T$ разрешается провести профилактический ремонт, при котором все $n$ элементов заменяются новыми. Пусть $D$ -- стоимость профилактического ремонта. Найдите оптимальное время $T$, минимизирующее средние убытки в единицу времени.
\end{problem}

\begin{remark} 
Примите в учет ``отсутствие последействия'' для показательного распределения из задачи \ref{exp_eps}.
Согласуется ли Ваш результат со здравым смыслом: не стоит проводить профилактический ремонт, если выход оборудования из строя не связан с его старением? Детали см. в книге Афанасьева~Л.Г. Очерк исследования операций. М.: Изд-во МГУ, 2007.  176~с. 
\end{remark}

\begin{problem}
Приведите пример таких с.в. $X$ и $Y$, что $X=-Y$, но $F_X=F_Y$  (имеют одинаковые  функции распределения).
\end{problem}

\begin{problem}
Случайная величина $X$ имеет функцию распределения $F_X$ и функцию плотности распределения $f_X$. Найти функции распределения и плотности распределения (если последние существуют) для случайных величин:\\
\indent а) $Y = aX+b$,\\
\indentб) $Z =e^X$,\\
\indent в) $V = X^2$ ($a$ и $b$ -- неслучайные величины).\\
Конкретизировать решения при $X \in \N(0,1)$.
\end{problem}

\begin{problem}
\label{inverse}
Пусть $F$ -- непрерывная функция распределения скалярной с.в. Пусть 
$$F^{-1}\left( y \right) = \inf \left\{ x: F\left( x \right) = u \right\}, u \in \left[ 0, 1 \right].$$
Тогда если с.в. $y$ -- равномерно распределена на отрезке $\left[ 0, 1 \right]$, то с.в. $F^{-1}\left( u \right)$ имеет распределение с функцией распределения $F$.
\end{problem}
\begin{remark}
Эта задача играет важную роль в методах Монте-Карло (см. главу 7 и \cite{24}).
\end{remark}

\begin{problem}
С.в. $X$, $Y$, $Z$ независимы и равномерно распределены на $(0,1)$. Докажите равенство распределений  $(XY)^{Z}$  и $X$.
\end{problem}

\begin{ordre}
Докажите, что величина $-Z(\ln X + \ln Y)$ показательно распределена на $[0,\infty)$. Для этого удобно воспользоваться заменой $U = Z W$, $V = W/Z$, где $W = - \ln X - \ln Y$, $f_W(w) = w e^{-w} [w>0]$. 
\end{ordre}

\begin{problem}
Случайные величины $X_1,\ldots,X_n$ независимы и имеют показательное распределение с параметром $\lambda$.
Докажите равенство распределений $(X_1+  \frac{1}{2} X_2 + \ldots+ \frac{1}{n} X_n)$ и  $\max\{X_1,\ldots,X_n\}$.
\end{problem}

\begin{problem}
Рассмотрим прибор, который может ломаться, и устройство, которое чинит этот прибор мгновенно, когда тот выходит из строя. При этом считается, что времена работы прибора от поломки до поломки распределены как $\xi \in \text{Exp}(\lambda)$ и не зависят друг от друга. Чинящее устройство, в свою очередь, тоже может ломаться. Распределение времени работы этого устройства $\eta \in \text{Exp}(\mu)$ и не зависит от времени работы прибора. Найти распределение времени работы всей системы в целом. 
\end{problem}
\begin{remark}
См. задачу \ref{exp_eps}. 
\end{remark}

\begin{problem}(А.Н. Соболевский.)
Пусть плотность $p(x,y)$ задает равномерное распределение на единичном квадрате $(0, 1)^2$. Проверьте, что $p(x | X =$\linebreak $= Y )$, понимаемая как предел при $\delta \to 0$ условной плотности $p(x | -\delta <$\linebreak $< X - Y <  \delta)$, постоянна и равна~1, в то время как аналогичный предел для $p (x | 1 - \delta < X/Y < 1 + \delta)$ равен $2x$. 
\end{problem}

\begin{problem}
Пусть $X$ и $Y$ --- независимые случайные величины, равномерно распределенные на $(-b,b)$. 
Найдите вероятность $q_b$ того, что уравнение $t^2+tX+Y=0$ имеет действительные корни. Докажите, что 
существует $\lim\limits_{b\to\infty} q_b=q$. Найдите $q$. 
\end{problem}

\begin{problem}
Пусть $X=(X_1,...,X_n)$ -- случайный вектор, компоненты которого имеют нормальное распределение, т.е. $X_i \sim \mathcal{N}(\mu_i,\sigma^2_i)$. Следует ли из этого, что вектор $X$ имеет многомерное нормальное распределение?
\end{problem}

\begin{problem}
Двумерный случайный вектор $X=(X_1,X_2)$  имеет следующую функцию плотности распределения: 
$$
f(x_1,x_2)=\begin{cases}
\dfrac{c}{\sqrt{x_1^2+x_2^2}} \text{ при } x_1^2+x_2^2\leqslant 1 , \\
\quad 0 \quad\text{ иначе}. 
\end{cases}
$$
\begin{enumerate}
\item Найдите $c$.
\item Найдите частные и условные распределения его компонент. 
\item Являются ли они 
\begin{enumerate}
\item стохастически зависимыми; 
\item коррелированными? 
\end{enumerate}
\end{enumerate}
\end{problem}

\begin{problem}
Попарно некоррелированные с.в. $\xi, \eta, \zeta$ обладают одинаковыми математическими ожиданиями и дисперсиями $\Exp\xi = \Exp\eta = \Exp\zeta = 0$, $\Var\xi = \Var\eta = \Var\zeta = \sigma^2$. Найдите $\max \Exp(\xi\eta\zeta)$ и $\min \Exp(\xi\eta\zeta)$.
\end{problem}

\begin{problem}
\label{scalar_prod}
Будем говорить, что  с.в. $\xi$ принадлежит классу $L_2$, если  выполнено ${\mathbb E}\xi^2<\infty$. $L_2$ является линейным пространством, в котором можно задать скалярное произведение следующим образом:
\[
\langle \xi_1, \xi_2 \rangle = \Exp\left[\xi_1 \xi_2\right].
\]
Для величин в $L_2$ введем среднеквадратичное расстояние
\[
\Vert \xi_1 - \xi_2 \Vert_2 = (\Exp{|\xi_1 - \xi_2|^2})^{1/2}.
\]
Пусть $H$ -- конечномерное  подпространство пространства  $L_2$ и $\widehat{\xi}$ -- проекция с.в. $\xi$ на подпространство $H$.

\begin{enumerate}

\item Покажите, что если $H$ является совокупностью всех постоянных, то 
$\widehat{\xi}  = \Exp{\xi}$ и как следствие 
\[
\Vert \xi - \Exp{\xi} \Vert_2 = \underset{\lambda \in H}{\min} \Vert \xi - \lambda \Vert_2.
\] 

\item Рассмотрим  последовательности с.в. $\xi_n \overset{L_2}{\longrightarrow} \xi$, $\eta_n \overset{L_2}{\longrightarrow} \eta$. Докажите, что в этом случае  $\Exp(\xi_n \eta_n) \rightarrow \Exp(\xi \eta)$. 

\end{enumerate}

\end{problem}

\begin{remark}
Рекомендуем также ознакомиться с аналогичными задачами \ref{condExp1}--\ref{condExp3} из раздела \ref{hard}. 
\end{remark}

\begin{problem}
\label{scalar_prod_1}
Пусть  $\xi_1$, $\xi_2$ $\in L_2$ (см. предыдущую задачу и задачу \ref{condExp1} раздела \ref{hard}) со средними значениями $a_1$ и $a_2$, дисперсиями $\sigma_1^2$ и $\sigma_2^2$ 
и коэффициентом корреляции~$r$. Покажите, что наилучшая линейная оценка
$\widehat{\xi_1}  = c_1 + c_2 \xi_2$ для~$\xi_1$ является 
\[
\widehat{\xi_1}  = a_1 + r \frac{\sigma_1}{\sigma_2} (\xi_2 - a_2).
\]
\end{problem}
\begin{remark}
Коэффициент корреляции определяется как  
$$
r = \frac{\Exp (\xi_1 - a_1) (\xi_2 - a_2)}{\sigma_1 \sigma_2}.
$$
Под наилучшей оценкой $\eta$ понимается такая оценка, которая доставляет решение задачи
$$\Exp\left[\| \eta - \xi_1 \|_2^2\right] \to \min.$$
\end{remark}

\begin{problem} 
\label{normal}
Пусть
$$
X=\begin{pmatrix}
X_1\\
X_2\\
X_3
\end{pmatrix}
\in \N\left(
\begin{pmatrix}
2\\
3\\
1
\end{pmatrix}, 
\begin{Vmatrix}
5 & 2 & 7\\
2 & 5 & 7\\
7 & 7 & 14
\end{Vmatrix}
\right) . 
$$
\begin{enumerate}
\item Найдите распределение случайной величины $Y_1=X_1+X_2-X_3$. 
\item Найдите распределение случайной величины $Y_2=X_1+X_2+X_3$. 
\item Найдите ${\mathbb E}(Y_2\, |\, X_1=5, X_2=3)$. 
\item Найдите ${\mathbb E}(Y_2\, |\, X_1=5, X_2<3)$. 
\item Найдите ${\mathbb P}(Y_2<10\, |\, X_1=5, X_2<3)$.

\item Пусть $X \in \N\left(m, R\right)$, где матрица $R$ -- неотрицательно определена (это означает, что характеристическая функция вектора $X$ представима в виде \[
\phi_{X}(t) = \Exp e^{it^TX} = \exp\left( it^{T}m-(1/2)t^{T}Rt \right)),
\] где верхний индекс $T$ означает транспонирование. Найдите распределение случайного вектора $Y = AX + b$. 

\end{enumerate}
\end{problem}

\begin{ordre}
Покажите, что
$\phi_{AX+b}(t) = e^{ib}\phi_{X}(A^{T}t)$.
\end{ordre}

\begin{problem} (Доли выборки.)
\label{sec:ordered_seq}
Покажите, что случайный вектор $W=(W_1, W_2, \ldots, W_n)^T$, который отражает упорядоченное расположение случайно (равномерно и независимо) брошенных $n$ точек на отрезок $[0,1]$, имеет плотность распределения 
$$
p_W(w)=\begin{cases}
n! , & \text{ если } 0\leqslant w_1\leqslant w_2\leqslant \ldots \leqslant w_n\leqslant 1, \\
0 & \text{в остальных случаях}. 
\end{cases}
$$
Покажите, что распределение вектора $U=(U_1, U_2, \ldots, U_n)^T$, если 
$$
U_1=W_1, \quad U_2=W_2-W_1, \ldots, U_n=W_n-W_{n-1} 
$$
имеет вид
$$
p_U(u)=\begin{cases}
n! , & \text{ если } \sum_{k=1}^n{u_k}\le 1, u_k\ge 0, k=1,...,n, \\
0 & \text{в остальных случаях}. 
\end{cases}
$$
\begin{ordre}
В общем случае: 
$$
U=\varphi(W) \,\Rightarrow\, p_U(u)=p_W\bigl(\varphi^{-1}(u)\bigr) \Bigl| J\Bigl( \frac{\partial W}{\partial U}\Bigr)\Bigr|,$$
$$
J(\cdot)=\det\Bigl( \frac{\partial \varphi_i^{-1}(u)}{\partial u_j}\Bigr) . 
$$
\end{ordre}
\begin{remark}
Эта задача тесно связана с задаче 13 раздела 7 (теорема Пайка): $W_k = X_{\left(k\right)}$. 
\end{remark}
\end{problem}

\begin{problem}(См. \cite{4}.)
Спортсмен стреляет по круговой мишени. Вертикальная и горизонтальная 
координаты точки попадания пули (при условии, что центр мишени -- начало 
координат) -- независимые случайные величины, каждая с распределением 
$\N(0,1)$. Покажите, что расстояние от точки попадания до центра имеет 
плотность распределения вероятностей $r\exp \left( {-{r^2} \mathord{\left/ 
{\vphantom {{r^2} 2}} \right. \kern-\nulldelimiterspace} 2} \right)$ для 
$r\ge 0$. Найдите медиану этого распределения.
\end{problem}

\begin{problem}
Пусть $X_n \overset{\text{ п.н. }}{\longrightarrow} X$, $Y_n \overset{\text{ п.н. }}{\longrightarrow} Y$. Докажите справедливость соотношений
\begin{enumerate}
\item $a X_n + b Y_n \overset{\text{ п.н. }}{\longrightarrow} a X + b Y$, где  $a, b = {\rm const}$;
\item $|X_n| \overset{\text{ п.н. }}{\longrightarrow} |X|$;
\item $X_n Y_n \overset{\text{ п.н. }}{\longrightarrow} XY$;
\item верны ли предыдущие пункты для других типов сходимостей (см. \cite{stoianov})? 
\end{enumerate}

\end{problem}

\begin{problem}
Докажите, что из сходимости по вероятности всегда следует сходимость по распределению. Обратное -- не всегда верно. Приведите контрпример.
\end{problem}

\begin{problem}
Пусть $X_n \overset{d}{\longrightarrow} a$, где $a = {\rm const}$. Докажите, что $X_n \overset{p}{\longrightarrow} a$.
\end{problem}

\begin{problem}
Приведите пример, показывающий, что из сходимости по вероятности не следует сходимость в среднем порядка $p>0$ $(\Exp (\xi_n - \xi)^p \to$\linebreak $\to 0)$ (см. \cite{stoianov}).
\end{problem}

\begin{problem}
Пусть в вероятностном пространстве  $(\Omega, \mathcal{F},\PR)$  $\Omega = [0,1)$, $\mathcal{F}$  -- $\sigma$-алгебра, содержащая полуинтервалы вида $\Omega _{in} =[(i-1)/n,\; i/n)$,  $i = \overline{1,n}$, $n \in \mathbb{N}$  и $\PR$ -- мера  Лебега $(\forall i,n:\PR\{ \omega \in \Omega _{in} \} =1/n)$.

Исследуйте сходимость следующих последовательностей случайных величин: 

$X_{1}^{(1)}  (\omega ),\; \; X_{2}^{(1)}  (\omega ),\; \; X_{2}^{(2)}  (\omega ),\; \; X_{3}^{(1)}  (\omega ),\; \; X_{3}^{(2)}  (\omega ),\; \; X_{3}^{(3)}  (\omega ),\ldots $

в случаях:

\textit{а}) $X_{n}^{(i)} \; (\omega )=n\; \; \text{при}\; \; \omega \in \Omega _{in} ,\; X_{n}^{(i)} (\omega )=0\; \; \text{при}\; \; \omega \in \Omega \backslash \Omega _{in} ;$

\textit{б}) $X_{n}^{(i)} \; (\omega )=n^{-1} \; \; \text{при} \; \omega \in \Omega _{in} ,\; X_{n}^{(i)} (\omega )=0\; \; \text{при} \; \; \omega \in \Omega \backslash \Omega _{in} ;$

\textit{в}) $X_{n}^{(i)} (\omega )=n\; \; \text{при}\; \; \omega \in \Omega _{in} ,\; \, X_{n}^{(i)} \; (\omega )=n^{-1} \; \text{при}\; \; \omega \in \Omega \backslash \Omega _{in} ;$

\textit{г}) $X_{n}^{(i)} \; (\omega )=n^{-1} \; \; \text{при}\; \; \omega \in \Omega _{in} ,\; X_{n}^{(i)} \; (\omega )=1-n^{-1} \; \text{при}$ $\omega \in \Omega \backslash \Omega _{in} ;$

\textit{д}) $X_{n}^{(i)} \; (\omega )\in N(m_{in} ,\; \sigma _{in}^{2} ),\; \, \text{где}\; \mathop{\lim }\limits_{in\to \infty } \; m_{in} =3,\; \mathop{\lim }\limits_{in\to \infty } \sigma _{in}^{2} =1.$
\end{problem}

\begin{problem}(Лемма Бореля--Кантелли.)
\label{bor_kant}
Известно, что для последовательности с.в. $\{X_k\}_{k=1}^\infty$ $\exists X$, что $\forall\,\varepsilon >0$: 
\[
\sum \limits_{k=1}^{\infty} \PR(|X_k - X| > \varepsilon) < \infty.
\] 
Докажите, что в этом случае $X_1, X_2, \ldots$ сходится к $X$ с вероятностью единица (п.н.). Справедливо также обратное утверждение для последовательности независимых в совокупности с.в. $\{X_k\}_{k=1}^{\infty}$.
\end{problem}

\begin{ordre}
Введите вспомогательное событие $A_n^m$ -- найдется такоe $k > n$, что $|X_k - X| > 1/m$. Покажите, что 
\[
\PR(A^m) = \PR(\mathop{\cap} \limits_{n=1}^{\infty} A_n^m) = 0.
\] 
Пусть событие  $A$ заключается в том, что  
$$
\exists m: \forall\,n \exists k(n): \; |X_{n+k(n)} - X| > 1/m \, \text{(отрицание сходимости п.н.)}.
$$
Убедитесь в равенстве событий  $A$ и $\cup_{m=1}^{\infty} A^m$.
\end{ordre}

\begin{problem}
\label{limsubseq}
Докажите следующее утверждение: если последовательность $\{\xi_n\}$ сходится  по вероятности, то из нее можно извлечь подпоследовательность  $\{\xi_{n_k}\}$, сходящуюся с вероятностью единица.   
\end{problem}

\begin{ordre}
Выберите $n_k$ в соответствии со свойством
\[
\PR(|\xi_{n_{k+1}} - \xi_{n_k}| > 2^{-k}) < 2^{-k}. 
\]
При помощи леммы  Бореля--Кантелли установите, что ряд \[\sum_k |\xi_{n_{k+1}} - \xi_{n_k}| \] сходится на множестве с вероятностной мерой 1. Предельное значение случайной величины в точках расхождения ряда можно определить равным $0$.
Детали доказательства можно найти в книге Ширяева~А.Н. Т.~1 [12].  
\end{ordre}

\begin{problem}\Star
\label{limpnnero}
Докажите, что не существует метрики, определенной на парах случайных величин, сходимость в которой равносильна сходимости почти наверное.
\end{problem}

\begin{ordre}
Используя результат задачи \ref{limsubseq}, придите к противоречию для последовательности, сходящейся по вероятности, но расходящейся с вероятностью 1. См. также \cite{stoianov}.
\end{ordre}

\begin{problem}
Число $\alpha$ из отрезка $[0, 1]$ назовем нормально приближаемым рациональными числами, если найдутся $c,\varepsilon>0$ такие, что 
при любом натуральном $q$ 
\begin{equation*}
\label{BorelKantel}
\min\limits_{p\in {\mathbb Z}} \Bigl|\alpha-\frac{p}{q} \Bigr|\geqslant \frac{c}{q^{2+\varepsilon}} . 
\end{equation*}
Используя лемму Бореля--Кантелли (см. задачу \ref{bor_kant}), докажите, что множество нормально приближаемых чисел на отрезке $[0, 1]$ имеет лебегову меру единица. 

\end{problem}
\begin{ordre}

Зафиксируем $c$, $\varepsilon >0$ и рассмотрим множество
\[A_{q} =\left\{\left. \alpha \in \left[0,\; 1\right]\; \right|\; \mathop{\min }\limits_{p\in {\mathbb Z}} \left|\alpha -\frac{p}{q} \right|<\frac{c}{q^{2+\varepsilon } } \right\}.\] 
Покажите, что $\mu \left(A_{q} \right)\le {2c\mathord{\left/ {\vphantom {2c q^{1+\varepsilon } }} \right. \kern-\nulldelimiterspace} q^{1+\varepsilon } } $. Таким образом, ряд $\sum_{q=1}^{\infty} \mu \left(A_{q} \right) $ сходится. В силу леммы Бореля--Кантелли (см. задачу \ref{bor_kant}) отсюда следует нужное утверждение.

\end{ordre}

\begin{remark}

В связи с полученным результатом будет интересно заметить, что $\forall\,\alpha \in [0, 1]$ существует такая бесконечная последовательность $q_{k} $ и соответствующая ей последовательность $p_{k} $, что
\[\left|\alpha -\frac{p_{k} }{q_{k} } \right|<\frac{1}{\sqrt{5} } \frac{1}{q_{k} ^{2} } .\] 
В теории цепных дробей показывается, что последовательность ${p_{k} \mathord{\left/ {\vphantom {p_{k}  q_{k} }} \right. \kern-\nulldelimiterspace} q_{k} } $  будет подпоследовательностью последовательности подходящих дробей для числа $\alpha $. Заметим также, что константу ${1\mathord{\left/ {\vphantom {1 \sqrt{5} }} \right. \kern-\nulldelimiterspace} \sqrt{5} } $ в неравенстве уменьшить нельзя.

Отметим также, что данная задача возникает в КАМ-теории (см.~Синай~Я.Г. Введение в эргодическую теорию.  М.: Фазис, 1996.  144~с.; Treschev~D., Zubelevich~O. Introduction to the perturbation theory of Hamiltonian systems. Berlin: Springer-Verlag,  2010.  211~p.).
\end{remark}

\begin{problem}
Пусть случайная величина $X_n$ принимает значения 
\begin{enumerate}
\item $(2^n, -2^n)$ с вероятностями $(1/2, 1/2)$;
\item $(-\sqrt{n}$, $\sqrt{n})$ с вероятностями $(1/2, 1/2)$;
\item $(n, 0, -n)$ с вероятностями $(1/4, 1/2, 1/4)$.
\end{enumerate}
Выполняется ли для последовательности независимых случайных величин 
$X_1$, $X_2$, $\ldots$ \textit{закон больших чисел}? 
\end{problem}

\begin{ordre}
 Воспользуйтесь тем, что при $n\to\infty$
\[
S_n\xrightarrow{p}0 \,\Leftrightarrow\, S_n\xrightarrow{d}0 \,\Leftrightarrow\, \forall t \; \varphi_{S_n}(t) \to 1,   
\]
\noindent где $S_n=\frac{X_1+\ldots +X_n}{n}$, $\varphi_{S_n}(t)$ -- характеристическая функция. В общем случае соответствие между сходимостью функций распределения $F_n(x)$ и сходимостью соответствующих $\varphi_n(t)$ такое (см.~также задачу~\ref{levi_kramer} раздела~\ref{zb4}):
\begin{enumerate}
\item Если $F_n(x) \to F(x)$ в точках непрерывности $F(x)$ и $F(x)$ удовлетворяет свойствам функции распределения, то $\forall\,t \; \varphi_n(t) \to \varphi(t)$, где $\varphi(t)$ х.ф. $F(x)$.
\item Если $\forall\,t \; \varphi_n(t) \to \varphi(t)$ и $\varphi(t)$ непрерывна в точке $t = 0$, то $\varphi(t)$ является х.ф. некоторого распределения $F(x)$  и $F_n(x) \to F(x)$ в точках непрерывности $F(x)$.

\end{enumerate}

\end{ordre}

\begin{problem}
При каких значениях $\alpha > 0$ к последовательности независимых случайных величин $\{ X_n\}_{n=1}^{\infty}$, 
таких что ${\mathbb P}\left( X_n=n^{\alpha}\right)={\mathbb P}\left( X_n=\right.$\linebreak $\left.=-n^{\alpha}\right)=1/2$, применим усиленный закон больших чисел? 
\end{problem}

\begin{problem}
Пусть $\{ X_n\}_{n=1}^{\infty}$ -- последовательность случайных величин с дисперсиями $\sigma_i^2$. Докажите, что если все 
корреляционные моменты (корреляции) $R_{ij}$ случайных величин $X_i$ и $X_j$ неположительны и при  
$\frac{1}{n^2}\sum\limits_{i=1}^{n} \sigma_i^2\to 0$, $n\to\infty$, то для последовательности $\{ X_n\}_{n=1}^{\infty}$ выполняется закон больших чисел. 
\end{problem}

\begin{problem}
Пусть $\{ X_n\}_{n=1}^{\infty}$ -- последовательность случайных величин с равномерно ограниченными дисперсиями, причем каждая 
случайная величина $X_n$ зависит только от $X_{n-1}$ и $X_{n+1}$, но не зависит от остальных $X_i$. Докажите выполнение для этой 
последовательности закона больших чисел.
\end{problem}

\begin{problem} Следует ли из УЗБЧ для схемы испытаний Бернулли, что с 
вероятностью 1 для любых, сколь угодно маленьких $\varepsilon >0$ и $\delta 
>0$, можно подобрать такое $n\left( {\varepsilon ,\delta } \right)\in {\mathbb N}$, что с вероятностью, не меньшей $1-\delta $ для всех $n\ge n\left( 
{\varepsilon ,\delta } \right)$, частота выпадения герба $\nu _n $ отличается 
от вероятности выпадения герба $p$ не больше, чем на $\varepsilon $? Следует 
ли это из обычного ЗБЧ?
\end{problem}

\begin{problem} Приведите пример последовательности независимых с.в. 
$\left\{ {\xi _n } \right\}_{n\in {\mathbb N}} $ таких, что предел $\mathop 
{\lim }\limits_{n\to \infty } \frac{\xi _1 +...+\xi _n }{n}$ существует по 
вероятности, но не существует с вероятностью 1 (cм. \cite{stoianov}).
\end{problem}

\begin{problem} Докажите \textit{локальную теорему Муавра--Лапласа} в общем случае:
\[
C_n^k p^k\left( {1-p} \right)^{n-k}\simeq \left( {2\pi np\left( {1-p} \right)} 
\right)^{-1 \mathord{\left/ {\vphantom {1 2}} \right. 
\kern-\nulldelimiterspace} 2}\exp \left( {-\frac{1}{2}\frac{\left( {k-np} 
\right)^2}{np\left( {1-p} \right)}} \right)
\]
равномерно по таким $k$, что $\left| {k-np} \right|\le \varepsilon \left( n 
\right)$, $\varepsilon \left( n \right)=o\left( {n^{2 \mathord{\left/ 
{\vphantom {2 3}} \right. \kern-\nulldelimiterspace} 3}} \right)$, $n\to 
\infty $, $p$ -- фиксировано.
\end{problem}

\begin{ordre}
См. книгу  Гнеденко~Б.В. \cite{2}. 
\end{ordre}

\begin{problem}
Докажите, что для с.в. $X_i \in \text{Be}(p)$, $p > 0$ выполнена локальная теорема Муавра--Лапласа ($n \gg 1$):
\[
\PR\left(\sum_{i = 1}^n X_i = k \right) \simeq
\left( {2\pi np\left( {1-p} \right)} 
\right)^{-1 \mathord{\left/ {\vphantom {1 2}} \right. 
\kern-\nulldelimiterspace} 2}\exp \left( {-\frac{1}{2}\frac{\left( {k-np} 
\right)^2}{np\left( {1-p} \right)}} \right).
\]
При этом, если $X_1,\ldots,X_n \in \text{Be}(1/n)$, то $\PR\left(\sum_{i = 1}^n X_i = k \right) \simeq e^{-1}/k!$, т.е. $\sum_{i = 1}^n X_i \mathop{\longrightarrow}\limits^{d} \text{Po}(1)$ при $n \to \infty$, где $\text{Po}\left(\lambda\right)$ -- распределение Пуассона (Poisson) с параметром $\lambda$. Будет ли иметь место сходимость по вероятности  $\sum_{i = 1}^n X_i \mathop{\longrightarrow}\limits^{p} \text{Po}(1)$ при $n \to \infty$?

\end{problem}

\begin{problem}
Рассмотрим прибор, называемый ``доска Гальтона'' (см.~рис.~\ref{Fig:gam_scheme.JPG}).
\imgh{30mm}{gam_scheme.JPG}{\smallДоска Гальтона}    
Принцип работы этого устройства таков: металлические шарики поступают в самый верхний канал. Наткнувшись на первое острие, они «выбирают» путь направо или налево. Затем происходит второй такой выбор и т.д. При хорошей подгонке деталей выбор оказывается случайным. Как видно, попадание шариков в нижние лунки не равновероятно. В этом случае мы имеем дело с нормальным распределением. Объясните почему?
\end{problem}

\begin{problem}
Книга объемом $500$ страниц содержит $50$ опечаток. Оцените вероятность того, что на случайно выбранной странице 
имеется не менее трех опечаток, используя нормальное и пуассоновское приближения, сравнить результаты. 
\end{problem}

\begin{problem}
В тесто для выпечки булок с изюмом замешано $N$ изюмин. Всего из данного теста выпечено $K$ булок. Оцените вероятность того, 
что в случайно выбранной булке число изюмин находится в пределах от $a$ до $b$. 
\end{problem}

\begin{problem}
В поселке $N$ жителей, каждый из которых в среднем $n$ раз в месяц ездит в город, выбирая дни поездки независимо от остальных. 
Поезд из поселка в город идет один раз в сутки. Какова должна быть вместимость поезда для того, чтобы он переполнился с вероятностью, 
не превышающей заданного числа $\beta$? 
\end{problem}

\begin{problem}\label{cauchy_gen}
Случайные величины $X_1,\ldots,X_n$ независимы и имеют распределение Коши, т.е. 
\[
f(x) = \frac{d}{\pi(d^2 + x^2)}, 
\quad
x \in \mathbb{R}.
\]
Докажите равенство распределений $(X_1+\ldots+X_n)/n$ и  $X_1$. Противоречит ли это ЗБЧ или ЦПТ?
\end{problem}


\section{Задачи повышенной сложности}
\label{hard}

\begin{problem}
Докажите, что при $n\to\infty$ ($p \in \left[1,\infty\right]$)
$$
X_n\xrightarrow{L_2} X \,\Rightarrow\, X_n\xrightarrow{L_1}X \, \Rightarrow\, X_n\xrightarrow{p}X 
\, \Leftarrow\, X_n\xrightarrow{\text{ п.н. }}X . 
$$
Также покажите, что 
$$
X_n\xrightarrow{p}X \, \Rightarrow\, X_n\xrightarrow{d}X . 
$$
С помощью контрпримеров покажите, что никакие другие стрелки импликации в эту схему в общем случае добавить нельзя. 
При каких дополнительных условиях можно утверждать, что 
$$
X_n\xrightarrow{\text{ п.н. }}X  \, \Rightarrow\, X_n\xrightarrow{L_1}X ?
$$
Кроме того, покажите, что при $n \rightarrow \infty$
$$
X_n\xrightarrow{p} X \; \Leftrightarrow\, \rho_P(X_n,X)={\mathbb E}\Bigl( \frac{|X_n-X|}{1+|X_n-X|}\Bigr)
\rightarrow 0.
$$
Это означает, что сходимость по вероятности метризуема.
\end{problem}

\begin{remark} Из книги А.Н.~Ширяева Т.~1  [\ref{chiraiev}] (см. также также главу 6 из \cite{Gupta} и \cite{stoianov}) имеем 
$$
X_n\xrightarrow{\text{ п.н. }}X  \, \Rightarrow\, X_n\xrightarrow{L_1}X , 
$$
т.е. возможен предельный переход под знаком математического ожидания, если семейство с.в. $\{ X_n\}$ является равномерно интегрируемым: 
$$
\sup\limits_n {\mathbb E}\bigl[ |X_n|\cdot {\mathbb I}_{\{ |X_n|>c\}} \bigr] \rightarrow 0, \quad c \to +\infty. 
$$

Отметим также, что сходимость по распределению, так же как и по вероятности, метризуема в отличие от сходимости п.н. (см. задачу \ref{limpnnero} из  раздела \ref{standart}). 

\end{remark}

\begin{problem}\Star(Лоэв--Колмогоров) Пусть $\left\{ {\xi _n } \right\}_{n\in 
{\mathbb N}} $ -- последовательность  независимых с.в. таких, что 
$\sum\limits_{n=1}^\infty {\frac{\Exp\left| {\xi _n } \right|^{\alpha _n 
}}{n^{\alpha _n }}} <\infty $, где $0<\alpha _n \le 2$, причем $\Exp\xi _n =0$ 
в случае $1\le \alpha _n \le 2$. Покажите, что тогда 
$\frac{1}{n}\sum\limits_{k=1}^n {\xi_k } \buildrel \text{п.н.} \over 
\longrightarrow 0$.
\end{problem}

\begin{problem}\Star(Теорема Колмогорова о трех рядах.) 
Пусть $\left\{ {\xi _n } 
\right\}_{n\in {\rm N}} $ -- последовательность независимых с.в. 
Покажите, что для сходимости ряда $\sum\limits_{n=1}^\infty {\xi _n } $ с 
вероятностью 1 необходимо, чтобы для любого $c>0$ сходились ряды 
$\sum\limits_{n=1}^\infty {\Exp\xi _n^c } $, $\sum\limits_{n=1}^\infty {\Var\xi 
_n^c } $, $\sum\limits_{n=1}^\infty {\PR\left( {\left| {\xi _n } \right|\ge c} 
\right)} $, где $\xi _n^c =\xi _n I\left( {\xi _n \le c} \right)$, и 
достаточно, чтобы эти ряды сходились при некотором $c>0$.
\end{problem}

\begin{problem}(Интеграл Лебега--Стилтьеса.)
Для неотрицательной с.в. $\xi$ определим ее математическое  ожидание $\Exp \xi$ посредством интеграла Лебега $\int \xi(\omega) d\PR(\omega)$, определяемого как
$$
\lim \limits_{n \to \infty} 
  \sum \limits_{k = 1}^{n 2^n} \frac{k-1}{2^n} \PR\left(  \frac{k-1}{2^n} \leq \xi <  \frac{k}{2^n}  \right) + n \PR(\xi \geq n) 
$$ 

\noindent Пусть $\xi: \Omega \mapsto \mathbb{R}$ -- произвольная с.в. и 
\[ \xi^{+} = \max \{\xi, 0 \}, \quad \xi^{-} = - \min \{\xi, 0\}.\] 
Тогда при условии $\Exp \xi^{+} < \infty$ или $\Exp \xi^{-} < \infty$ оператор $\Exp$ определяется как 
$\Exp \xi = \Exp \xi^{+} - \Exp \xi^{-}$. 

Докажите следующие утверждения:

\begin{enumerate}
\item Для с.в. $\xi \geq 0$ существует последовательность  \textit{простых} с.в. (принимают конечное число значений)  таких, что $\xi_{n+1} \geq \xi_n$ и $\xi_n \mathop{\longrightarrow} \limits^{\text{п.н.}} \xi$. 

Следствие: $\forall\,\xi \geq 0 $ определено $ \Exp \xi$.
\item Если  $\xi_{n+1} \geq \xi_n \geq 0$ и $\xi_n \mathop{\longrightarrow} \limits^{\text{п.н.}} \xi$, то $\Exp \xi_n \to \Exp \xi$. 
\begin{enumerate}
\item[Следствие 1:] Пусть $\xi_n \geq 0$, тогда 
\[ \Exp\left(\sum \limits_{n=1}^{\infty} \xi_n \right) = \sum \limits_{n=1}^{\infty} \Exp \xi_n. \]
\item[Следствие 2:] Пусть $\Exp \vert \xi \vert < \infty$ и $\PR(A) \to 0$, тогда $\Exp (\vert \xi \vert I_A) \to 0$. 
\end{enumerate}
\item Пусть $\xi_n \geq 0$, тогда $\Exp (\lim \limits_{n\to\infty} \inf \xi_n) \leq \lim \limits_{n\to\infty} \inf \Exp \xi_n$. 
\item Пусть $\vert \xi_n \vert \leq Y$ и $\Exp Y < \infty$, тогда $\Exp (\lim \limits_{n\to\infty} \inf \xi_n) \leq \lim \limits_{n\to\infty} \inf \Exp \xi_n  \leq$\linebreak $\leq \lim \limits_{n\to\infty} \sup \Exp \xi_n  \leq \Exp (\lim \limits_{n\to\infty} \sup \xi_n)$.
\begin{enumerate}
\item[Следствие 1:] Пусть $\xi_n \geq 0$, $\vert \xi_n \vert \leq Y$, $\Exp Y < \infty$ и $\xi_n \mathop{\longrightarrow} \limits^{\text{п.н.}} \xi$, тогда $\Exp \vert \xi \vert < \infty$, $\Exp \xi_n \to \Exp \xi$ и $\Exp \vert \xi_n - \xi \vert \to 0$.
\item[Следствие 2:] Пусть $\vert \xi_n \vert \leq Y$, $\Exp Y^p < \infty$, $p > 0$ и $\xi_n \mathop{\longrightarrow} \limits^{\text{п.н.}} \xi$, тогда $\Exp \vert \xi \vert^p < \infty$ и $\Exp \vert \xi_n - \xi \vert^p \to 0$. 
\end{enumerate}
\end{enumerate} 

\end{problem}

\begin{ordre}
См. главу~4 в книге Константинов~Р.В. Лекции по функциональному анализу. Долгопрудный: МФТИ, 2007.  156 с. 
\end{ordre}

\begin{remark}
Положим $g : \mathbb{R} \to \mathbb{R}$ -- борелевская функция ($\forall\,A \in$\linebreak $\in \mathcal{B}(\mathbb{R}): \;  g^{-1}(A) \in \mathcal{B}(\mathbb{R}) $). Если существует один из интегралов 
$$
\int \limits_{A} g(x) dF_{\xi}(x) \quad \text{или} \quad \int \limits_{\xi^{-1}(A)} g(\xi(\omega)) d\PR(\omega),
$$
то существует и другой, и они совпадают. Интеграл  $\int \limits_{\mathbb{R}} g(x) dF_{\xi}(x)$ называется интегралом \textit{Лебега--Стилтьеса}.\\
\end{remark}


\begin{problem}
Симметричную монету независимо бросили $n$ раз. Результат бросания записали в виде последовательности нулей и единиц. 
\begin{enumerate}
\item Покажите, что с вероятностью, стремящейся к единице при $n\to \infty $, длина максимальной подпоследовательности  из подряд идущих единиц $l_n$ лежит в промежутке $(\log_2 \sqrt{n} ,\; \log_2 n^{2} )$. 

\item\Star Покажите, что серия из гербов длины $\log_{2} n$ наблюдается с вероятностью, стремящейся к единице при $n\to\infty$.

\item\DStar Верно ли, что  $\frac{ l_n}{ \log_2 n } \overset{\text{п.н.}}{\longrightarrow} 1$?
\end{enumerate}
\end{problem}

 \begin{ordre}
 В пункте а) для нахождения нижней оценки разобьем всю последовательность бросаний на участки длиной $\log_2 \sqrt{n}$. Оцените вероятность того, что хотя бы один из участков состоит полностью из единиц.    
 \end{ordre}

\begin{remark}
 См.  Erdos~P.,  Renyi~A. On a new law of large numbers // Journal Analyse Mathematique. 1970. V. 23. P.~103--111.
 
\end{remark}

\begin{problem}(Закон 0 и 1.)
\label{1and0law}
Пусть $(\Omega,\mathcal{F},{\mathbb P})$ --- вероятностное пространство, $\xi_1,\xi_2,\ldots$ -- некоторая последовательность независимых с.в. 
Обозначим $\mathcal{F}_n^{\infty}=\sigma(\xi_{n},\xi_{n+1},\ldots)$ -- $\sigma$-алгебру, порожденную с.в. $\xi_{n},\xi_{n+1},\ldots$ и пусть 
$$
{\mathcal X}=\bigcap\limits_{n=1}^{\infty} \mathcal{F}_{n}^{\infty} . 
$$
Поскольку пересечение $\sigma$-алгебр есть снова $\sigma$-алгебра, то ${\mathcal X}$  есть \linebreak$\sigma$-алгебра. Эту $\sigma$-алгебру 
будем называть ``хвостовой'' или ``остаточной'', в связи с тем, что всякое событие $A\in{\mathcal X}$ не зависит от значений с.в. 
$\xi_1,\xi_2,\ldots,\xi_n$ при любом конечном $n$, а определяется лишь ``поведением бесконечно далеких значений последовательности 
$\xi_1,\xi_2,\ldots$''. \\
\noindent С помощью задачи $\ref{SigmaAlgebra}$ раздела \ref{standart} докажите справедливость следующего утверждения (``закон 0 и 1''): 
\par Пусть $\xi_1,\xi_2,\ldots$ -- последовательность независимых в совокупности с.в. и $A\in{\mathcal X}$ 
(событие $A$ принадлежит ``хвостовой'' $\sigma$-алгебре). Тогда ${\mathbb P}(A)$ может принимать лишь два значения, $0$ или $1$. 
\end{problem}

\begin{ordre}
Идея доказательства состоит в том, чтобы показать, что каждое <<хвостовое>> событие $A$ не зависит от самого себя и, значит, 
${\mathbb P}(A\cap A)={\mathbb P}(A)\cdot {\mathbb P}(A)$, т.е. ${\mathbb P}(A)={\mathbb P}^2(A)$, откуда 
${\mathbb P}(A)=0$ или $1$. 

Полученный результат, в частности, означает, что $\sum  \xi_k$ сходится или расходится с вероятностью 1.  

Существует также более общий результат -- закон Хьюитта и Сэвиджа (см.~Ширяев~А.Н. Т.~2 \cite{21}), где под ${\mathcal X}$ подразумевается $\sigma$-алгебра перестановочных событий:
\[
\mathcal{X} = \{ \xi^{-1}(B), \; B \in \mathcal{B}(\mathbb{R}^{\infty}): \;
\PR (\xi^{-1}(B) \triangle (\pi\xi)^{-1}(B)) = 0, \; \forall \pi \}, 
\] 
a $\pi$ -- перестановка конечного числа элементов.
\end{ordre}

\begin{problem}
Сто паровозов выехали из города по однополосной линии, каждый с постоянной скоростью. Когда движение установилось, то из-за того, что быстрые догнали идущих впереди более медленных, образовались караваны (группы, движущиеся со скоростью лидера). Найдите математическое ожидание и дисперсию числа караванов. Скорости различных паровозов независимы и одинаково распределены, а функция распределения скорости непрерывна.
\end{problem}

\begin{problem}
Согласно законам о трудоустройстве в городе \textit{М}, наниматели обязаны предоставить всем рабочим выходной, если хотя бы у одного из них день рождения, и принимать на службу рабочих независимо от их дня рождения. За исключением этих выходных, рабочие трудятся весь год из 365 дней. Условимся что день рождения рабочего выбирается равновероятно из 365 дней. Работодатель максимизировал среднее число трудовых человеко-дней в году (то есть произведение числа рабочих на число трудовых дней). Сколько рабочих трудятся на фабрике в городе~\textit{М}?

\end{problem}

\begin{problem}
В каждую $i$-ю единицу времени живая клетка получает случайную дозу облучения $X_i$, причем $\{ X_i\}_{i=1}^{n}$ имеют 
одинаковую функцию распределения $F_X(x)$ и независимы в совокупности $\forall\,n \in \mathbb{N}$. Получив интегральную дозу облучения, 
равную $\nu$ ($\nu \gg \Exp X_1$), клетка погибает. Обозначим через $T$ -- с.в., равную времени жизни клетки.
Оцените среднее время жизни клетки ${\mathbb E}T$. 
\end{problem}

\begin{ordre}

Докажите \textit{тождество Вальда}: 
$$
{\mathbb E}S_T={\mathbb E}X_1\cdot {\mathbb E}T, \; \text{где} \; S_T = \sum_{i=1}^{T}X_i,    
$$

\noindent введя вспомогательную с.в.

$$
Y_j=\begin{cases}
1, &\text{ если }\quad X_1+\ldots +X_{j-1}=S_{j-1}<\nu, \\
0 &\text{ в остальных случаях} 
\end{cases}
$$
 
\noindent и записав $S_{T}$ в следующем виде: $S_T = \sum_{i=1}^{\infty} Y_i X_i$. 

\end{ordre}

\begin{problem}\Star(Теорема Дуба.)
\label{sec:doob}
Пусть $Y_0,\dots,Y_n$ -- последовательность с.в., являющаяся мартингалом (см.~замечание).
Показать, что для  \textit{момента остановки}
$\tau = \inf \{k\leq n: Y_k\geq \lambda \}$  ($\tau = n$, если $\max_{k\leq n}Y_k <\lambda$), верно
\begin{equation*}
\mathbb{E}{Y_{\tau}} = \mathbb{E}Y_n.
\end{equation*}
\end{problem}
\begin{remark}
Определение мартингала (см. Ширяев~А.Н. Т.~2 [\ref{chiraiev}]). Пусть задано вероятностное пространство $(\Omega,\mathcal{F},\mathbb{P})$  с семейством ($\mathcal{F}_n$) $\sigma$-алгебр $\mathcal{F}_n$, $n\geq 0$ таких, что $\mathcal{F}_0\subseteq\dots\subseteq\mathcal{F}_n\subseteq \mathcal{F}$. 
Пусть $Y_0,\dots,Y_n$ -- последовательность с.в., заданных на $(\Omega,\mathcal{F},\mathbb{P})$, для каждого $n\geq 0 $ величины $Y_n$ являются $\mathcal{F}_n$-измеримыми. Мартингалом называется  такая  последовательность $Y = (Y_n,\mathcal{F}_n)$, что 
\begin{equation*}
\mathbb{E}|Y_n|\leq\infty,
\end{equation*}
\begin{equation*}
\mathbb{E}(Y_{n+1}|\mathcal{F}_n)=Y_n \quad \text{почти наверное}.
\end{equation*}

Для доказательства теоремы Дуба необходимо показать, что 
 для любого $A\in \mathcal{F}_n$ 
\begin{equation*}
\int_{A} Y_{\tau} d\mathbb{P} = \int_{A}Y_{n}d\mathbb{P}.
\end{equation*}
\end{remark}

\begin{problem}(Мартингалы и теорема о баллотировке.)
Пусть $S_n =\xi _1 +....+\xi _n $, где с.в. $\left\{ {\xi _k } \right\}_{k\in {\mathbb N}} $ -- 
независимы и одинаково распределены \[\xi _k =\left\{ {\begin{array}{l}
 1,\quad p=1 \mathord{\left/ {\vphantom {1 2}} \right. 
\kern-\nulldelimiterspace} 2, \\ 
 -1,\;\,p=1 \mathord{\left/ {\vphantom {1 2}} \right. 
\kern-\nulldelimiterspace} 2. \\ 
 \end{array}} \right.\] Покажите, что тогда

$\PR\left( {S_1 >0,...,S_n >0\left| {S_n =a-b} \right.} 
\right)=\frac{a-b}{a+b},$ где $a>b$ и $a+b=n$.
\end{problem}

\begin{remark}
Проинтерпретируем этот результат: $\xi_k =1$ будем 
интерпретировать как голос, поданный на выборах за кандидата $A$, 
$\xi _k =-1$ -- за кандидата $B$. Тогда $S_n $ есть разность числа 
голосов, поданных за кандидатов $A$ и $B$, если в голосовании 
приняло участие $n$ избирателей, а $\PR\left( {S_1 >0,...,S_n >0\left| {S_n 
=a-b} \right.} \right)$ есть вероятность того, что кандидат $A$ все 
время был впереди кандидата $B$, при условии, что $A$ в общей 
сложности собрал $a$ голосов, а $B$ собрал $b$ голосов и $a>b$, 
$a+b=n$. См.~также Ширяев~А.Н. T.~2 [\ref{chiraiev}]. 
\end{remark}

\begin{problem}(Закон арксинуса, см. Ширяев~А.Н. Т.~2 [\ref{chiraiev}].)
В условиях замечания к предыдущей задаче 
найдите при 20 бросаниях, с какой вероятностью один из игроков
\begin{enumerate}
\item[а)] никогда не будет впереди,
\item[б)] будет впереди не более одного раза.
\end{enumerate}
\end{problem}

\begin{remark} 
Воспользуйтесь \textit{законом арксинуса}

$\PR\left( {k_n <xn} \right)\simeq 
\frac{2}{\pi }\arcsin \sqrt x $, где с.в. $k_n =\left| {\left\{ 
{k=1,...,n:\;\;S_k \ge 0} \right\}} \right|$.
\end{remark}

\begin{problem}\Star(Случайные блуждания.)
Заблудившийся грибник оказался в центре леса, имеющего форму круга с радиусом 5 км. Грибник движется по следующим правилам: пройдя 100 метров в случайно выбранном направлении север-юг-запад-восток (на это у грибника уходит 2 минуты), он решает снова случайно выбрать направление движения и~т.д. Оцените математическое ожидание времени, через которое грибник выйдет из леса. Как изменится ответ, если выбор направления движения грибника осуществляется равновероятно на $[0,2\pi)$? 
\end{problem}

\begin{problem} 
Покажите, что последовательность дискретных с.в. $\left\{ 
{\xi _n } \right\}_{n\in {\mathbb N}} $, принимающих значения в ${\mathbb N}$, 
сходится по распределению к дискретной с.в. $\xi $ тогда и только тогда, 
когда при $n\to \infty $ для любого $k\in {\mathbb N}$ $\PR\left( {\xi _n =k} 
\right)\to \PR\left( {\xi =k} \right)$.
\end{problem}

\begin{problem}
На множестве действительных положительных
чисел задано некоторое вероятностное распределение $P$.
Из одной и той же точки плоскости начинают прыгать две лягушки. Направление прыжка лягушки выбирают случайно и равномерно (на окружности), а длину -- согласно распределению $P$. Направление и длина каждого прыжка независимы,
также как и прыжки лягушек.
Первая лягушка сделала $n$ прыжков, а вторая $m$
(где $m,n>0$ и $m+n>2$).
Докажите, что вероятность того, что первая лягушка удалилась на большее расстояние
от исходной точки, чем вторая, равна $n/(n+m)$.

\begin{ordre}
Для независимых случайных векторов $A$, $B$ и $C$, имеющих произвольное непрерывное распределение на плоскости, справедливо тождество
$$
\mathbb P \biggl( \| A \|  > \|B + C\| \biggr) + \mathbb P \biggl( \| B \|  > \|A + C\| \biggr) + \mathbb P \biggl( \| C \|  > \|B + A\| \biggr) = 1.
$$
\end{ordre}

\begin{remark}
Bernardi~O. A short proof of Rayleigh’s Theorem with extensions // The American Mathematical Monthly. 2013. V.\;120(4). P.~362--364.
\end{remark}
\end{problem}

\begin{problem}(Сходимость по моментам.) 
\label{moments_conv}
Пусть $\left\{ {F_n \left( x \right)} \right\}_{n\in {\mathbb N}}$ -- последовательность функций  распределения, 
имеющих все моменты: 
\[M_{n,k} =\int {x^kdF_n \left( x \right)} <\infty. \] 
Пусть для всех $k\in {\mathbb N}$ имеют место следующие сходимости: $M_{n,k} \to$\linebreak $\to 
M_k \ne \pm \infty $. Тогда (по теореме Хелли) существуют такая подпоследовательность $\left\{ 
{F_{n_m } \left( x \right)} \right\}_{m\in {\mathbb N}} $ и функция 
распределения $F\left( x \right)$\linebreak (с моментами $\left\{ {M_k } 
\right\}_{k\in {\mathbb N}} )$, что $F_{n_m } \left( x \right)\to F\left( x 
\right)$ в точках непрерывности~$F\left( x \right)$.\\

\begin{enumerate}
\item Покажите, что если моменты $\left\{ {M_k } \right\}_{k\in {\mathbb N}} $ 
однозначно определяют функцию $F\left( x \right)$ (достаточным условием для 
этого будет существование такого $c >0$, что 
$ M_k^{1/k} / k \leq c $, $k\in {\mathbb N})$, то в качестве подпоследовательности можно брать саму последовательность \linebreak (теорема Шохата--Тамаркина).

\item  Проинтерпретируйте полученный результат с точки зрения сходимости по 
распределению соответствующих с.в.
\end{enumerate}
\end{problem}
\begin{remark}
См., например, Сачков В.Н. Вероятностные методы в комбинаторном анализе.  М.: Наука, 1978.  288 с.
\end{remark}

\begin{problem}

Покажите, что все моменты распределений с плотностями
$$p_{\lambda } \left(x\right)=\frac{1}{24} e^{-x^{{1\mathord{\left/ {\vphantom {1 4}} \right. \kern-\nulldelimiterspace} 4} } } \left(1-\lambda \sin x^{{1\mathord{\left/ {\vphantom {1 4}} \right. \kern-\nulldelimiterspace} 4} } \right),\; x\ge 0$$
при любом значении параметра $\lambda \in \left[0,1\right]$ совпадают.

\begin{remark}

Достаточное \textit{условие Карлемана} того, что моменты однозначно определяют распределение с.в. $\xi$ на всей прямой, имеет вид (см. \cite{stoianov}):
$$
\sum _{n=1}^{\infty }\left(\Exp\left[|\xi|^{2n}\right] \right)^{{-1\mathord{\left/ {\vphantom {-1 \left(2n\right)}} \right. \kern-\nulldelimiterspace} \left(2n\right)} } =\infty.
$$
Набор (всех) моментов комплексной с.в. однозначно определяет ее характеристическую функцию (а, следовательно, и ее распределение вероятностей) тогда и только тогда, когда
$$\varlimsup_{n\to\infty}n^{-1}\left({\Exp\left[|\xi|^{n}\right]}\right)^{1/n}<\infty.$$
Причем характеристическая функция с.в. будет аналитической в этом случае \cite{28}.
\end{remark}

\end{problem} 

\begin{problem}
Приведите примеры таких с.в. $X$, $Y$ и $Z$, что вероятностные распределения сумм $X+Y$ и $X+Z$ совпадают, но распределения с.в. различны.
\end{problem}
\begin{ordre}
Удобнее сначала подобрать соответствующие характеристические функции.
\end{ordre}
\begin{remark}
См. книгу Секея~Г. [\ref{sekei}]. Эту же книгу можно рекомендовать и по двум следующим задачам.
\end{remark}

\begin{problem}
Приведите примеры независимых одинаково распределенных с.в. $X$ и $Y$, для которых распределение суммы  $X+Y$ не однозначно определяет распределение $X$ и $Y$.
\end{problem}

\begin{remark}
Согласно теоремам Крамера и Райкова распределение суммы $X+Y$  однозначно определяет распределение слагаемых, коль скоро  $X+Y$ имеет нормальное или Пуассоновское распределение. См. также замечание к задаче \ref{sec:infdiv} раздела \ref{zb4}.
\end{remark}

\begin{problem}
Пусть 
  $X$ стохастически меньше, чем $Y$, т. е. 
 $$\forall\,t \; F_X(t) \leq F_Y(t), \exists t_0: \; F_X(t_0) < F_Y(t_0).$$
 Парадоксально, но может так случиться, что 
  $\PR(X > Y) \geq 0.99$. Приведите пример таких с.в. $X$ и $Y$.
\end{problem}

\begin{problem}(Парадокс транзитивности \cite{2013}, \cite{book12}.)
Будем говорить, что с.в. $X$ больше по вероятности с.в. $Y$, если ${\mathbb P}(X>Y)>{\mathbb P}(X\le Y)$. 
Пусть известно, что для с.в. $X$, $Y$, $Z$ выполнена следующая цепочка равенств: 
$$
{\mathbb P}(X>Y)={\mathbb P}(Y>Z)=\alpha>\frac{1}{2} . 
$$
Верно ли, что $X$ больше по вероятности $Z$ и почему? 
\end{problem}

\begin{problem}

Требуется определить, начиная с какого этажа брошенный с балкона 100-этажного здания стеклянный шар разбивается. В наличии имеется два таких шара. Предложите метод нахождения граничного этажа, минимизирующий математическое ожидание числа бросков. Рассмотреть случай большего числа шаров.  

\end{problem}

\begin{problem}\Star
На подоконнике лежит $N$ помидоров. Вечером $i$-го дня ($1 \leqslant i \leqslant$\linebreak $\leqslant N$) портится один помидор. Каждое утро человек съедает один (случайно выбранный) свежий помидор из оставшихся. Таким образом, каждый помидор либо испортился, либо был съеден.

\begin{enumerate}
\item Получите рекуррентную формулу для математического ожидания количества съеденных помидоров от числа $N$.
\item Найдите асимптотическую оценку количества съеденных помидоров при $N \rightarrow \infty$.
\end{enumerate}

\end{problem}

\begin{problem}
\begin{enumerate}
\item Имеется монетка (несимметричная). Несимметричность монетки заключается в том, что либо орел выпадает в два раза чаще решки, 
либо наоборот (априорно, до проведения опытов, оба варианта считаются равновероятными). Монетку бросили $10$ раз. Орел выпал $7$ раз. 
Определите апостериорную вероятность того, что орел выпадает в два раза чаще решки (апостериорная вероятность считается с учетом 
проведенных опытов; иначе говоря, это просто условная вероятность). 

\item Определите апостериорную вероятность того, что орел выпадает не менее чем в два раза чаще решки. Если несимметричность 
монетки заключается в том, что либо орел выпадает не менее чем в два раза чаще решки, либо наоборот (априорно оба варианта считаются 
равновероятными). 
\end{enumerate}
\end{problem}

\begin{problem}(Распределение Гумбеля или двойное экспоненциальное распределение.)
\label{gumbel}
Пусть $\xi _{1} ,\ldots,\xi _{n} $ -- независимые одинаково распределенные с.в., и существуют такие константы $\alpha, T>0$, что
\[\mathop{\lim }\limits_{y\to \infty } e^{{y\mathord{\left/ {\vphantom {y T}} \right. \kern-\nulldelimiterspace} T} } \left[1-\PR\left(\xi _{1} <y\right)\right]=\alpha. \] 
Покажите, что при $n\to \infty$
\[
\max \left\{\xi _{1},\ldots,\xi _{n} \right\}-T\ln \left(\alpha n\right)\xrightarrow[ ]{d} \zeta,  
\]
\noindent где $\PR\left(\zeta < t \right)=\exp \left\{-e^{-{t \mathord{\left/ {\vphantom {t /  T}} \right. \kern-\nulldelimiterspace} T} } \right\}.$ Исследуйте поведение $\PR\left(\zeta \ge t \right)$ при $t \to \infty$.

\end{problem}

\begin{problem}(Logit-распределение или распределение Гиббса.)
\label{gibbs}
Пусть известно, что $\xi_{1} ,\ldots,\xi_{n}$~--- независимые одинаково распределенные по \textit{закону Гумбеля} с.в. (см. предыдущую задачу). Пусть  $X_{k} =C_{k} +\xi _{k} ,$ $k=1,\ldots,n$, где $C_{k} $ -- некоторые константы. Положим $\kappa=\arg \max\limits_{k=1,\dots,n} \left\{ X _{k} \right\}$. Покажите, что с.в. $\kappa$ имеет распределение:
\[\PR\left(\kappa=k\right)=\frac{e^{{C_{k} \mathord{\left/ {\vphantom {C_{k}  T}} \right. \kern-\nulldelimiterspace} T} } }{\sum _{l=1}^{n}e^{{C_{l} \mathord{\left/ {\vphantom {C_{l}  T}} \right. \kern-\nulldelimiterspace} T} }  } , \quad k=1,\ldots,n.\] 
\end{problem}
\begin{remark}
См. монографию Andersen~S.P., de Palma~A., Thisse~J.-F.  Discrete choice theory of product differentiation.   Cambridge: MIT Press,  1992.  448 p.; см.~также \cite{222}. Данная задача и предыдущая играют важную роль при поиске стохастических равновесий в транспортных сетях arXiv:1701.02473 (см. таже задачу 15 раздела 6).
\end{remark}

\begin{problem}(Рекорды \cite{4}.)
\label{records}
Пусть $X_1 ,X_2 ,\ldots $ -- независимые 
с.в. с одной и той же плотностью распределения вероятностей 
$p(x)$. Будем говорить, что наблюдается рекордное значение в момент времени 
$n>1$, если $X_n >\max \left \{X_1,\ldots,X_{n-1}  \right\}$. Докажите 
следующие утверждения:

\begin{enumerate}
\item вероятность того, что рекорд зафиксирован в момент времени $n$, 
равна $1/n$;

\item математическое ожидание числа рекордов до момента времени $n$ 
равно 
\[
\sum\limits_{1<k\le n} {\frac{1}{k}} \approx \ln n;
\]

\item пусть $Y_n $ -- с.в., принимающая значение $1$, если 
в момент времени $n$ зафиксирован рекорд, и значение $0$ -- в противном случае. 
Тогда с.в. $Y_1 ,Y_2 ,\ldots$ независимы в совокупности;

\item дисперсия числа рекордов до момента времени $n$ равна
\[
\sum\limits_{1<k\le n} {\frac{k-1}{k^2}} \approx \ln n;
\]

\item если $T$ -- момент появления первого рекорда после момента времени $1$, то $\Exp T= \infty$.
\end{enumerate}
\end{problem}

\begin{remark}
Приведем для справки следующую теорему. Пусть $\eta_0,\eta_1,\dots$ -- последовательность независимых с.в. с одной и той же непрерывной функцией распределения. Для каждого $n\in \mathbb{Z}_{+}$ по с.в. $\eta_0,\eta_1,\dots,\eta_n$ построим вариационный ряд 
$$\eta_{0,n}\leq \eta_{1,n}\leq\dots\leq\eta_{n,n}.$$
Рекордные моменты $\{\nu(n),n\in\mathbb{Z}_{+}\}$ определяются следующим образом: $\nu(0) = 0$ и
$$\nu(n+1)=\min\{j>\nu(n): \eta_j>\eta_{j-1,j-1}\},\quad n\in \mathbb{Z}_{+}.$$
Верен следующий результат (см. Якымив А.Л. Вероятностные приложения тауберовых теорем.  М.: Физматлит; Наука, 2005.  256 c.):
$$\mathbb{P}\left(\nu(n)>t\right)\equiv \frac{t^{-1}\ln^{n-1}(t)}{(n-1)!}.$$

\medskip

Для доказательства этого результата существенно используются тауберовы теоремы.
Тауберовыми теоремами называют теоремы, которые позволяют из асимптотических свойств п.ф. и преобразований Лапласа функций и последовательностей  (а также других интегральных преобразований) вывести асимптотики этих функций и последовательностей (то есть эти теоремы обратные к абелевым). 
\end{remark}

\begin{problem}
Зрелый индивидуум производит потомков согласно производящей функции (п.ф.) $h(s)$.
Предположим, что все начинается с $k$ индивидуумов, достигших зрелости. Каждый потомок достигает зрелости с вероятностью $p$. Докажите, что п.ф. числа зрелых индивидуумов в следующем поколении равна $H_k(s) = [h(ps + 1 - p)]^{k}$.   
\end{problem}

\begin{problem}(Распределение Коши.)\label{cauchy_rad_emitt}
Радиоактивный источник испускает частицы в случайном направлении (при этом все направления равновероятны). Рассмотрим плоскость, в которой находится источник. Введем систему координат так, что источник имеет координаты $(\theta, d)$, где $d>0$, а (бесконечная) фотопластина совпадает с осью абцисс.
\begin{enumerate}
\item 
Покажите, что горизонтальная координата точки попадания частицы в плоскость имеет  {\it распределение Коши} с плотностью:
$$
p(x)=\frac{d}{\pi(d^2 + (x - \theta)^2)}.
$$
\item 
Можно ли вычислить математическое ожидание этого распределения?
\item 
Предположим, что параметр $\theta$ неизвестен, но имеется $n\gg1$ независимых реализаций рассматриваемой с.в. $X_1,\ldots,X_n$. Предложите способ оценивания параметра $\theta$. То есть необходимо указать измеримую функцию от выборки $X_1,\ldots,X_n$ $\hat{\theta}_n = \hat \theta_n(X_1,\ldots,X_n)$, значение которой будет расцениваться в качестве приближения к неизвестному истинному значению $\theta$.

\begin{remark}
К такого вида оценке целесообразно предъявить, например, следующие требования:
\begin{enumerate}
\item состоятельность: $\widehat \theta_n \xrightarrow{\text{п.н.}} \theta$ или $\widehat \theta_n \xrightarrow{p} \theta$;
\item несмещенность: $\Exp(\widehat \theta_n) = \theta$?
\end{enumerate}

Для сравнения между собой различных оценок одного и того же параметра выбирают некоторую \textit{функцию риска}, которая измеряет отклонение оценки от истинного значения параметра. В качестве функции риска часто выбирают дисперсию оценки. 
См., например, Косарев Е.Л. Методы обработки экспериментальных данных.  М.: Физматлит, 2008.  208 с, а также выступление Голубева Г.К. ``Вероятностные методы классической математической статистики'' на www.mathnet.ru. 

Интересная интерпретация распределения Коши возникает при изучении стрельбы лучником с завязанными глазами (см. популярные книги Б.~Мандельброта).
\end{remark}

\end{enumerate}
\end{problem}

\begin{problem}
\label{condExp1}
Предположим, что с.в. $X\in L_2$, это означает ${\mathbb E}(X^2)<\infty$. Докажите, что 
\begin{equation*}
\label{UMO}
\| X-{\mathbb E}(X|Y_1,\ldots,Y_n)\|_{2}=\min\limits_{\varphi\in H} \| X-\varphi(Y_1,\ldots,Y_n)\|_{2} , 
\end{equation*}
где $H$ -- подпространство пространства $L_2$ всевозможных борелевских функций $\varphi(Y_1,\ldots,Y_n)\in L_2$; 
${\mathbb E}(X|Y_1,\ldots,Y_n)$ -- условное математическое ожидание с.в. $X$ относительно $\sigma$-алгебры, порожденной с.в. 
$Y_1,\ldots,Y_n$, часто говорят просто относительно с.в. $Y_1,\ldots,Y_n$; 
$$
\| X\|_{2}=\sqrt{\langle X,X\rangle}=\sqrt{{\mathbb E}(X\cdot X)}=\sqrt{{\mathbb E}(X^2)} . 
$$
Обобщите утверждение задачи на случай, когда $X$ -- случайный вектор.
\end{problem}

\begin{ordre}
Покажите, что $(X-{\mathbb E}(X|H)) \bot \xi,\; \forall\,\xi\in H$, т.е. ${\mathbb E}(\cdot|H)$ 
является проектором на подпространство $H$ в $L_2$. Детали см., например, в учебнике: Розанов Ю.А. Теория вероятностей, случайные процессы и математическая статистика. М.: Наука, 1985.  320 с.
\end{ordre}

\begin{remark}
Данную задачу и последующие следует сравнить с задачами \ref{scalar_prod} и \ref{scalar_prod_1} 
из раздела \ref{standart}.
\end{remark}

\begin{problem}
\label{cond}
Используя соотношение из задачи \ref{condExp1}  в качестве определения условного математического ожидания, докажите справедливость основных свойств математического ожидания, в частности,
\[
\Exp(\Exp(X|Y)) = \Exp(X). 
\]
\end{problem}

\begin{ordre}
$\langle 1, X \rangle = \langle 1, \Exp(X|Y) \rangle$.
\end{ordre}

\begin{problem}
\label{condExp3}
Докажите, что если в условиях предыдущей задачи вектор $(X,Y_1,\ldots,Y_n)^T$  является нормальным случайным вектором (без ограничения 
общности можно также считать, что $(Y_1,\ldots,Y_n)^T$  -- невырожденный нормальный случайный вектор), то в качестве $H$ можно взять 
подпространство всевозможных линейных комбинаций с.в. $Y_1,\ldots,Y_n$. То есть  можно более конкретно сказать, на каком именно 
классе борелевских функций достигается минимум в задаче $\ref{UMO}$. 
\end{problem}

\begin{ordre}
Будем искать 
${\mathbb E}(X|Y_1,\ldots,Y_n)$ в виде 
$$
\label{Gauss}
{\mathbb E}(X|Y_1,\ldots,Y_n)=c_1 Y_1+\ldots +c_n Y_n . 
$$
Докажите следующие утверждения:

\begin{enumerate}
\item $X-c_1 Y_1-\ldots-c_n Y_n, Y_1,\ldots, Y_n$ -- независимы;
\item $X-c_1 Y_1-\ldots-c_n Y_n$ ортогонален подпространству $H$ пространства $L_2$ всевозможных борелевских функций $\varphi(Y_1,\ldots,Y_n)\in L_2$.
\end{enumerate}
 
\end{ordre}

\begin{problem}
Ведущий приносит два одинаковых конверта и говорит, что в них лежат деньги, причем в одном вдвое больше, чем в другом. Двое участников берут конверты и тайком друг от друга смотрят, сколько в них денег. Затем один говорит другому: ``Махнемся не глядя?'' (предлагая поменяться конвертами). Стоит ли соглашаться?
\end{problem}
\begin{remark}
См. \cite{book2012} и книгу \cite{book12} (этот материал полезно посмотреть и в связи со следующей задачей). В последней книге затрагиваются  вопросы о аксиоматике теории вероятностей. Отметим, что аксиоматика А.Н.~Колмогорова, построенная на теории меры, не единственный способ ввести с.в. 
\par Интересные материалы на эту тему имеются в статье Давида Мамфорда ``На заре эры стохастичности'', в сборнике ``Математика: границы и перспективы'' / под ред. Д.В.~Аносова и А.Н.~Паршина. М.: Фазис, 2005. 
Также отметим в этом сборнике статью W.T.~Gowers’а и аналогичный сборник ``Математические события XX века''. М.: Фазис, 2003.
По данной тематике полезны следующие материалы: Jaynes~E.T. Probability theory: the logic of science.  Cambridge University Press, 2003.; Кановей~В.Г., Любецкий~В.А. Современная теория множеств: абсолютно неразрешимые классические проблемы. М.: \mbox{МЦНМО}, 2013.  320~с.  См. также задачу \ref{banah_tar}.
\end{remark}

\begin{problem}(Спящая красавица.)
В воскресенье с красавицей обговаривается схема эксперимента, согласно которой  вечером в воскресенье красавица засыпает. Далее подкидывается симметричная монетка. Если монетка выпадает орлом, то красавицу будят в понедельник (потом снова дают снотворное), затем будят еще раз во вторник (потом снова дают снотворное). Если решкой, то будят только в понедельник (потом снова дают снотворное). В среду красавицу пробуждают окончательно в любом случае. Снотворное стирает красавице память в том смысле, что она помнит правила, оговоренные с ней в воскресенье, но не помнит, сколько раз уже ее будили, и не знает, какой сегодня день недели. Каждый раз, когда красавицу будят ей предлагают оценить вероятность того, что монетка выпала решкой. Что может ответить красавица?
\end{problem}

\begin{problem}(Задача о лабиринте, А.Н. Соболевский.)
``В лесах дремучих стоит дом не дом, чертог не чертог, а дворец зверя лесного, чуда морского, весь в огне, в серебре и золоте и каменьях самоцветных''. Красная девица входит на широкий двор, в ворота растворенные, и находит там три двери, а за ними три горницы красоты несказанной, а в каждой из тех горниц еще по три двери, ведущие в горницы краше прежних.
Походив по горницам, красная девица начинает догадываться, что дворец построен ярусами: двери со двора ведут в горницы первого яруса, из тех --- в горницы второго яруса, и так далее. В каждой горнице есть вход и три выхода, ведущие в три горницы следующего яруса. В горницах последнего $n$-го яруса растворены окна широкие во сады диковинные, плодовитые, а в садах птицы поют и цветы растут.
Вернувшись на широкий двор и отдохнув, красная девица видит, что произошла перемена: ворота, через которые она вошла, и часть дверей внутри дворца сами собой затворились,
да не просто так, а каждая дверь с вероятностью 1/3 независимо от других. Немного обеспокоенная, красная девица начинает метаться из горницы в горницу сквозь оставшиеся незатворенными двери в поисках выхода. Покажите, что при больших $n$ вероятность того, что
она сможет добраться до окон, растворенных в сады, близка к $(9 - \sqrt{27}) / 4 \approx 0.95$.
\end{problem}

\begin{remark}
См.  Соболевский А.Н. Конкретная теория вероятностей.  
\noindent \verb|http://ium.mccme.ru/postscript/f10/sobolevskii-main.pdf |.

Задача свзяна с изучением вырождением ветвящегося процесса (см. задачу 2 раздела 6).

\end{remark}

 \begin{problem}(Задача М. Гарднера о разборчивой невесте.)
 В аудитории находится невеста, которая хочет выбрать себе жениха. За дверью выстроилась очередь из $N$ женихов. Относительно любых 
 двух женихов невеста может сделать вывод, какой из них для неё предпочтительнее. Таким образом, невеста задает на множестве женихов 
 отношение порядка (естественно считать, что если $A$ предпочтительнее $B$, а $B$ предпочтительнее $C$, то $A$ предпочтительнее $C$). 
 Предположим, что все $N!$ вариантов очередей равновероятны и невеста об этом знает (равно, как и число $N$). Женихи запускаются 
 в аудиторию по очереди. Невеста видит каждого из них в первый раз! Если на каком-то женихе невеста остановится (сделает свой выбор), 
 то оставшаяся очередь расходится. Невеста хочет выбрать наилучшего жениха (исследуя $k$-го по очереди жениха, невеста лишь может 
 сравнить его со всеми предыдущими, которых она уже просмотрела и пропустила). Оцените (при $N\to\infty$) вероятность того, что невесте 
 удастся выбрать наилучшего жениха, если она придерживается следующей стратегии: просмотреть (пропустить) первых по очереди $[N/e]$ 
 кандидатов и затем выбрать первого кандидата, который лучше всех предыдущих (впрочем, такого кандидата может и не оказаться, тогда, 
 очевидно, невеста не смогла выбрать наилучшего жениха). 
 \end{problem}
 \begin{remark}
 Можно показать, что описанная стратегия будет асимптотически оптимальной (подробнее об этом будет написано в Части 2). Популярное изложение имеется у С.М. Гу\-сейн-Заде. Однако полезно обратить внимание, что эта задача является ярким примером целого направления: оптимальной остановки стохастических процессов, перетекающих в оптимальное управление процессами Маркова, а точнее еще более общего направления управляемых процессов Маркова. Наиболее полезным для приложений (особенно в области ``Исследование операций'') во всем этом является распространение принципа динамического программирования Беллмана на стохастический случай ({\itпринцип Вальда--Беллмана}). Об этом можно прочитать, например, у Е.Б.~Дынкина~-- А.А.~Юшкевича, А.Н. Ширяева, В.И.~Аркина~-- И.В.~Евстигнеева. Для введения можно посмотреть книги Е.С.~Вентцель или Ю.А.~Розанова \cite{6}, см. также задачу~\ref{m_field_games} из раздела~\ref{macrosystems}. Также об этом направлении подробнее планируется нписать в Части 2. Студенты ФУПМ МФТИ также столкнутся с этим направлением в курсе стохастических дифференциальных уравнений. В таком контексте полезно будет посмотреть книгу: Оксендаль~Б.  Стохастические дифференциальные уравнения: Введение в теорию и приложения.  М.: Мир, 2003.  408~с.
 \end{remark}
 
 \begin{problem}(Биржевой парадокс.)
Рассмотрим любопытный экономический пример. Пусть имеется начальный капитал $K_1$, который требуется увеличить. Для этого имеются две возможности: вкладывать деньги в надежный банк и покупать на бирже акции некоторой компании. Пусть $u$ -- доля капитала, вкладываемая в банк, а $v$ -- доля капитала, расходуемая на приобретение акций ($0 \leq u + v \leq 1$). Предположим, что банк гарантирует $b \times 100 \%$ годовых, а акции приносят $X \times 100 \%$ годовых, где $X$ -- с.в. с математическим ожиданием $m_X > b > 0$.  Таким образом, через год капитал составит величину $K_2 = K_1 (1 + b u + Xv)$. Очевидно, что, если придерживаться стратегии, максимизирующей средний доход за год, то выгодно присвоить следующие значения: $u = 0$, $v = 1$. 

Рассмотрите прирост капитала $K_{t+1}$   за  $t$ лет, считая $X_1, \ldots,  X_t$ независимыми с.в.  Покажите, что при ежегодном вложении капитала в акции  
\[
\Exp(K_t) \rightarrow \infty \; \text{при} \; t  \rightarrow \infty,
\]
\noindent но при этом в случае $\Exp\left[ \log (1 + X) \right] < 0$     

\[
K_t \overset{\text{п.н.}}{\longrightarrow}  0 \; \text{при} \; t  \rightarrow \infty.
\]
Приведите пример такой с.в. $X$.

\end{problem}

\begin{ordre}
Воспользуйтесь усиленным законом больших чисел, введя замену $K_{t+1} = K_1 e^{t Y_t}$, где  $Y_t = \frac{1}{t} \sum \limits_{i=1}^{t}\log(1+X_i)$. 
\end{ordre}

\begin{remark}
Рассмотренный парадокс хорошо иллюстрирует особенность поведения случайной последовательности, которая в отличие от детерминированной может сходиться в разных смыслах к разным значениям. 

Причина рассмотренного биржевого парадокса -- неудачный выбор критерия оптимальности. В качестве альтернативы могут быть выбраны следующие критерии:

\begin{enumerate}
\item \textit{Логарифмическая стратегия} или \textit{стратегия Келли}:
\[
Y_t \rightarrow \lambda(u) = \Exp\left[ \log (1 + b u + X (1-u)) \right], t \to \infty,
\]
\[
\lambda(u) \rightarrow \max \limits_{ 0 \leq u \leq 1},
\]

\[
\Exp\left[ \log (K_{t+1}) \right] \rightarrow \max \limits_{ u_1 \ldots u_t }.
\]

\item {\itВероятностный критерий}:
\[
\mathbb{P}_\phi (u_1 \ldots u_t ) = \mathbb{P} (K_{t+1} \geq -\phi) \rightarrow \max \limits_{ u_1, \ldots, u_t }.
\]
Используя метод \textit{динамического программирования} (см. Часть 2, а также Bertsekas D.P. Dynamic Programming and Optimal Control. V.~1.  Belmont, Massa\-chu\-setts: Athena Scientific,  1995.  558 p.), можно показать, что оптимальное управление удовлетворяет следующему соотношению: 

\[
u_t(K_t) = \begin{cases}
\begin{array}{cc}
0, & \phi \geq - K_t(1+b), \\
1, & \phi < - K_t(1+b).
\end{array}\end{cases}
\]

\end{enumerate}

Подобный тип управления (стратегии) созвучен исторической практике накопления капитала: в эпоху первичного накопления люди зачастую серьезно рисковали ради денег, но как только они накапливали сумму, достаточную для безбедного существования при ее вложении хотя бы в банк, необходимость в риске для них отпадала. Более детальное рассмотрение задачи см. в Кибзун А.И., Кан Ю.С. Задачи стохастического программирования с вероятностными критериями.  М.: Физматлит, 2009.  347 с., см. также книгу \cite{book12}, которая будет полезна и для решения следующей задачи.
\end{remark} 
 
 \begin{problem}(Парадокс страхования.)
 В некотором городе имеется много жителей и одна большая страховая компания. Каждый год каждый житель (независимо от остальных) с вероятностью $p = 0,1$ может потерять половину ($b = 0,5$) своего состояния. Страховая компания предлагает жителям страховать эту половину состояния под годовую ставку $c = 5,5 \%$. То есть житель, владеющий собственностью $V$, за год должен заплатить страховой компании $cV>pbV$. Если этот (застрахованный) житель теряет половину совего состояния, то страховая компания полностью компенсирует ему эти потери. Почему при такой ставке $с$ страховая компания может рассчитывать на прибыль? Почему в таком случае страхование может быть выгодно и жителям города?
 \begin{ordre}
 Для ответа на второй вопрос, считайте, что каждый житель довольно долго живет ($n\gg 1$ лет) и сравнивает свои потери в случае использования услуг страховой компании $= V_0(1-c)^n$ с ожидаемыми потерями в случае отказа от услуг страховой компании $\approx V_0(1-b)^{np}$.
 \end{ordre}
 \end{problem}

\begin{problem}(Равновесие Нэша \cite{27}.)
Один игрок прячет (зажимает в кулаке) одну или две монеты достоинством 10 рублей. Другой игрок должен отгадать, сколько денег у первого спрятано. Если отгадывает, то получает деньги, если нет -- платит 15 рублей. Каковы  должны быть стратегии игроков при многократном повторении игры?

\end{problem}

\begin{problem} (Равновесие Байеса--Нэша.)
В аудитории находится 100 человек (игроков). Каждого просят написать целое число от 1 до 100. Победителем окажется тот участник, который написал число, наиболее близкое к 2/3 от среднего арифметического всех чисел. Требуется найти  \textit{наилучший ответ} при фиксированных стратегиях соперников: каждый соперник, просчитав на $X \in \mathrm{Po}(2)$  хода  вперед, выбирает наугад (равномерно) число от 1 до $100 \cdot(2/3)^X$. К примеру, возможен следующий ход мыслей соперника: ``поскольку все догадались, что не стоит писать число большее $100 \cdot(2/3)$, то не стоит писать число большее $100 \cdot(2/3)^2$''.
\end{problem}

\begin{remark}
Стратегия (действие игрока $i$ в зависимости от своего типа $t_i \in T_i$) $s_i^*: T_{i} \rightarrow A_i$ называется \textit{наилучшим ответом} на заданные стратегии соперников $s_{-i}(t_{-i})$, если она является решением задачи максимизации ожидаемого выигрыша рассматриваемого игрока. При этом усреднение производится по всем неизвестным переменным, относящимся к соперникам (типам соперников $t_{-i}$): 
 \[
 \underset{t_{-i}}{\sum} u_i(s_i^*, s_{-i}(t_{-i}), t_i, t_{-i}) \cdot  \PR(t_{-i} | t_i) = \]\[ \underset{a_i \in A_i}{\max} \underset{t_{-i}}{\sum} u_i(a_i, s_{-i}(t_{-i}), t_i, t_{-i}) \cdot \PR(t_{-i} | t_i), 
 \]
 \noindent где $u_i$ -- выигрыш игрока $i$  при заданных действиях и типах всех игроков, $A_i$ -- множество действий игрока $i$, $\PR(t_{-i} | t_i)$ -- представление о типах остальных игроков при известном своем типе. Таким образом можно определить, вообще говоря, многозначное отображение наилучших ответов $S^*(s) = \left\{s_i^*\left(s_{-i}\right)\right\}_i$. Неподвижная точка этого отображения $S^*(s^*) = s^*$ будет \textit{равновесием Байеса--Нэша}. Такая точка не всегда существует, но всегда существует ее аналог в случае, когда игрок может смешивать несколько стратегий (т.е. смешанной стратегией будет дискретное распределение $(p_1,\ldots,p_{|A_i|}),\; \PR(s_i = k) = p_k$). 

См. также книгу R.B. Myerson. Game Theory, Harvard University Press, 1997. 
\end{remark}

\begin{problem}(Аукцион первой цены.)
Два участника аукциона конкурируют за покупку некоторого объекта. Ценности объекта    $v_1$ и $v_2$ для участников являются независимыми с.в., равномерно распределенными на отрезке [0, 1]. Участник  имеет точную информацию о своем значении $v_i$, но не знает $v_j$. Участники делают ставки из диапазона [0, 1] одновременно и независимо друг от друга. В данном аукционе побеждает тот, кто поставил большую ставку. При равенстве ставок бросается жребий. Каждый обязан заплатить по средней ставке, даже если ему объект не достается! Отказаться от участия в этом аукционе нельзя.
\begin{enumerate}
\item Выпишите функции выигрыша игроков в данной игре.
\item Найдите наилучший ответ игрока в данной игре в классе квадратичных стратегий: $b_i(v_i) = cv_i^2$, где $c > 0$. Зависит ли наилучший ответ от стратегии соперника в данном случае?
\item Покажите, что в этой игре нет других (кроме найденных в пункте б)) решений с гладкими  возрастающими симметричными стратегиями. \textit{Симметричными} стратегиями является пара $(b_1(v_1), b_2(v_2)) =$\linebreak $= (b(v_1), b(v_2))$.
\end{enumerate}
\end{problem}

\begin{ordre}
Стратегией игрока $i$ в данной игре является функция $b_i: [0, 1] \rightarrow [0, 1],$ равная величине ставки при ценности объекта  $v_i$. См. книгу Myerson R.B. Game Theory. Harvard University Press,  1997.  558~p. 
\end{ordre}

\begin{problem}(Распределения канторовского типа.)
Пусть $X_{k} $ -- взаимно независимые с.в. с распределением $\Be(1/2)$. 
 В сумме $\sum_{k=1}^\infty 2^{-k} X_{k}$ рассмотрим только слагаемые с четными номерами, или, что с точностью до множителя 3 (в дальнейшем потребуется для удобства), есть $Y=3\sum _{s=1}^{\infty }4^{-s} X_{s}  $. Покажите, что функция распределения $F(x)=$\linebreak $=\PR\left(Y\le x\right)$ является сингулярной (когда не оговаривается относительно какой меры, подразумевается, что относительно \textit{меры Лебега}, т.е. равномерной).

\begin{ordre}
Можно рассматривать $Y$ как выигрыш игрока, который получает $3\cdot 4^{-k} $, когда $k$-е бросание симметричной монеты дает в результате решку. Ясно, что полный выигрыш лежит между 0 и $3\left(4^{-1} +4^{-2} +\ldots \right)=1$. Если первое подбрасывание монеты привело к решке, то полный выигрыш $\ge {3\mathord{\left/ {\vphantom {3 4}} \right. \kern-\nulldelimiterspace} 4} $, тогда как в противоположном случае $Y\le 3\left(4^{-2} +4^{-3} +\ldots \right)=4^{-1} $. То есть неравенство ${1\mathord{\left/ {\vphantom {1 4}} \right. \kern-\nulldelimiterspace} 4} <Y<{3\mathord{\left/ {\vphantom {3 4}} \right. \kern-\nulldelimiterspace} 4} $ не может быть осуществлено ни при каких обстоятельствах, значит, $F(x)={1\mathord{\left/ {\vphantom {1 2}} \right. \kern-\nulldelimiterspace} 2} $ в интервале $x\in \left({1\mathord{\left/ {\vphantom {1 4}} \right. \kern-\nulldelimiterspace} 4} ,{3\mathord{\left/ {\vphantom {3 4}} \right. \kern-\nulldelimiterspace} 4} \right)$. Чтобы определить, как ведет себя функция распределения на интервале $x\in \left(0,{1\mathord{\left/ {\vphantom {1 4}} \right. \kern-\nulldelimiterspace} 4} \right)$, покажите, что на этом интервале график отличается только преобразованием подобия $F(x)=(1/2)F(4x)$.

\end{ordre}

\begin{remark}
Пример, когда свертка двух сингулярных распределений имеет непрерывную плотность: с.в. $X=\sum _{k=1}^{\infty }2^{-k} X_{k}  $ имеет равномерное распределение на интервале $\left(0, 1\right)$. Обозначим сумму членов ряда с четными и нечетными номерами через $U$ и $V$ соответственно. Ясно, что $U$ и $2V$ имеют одинаковое распределение и их распределение относится к канторовскому типу.
\end{remark}

\end{problem}

\begin{problem}
Напомним, что \textit{сингулярными} мерами называются меры, функции распределения $F(x)$ которых непрерывны, но точки их роста ($x$ -- точка роста $F(x)$, если для любого $\varepsilon >0$ выполняется: $F(x+\varepsilon )-F(x-\varepsilon )>0$) образуют множество нулевой меры Лебега. Покажите, что мера, соответствующая функции Кантора, сингулярна по отношению к мере Лебега.

\end{problem}

\begin{remark} (Соболевский А.Н.  Конкретная теория вероятностей.)
Любая вероятностная мера может быть представлена в виде суммы абсолютно непрерывной, дискретной и сингулярной мер. 

Сингулярные распределения вероятности возникают в эргодической теории и математической статистической физике как инвариантные меры диссипативных динамических систем, обладающих т.н. ``странными аттракторами''. Количественное изучение таких мер относится к геометрической теории меры и известно под названием ``фрактальной геометрии''. Основные импульсы развития этой дисциплины исходили из работ К. Каратеодори, Ф. Хаусдорфа, А. Безиковича 1920-х годов, а позднее -- Б. Мандельброта и многочисленных физиков, которые занимались ``динамическим хаосом'' в 1980-х годах (П. Грассбергер, И. Прокачча, Дж. Паризи, У. Фриш). Подробнее о мультифрактальных мерах см., например, книги: Федер Е. Фракталы.  М.: Мир, 1991.  260 с.;  Falconer K. Fractal Geometry:  Mathematical Foundations and Applications. John Wiley \&  Sons, 1990;  Песин~Я.Б. Теория размерности и динамические системы.  М.--Ижевск: Ин-т компьютерных исследований, 2002.  404 с.

\end{remark}

\begin{problem}(Модель Эрдёша--Реньи.)
\label{sec:erdRenyi}
 Пусть есть конечное множество\linebreak (в дальнейшем множество вершин) $V$, $\xi _{vv'} $ -- независимые с.в., занумерованные парами $\left\{v,v'\right\}\in V\times V$, $\vert V \vert = N$, $\xi _{vv'} \in \Be(p)$. 
Таким образом, можно задать абстрактный случайный граф на фиксированном множестве вершин. Покажите, что 
 
\begin{enumerate}
\item  при $p=\frac{1}{N^{1+\varepsilon } } $, $\varepsilon >0$, среднее число не изолированных вершин в случайном графе есть $o\left(N\right)$;
 
\item  при $p=\frac{1}{N^{1-\varepsilon } } $, $\varepsilon >0$, с вероятностью близкой к 1 существует связная компонента порядка $N$ ($N \gg 1$).
\end{enumerate}
\end{problem}
 \begin{remark}
 См. брошюру  Малышева~В.А. \cite{27}, аналогично для следующей задачи.
 \end{remark}

\begin{problem} 
Рассматривается конфигурация спинов $\omega =\left\{x_{mn} \right\}$ (где $x_{mn} \in$\linebreak $\in Be(p) $ -- независимые с.в.) на двумерной решетке $\left\{(m,n) \in {\mathbb Z}^{2} \right\}$. Вершину $(m,n)$ назовем занятой, если $x_{mn} =1$. Соединим ребром все соседние (находящиеся на расстоянии 1) занятые вершины. Получится случайный граф $G=G\left(\omega \right)$. Назовем кластером графа $G$ максимальное подмножество $A$ вершин решетки такое, что для любых двух $v,v'\in A$ существует связывающий их путь по ребрам графа $G$. Докажите, что существует такое $0<\bar{p}<1$, что при $p<\bar{p}$ все кластеры конечны с вероятностью 1, а при $p>\bar{p}$ с положительной вероятностью есть хотя бы один бесконечный кластер.
 \end{problem}
 
\begin{ordre}
Покажите, что при достаточно малых значениях $p$ вероятность события, что все кластеры конечны, равна~1. Покажите, что вероятность того, что кластер, содержащий начало координат и имеющий не менее $N$ вершин, не превосходит $\left(Cp\right)^{N} \mathop{\to }\limits_{N\to \infty } 0$, где $C$ -- некоторая константа. А значит, и событие: бесконечный кластер содержит начало координат -- имеет нулевую вероятность.
\end{ordre}

\begin{problem}

На некоторой реке имеется 6 островов (см. рис. \ref{Fig:graphs_bridges.png}), соединенных между собой системой мостов. Во время летнего наводнения часть мостов была разрушена. При этом каждый мост разрушается с вероятностью ${1\mathord{\left/ {\vphantom {1 2}} \right. \kern-\nulldelimiterspace} 2} $, независимо от других мостов. Какова вероятность того, что после наводнения можно будет перейти с одного берега на другой, используя не разрушенные мосты? См.  Cherny A.S. The Kolmogorov student's competitions on probability theory. MSU.\\
\noindent \verb|http://www.newton.ac.uk/preprints/NI05043.pdf |

\end{problem}

\imgh{30mm}{graphs_bridges.png}{Схема мостов}

\imgh{70mm}{guk.jpg}{Иллюстрация к задаче Перколяция. Анастасия Ковылина, 2014}

\begin{problem}\Star(Перколяция.)
В квадратном пруду (со стороной, равной 1) 
выросли (случайным образом) $N\gg 1$ цветков лотоса, имеющих форму круга 
радиуса $r>0$. Назовем $r = r_N $  \textit{радиусом перколяции}, если с вероятностью 0.99, не 
любящий воду жук сможет переползти по цветкам лотоса с северного берега на 
южный, не замочившись. Покажите, что $r_N \sim \dfrac{C}{\sqrt{N}}$. Оцените $C$.

\end{problem}

\begin{remark} 
См. Кестен Х. Теория просачивания для математиков.  М.: Мир, 1986.  392 c.; Grimmett G. Percolation.  Springer,  1999.  447 p.
\end{remark}

\begin{problem} (Метод ренорм группы.) В квадратной таблице задан процесс окрашивания ячеек: с вероятностью $p$ ячейка независимо окрашивается в черный цвет, иначе остается белой. Исследуем в такой системе размеры связных компонент из черных ячеек. При небольших значениях $p$ образуется много отделимых клеточных областей, в то время как при $p \simeq 1$ c большой вероятностью получится несколько  компонент, соизмеримых со всей таблицей. Установлено, что значение $p_c = 0.5927462\ldots$ является критическим, при превышении которого в результате раскраски  среднее значение размера компоненты растет вместе с увеличением размера таблицы (таблица достаточно большая). Определим $\pi(s)$ как плотность распределения площадей компонент, где площадь одной клетки равна $a$. Допустим, что $\pi(s)$ зависит только от параметров $a$ и $\Exp s$. Так как исследуемая плотность не должна зависеть от масштаба, то 
\[
\pi(s) = f\left( \frac{s}{a}, \frac{a}{\Exp s} \right).
\]      
Изменим площадь одной клетки $a \to a/b$, тогда  ввиду независимости от масштаба в области ($p < p_c$) плотность должна принять аналогичный вид  
\[
\pi(s) = c(b) f\left( \frac{s}{a/b}, \frac{a/b}{\Exp s} \right) = c(b) f\left( \frac{sb}{a}, \frac{a}{\Exp (sb)} \right) .
\]
Если $p \to p_c$, то $\Exp s \to \infty$ (неограниченный рост компонент c большой площадью) и, следовательно,
\[
\pi(s) = c(b) f\left( \frac{sb}{a}, 0 \right) = c(b)\pi(sb).
\]
Докажите, что распределение со свойством $\pi(s) = c(b) \pi(sb)$  является степенным ($\pi(s) \propto s^{-\beta}$).
\end{problem}

\begin{remark}
В контексте этой и последующих трех задач рекомендуем ознакомиться с работами Newman~M.E.J. Power laws, Pareto distributions and Zipf's law // Contemporary physics.  2005. V.~46, N~5. P.~323--351;
 Mitzenmacher\;M.\ \;A brief history of generative models for power law and lognormal distributions // Internet mathematics. 2004. V.~1, N~2. P.~226--251.
 
 Степенные законы распределения с.в. уже встречались ранее в задачах \ref{records} и \ref{cauchy_rad_emitt}. Также эти законы будут встречаться в следующих четырех задачах, а также в следующих разделах, см. задачи \ref{lognormal}, \ref{bluzd_ust}, \ref{Holtsmark}, \ref{Cramer_Zone} раздела \ref{zb4} и задачи 4, 5 раздела 6.
\end{remark}

\begin{problem} (Масштабируемость степенного распределения, Б. Мандельброт.)
Закон А. Лотки гласит, что ``на каждого ученого, который написал за всю жизнь не менее $k$ научных работ, приходится $k^2$ ученых, творчество которых свелось к их первой работе (т.е. они сделали всего одну научную работу)''. Покажите, что отсюда можно сделать следующий довольно неожиданный вывод: что если про (случайно выбранного) ученого известно, что он уже написал $n$ работ, то можно надеяться, что в будущем он напишет еще в среднем столько же работ -- $n$.
\begin{ordre}
Воспользуйтесь следующим свойством степенного распределения ($x>h>0, \alpha > 1$): если $\PR\left( X \ge x \right) \propto x^{-\alpha}$, то $$\PR\left( X \ge x | X \ge h \right) = \left(\frac{x}{h}\right)^{-\alpha}.$$  Следовательно, $$\Exp\left[X-h|X \ge h\right] = \frac{h}{\alpha - 1}.$$ Из условий задачи следует, что $\alpha = 2$. 
\end{ordre}
\end{problem}

\begin{problem}(Логнормальное распределение и мультипликативный процесс.)
\label{lognorm}
Рассмотрим частицу, которая при перемещении из одного места
в другое может разделиться на несколько меньших частиц вследствие соударения или другого воздействия. Обозначим через $K(t) \sim \Po(\lambda t)$ число отделившихся частей к моменту времени $t$. Первоначальный размер частицы равен $s_0$. Пусть $D_i$  -- доля частицы, отделившаяся при $i$-м соударении. Тогда размер частицы в момент времени $t$ имеет вид (следует сравнить с задачей \ref{sec:cpoisson} раздела \ref{zb4})
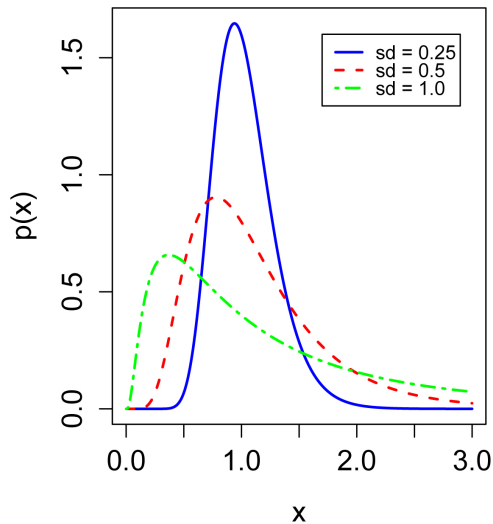
\imgh{70mm}{log_norm_dist.jpg}{График функции плотности лог-нормального распределения с $\mu = 0$ и различных значениях среднеквадратического отклонения $\sigma \, (0.25;\, 0.5;\, 1.0)$  (standard deviation)}
\[
Z(t) = s_0 \prod \limits_{i=1}^{K(t)} (1 - D_i),
\]
\[
S(t) = \ln Z(t) =  \mu +  \sum \limits_{i=1}^{K(t)} Y_i,
\]
\noindent где $\mu = \ln s_0$, $Y_i = \ln (1 - D_i)$ -- независимые с.в. (по~пред\-положению), причем $\mathbb{E} Y = a < 0$, $\mathbb{D} Y = \sigma^2$. 

Покажите, что при $t\to\infty$ с.в. $Z(t)$ стремится к лог-нормальному распределению (см. замечание):
\[
Z(t) \sim\mathrm{Log}\mathcal{N} \left(\mu + \lambda t a, \lambda t (a^2 + \sigma^2) \right).
\] 
\end{problem}

\begin{ordre}
Следует использовать Ц.П.Т. и задачу \ref{sec:ordered_seq} раздела \ref{standart}.
\end{ordre}

\begin{remark}
Cм. Mitzenmacher~M. A Brief History of Generative Models for Power Law and Lognormal Distributions // Internet Math. 2004.  V.~1, N~2.  P. 226--251.

Говорят, что $X$  имеет лог-нормальное распределение $\mathrm{Log}\mathcal{N}(\mu,  \sigma^2)$ с параметрами $\mu$ и $\sigma$, если распределение задаётся плотностью вероятности, имеющей вид (см. рис. \ref{Fig:log_norm_dist.jpg}):
\[
f_X(x) = \frac{1}{x \sigma \sqrt{2 \pi}} e^{- \frac{(\ln x - \mu)^2}{2 \sigma^2}}, \; x > 0.
\]
\end{remark}

Если после логарифмирования с.в. нормально распределена, то исходная с.в. логарифмически нормально распределена.



\begin{problem}(Growth model with preferential attachment.)
\label{prefattach}
Рассмотрим следующую модель роста неориентированного графа, состоящего изначально из одной вершины и петли: пусть в единицу времени  в 
графе добавляются новая вершина и новое 
ребро, соединяющее новую вершину со случайно выбранной имеющейся вершиной 
(вероятность выбора пропорциональна степени вершины -- принцип предпочтительного присоединения).

Покажите, что в этой модели применим  степенной закон для степеней вершин, а 
именно, что доля вершин со степенью не меньше $d$ приближенно равна $1/d^2$.
\begin{remark}
Обозначим через $d_k (t)$  ожидаемое (среднее) значение степени $k$-й вершины в момент 
времени $t$. Покажите, что для $d_k (t)$ справедливо следующее 
дифференциальное уравнение:
\[
\frac{d }{d t}d_k (t)=\frac{d_k (t)}{2t}
\]
с условием $d_k (k)=1 $. Откуда следует, что $d_k (t)= \sqrt 
{t/k} $.

Детали см.~в пишущейся книге Blum A., Hopcroft J., Kannan R. Foundations of Data Science // e-print. 2016.\\ 
\noindent\verb|http://www.cs.cornell.edu/jeh/book.pdf |
\end{remark}
\end{problem}

\begin{problem}(Модель Бакли--Остгуса и предпочтительное присоединение.)
\label{pref_attach}
Рассмотрим следующую модель роста web-графа, которую, с точностью до 
небольшой оговорки, можно называть моделью Бакли--Остгуса (см.~задачу \ref{hnmgraph} из раздела \ref{combinatorics}). 

\textit{База индукции.} Сначала имеется всего одна вершина, которая ссылается сама на себя (вершина 
с петлей). 

\textit{Шаг Индукции.} Предположим, что уже имеется некоторый (ориентированный граф). Пусть 
появляется новая вершина. Тогда с вероятностью $\beta >0$ из этой вершины 
проводится ребро равновероятно в одну из существующих вершин, а с 
вероятностью $1-\beta $ из этой вершины проводится ребро в одну из 
существующих вершин не равновероятно, а с вероятностями пропорциональными 
входящим степеням вершин (выходящая степень всех вершин одинакова и равна 1) 
-- \textit{правило предпочтительного присоединения} (от англ. preferential attachment).

Другими словами, если уже построен граф из $n-1$ вершины, то новая $n-$я вершина 
сошлется на вершину $i=1,...,n-1$ с вероятностью
\[
\frac{\mbox{in}\deg _{n-1} \left( i \right)+a}{\left( {n-1} \right)\left( 
{a+1} \right)},
\]
где $\mbox{in}\deg _{n-1} \left( i \right)$ -- входящая степень вершины $i$ в 
графе, построенном на шаге $n-1$. Параметры $\beta $ и $a$ связаны следующим 
образом
\[
a=\frac{\beta }{1-\beta }.
\]
Интернету (host-графу) наилучшим образом соответствует значение $a=0.277$.

Далее вводится число $m$ -- среднее число web-страниц на одном сайте, и 
каждая группа web-страниц с номерами $km,km+1,...,\left( {k+1} \right)m$ 
объединяется в один сайт. При этом все ссылки (имеющиеся между 
web-страницами) ``наследуются'' содержащими их сайтами -- получается, что с 
одного сайта на другой может быть несколько ссылок. Впрочем, сайты могут 
совпадать -- внутри одного сайта web-страницы также могут друг на друга 
сослаться. Пусть, скажем, получилось, что для заданной пары сайтов таких 
(одинаковых) ссылок оказалось $l\le m$, тогда мы превращаем их в одну 
ссылку, но с весом (вероятностью перехода) $l \mathord{\left/ {\vphantom {l 
m}} \right. \kern-\nulldelimiterspace} m$. 

Установите степенной закон распределения входящих степеней вершин в модели 
Бакли--Остгуса. При этом ограничьтесь случаем $m=1$ (этот закон не 
будет зависеть от $m)$.

\begin{ordre}
Обозначим через $X_k \left( t \right)$ число вершин с входящей степенью $k$ 
в момент времени $t$, т.е. когда в графе имеется всего $t$ вершин. Заметим, 
что по определению
\[
t=\sum\limits_{k\ge 0} {X_k \left( t \right)} =\sum\limits_{k\ge 1} {kX_k 
\left( t \right)} =\sum\limits_{k\ge 0} {kX_k \left( t \right)} .
\]
Поэтому, для $k\ge 1$ вероятность того, что $X_k \left( t \right)$ 
увеличится на единицу при переходе на следующий шаг $t\to t+1$ по формуле 
полной вероятности равна
\[
\beta \frac{X_{k-1} \left( t \right)}{t}+\left( {1-\beta } 
\right)\frac{\left( {k-1} \right)X_{k-1} \left( t \right)}{t}.
\]
Аналогично, для $k\ge 1$ вероятность того, что $X_k \left( t \right)$ 
уменьшится на единицу при переходе на следующий шаг $t\to t+1$
\[
\beta \frac{X_k \left( t \right)}{t}+\left( {1-\beta } \right)\frac{kX_k 
\left( t \right)}{t}.
\]
Таким образом ``ожидаемое'' приращение $\Delta X_k \left( t \right)=X_k \left( {t+1} \right)-X_k \left( t \right)$ за 
$\Delta t=\left( {t+1} \right)-t=1$
будет
\[
\frac{\Delta X_k \left( t \right)}{\Delta t}=\beta \frac{\left( {X_{k-1} 
\left( t \right)-X_k \left( t \right)} \right)}{t}+\left( {1-\beta } 
\right)\frac{\left( {k-1} \right)X_{k-1} \left( t \right)-kX_k \left( t 
\right)}{t}.
\]
Для $X_0 \left( t \right)$ аналогичное уравнение будет иметь вид
\[
\frac{\Delta X_0 \left( t \right)}{\Delta t}=1-\beta \frac{X_0 \left( t 
\right)}{t}.
\]
К сожалению, соотношения полученные соотношения -- не есть точные 
(разностные) уравнения, описывающие то, как меняется $X_k \left( t \right)$, 
хотя бы потому, что изменение $X_k \left( t \right)$ происходит случайно. 
Предложенная выше динамика полностью детерминированная. Однако, для больших 
значений $t$, когда наблюдается концентрация случайных величин $X_k \left( t 
\right)$ вокруг своих математических ожиданий (средних значений), реальная 
динамика поведения $X_k \left( t \right)$ и динамика поведения средних 
значений $X_k \left( t \right)$ становятся близкими -- вариация на тему 
\textit{теоремы Куртца} -- см. задачу 3 раздела 6. Таким образом, на полученную систему можно смотреть, как на динамику 
средних значений, вокруг которых плотно сконцентрированы реальные значения. 
Под плотной концентрацией имеется в виду, что разброс значений величины 
контролируется квадратным корнем из ее среднего значения.

Будем искать решение этой системы на больших временах ($t\to \infty )$ в 
виде $X_k \left( t \right)\sim c_k \cdot t$ (иногда такого вида режимы 
называют \textit{промежуточными асимптотиками}):
\[
c_0 =\frac{1}{1+\beta },
\quad
\frac{c_k }{c_{k-1} }=1-\frac{2-\beta }{1+\beta +k\cdot \left( {1-\beta } 
\right)}\simeq 1-\left( {\frac{2-\beta }{1-\beta }} \right)\frac{1}{k}.
\]
Откуда получаем следующий \textbf{\textit{степенной закон}}
\[
c_k \sim k^{-\frac{2-\beta }{1-\beta }}=k^{-2-a}.
\]
Заметим, что если построить на основе этого соотношения \textit{ранговый закон} распределения 
вершин по их входящим степеням, т.е. отранжировать вершины по входящей 
степени, начиная с вершины с самой высокой входящей степенью, то также 
получим степенной закон 
\[
\mbox{in}\deg _r \sim r^{-1-\beta }.
\]
Действительно, обозначив для краткости $\mbox{in}\deg _r $ через $x$, 
получим, что нам нужно найти зависимость $x\left( r \right)$, если известно, 
что 
\[
\frac{dr\left( x \right)}{dx}\sim -x^{-\frac{2-\beta }{1-\beta 
}}=-x^{-1-\frac{1}{1-\beta }},
\]
где зависимость $r\left( x \right)$ получается из зависимости $x\left( r 
\right)$ как решение уравнения $x\left( r \right)=x$.
\end{ordre}

}

\begin{remark}
Для более глубокого погружения в тематику задачи можно рекомендовать следующие обзорные работы Mitzenmacher M. 
A Brief History of Generative Models for Power Law and Lognormal Distributions //
Internet Mathematics. 2004. V. 1, N 2. P. 226--251; 39.	Подлазов А.В. Закон Ципфа и модели конкурентного роста. Новое в синергетике. Нелинейность в современном естествознании. М.: ЛИБРОКОМ, 2009. – С. 229–256; Newman М.Е.J. Networks: An Introduction. Oxford University Press, 2010 и популярную статью arXiv:1701.02595.
\end{remark}

\end{problem}
 
\begin{problem}(Обобщенная схема размещений)
\noindent
\begin{enumerate}
\item  Пусть для целочисленных  неотрицательных с.в. $\eta _1,\ldots,\eta _N$ существуют независимые одинаково распределенные с.в. 
$\xi_1,\ldots,\xi_N$ такие, что
\[ 
\PR\left( {\eta _1 =k_1,\ldots,\eta _N =k_N } \right) =  \]\[ =\PR\left( {\left. {\xi _1 
=k_1 ,\ldots,\xi _N =k_N } \right|\xi _1 +\ldots+\xi _N =n} \right). \;\;\;\; (*) \]
 
Введем независимые одинаково распределенные с.в. $\xi _1^{\left( r \right)} 
$,{\ldots},$\xi _N^{\left( r \right)} $, где $r$ -- целое неотрицательное число 
и
\[
\PR\left( {\xi _1^{\left( r \right)} =k} \right)=\PR\left( {\left. {\xi _1 =k} 
\right|\xi _1 \ne r} \right),
\quad
k=0,1,\ldots ,
n\]
Пусть $p_r = \PR\left( {\xi _1 =r} \right)$ и $S_N =\xi _1 +\ldots+\xi _N$,
$S_N^{\left( r \right)} =\xi _1^{\left( r \right)} +\ldots+\xi _N^{\left( r \right)} $. Пусть $\mu _r \left( {n,N} \right)$ -- число с.в. $\eta _1$,{\ldots},$\eta_N $, принявших значение $r$. Покажите, что с.в. типа $\mu_r \left( {n,N} \right)$ можно изучать с помощью \textit{обобщенной схемы размещений}: для любого $k=0,\ldots,N$
\[
\PR\left( {\mu _r \left( {n,N} \right)=k} \right)=C_n^k p_r^k \left( {1-p_r } 
\right)^{N-k}\frac{\PR\left( {S_{N-k}^{\left( r \right)} =n-kr} 
\right)}{\PR\left( {S_N =n} \right)}.
\]
Напомним, что в классической схеме размещений $n$ различных частиц по $N$ 
различным ячейкам было доказано, что распределение числа заполнений ячеек $\eta _1 
$,{\ldots},$\eta _N $ имеет вид (\textit{мультиномиальное распределение}):
\[
\PR\left( {\eta _1 =k_1 ,\ldots,\eta _N =k_N } \right)=\frac{n!}{k_1!\ldots k_N! N^n},
\]
где $k_1,\ldots,k_N$ -- неотрицательные целые числа, $\sum_{i} k_i =n$. Если положить $\xi_1,\ldots,\xi_N \in \Po \left( \lambda 
\right)$ -- независимые с.в. ($\lambda >0$ -- произвольно), то получим (*).
 
\item * Дан случайный граф (модель Эрдеша--Реньи, см. задачу \ref{sec:erdRenyi}) $G\left( {n,\;p} \right)$. Пусть 
$p=c\frac{\ln n}{n}$. Покажите, что при $n\to\infty$ и $c>1$ граф $G\left( {n,\;p} \right)$ 
почти наверное связен, а при $c>1$ почти наверное не связен.

\end{enumerate}
\begin{remark}
См. монографию Колчин В.Ф. Случайные графы.  М.: Физматлит, 2004.  256 с. Пункт б) подробно разобран в монографии Алона--Спенсера \cite{15}. Мы также рекомендуем смотреть популярные тексты А.М.~Райгородского по тематике случайных графов и их приложений, см. также раздел
\ref{combinatorics}.
\end{remark}
\end{problem}

\begin{problem}(Вероятностное доказательство формулы Эйлера \cite{7}.)
Пусть $X$ -- 
целочисленная с.в. с распределением

\[\PR\left( {X=n} \right)=\frac{1}{\zeta (s)n^s},\;\]
\noindent где $\zeta 
(s)=\sum\limits_{n\in {\mathbb N}} {n^{-s}} ,\quad s>1$.

Пусть $1<p_1 <p_2 <p_3 <\ldots $ -- простые числа, и пусть $A_k $ -- событие 
= \{$X$ делится на $p_k ${\}}.

\begin{enumerate}

\item Найдите $\PR(A_k)$ и покажите, что события $A_1 ,A_2 ,\ldots $ 
независимы.

\item Покажите, что
\begin{center}
$\prod\limits_{k=1}^\infty {(1-p_k ^{-s})} =\frac{1}{\zeta (s)}$ (формула 
Эйлера).
\end{center}

\end{enumerate}
\end{problem}

\begin{problem}\Star(Статистика теоретико-числовых функций.)
\label{sec:z_func_riman}
Довольно часто 
вероятностные соображения (например, независимость) используются в теории 
чисел не совсем строго, но зато весьма часто они позволяют угадать 
правильный ответ. Поясним сказанное, пожалуй, наиболее популярным примером (Дирихле, 1849) 
из книги Кац М. 
Статистическая независимость в теории вероятностей, анализе и теории 
чисел.  М.: ИЛ, 1963.  153 с.

Пусть $A$ -- некоторое множество положительных целых чисел. Обозначим через 
$A\left( n \right)$ количество тех его элементов, которые содержатся среди 
первых $n$ чисел натурального ряда. Если существует предел $\mathop {\lim 
}\limits_{n\to \infty } {A\left( n \right)} \mathord{\left/ {\vphantom 
{{A\left( n \right)} n}} \right. \kern-\nulldelimiterspace} n=\PR\left( A 
\right)$, то он называется плотностью $A$. К сожалению, вероятностная мера 
$\PR\left( A \right)$ не является вполне аддитивной (счетно-аддитивной).

Рассмотрим целые числа, делящиеся на простое число $p$. Плотность множества 
таких чисел, очевидно, равна $1 \mathord{\left/ {\vphantom {1 p}} \right. 
\kern-\nulldelimiterspace} p$. Возьмем теперь множество целых чисел, которые 
делятся одновременно на $p$ и $q$ ($q$ -- другое простое число). Делимость 
на $p$ и $q$ эквивалентна делимости на $pq$, и, следовательно, плотность 
нового множества равна $1/(pq)$. Так как $1/(pq)=( 1 /p)\cdot (1/q)$, 
то мы можем истолковать это так: ``события'', заключающиеся в 
делимости на $p$ и $q$, независимы. Это, конечно, выполняется для любого 
количества простых чисел.

Поставим теперь задачу посчитать долю несократимых дробей или, другими 
словами, ``вероятность'' несократимости дроби (фиксируется знаменатель дроби 
$n$, а затем случайно с равной вероятностью $1 \mathord{\left/ {\vphantom 
{1 n}} \right. \kern-\nulldelimiterspace} n$ выбирается любое число от 1 до 
$n$ в качестве числителя и подсчитывается доля случаев, в которых 
полученная дробь оказывалась несократимой) в следующем (чезаровском) смысле 
(здесь и далее индекс $p$ может пробегать только простые числа):
\[
\mathop {\lim }\limits_{N\to \infty } \frac{1}{N}\sum\limits_{n=1}^N 
{\frac{\# \left\{ {k<n:\;\;\text{Н.О.Д.}\left( {n,k} \right)=1} \right\}}{n}} 
=\mathop {\lim }\limits_{N\to \infty } \frac{1}{N}\sum\limits_{n=1}^N 
{\frac{\phi \left( n \right)}{n}} = \]\[=\mathop {\lim }\limits_{N\to \infty } 
\frac{1}{N}\sum\limits_{n=1}^N {\prod\limits_p {\left( {1-\frac{\rho _p 
\left( n \right)}{p}} \right)} } ,
\]
где $\phi \left( n \right)$ -- функция Эйлера, $$\rho _p \left( n 
\right)=\left\{ {\begin{array}{l}
 1,\quad n\mbox{ делится на }p, \\ 
 0\quad \mbox{иначе}. \\ 
 \end{array}} \right. $$

\noindent Согласно введенному выше определению плотности:

\[
\Exp\left\{ {\prod\limits_{p\le p_k } {\left( {1-\frac{\rho _p \left( n 
\right)}{p}} \right)} } \right\}=\prod\limits_{p\le p_k } {\Exp\left\{ {\left( 
{1-\frac{\rho _p \left( n \right)}{p}} \right)} \right\}} 
=\prod\limits_{p\le p_k } {\left( {1-\frac{1}{p^2}} \right)} .
\]
С учетом этого хочется написать следующее:
\[
\mathop {\lim }\limits_{N\to \infty } \frac{1}{N}\sum\limits_{n=1}^N 
{\frac{\phi \left( n \right)}{n}} =\Exp\left\{ {\frac{\phi \left( n 
\right)}{n}} \right\}=\Exp\left\{ {\prod\limits_p {\left( {1-\frac{\rho _p 
\left( n \right)}{p}} \right)} } \right\}\mathop =\limits^? 
\]
\[
\mathop =\limits^? \prod\limits_p {\Exp\left\{ {\left( {1-\frac{\rho _p \left( 
n \right)}{p}} \right)} \right\}} =\prod\limits_p {\left( {1-\frac{1}{p^2}} 
\right)} =\frac{1}{\zeta \left( 2 \right)}=\frac{6}{\pi ^2}.
\]
Будь введенная вероятностная мера, по которой считается это математическое 
ожидание, счетно-аддитивной, то можно было бы поставить точку, получив 
ответ. Однако это не так. Несмотря на правильность ответа, 
приведенное выше рассуждение не может считаться доказательством. Впрочем, 
часто вероятностные рассуждения удается пополнить, используя их в качестве 
основы. Так, в разобранном нами примере все сводится к обоснованию равенства 
``$?$''.

Легко понять, что полученный ответ несет определенную информацию о 
статистических свойствах функции Эйлера.

В теории чисел такого типа задачи занимают крайне важное место. Достаточно 
сказать, что гипотеза Римана ``на миллион'' о распределении нетривиальных 
нулей дзета-функции Римана $\zeta \left( z \right)$ равносильна
следующему свойству функции Мёбиуса $\mu(n)$: для любого $\epsilon > 0$

\begin{center}
$\left| {\sum\limits_{n=1}^N {\mu \left( n \right)} } \right| = O\left(N^{0.5 + \epsilon} \right) $ 
(Дж. Литлвуд),
\end{center}
где $\mu(n)$ определяется уравнением
\[
\sum\limits_{d \vert n} \mu(d) = 
\begin{cases}
1, & n = 1; \\
0, & n > 1.
\end{cases}
\quad
\Leftrightarrow
\quad
\mu(n) = 
\begin{cases}
(-1)^r, & n = p_1\ldots p_r; \\
0, & n \;  \mod \; p^2 = 0.
\end{cases}
\]
Что в свою очередь (Х.М. Эдвардс) в определенном смысле ``завязано'' на 
случайности последовательности $\left\{ {\mu \left( n \right)} 
\right\}_{n\in {\mathbb N}} $.\\

\noindent Используя указанный выше формализм, найдите долю чисел натурального ряда, 
``свободных от квадратов'', т.е. не делящихся на квадрат любого простого 
числа.

\end{problem}

\begin{remark} 
Известный российский математик Владимир Игоревич Арнольд 
последние десять лет жизни активно развивал описанное направление, которое 
он называл {\itЭкспериментальной математикой} (помимо популярных книжек и 
статей, осталось и несколько видеолекций на эту тему на mathnet.ru
с выступлениями на семинаре МИАН, в летней школе <<Современная математика>> и 
на мехмате). Получая схожим образом ``ответы'', их далее можно проверять, 
ставя численные эксперименты. При современных возможностях вычислительных 
машин можно отслеживать логарифмические функции в асимптотике (т.е. 
проверять гипотезы с логарифмами), но не с повторными логарифмами, которые 
так же, как и в случайных процессах, встречаются в теории чисел. Таким образом, 
у В.И. Арнольда получалось довольно много теорем (десятки, а возможно, даже сотни). 
Часть теорем, конечно, была известна ранее (см., например, книги Карацубы~А.А. Основы аналитической теории чисел. М.: Наука, 1975.  240 с; Иванец~Х., Ковальский~Э. Аналитическая теория чисел. М.: МЦНМО, 2014. 712 с.), но удавалось 
получать и новые формулировки.

Следующий пример, взятый из другой книги М. Каца \cite{20}, демонстрирует, что отмеченным выше 
способом можно получить и неверный результат.

\begin{example}(Из журнала Nature, 1940.) Из изложенного выше следует, что количество целых чисел, не превосходящих $N$ и не делящихся ни на одно из простых 
чисел $p_1 $, $p_2 $, {\ldots}, $p_k $, равно приблизительно 
$N\prod\limits_{j=1}^k {\left( {1-\frac{1}{p_j }} \right)} $. Рассмотрим 
теперь количество целых чисел, не превосходящих $N$ и не делящихся ни на одно из 
простых чисел, меньших $\sqrt N $. Такими числами могут быть только простые 
числа, лежащие между $\sqrt N $ и $N$, число которых
\[
\pi \left( N \right)-\pi \left( {\sqrt N } \right)\simeq N\prod\limits_{p_j 
<\sqrt N } {\left( {1-\frac{1}{p_j }} \right)} .
\]
Но из теории чисел известно, что $\pi \left( N \right)\simeq N \mathord{\left/ 
{\vphantom {N {\ln N}}} \right. \kern-\nulldelimiterspace} {\ln N}$ и
\[
\prod\limits_{p_j <\sqrt N } {\left( {1-\frac{1}{p_j }} \right)} \simeq {\exp 
\left( {-\gamma } \right)} \mathord{\left/ {\vphantom {{\exp \left( {-\gamma 
} \right)} {\ln \sqrt N }}} \right. \kern-\nulldelimiterspace} {\ln \sqrt N 
}={2\exp \left( {-\gamma } \right)} \mathord{\left/ {\vphantom {{2\exp 
\left( {-\gamma } \right)} {\ln N}}} \right. \kern-\nulldelimiterspace} {\ln 
N},
\]
где $\gamma $ -- константа Эйлера. Следовательно, $2\exp \left( {-\gamma } 
\right)=1$. Пришли к неверному соотношению!
\end{example}

Интересно в этой связи отметить также вероятностный способ получения 
правильной асимптотической формулы $\pi \left( N \right)\simeq N 
\mathord{\left/ {\vphantom {N {\ln N}}} \right. \kern-\nulldelimiterspace} 
{\ln N}$ для количества простых чисел, не превосходящих $N$, приведенный в книге 
Куранта Р.,  Роббинса Г.  Что такое математика.  М.: МЦНМО, 2007.  568 с., см. также \cite{2013}.

В заключение заметим, что применение вероятностных соображений в теории 
чисел продолжает привлекать ведущих математиков и по сей день (см., например, задачи \ref{permutation}, \ref{permutation1} раздела \ref{zb4}, статью
 Синай Я.Г. Статистические свойства функции Мебиуса // Автомат. и телемех. 2013. №~10. С. 6--14, а также книгу Тао Т.  Структура и случайность.  М.: МЦНМО, 2013. 360~с.).
\end{remark}

\begin{problem}\Star(Парадокс Банаха--Тарского.)
\label{banah_tar}
С. Банах показал, что если предполагать только аддитивность меры, то в одномерном и двумерном пространствах любое ограниченное множество становится измеримым (имеет длину и площадь). Таким образом, в одно- и двумерном случаях равномерное распределение можно задать на любом (ограниченном) множестве, если от вероятности требовать только аддитивность. Приведите пример, показывающий, что в трехмерном пространстве это сделать невозможно. 
\end{problem}
\begin{remark}
Базируясь на аксиоме выбора, шар в трехмерном пространстве допускает такое разбиение на конечное число непересекающихся множеств, из которых можно составить передвижением (как твердых тел = перенос + поворот) два шара того же радиуса (см., например,  Босс~В. Т.~6, 12, 16 \cite{2013}).
\end{remark}

\begin{problem}
Найдите вероятность $q_n$ того, что случайная $(0,1)$-матрица размера $n\times n$ является невырожденной над полем 
$GF_2=\{ 0,1\}$. Доказать, что существует $\lim\limits_{n\to\infty} q_n=q>0$. 
\end{problem}

\begin{problem}\Star(Закон Вигнера \cite{7}.) 
Пусть $\xi _{ij}^{n} $, $1\le i,j\le n$, $n=1,2,\ldots $ -- совокупность одинаково распределенных с.в., удовлетворяющих условиям:

\begin{enumerate}
\item \label{zero}  при каждом $n$ с.в. $\xi _{ij}^{n} $, $1\le i\le j\le n$ независимы;

\item \label{first}  матрица $\left(A^{n} \right)_{ij} =\frac{1}{2\sqrt{n} } \left(\xi _{ij}^{n} \right)$ симметричная, т.е. $\xi _{ij}^{n} =\xi _{ji}^{n} $;

\item \label{second} с.в. $\xi _{ij}^{n} $ имеет симметричное распределение, т.е. для всякого борелевского множества $B\in B({\mathbb R})$ выполнено равенство $\PR\left(\xi _{ij}^{n} \in B\right)=$\linebreak $=\PR\left(\xi _{ij}^{n} \in -B\right)$;

\item \label{third} все моменты с.в. $\xi _{ij}^{n} $ конечны, т.е. $\Exp\left[\left(\xi _{ij}^{n} \right)^{k} \right]<\infty $ при всех $k\ge~1$, причем дисперсия равна единице: $\Var\left[\xi _{ij}^{n} \right]=1$.
\end{enumerate}

Рассмотрим дискретную вероятностную меру $\mu ^{n} $ собственных значений $\lambda _{1}^{(n)} ,\ldots ,\lambda _{n}^{(n)} $ случайной матрицы $A^{n} $: для всякого борелевского множества $B\in B({\mathbb R})$:
\begin{center}
$\mu ^{n} (B)=\frac{1}{n} \sum _{i=1}^{n}I\left\{\lambda _{i}^{(n)} \in B\right\} $. 
\end{center}
Ясно, что такая мера сама по себе случайна, так как зависит от собственных значений случайной матрицы. Пусть $L_{k}^{n} $ -- $k$-й момент меры $\mu ^{n} $ матрицы $A^{n} $, т.е. 
\begin{center}
$L_{k}^{n} =\int _{-\infty }^{\infty }x^{k} \mu ^{n} (dx) =\frac{1}{n} \sum _{i=1}^{n}\left(\lambda _{i}^{(n)} \right)^{k}  $ (также с.в.).
\end{center}
Докажите, что $\mathop{\lim }\limits_{n\to \infty } \Exp\left[L_{k}^{n} \right]=m_{k} $, $\mathop{\lim }\limits_{n\to \infty } \Var\left[L_{k}^{n} \right]=0$, где
\begin{center}
 $m_{k} =\frac{2}{\pi } \int _{-1}^{1}\lambda ^{k} \sqrt{1-\lambda ^{2} } d\lambda  $,
\end{center}
то есть случайные меры собственных значений в некотором смысле сходятся к неслучайной мере на действительной прямой с плотностью, задаваемой \textit{полукруговым законом Вигнера}: 

\[p(\lambda )=\left\{\begin{array}{cc} {\frac{2}{\pi } \sqrt{1-\lambda ^{2} } d\lambda , } & {-1\le \lambda \le 1;} \\ {0} & {\text{иначе.}} \end{array}\right. \] 

\noindent Применив неравенство Чебышёва, покажите, что $L_{k}^{n} \mathop{\to }\limits_{}^{p} m_{k}$, $n\to \infty$.

\end{problem}

\begin{remark}
Если $\eta_{n} $ -- последовательность мер, моменты которых сходятся к соответствующим моментам меры $\eta $, то при дополнительных условиях на рост моментов мер $\eta_{n} $ (см. задачу \ref{moments_conv}) сами эти меры слабо сходятся к мере $\eta $.

См. по данной тематике Tao T. Topics in random matrix theory.  V. 132. American Mathematical Soc., 2012;  Tropp J.A. An introduction to matrix concentration inequalities //  arXiv:1501.01571 (2015) (приложения); van Handel R. Structured random matrices // arXiv:1610.05200 (2016),  Ver\-shy\-nin~R. High Dimensional Probability for Mathematicians and Data Scientists. \& Four lectures on probabilistic methods for data science //  arXiv:1612.06661 (2016).

\end{remark}

\begin{ordre} (См. Коралов--Синай \cite{7}.) Покажите, что 

\[\Exp\left[L_{k}^{n} \right]=\frac{1}{n} \Exp\left[\sum _{i=1}^{n}\left(\lambda _{i}^{(n)} \right)^{k}  \right]=\frac{1}{n} \Exp\left[\mathrm{tr}\left[\left(A^{n} \right)^{k} \right]\right] = \]\[ =\frac{1}{n} \left(\frac{1}{2\sqrt{n} } \right)^{k} \Exp\left[\sum _{i_{1} ,\ldots ,i_{k} =1}^{n}\xi _{i_{1} i_{2} }^{n} \xi _{i_{2} i_{3} }^{n} \cdots \xi _{i_{k} i_{1} }^{n}  \right],\] 
\noindent где $\lambda _{i}^{(n)} $ -- собственные значения матрицы $A^{n} $ и $\mathrm{tr}\left[\left(A^{n} \right)^{k} \right]$ -- след ее $k$-й степени.
С учетом условий \[\Exp\left[\xi _{i_{1} i_{2} }^{n} \xi _{i_{2} i_{3} }^{n} \cdots \xi _{i_{k} i_{1} }^{n} \right]=\prod_{J_k} \Exp\left[\left(\xi _{ij}^{n} \right)^{p(i.j)} \right], \]\[J_k = \{ i,j: i \leq j, \;
\sum p(i,j) = k\}, \] причем в силу условия \ref{second} нечетные моменты равны нулю: $$\Exp\left[\left(\xi _{ij}^{n} \right)^{p(i.j)} \right]=0,$$ если $p(i,j)$ -- нечетное. Таким образом, для четных моментов $k=2s$ имеем

\[\Exp\left[L_{2r}^{n} \right]=\frac{1}{2^{2r} n^{r+1} } \Exp\left[\sum _{i_{1} ,\ldots ,i_{2r} =1}^{n}\xi _{i_{1} i_{2} }^{n} \xi _{i_{2} i_{3} }^{n} \cdots \xi _{i_{2r} i_{1} }^{n}  \right]= \]\[ =\frac{1}{2^{2r} n^{r+1} } \sum _{i_{1} ,\ldots ,i_{2r} =1}^{n}\prod_{J_r} \Exp\left[\left(\xi _{ij}^{n} \right)^{2p(i.j)} \right]  .\] 

 Задачу можно свести к комбинаторному подсчету соответствующих путей (сопоставленных ненулевым слагаемым $\xi _{i_{1} i_{2} }^{n} \xi _{i_{2} i_{3} }^{n}, \ldots, \xi _{i_{2r} i_{1} }^{n} $) на множестве $\left\{1,2,\ldots ,n\right\}$, где каждое ребро (петли при этом не запрещаются) проходится четное число раз (без учета направления). Совокупность таких путей можно разбить на два класса: пути, накрывающие дерево из $r$ ребер, каждое из которых проходится дважды, и остальные пути, для которых либо есть петли либо циклы, либо есть ребра, через которые путь проходит по крайней мере четыре раза. Для путей второго класса (с учетом условия \ref{third}) получите оценку: \[
 \Exp\left[
    \sum  \limits_{ i_{1} < \ldots < i_{2r} } \xi _{i_{1} i_{2} }^{n} \xi _{i_{2} i_{3} }^{n} \cdots           \xi _{i_{2r} i_{1} }^{n}  
\right]\le C_{r} n^{r} ,
 \] \noindent где $C_{r} $ -- некоторая константа.

Как следствие, вклад таких путей в искомое математическое ожидание стремится к нулю при $n\to \infty $: \[\mathop{\lim }\limits_{n\to \infty } \frac{1}{2^{2r} n^{r+1} } C_{r} n^{r} k_{r} =0.\] Для путей первого класса в силу условия \ref{third} на дисперсию получите $\Exp \left(\xi _{i_{1} i_{2} }^{n} \xi _{i_{2} i_{3} }^{n} \cdots \xi _{i_{2r} i_{1} }^{n}\right) =1$. Покажите, что число путей первого класса равно \[\frac{n!(2r)!}{(n-r-1)!r!(r+1)!}. \] Для этого сопоставьте таким путям неотрицательные траектории одномерного симметричного случайного блуждания: $\left(\omega _{0} \omega _{1} \ldots \omega _{2r} \right)$, где $\omega _{0} =\omega _{2r} =0$, $\omega _{i} \ge 0$ при всех $i=1,\ldots ,2r$. Каждой фиксированной такой траектории $\left(\omega _{0} \omega _{1} \ldots \omega _{2r} \right)$ соответствует $n(n-1)\ldots (n-r)$ допустимых путей первого класса. Число неотрицательных траекторий с начальной и конечной нулевой точкой равно $\frac{(2r)!}{r!(r+1)!} $ (см. задачу \ref{sec:katalan} из раздела \ref{genF}).
Итак, 

\[\mathop{\lim }\limits_{n\to \infty } \Exp\left[L_{2r}^{n} \right]=\mathop{\lim }\limits_{n\to \infty } \frac{1}{2^{2r} n^{r+1} } \frac{n!(2r)!}{(n-r-1)!r!(r+1)!} =\frac{(2r)!}{2^{2r} r!(r+1)!} =m_{2r} .\] 
Утверждение касательно дисперсии $L_{k}^{n} $ доказывается аналогично.

\end{ordre}

\begin{problem}
На автобусных остановках обычно указывается интервал движения автобуса $m$, т.е. среднее время между двумя последовательными прибытиями автобусов. Предположим, что эти интервалы независимы и одинаково распределены с математическим ожиданием $m$ и дисперсией $\sigma^2$. Покажите, что если пассажир приходит на остановку в случайный момент времени, то математическое ожидание  времени пребывания на остановке до прихода  автобуса будет равно $\left(m^2+\sigma^2\right)/2m$.

Если известна функция распределения с.в., равной интервалу между последовательными прибытиями автобусов, то как стоит автобусу оптимально выбрать продолжительность остановки? Оптимальность необходимо понимать в смысле минимизации среднего времени ожидания автобуса. Решите задачу для случая показательного распределения.
\end{problem}
\begin{remark}
См. [\ref{sekei}], аналогично для следующей задачи.
\end{remark}

\begin{problem}
Продавец газет заказывает ежедневно $N$ газет. С каждой проданной газеты он получает прибыль $b$ и теряет $c$ на каждой газете, оставшейся непроданной.  Каким нужно выбрать  $N$, чтобы максимизировать ожидаемую прибыль, если число покупателей за день имеет распределение $\Po(\lambda)$? 
\end{problem}

\section{Производящие и характеристические \\ функции}
\label{genF}

Помимо источников литературы, указанных в тексте задач, рекомендуется ознакомиться со следующими книгами:
\begin{itemize}
 \item Айгнер М. Комбинаторная теория.  М.: Мир, 1982.  561 с.
 \item Гульден Я., Джексон Д. Перечислительная комбинаторика.  М.: Наука, 1990.  504 с.
 \item Стенли Р. Перечислительная комбинаторика.   М.: Мир, 1990.  440~с.
  \item Рыбников К.А. Комбинаторный анализ. Очерк истории.  М.: \linebreakМехМат, 1996.  125 с.
  \item Леонтьев В.К. Комбинаторика и информация. Часть 1. Комбинаторный анализ.  М.: МФТИ, 2015.  174 с.
\end{itemize}

\begin{problem}(Счастливые билеты.)
Трамвайные билеты имеют шестизначные номера. Билет называют счастливым, если 
сумма его первых трех цифр равна сумме трех последних. Вычислите приближенно 
вероятность того, что Вам достанется счастливый билет (предполагается, что 
появление билетов равновероятно).

\end{problem}

\begin{ordre}
Покажите, что число счастливых билетов совпадает с числом билетов, 
сумма цифр у которых равна 27. Запишите \textit{производящую функцию} (\textit{ПФ}) для 
последовательности $\left\{ {a_k } \right\}_{k=0}^{54} $, где $a_k $ -- число 
билетов с суммой цифр, равной $k$. Для нахождения коэффициента ПФ $a_{27} $ 
воспользуйтесь \textit{теоремой Коши} (ТФКП). Для оценки полученного интеграла примените 
\textit{метод стационарной фазы}. (см. книгу Федорюк~М.В. Метод перевала.  М.: Наука, 1977, а также книгу \cite{lando}. Последнюю книгу можно также рекомендовать и для решения последующих двух задач).
\end{ordre}

\begin{problem}(Числа Каталана.)
\label{sec:katalan}
Пусть в очереди в столовой МФТИ за булочками по цене 10 рублей 50 копеек стоят $2n$ студентов. Пусть у $n$ человек нет пятидесятикопеечной монеты, но есть рублевая, 
а у $n$ человек  есть монета в 50 копеек. Пусть изначально касса пуста. 
Найдите вероятность события, что никто из студентов в очереди не будет ждать 
свою сдачу. Считайте, что все способы расстановки студентов в очереди равновероятны.
\end{problem}

\begin{ordre}
Задачу можно проинтерпретировать в терминах правильных скобочных структур, описываемых числами Каталана. Если обозначить левую скобку буквой $a$, а правую -- $b$, то можно переписать правильные скобочные структуры в виде «слов» в алфавите $\left\{a,b\right\}$ (\textit{язык Дика}). Несложно показать, что «некоммутативный производящий ряд», перечисляющий слова языка (этот ряд представляет собой просто формальную сумму всех слов языка, вкдючая пустое слово $\lambda$, выписанных в порядке возрастания длины):

\[D(a,b)=\lambda +ab+aabb+abab+aaabbb+aababb+\ldots \] 
удовлетворяет уравнению

\[D(a,b)=\lambda +aD(a,b)bD(a,b).\] 

Перейдите от некоммутативного производящего ряда к обычному, сделав подстановку $a=x$, $b=x$, $\lambda =x^{0} =1$.

Можно искать ПФ языка Дика с помощью \textit{теоремы Лагранжа}, связывающую ее с ПФ подъязыка его неразложимых слов.

\textbf{Определение. }Слово $w=\beta _{1} \ldots \beta _{m} $ языка $L$ называется \textit{неразложимым} в $L$, если никакое его непустое подслово 
$$\beta _{i}\beta _{i+1}\dots \beta _{i+l} , 1\le i,\; i+l\le m,\; l\ge 0,$$ отличное от самого слова $w$, не принадлежит языку $L$.

В частности, пустое слово в любом языке, содержащем его, неразложимо.

Несложно проверить, что язык Дика удовлетворяет нижеперечисленным свойствам:

1) пустое слово входит в язык $L$;

2) начало всякого неразложимого слова не совпадает с концом другого или того же самого неразложимого слова;

3) если между любыми двумя буквами любого слова языка $L$ вставить слово языка $L$, то получится слово языка $L$;

4) если из любого слова языка $L$ выкинуть подслово, входящее в язык $L$, то получится слово языка $L$.

Обозначим через $n(y)=n_{0} +n_{1} y+n_{2} y^{2} +\ldots $ ПФ для числа неразложимых слов языка $L$.

\textbf{Теорема.} {\itПФ $l(x)$ для языка $L$, удовлетворяющего свойствам {\rm 1)--4)}, и ПФ $n(x)$ для подъязыка неразложимых слов в нем связаны между собой} \textbf{\textit{уравнением Лагранжа}}:

$$l(x)=n\left(xl(x)\right).$$

Приведем уравнение к классическому виду. Положим $xl(x)=$\linebreak $=\tilde{l}(x)$. Тогда уравнение Лагранжа примет вид

$$\tilde{l}(x)=xn(\tilde{l}(x)).$$

Неразложимые слова в языке Дика -- это $\lambda $ и $ab$. Отсюда немедленно получаем уравнение $l(x)=1+\left(xl(x)\right)^{2} $ на ПФ для языка Дика.

\textbf{Замечание:} Уравнение Лагранжа -- функциональное уравнение, связывающее между собой  ПФ для числа слов в языке и числа неразложимых слов в нем. Оказывается, если одна из функций известна, то оно всегда разрешимо (см.  указания к следующей задаче об уточнении приведенной здесь теоремы Лагранжа).

\end{ordre}

\begin{problem}(Остовные деревья в полном графе.) 
Пусть имеется полный граф с $n$ вершинами $\{ 1,2,\ldots ,n\} $. Каждое из $C_n^2=n(n-1)/2 $ ребер графа с вероятностью $1/2 $ удаляется. Найдите вероятность того, что полученный после удаления ребер граф будет остовным деревом.

\end{problem}

\begin{ordre}

Обозначим $t_{n} $ -- число остовных деревьев на пронумерованных вершинах $\{ 1,2,\ldots ,n\} $. Ясно, что искомая в задаче вероятность есть $t_n2^{-C_n^2}$.

Выделим одну вершину и посмотрим на те связные компоненты или блоки, на которое разобьется остовное дерево, если проигнорировать все ребра, проходящие через выделенную вершину. Если невыделенные вершины образуют $m$ компонент размеров $k_{1} ,k_{2} ,\ldots ,k_{m} $, то их можно соединить с выделенной вершиной $k_{1} k_{2} \cdot\ldots\cdot k_{m} $ способами.

Такие рассуждения приводят к рекуррентному соотношению

$$ t_{n} = \sum _{m>0}\frac{1}{m!}  \; \sum _{\sum_{i=1}^m k_i = n-1}\left(\begin{array}{c} {n-1} \\ {k_{1} ,k_{2} ,\ldots ,k_{m} } \end{array}\right) \, k_{1} k_{2} \cdot\ldots\cdot k_{m} \, t_{k_{1} } t_{k_{2} }\cdot\ldots\cdot t_{k_{m} } $$

\noindent при любом $n>1$. 

Теперь обозначим $u_{n} =nt_{n} $, тогда рекуррентное соотношение примет следующий вид:

$$ \frac{u_{n} }{n!} =\sum _{m>0}\frac{1}{m!}  \; \sum _{\sum_{i=1}^m k_i = n-1}\frac{u_{k_{1} } }{k_{1} !} \frac{u_{k_{2} } }{k_{2} !} \cdots \frac{u_{k_{m} } }{k_{m} !}  ,\quad n>1.$$

Обозначим через $U(x)$ \textit{экспоненциальную производящую функцию} ({\itЭПФ}) для последовательности $\left\{u_{n} \right\}$ (то есть $U(x)=\sum _{n=0}^{\infty }\frac{u_{n} }{n!} x^{n}  $). Таким образом,

$$U(x)=xe^{U(x)}.$$

Для нахождения явной формулы для этой последовательности можно воспользоваться следующим уточнением \textit{теоремы Лагранжа}.

\textbf{Теорема.} Пусть функции $\varphi =\varphi (x)$ ($\varphi (0)=0$) и $\psi =\psi (z)$ связаны между собой уравнением Лагранжа: 

$$\varphi (x)=x\psi \left(\varphi (x)\right).$$ 

Тогда коэффициент при $x^{n} $ в функции $\varphi $ равен коэффициенту при $z^{n-1} $ в разложении $\frac{1}{n} \psi ^{n} (z)$: $[x^{n}]\varphi(x) = [z^{n-1}]\frac{1}{n} \psi ^{n} (z)$.

Интересный пример использования аппарата производящих функций был описан в лекции Станислава Смирнова ``Сколько простых кривых на решетке?'', прочитанной в ЛШСМ-2016 (видео доступно на сайте mathnet.ru).

\end{ordre}

\begin{problem}(Задача с марками / coupon collector's problem \cite{13}.)
\label{laplas}
Пусть Вы хотите собрать коллекцию из $N$ 
марок. Для этого Вы каждый день покупаете конверт со случайной маркой 
(марки появляются на купленных конвертах равновероятно).
\begin{enumerate}
\item Введем дискретную с.в. $X$, равную номеру впервые купленной Вами 
повторной марки. Найдите математическое ожидание с.в. $X$.

\item Покажите, что распределение случайной величины $X$ имеет асимптотически 
\textit{ распределение Релея } при $n=t\sqrt N $ ($N\gg 1)$:
$$ \PR \left ( {X>t\sqrt N } \right ) \sim e^{-\frac{t^2}{2}},
\quad \PR\left ( {X=t\sqrt N } \right ) \sim \frac{1}{\sqrt N }te^{-\frac{t^2}{2}}.$$
\item Определите математическое ожидание номера купленной Вами марки, которая 
станет недостающей в собранной Вами коллекции марок. Решение этой задачи используется в Compressed Sensing (см. Часть 2).
\end{enumerate}
\end{problem}

\begin{ordre}

а) Покажите, что  
$$\PR \left ( X>n \right ) =\frac{n!}{N^n} \left[ {z^n} \right]\left( {1+z} 
\right)^N=n! \; \left[ {z^n} \right] \left( {1+\frac{z}{N}} \right)^N,$$
$$\Exp X
=
\sum\limits_{n=0}^{\infty} 
n! \; \left[ z^n \right] \left( 1+\frac{z}{N} \right)^N
=
\int\limits_0^{\infty} e^{-t} \left( 1+\frac{t}{N} \right)^N dt. $$

Далее для приближенного вычисления интеграла  
$$\int\limits_0^\infty 
{e^{-t}\left( {1+\frac{t}{N}} \right)^Ndt} =N\int\limits_0^\infty 
{e^{N\left( {\ln (1+u)-u} \right)}du} $$
 воспользуемся \textit{методом Лапласа}:
$$
 \int\limits_0^\infty {f(u)e^{NS(u)}du} \approx f(u_0 ) e^{NS(u_0 
)}\int\limits_{u_0 -\delta }^{u_0 +\delta } {e^{\frac{NS''(u_0 )(u-u_0 
)^2}{2}}} du\approx $$
$$\approx \sqrt {2\pi } \frac{f(u_0 )e^{NS(u_0 )}}{\sqrt {-NS''(u_0 )} }, $$
где $u_0 $ -- единственная точка максимума вещественнозначной функции 
$S(u)$ на полубесконечном интервале $(0, +\infty )$. Основная идея 
асимптотического представления интеграла Лапласа заключается в представлении 
функции $S(u)$ в окрестности точки максимума $u_0 $ в виде ряда Тейлора.

б) Согласно \textit{теореме Коши} (из курса ТФКП): 
$$ \PR(X>n)=n!\left[ {z^n} \right]\left( {1+\frac{z}{N}} 
\right)^N=\frac{1}{n!}\frac{1}{2\pi i }\oint\limits_{\vert z\vert =\rho } 
{\left( {1+\frac{z}{N}} \right)^N\frac{dz}{z^{n+1}}}. $$
Перейдя к полярным координатам, получите, что
$$
\frac{1}{2i\pi }\oint\limits_{\vert z\vert =\rho } {\left( {1+\frac{z}{N}} 
\right)^N \frac{dz}{z^{n+1}}}  =\frac{1}{2\pi }\int\limits_{-\pi }^\pi {\left( 
{1+\frac{\rho e^{i\theta }}{N}} \right)^N\frac{d\theta }{\rho ^ne^{in\theta 
}}}= $$
$$=\frac{1}{2\pi }\int\limits_{-\pi }^\pi {e^{f(\rho e^{i\theta })}d\theta 
} ,
$$
где $f(z)=N\ln \left( {1+\frac{z}{N}} \right)-n\ln z$, $z=\rho 
e^{i\theta }$.

Воспользуйтесь \textit{методом перевала}: разложите в ряд Тейлора функцию $f$ в окрестности седловой 
точки: $$f(\rho e^{i\theta })=f(\rho )-\frac{1}{2}\beta (\rho )\theta 
^2+O(\theta ^3)$$ для $\left| \theta \right|<\delta $ (разложение Тейлора), 
где $$\beta (\rho )=\rho ^2\left. {\left[ {\left( {\frac{d}{dz}} 
\right)^2f(z)} \right]} \right|_{z=\rho }. $$

Замените интеграл по всей окружности $\vert z\vert =\rho $ на интеграл по ее 
части: 
$$
\int\limits_{-\delta }^\delta {e^{f(\rho e^{i\theta })}d\theta } 
\approx e^{f(\rho )}\int\limits_{-\delta }^\delta {e^{-\frac{1}{2}\beta 
(\rho )\theta ^2}d\theta }= $$
$$=\frac{e^{f(\rho )}}{\sqrt {\beta (\delta )} 
}\int\limits_{-\delta \sqrt {\beta (\rho )} }^{\delta \sqrt {\beta (\rho )} 
} {e^{-\frac{1}{2}u^2} du}\mathop \to \limits_{\beta (\rho )\to \infty } 
\frac{e^{f(\rho )}}{\sqrt {\beta (\delta )} }\int\limits_{-\infty }^\infty 
{e^{-\frac{1}{2}u^2}du} = $$
$$=\sqrt {2\pi }\frac{e^{f(\rho )}}{\sqrt {\beta 
(\delta )} }.
$$

в) Пусть $Y$ -- дискретная с.в., равная номеру купленной Вами 
марки, которая станет недостающей в собранной Вами коллекции марок.

Тогда $$\PR \left ( Y\le n \right ) =\frac{n!}{N^n}\left[ {z^n} \right](e^z-1)^N=n!\left[ 
{z^n} \right](e^{\frac{z}{N}}-1)^N,$$ а следовательно, 
$$\Exp Y=\sum\limits_{n=0}^\infty {n!\left[ {z^n} \right]\left( 
{e^z-(e^{\frac{z}{N}}-1)^N} \right)} =\int\limits_0^\infty {\left[ 
{1-(1-e^{-\frac{t}{N}})^N} \right]} dt.$$
Сделав замену переменных 
$y=1-e^{-\frac{t}{N}}$, получите, что $$\Exp Y=N\left( {1+\frac{1}{2}+\cdots 
+\frac{1}{N}} \right). $$

Заметим, что этот пункт задачи можно также решить с помощью формулы полной вероятности, составив рекуррентные формулы для $x_k$ -- среднего числа купленных марок, если до получения полной коллекция марок не хватает ровно  $k$ марок.

\end{ordre}

\begin{problem}(Урновая схема \cite{13}.) 
Рассмотрим случайное размещение $n$ различных шаров по $m$ различным урнам. 
Пусть случайные величины MIN, MAX -- размер наименее или наиболее заполненной 
урны в случайном размещении. Получите функции распределения этих случайных 
величин, а именно: $$\PR \left ( {\rm MAX}\le l\right ) =\frac{n!\left[ {z^n} \right]e_l 
(z)^m}{m^n}=n!\left[ {z^n} \right]e_l \left( {\frac{z}{m}} \right)^m,$$

$$\PR\left ( {\rm MIN}>l\right ) =n!\left[ {z^n} \right]\left( {e^{\frac{z}{m}}-e_l \left( 
{\frac{z}{m}} \right)} \right)^m,
$$
где $e_l (z)=1+z+\frac{z^2}{2!}+\cdots 
+\frac{z^l}{l!}$.

\end{problem}

\begin{problem}(Задача о циклах в случайной перестановке \cite{13}.)
\label{cycle}

\begin{enumerate}
\item Найдите математическое ожидание числа циклов длины $r$ в 
случайной перестановке длины $n$.

\item Найдите математическое ожидание числа циклов в 
случайной перестановке длины $n$.

\item (Сто заключенных.) 
В коридоре находятся 100 человек, у каждого свой номер (от 1 до 100). Их по одному заводят в комнату, в которой 
находится комод со 100 выдвижными ящиками. В ящики случайным образом 
разложены карточки с номерами (от 1 до 100). Каждому разрешается заглянуть в 
не более чем 50 ящиков. Цель каждого -- определить, в каком ящике находится 
его номер. Общаться и передавать друг другу информацию запрещается. 
Предложите стратегию, которая с вероятностью, не меньшей $0.3$ (в 
предположении, что все $100!$ способов распределения карточек по ящикам 
равновероятны), приведет к выигрышу всей команды. Команда выигрывает, если 
каждый из 100 участников верно определил ящик с карточкой своего номера.

\end{enumerate}
\end{problem}
\begin{ordre}
 ЭПФ для последовательностей: 
 
--  
числа циклов длины $n$:
\[C(z)=\sum _{n=0}^{\infty }\frac{(n-1)!}{n!} z^{n}  =\log \frac{1}{1-z}; \]

-- числа множеств на $n$ элементах:
\[S(z)=\sum _{n=0}^{\infty }\frac{1}{n!} z^{n}  =e^{z}; \]

-- числа перестановок на $n$ элементах:
\[P(z)=\sum _{n=0}^{\infty }\frac{n!}{n!} z^{n}  =\frac{1}{1-z}. \]

Перестановка есть не что иное, как совокупность циклов:
\[P(z)=S\left(C(z)\right).\]

Если в такой функции ставить метку (переменную $u$) для циклов длины $r$, получим, что ``двойная''  ЭПФ, перечисляющая число циклов длины $r$ в перестановке длины $n$, равна

$$\exp \left\{(u-1)\frac{z^{r} }{r} +\log \frac{1}{1-z} \right\}=\frac{\exp \left\{(u-1)\frac{z^{r} }{r} \right\}}{1-z} .$$

Заметим, что $\left[z^{n} \right]\left. \frac{\partial }{\partial u} \left(\frac{\exp \left\{(u-1)\frac{z^{r} }{r} \right\}}{1-z} \right)\right|_{u=1} $ -- среднее (математическое ожидание) число циклов длины $r$ в перестановке длины $n$.

\textbf{Стратегия для ста заключенных.} Каждый человек вначале открывает ящик под номером, равным его собственному номеру, затем -- под номером, который указан на карточке, лежащей в ящике, и т.д. Среднее число циклов длины $r$ в случайной 
перестановке -- $1/r$. Тогда среднее число циклов длины, большей $n/2$,
есть $\sum\limits_{k=n \mathord{\left/ {\vphantom {n 2}} \right. 
\kern-\nulldelimiterspace} 2}^n {(1/k)} $. Это и есть вероятность 
существования цикла длины, большей $n/2$. Поэтому вероятность успеха команды  
есть $1-\sum\limits_{k=51}^{100} {(1/k)} \approx 0,31$ (для сравнения, если 
произвольно открывать ящики, то вероятность успеха будет 
$2^{-100}\approx 10^{-30}$). В случае, когда карточки 
разложены не случайным образом, то следует сделать случайной нумерацию 
ящиков, и далее следовать описанной выше стратегии.
\end{ordre}

\begin{problem}(Задача о беспорядках \cite{lando}.)
\label{permloop}
Группа из $n$ фанатов выигрывающей футбольной команды на радостях 
подбрасывают в воздух свои шляпы. Шляпы возвращаются в случайном порядке --
по одной к каждому болельщику. Какова вероятность того, что никому из 
фанатов не вернется своя шляпа? Найдите математическое ожидание и дисперсию числа шляп, вернувшихся к своим хозяевам.
\end{problem}

\begin{ordre}
Формально задача сводится к подсчету числа беспорядков $d_{n} $ на множестве из $n$ элементов (перестановка $\pi $ элементов множества $\left\{1,2,\ldots ,n\right\}$ называется  \textit{беспорядком}, если $\pi (k)\ne k$ ни при каких\linebreak $k=1,\ldots ,n$).

Получите ЭПФ для числа беспорядков: $$D(x)=\sum _{n=0}^{\infty }\frac{d_{n} }{n!} x^{n}  =\frac{e^{-x} }{1-x}. $$

Пусть $d_{n,k} $ -- число перестановок на множестве из $n$ элементов, оставляющих на месте ровно $k$ элементов (то есть число неподвижных точек равно $k$), тогда $d_{n,0} =d_{n} $. Более того, $d_{n,k} =C_{n}^{k} d_{n-k} $. Из правила суммы несложно получить

$$ n!=\sum _{k=0}^{n}d_{n,k}  =\sum _{k=0}^{n}C_{n}^{k} d_{n-k}  =\sum _{k=0}^{n}C_{n}^{n-k} d_{n-k}  =\sum _{k=0}^{n}C_{n}^{k} d_{k}.$$

То есть получилась биномиальная свертка двух рядов:
$$
\sum _{n=0}^{\infty }\frac{n!}{n!} x^{n}  =\sum _{n=0}^{\infty }\frac{d_{n} }{n!} x^{n}  \sum _{n=0}^{\infty }\frac{1}{n!} x^{n}.
$$
Для нахождения среднего и дисперсии числа шляп, вернувшихся к своим хозяевам, удобно воспользоваться указанием к предыдущей задаче о подсчете среднего и дисперсии циклов длины 1 в случайной перестановке длины $n$.

Сравните подход с помощью производящих функций для этой задачи с вашими решениями задач~17 и \ref{sec:latters} из раздела стандартных задач (раздел \ref{standart}).
\end{ordre}

\begin{problem}(Ожерелья.) 
Найдите вероятность того, что случайная раскраска ожерелья из $n$ бусин в $k$ различных цветов имеет ровно $m\le n$ бусин первого цвета. Ожерелья, получающиеся одно из другого с помощью плоского поворота, считаются эквивалентными. Положите $n=7$,\linebreak $k=2$, $m=3$.

\end{problem}

\begin{remark}
Ознакомиться с классическими работами по теории перечисления (Дж. Пойа, Дж.-К. Рота) можно по книге:

Перечислительные задачи комбинаторного анализа: сб. переводов / под ред. Г.П.~Гаврилова.   М.: Мир, 1979.  362~с.

Популярное изложение имеется, например, в учебном пособии: Краснов М.Л. и др. Вся высшая математика: учебник. Т.~7.  М.:\linebreak КомКнига, 2006.  208 с.
\end{remark}

\begin{ordre}

Сопоставим каждой раскраске функцию $f$ как отображение из множества пронумерованных бусин в нераскрашенном ожерелье в множество пронумерованных бусин в раскрашенном ожерелье. Основной нюанс в решении комбинаторных задач такого типа заключается в том, что некоторые функции (раскраски) оказываются эквивалентными, так как получаются одна из другой с помощью некоторой подстановки (в данном случае задающей поворот в плоскости): $f_{1} \sim f_{2} $, если найдется подстановка $g\in G$, что $f_{1} (g)=f_{2} $, т.е. $\forall\,i=1,\ldots ,n$   $f_{1} \left(g\left(d_{i} \right)\right)=f_{2} \left(d_{i} \right).$

Для решения этой задачи воспользуемся подходом Пойа.
\textit{Цикловым индексом подстановки} называют одночлен
\[x_{1} ^{k_{1} } x_{2} ^{k_{2} } \ldots x_{n} ^{k_{n} } ,\] 
где $\left(k_{1} ,\; k_{2} ,\; \ldots ,\; k_{n} \right)$ -- тип подстановки, т.е. подстановка представима в виде $k_{1} $ цикла длины 1, $k_{2} $ циклов длины 2  и т.д.

\textit{Цикловым индексом группы подстановок} $G$ называют среднее арифметическое цикловых индексов ее элементов:
\[P_{G} (x_{1} ,x_{2} ,\ldots ,x_{n} )=\frac{1}{|G|} \sum _{g\in G}x_{1} ^{k_{1} } x_{2} ^{k_{2} } \ldots x_{n} ^{k_{n} }  .\] 
Покажите, что цикловой индекс группы поворотов $C_n$ равен $$P_{G} (x_{1} ,x_{2} ,\ldots ,x_{n} )=\frac{1}{n} \sum _{j=1}^{n}\left(x_{\frac{n}{(n,j)} } \right)^{(n,j)}  ,$$ где $(n,j)$ -- наибольший общий делитель $n$ и $j$. Каждому цвету $r_{i} $ $i=1,\ldots ,k$ придадим некоторый \textit{вес} $w(r_{i} )$. \textit{Весом функции} $f$ назовем произведения весов полученной раскраски:
\[W(f)=\prod _{i=1}^{n}w(f(d_{i} )).\] 
Ясно, что эквивалентные функции имеют одинаковый вес:
\[f_{1} \sim f_{2} \Rightarrow W\left(f_{1} \right)=W\left(f_{2} \right).\] 
Весом класса эквивалентности называется вес любой функции из этого класса; если $F$ -- класс эквивалентности и $f\in F$, то $W(F)=$\linebreak $=W(f)$. Заметим, что и у неэквивалентных функций могут совпадать веса.

\textbf{Теорема Пойа (1937).}
{\itСумма весов классов эквивалентности равна}
$$\sum _{F}W(F) =P_{G} \left(\sum _{i=1}^{k}w(r_{i} ) ,\; \sum _{i=1}^{k}w^{2} (r_{i} ) ,\; \sum _{i=1}^{k}w^{3} (r_{i} ) ,\; \ldots ,\sum _{i=1}^{k}w^{n} (r_{i} ) \right),$$ {\itгде $P_{G} $ -- цикловой индекс группы подстановок $G$.}

\end{ordre}

\textbf{Следствие}. 
Число классов эквивалентности равно 
$$P_{G} \left(k,\; k,\; k,\; \ldots ,k\right).$$ 

Согласно следствию теоремы Пойа число различных ожерелий равно $$P_{G} (k,k,\ldots ,k)=\frac{1}{n} \sum _{j=1}^{n}k^{(n,j)}  =\frac{1}{n} \sum _{s|n}k^{s} \phi \left(\frac{n}{s} \right), $$ где $\phi \left(l\right)$ -- число взаимно простых делителей числа $l$, не превосходящих $l$ (\textit{функция Эйлера}), обозначение $s|n$ -- $s$ является делителем~$n$.

Для подсчета числа ожерелий с ровно $m$ бусинами первого цвета положите вес этого цвета $w(r_{1} )=x$, веса остальных цветов -- единицей. В ПФ (полученной согласно теореме Пойа) возьмите коэффициент при $x^{m} $.

\begin{problem}(Изомеры органических молекул.) 
Рассмотрим математическую модель органической молекулы: в центре тетраэдра поместим атом углерода $C$, в вершинах тетраэдра равновероятно помещаются метил $\left({\rm CH}_{3} \right)$, этил $\left({\rm C}_{2}{\rm H}_{5} \right)$, водород $\left({\rm H}\right)$ и хлор $\left({\rm Cl}\right)$. Найдите вероятность того, что случайная молекула заданной структуры окажется метаном ${\rm CH}_{4} $. 

\end{problem}

\begin{ordre}
 Воспользуйтесь подходом Пойа, описанным в указаниях к предыдущей задаче. Покажите, что цикловой индекс группы вращения тетраэдра 
 $$P_{G} (x_{1} ,x_{2} ,x_{3} ,x_{4} )=\frac{1}{12} \left(8x_{1} x_{3} +3x_{2} ^{2} +x_{1} ^{4} \right).$$
\end{ordre}

\begin{remark}
Для решения последующих трех задач рекомендуется ознакомиться с подходом Егорычева вычисления комбинаторных сумм, изложенным в следующих книгах:

Леонтьев В.К. Избранные задачи комбинаторного анализа.  М.: \linebreakМГТУ,  2001.  184~с.

Егорычев Г.П. Интегральное представление и вычисление комбинаторных сумм.  Новосибирск: Наука,  1977.  284~с.
\end{remark}

\begin{problem}
Обозначим через $E^n$ множество бинарных последовательностей длины $n$, или множество вершин 
единичного $n$-мерного куба, а через $E_k^n $ -- $k$-й слой куба $E^n$, то есть 
подмножество точек $E^n$, имеющих ровно $k$ единичных координат. Пусть 
$X= \langle x, y \rangle$ -- случайная величина, где $ 
x \in E_p^n $, $ y \in E_q^n $ -- векторы, независимые и равномерно 
распределенные на $E_p^n $ и $E_q^n $, соответственно. Обозначим 
через $a_{p,q} (k)=\PR\left( {X=k} \right)$. Доказать следующие утверждения:

\begin{enumerate}
\item $\sum\limits_{k=0}^n {a_{p,q} (k)z^k} =\frac{1}{2\pi i}C_n^p 
\oint\limits_{\left| u \right|=\rho } 
{\frac{(1+zu)^p(1+u)^{n-p}}{u^{q+1}}du}; $

\item $a_{p,q} (k)=\frac{C_p^k C_{n-p}^{q-k} }{C_n^q };$

\item $\Exp X=\frac{pq}{n};$

\item $\Var X=\frac{pq}{n(n-1)}\left( {n+\frac{pq}{n}-(p+q)} \right) .$
\end{enumerate}

\end{problem}

\begin{problem}
Пусть $\xi$ -- с.в., равномерно распределенная на множестве всех пар векторов 
$(x,y)\in \{ 0,1\}^n\otimes \{ 0,1\}^n$, равная скалярному произведению $ \langle x,y \rangle=\sum\limits_{k=1}^{n} x_k y_k$. Найдите 
$$
{\mathbb P}(\xi=k), \; {\mathbb E}\,\xi, \; \Var \xi . 
$$
\end{problem}

\begin{problem}
Пусть $\xi$ -- с.в., равномерно распределенная на множестве бинарных матриц ($\{0,1\}^{m\times n}$) порядка 
$m\times n$ и равная числу нулевых столбцов матрицы. Доказать, что 
\begin{enumerate}
\item $P_k(m,n)={\mathbb P}(\xi=k)=C_n^k \bigl( 2^m -1\bigr)^{n-k} 2^{-mn}$; 

\item ${\mathbb E}\xi= n2^{-m}$; 
\item если $2^m-1=\alpha n$, где $\alpha$ не зависит от $n$, то 
$$
\lim\limits_{n\to\infty} P_k(m,n)=e^{-\lambda}\, \frac{\lambda^k}{k!} ,\; \text{ где } \lambda=\alpha^{-1} ;
$$

г) Что можно сказать об этом предельном распределении?

\end{enumerate}
\end{problem}

\begin{problem}
Может ли функция $\varphi(t)=\begin{cases}1,\quad t\in[-T,T]\\
0,\quad t\notin[-T,T] \end{cases}$  
быть \textit {характеристической функцией} (х.ф.) некоторой с.в.? Изменится ли ответ, если <<чуть-чуть>> размазать (сгладить) разрывы функции 
$\varphi(t)$ в точках $t=\pm T$? 
\end{problem}

\begin{problem}
Будут ли функции~$\cos(t^2)$, $\exp(-t^4)$, $\arcsin((\cos t)/2) / \arcsin(1/2)$ характеристическими для каких-нибудь случайных величин?
\end{problem}

\begin{problem}

Докажите, что выпуклая линейная комбинация х.ф. есть х.ф.  ``смеси'' слагаемых. То есть если $\phi_{k}(t) = \Exp \left[\exp (i\xi_k t)\right]$, $\sum_k a_k = 1$, $a_k \geq 0$,  то $\sum_k a_k \phi_{k}(t)$ есть х.ф. случайной величины
\[
\xi = \sum_k [\zeta = k] \xi_k,
\]
где $\zeta$ -- не зависящая от $\{\xi_k\}$ c.в. с распределением $\PR(\zeta=k) = a_k$. Здесь $[\zeta = k]$ есть индикаторная функция события $\zeta = k$.
\end{problem}

\begin{problem}
Случайные величины $X_1, X_2$ независимы и имеют распределение Коши (Cauchy) $\text{Ca}(x_1, d_1), \; \text{Ca}(x_2, d_2)$, где плотность распределения  $\text{Ca}(x_1, d_1)$ определяется как    
\[
f(x) = \frac{d_1}{\pi(d_1^2 + (x-x_1)^2)}, x\in \mathbb{R}^1
\]
Докажите, что распределением $X_1 + X_2$ является $\text{Ca}(x_1+x_2, d_1+d_2)$.
\end{problem}

\begin{ordre}
Убедитесь, что характеристическая функция $X_1$ имеет вид $\varphi(t) = e^{ix_1t - d_1|t|}$.
\end{ordre}

\begin{problem}(Квадратичные формы.)
Пусть $A$ -- симметричная ($n \times n$)-мат\-рица, а $X$ -- $n$-мерный нормально распределенный вектор, $X \in \N(0, S)$. Покажите, что для случайной квадратичной формы $Q = X^TAX$ справедливы равенства 
\[
\Exp e^{tQ} =\det (I - 2t A S)^{-1/2},
\] 
\[
\Exp Q = \text{tr}(AS), \quad \Var Q = 2 \text{tr}((AS)^2),
\]
\[
\Exp [ (WX +m)^{T} A (WX +m) ] = m^T A m + \text{tr}(W^TAWS),
\]
где $W$ и $m$ -- неслучайные матрицы.
\end{problem}

\begin{problem}
Пусть $\varphi_{\xi}$~-- х.ф. абсолютно непрерывной случайной величины~$\xi$ с~плотностью $p_{\xi}$. Рассмотрим $f_1 =\text{Re}(\varphi_{\xi})$ и~$f_2 = \text{Im}(\varphi_{\xi}) $. Существуют ли случайные величины $\eta_1,\,\eta_2$, для которых $f_1,\,f_2$~являются их х.ф.? 
\end{problem}

\begin{problem}
Реализуем $m$~раз схему из $n$~испытаний Бернулли с вероятностью успеха~$p$. Считаем, что все реализации схем взаимно независимы. На выходе получим $m$~случайных векторов $x_1,\,\dots,\,x_m$ с координатами $0$~и~$1$. Некоторые из этих векторов могут совпадать. Скажем, что векторы $x_i,\,x_j,\,x_k$ образуют прямой угол с вершиной в~$x_k$, если скалярное произведение ~$\langle x_i-x_k,\,x_j - x_k\rangle = 0$ (помимо обычных прямых углов, под это определение попадают и <<вырожденные>>, т.\,е. образованные совпадающими векторами). Найдите математическое ожидание числа прямых углов во множестве~$\{x_1,\,\dots,\,x_m\}$.
\end{problem}

\begin{problem} 
Найдите вероятность того, что пара случайно выбранных из $\{ 0,1\}^n$ векторов является ортогональной 
\begin{enumerate}
\item над полем $\mathbb{F}_2=\{ 0,1\}$; 

\item над полем действительных чисел. 
\end{enumerate}
\end{problem}

\begin{problem}(Об оценке хвостов.)
Пусть $X$ -- неотрицательная целочисленная с.в., а $\Psi (z)=\Exp z^X$ -- её ПФ. Докажите, что
\[
{\rm\PR}(X\le x)\le x^{-r}\Psi (x)\mbox{ для }0<x\le 1;
\]
\[
{\rm \PR}(X\ge x)\le x^{-r}\Psi (x)\mbox{ для }x\ge 1.
\]
\begin{remark}
Для решения этой и последующих шести задач можно рекомендовать книгу \cite{29}.
\end{remark}
\end{problem}

\begin{problem}(Игра У. Пенни, 1969.)
Авдотья и Евлампий играют в игру: они 
бросают монету до тех пор, пока не встретится РРО или РОО. Если первой 
появится последовательность РРО, выигрывает Авдотья, если РОО -- Евлампий. Будет ли игра честной (одинаковы ли вероятности выигрыша у обоих игроков)?
\end{problem}

\begin{ordre}

Справедливы следующие равенства:
\[
\mbox{1+N(Р+О)=N+А+В},
\]
\[
\mbox{NРРО=А},
\]
\[
\mbox{NРОО=В+АО},
\]
где $A$ -- конфигурации, выигрышные для Авдотьи, $B$ --  конфигурации, выигрышные для Евлампия, $N$ -- конфигурации последовательностей, для которых ни один из игроков не выиграл.
Заменяя Р и О на $1/2$, найдите значения $A$ и $B$ -- вероятностей выигрыша Авдотьи и Евлампия.
\end{ordre}

\begin{problem}
Теперь трое игроков: Авдотья, Евлампий и Компьютер. Играют пока 
не выпадет одна из следующих последовательностей: $A\!=\,$РРОР, $B\!=\,$РОРР, $C\!=\,$ОРРР. Каковы шансы каждого выиграть?
\end{problem}

\begin{problem}
Рассмотрим следующую игру: первый игрок называет одну из восьми комбинаций ООО, ООР, ОРО, РОО, ОРР, РОР, РРО, РРР, второй игрок называет другую, потом они бросают симметричную монету до тех пор, пока в последовательности орлов (О) и решек (Р) -- результатов бросания -- не появится одна из выбранных комбинаций. Тот, чья комбинация появится, выиграл. Будет ли игра честной? Найдите для каждой комбинации, выбранной первым игроком, выгодные комбинации второго игрока.
\end{problem}

\begin{ordre} 
Для решения этой задачи, см. также \cite{book12}. Например, если первый игрок выбрал ОРО, то выбираю ООР, второй игрок выигрывает с вероятностью ${2\mathord{\left/ {\vphantom {2 3}} \right. \kern-\nulldelimiterspace} 3} $. Действительно, пусть $p_{OO} $, $p_{OP} $, $p_{PO} $, $p_{PP} $ -- вероятности выигрыша первого игрока при последних значениях ОО, ОР, РО, РР (в предположении, что раньше никто не выиграл). Их можно найти из системы уравнений:

\[\left\{\begin{array}{l} {p_{OO} =\frac{1}{2} p_{OO} +\frac{1}{2} \cdot 0,} \\ \\ {p_{OP} =\frac{1}{2} \cdot 1+\frac{1}{2} p_{PP} ,} \\ \\ {p_{PO} =\frac{1}{2} p_{OO} +\frac{1}{2} p_{OP} ,} \\ \\ {p_{PP} =\frac{1}{2} p_{PO} +\frac{1}{2} p_{PP} .} \end{array}\right. \] 

Откуда получаем

\[p_{OO} =0, \; p_{OP} =\frac{2}{3}, \;  p_{PO} =\frac{1}{3}, \; p_{PP} =\frac{1}{3} .\] 

То есть общая вероятность выигрыша первого игрока равна ${1\mathord{\left/ {\vphantom {1 3}} \right. \kern-\nulldelimiterspace} 3}$.
\end{ordre}

\begin{problem}
В вершине пятиугольника $ABCDE$ 
находится яблоко, а на расстоянии двух ребер, в вершине $C$, находится 
червяк. Каждый день червяк переползает в одну из двух соседних вершин с 
равной вероятностью. Так, через один день червяк окажется в вершине $B$ или 
$D$ с вероятностью $1 \mathord{\left/ {\vphantom {1 2}} \right. 
\kern-\nulldelimiterspace} 2$. По прошествии двух дней червяк может снова 
оказаться в $C$, поскольку он не запоминает своих предыдущих положений. 
Достигнув вершины $A$, червячок останавливается пообедать.

\begin{enumerate}
\item Чему равны математическое ожидание и дисперсия числа дней, прошедших до обеда?
\item Какую оценку дает неравенство Чебышёва для вероятности  того, что это число дней будет 100 или больше?
\item Что можно сказать о величине $p$-оценки из задачи ``об оценке хвостов''.
\end{enumerate}
\end{problem}

\begin{problem}
Пять человек стоят в вершинах пятиугольника $ABCDE$ и 
бросают друг другу диски фрисби. У них имеется два диска, которые в 
начальный момент находятся в соседних вершинах. В очередной момент времени 
диски бросают либо налево, либо направо с одинаковой вероятностью. Процесс 
продолжается до тех пор, пока обе тарелки не окажутся в руках одного из игроков.

\begin{enumerate}
\item Найдите математическое ожидание и дисперсию числа пар бросков.
\item Найдите ``замкнутое'' выражение через числа Фибоначчи для вероятности того, что игра продлится более 100 шагов.
\end{enumerate}
\end{problem}

\begin{problem}
Обобщите предыдущую задачу на случай $m$-угольника и найдите 
математическое ожидание и дисперсию числа пар бросков до столкновения 
дисков. Докажите, что если $m$ нечетно, то ПФ случайной величины,  равной числу бросаний, 
представима в следующем виде:\footnote{ Воспользуйтесь подстановкой $z=1 
\mathord{\left/ {\vphantom {1 {\cos ^2\theta }}} \right. 
\kern-\nulldelimiterspace} {\cos ^2\theta }$.}
\[
G_m (z)=\prod\limits_{k=1}^{(m-1)/2} {\frac{p_k z}{1-q_k z}} ,
\]
где
\[
p_k =\sin ^2\frac{(2k-1)\pi }{2m},
\quad
q_k =\cos ^2\frac{(2k-1)\pi }{2m}.
\]
\end{problem}

\begin{problem}(Загадочный случайный суп \cite{4}.)
Студент, решивший отобедать в 
столовой, может обнаружить в своей тарелке с супом случайное число $N$ 
инородных частиц со средним $\mu $ и конечной дисперсией. С~вероятностью $p$ выбранная частица является мухой, иначе это таракан; типы 
разных частиц независимы. Пусть $X$ -- количество мух и $Y$ -- количество тараканов.

\begin{enumerate}
\item Покажите, что ПФ с.в. $X$ удовлетворяет равенству
\[
\psi _X (s)=\psi _N (ps+1-p).
\]

\item Предположим, что с.в. $N$ имеет пуассоновское 
распределение с параметром $\mu$:  $N\in \Po\left( \mu \right))$. Покажите, что $X$ имеет пуассоновское распределение с 
параметром $p\mu $, а с.в. $X$ и $Y$ независимы. Покажите, что
\[
\psi _N (s)=\psi _N^2 \Bigl( {\frac{1}{2}(1+s)} \Bigr).
\]
\end{enumerate}

\end{problem}

\begin{problem}
На каждой упаковке овсянки печатается купон одного из $k$ различных цветов. Считая, что при отдельной покупке купон каждого цвета может встретиться с равной вероятностью и различные покупки независимы, требуется найти производящую функцию распределения вероятности, математическое ожидание и дисперсию минимального числа упаковок, которые придется купить для того, чтобы собрать купоны всех $k$ цветов.
\end{problem}

\begin{ordre}
Воспользуйтесь вспомогательными вероятностями: $p(n)$ -- вероятность разделить $n$ одинаковых объектов на $k$ групп так, что в каждой группе будет не менее  одного объекта, $f(m)$ -- вероятность того, что минимальное число купленных  упаковок равно $m$.
\[
p(n) = \sum_{m = k}^{n} f(m), 
\]
\[
P(z) = \sum_{n=k}^{\infty} p(n) z^n = \sum_{n=k}^{\infty} \sum_{m = k}^{n} f(m) z^n = 
\sum_{n=m}^{\infty} \sum_{m = k}^{\infty} f(m) z^m z^{n-m} = 
\]
\[=
F(z) \frac{1}{1-z}.
\]    
\end{ordre}

\begin{remark}
См. также \cite{1} и задачу \ref{laplas} из этого раздела и указание к ней.
\end{remark}

\begin{problem}\DStar(Задача о предельной форме диаграмм Юнга.)
\label{limung}
\textit{Разбиением} $\lambda $ натурального числа $n$ называется набор натуральных 
чисел $(\lambda _1 ,\ldots ,\lambda _N )$, для которого $\lambda _1 \ge 
\ldots \ge \lambda _N >0$ и $\lambda _1 +\ldots +\lambda _N =n$. 
\textit{Диаграммой Юнга} разбиения $\lambda =(\lambda _1 ,\ldots ,\lambda _N )$ называется 
подмножество первого квадранта плоскости, состоящее из единичных 
квадратиков. Квадратики размещаются по строкам, выровненным по левому краю, 
причем число квадратиков в $i$-й строке равно $\lambda _i $.
Множество всех диаграмм Юнга, 
соответствующих натуральному числу $n$ (или иначе множество всех разбиений 
натурального числа $n$) обозначим через $\mathcal P_n $.
Введем равномерную меру $\mu _n $ на $\mathcal P_n$, то есть $$\mu _n (\lambda 
)=\frac{1}{p(n)}$$ для всех $\lambda \in \mathcal P_n$, где  $p(n)=|\mathcal P_n|$ -- число разбиений натурального $n$ (число диаграмм Юнга, 
соответствующих натуральному $n$).

Пусть $r_k (\lambda )$ -- число слагаемых в разбиении $\lambda$, равных $k$ (иначе говоря, выполнено равенство $n = \sum_k k r_k (\lambda )$). Число слагаемых в разбиении $\lambda$, больших или равных $\lceil t \rceil$ (где $t$ -- неотрицательное действительное число), обозначим через $\phi _\lambda (t)=\sum\limits_{k\ge t} {r_k (\lambda )} $. Заметим, что $\phi _\lambda (t)$, $t \ge 0$ -- 
ступенчатая функция, непрерывная справа, замыкание внутренности подграфика которой и есть диаграмма Юнга, соответствующая разбиению $\lambda $. Сделаем шкалирование (скейлинг) функции $\phi _\lambda (t)$ с множителем $a > 0$: 
$$\tilde {\phi }_\lambda^a (t)=\frac{a}{n}\sum\limits_{k\ge at} {r_k (\lambda 
)} =\frac{a}{n}\phi _\lambda (at).$$ Заметим, что после шкалирования 
$$\int\limits_0^\infty {\tilde {\phi }_\lambda^a (t)dt} =1.$$

\textbf{Теорема (А.М. Вершик).} {\itПусть $a_n =\sqrt n $, тогда для любых $\varepsilon >0$ и $\delta > 0$ 
существует $n_{\varepsilon, \delta} $, что для любого $n>n_{\varepsilon, \delta} $}
\[
\mu _n \left\{ {\lambda :\quad \mathop {\sup }\limits_t \left| {\tilde {\phi 
}_\lambda^{\sqrt n} (t)-C(t)} \right|<\varepsilon } \right\}>1-\delta ,
\]
{\itгде} $C(t)=\int\limits_t^\infty {\frac{e^{-\sqrt {\varsigma (2)} 
u}}{1-e^{-\sqrt {\varsigma (2)} u}}du} =-\frac{\sqrt 6 }{\pi }\ln 
(1-e^{-\left( {\pi \mathord{\left/ {\vphantom {\pi {\sqrt 6 }}} \right. 
\kern-\nulldelimiterspace} {\sqrt 6 }} \right)t})$.

Другими словами, типичная (в смысле равномерной меры) диаграмма Юнга после 
шкалирования $\tilde {\phi }_\lambda^{\sqrt n} (t)=\frac{1}{\sqrt n }\sum\limits_{k\ge 
\sqrt n t} {r_k (\lambda )} $ с ростом $n$ имеет \textit{предельную форму {\rm(}limit shape{\rm)}}, заданную 
с помощью $C(t)$, или в более симметричной форме: $e^{-\frac{\pi }{\sqrt 6 
}x}+e^{-\frac{\pi }{\sqrt 6 }y}=1$.

Интерпретация этой задачи с точки зрения статистической физики следующая.  $\mathcal P_n $ -- множество всех состояний системы (бозе-частиц) 
с постоянной энергией $n$, каждое состояние $\lambda \in \mathcal P_n $ однозначно 
характеризуется набором значений $r_k (\lambda )$ (чисел заполнений $k$-го 
уровня энергии). Заметим, что у бозе-частиц нет ограничений на значения $r_k 
(\lambda )$, в отличие от ферми-частиц, где $r_k (\lambda )\le 1$. $\mathcal P_n $ в 
статистической физике называется каноническим ансамблем (равномерная 
мера $\mu _n $ соответствует тому, что все состояния равновероятны). Теорема утверждает, что  шкалированная функция плотности распределения по уровням энергии у системы бозе-частиц (с полной энергией $n$)  имеет в пределе вид функции $C(t)$ (иначе теорема утверждает, что существует предельное распределение энергий частиц, позволяющее делать заключение типа: какая предельная
доля общей энергии приходится на частицы с данными энергиями).

\begin{remark}
С детальным доказательством теоремы можно ознакомиться по статье: Вершик А.М.  Статистическая механика комбинаторных разбиений и их предельные конфигурации // Функциональный анализ и его приложения.  1996.  Т. 30, вып. 2.  С. 19--39.
\end{remark}

Согласно приведенной выше статье, далее описаны только основные шаги и идеи доказательства. Читателю требуется с учетом приведенных ниже указаний качественно обосновать результат о предельной форме диаграмм Юнга.

\begin{ordre} Заметим, что для обоснования существования предельной формы шкалированных диаграмм Юнга необходимо для начала проверить выполнимость следующих пределов:
$$
\mathop {\lim }\limits_{n\to 
\infty } \mathbb E_{n}\tilde {\phi }_\lambda^{\sqrt n} (t) = C(t),
$$
$$\mathop {\lim }\limits_{n\to 
\infty }\mathbb E_{n } \left[ {\left( {\tilde {\phi }_\lambda^{\sqrt n} (t)-C(t)} 
\right)^2} \right]=0,$$
что на самомо деле не так просто.

Основная идея обоснования утверждения теоремы заключается, во-первых, в переходе к пространству всех диаграмм Юнга $\mathcal P=\bigcup\limits_n {\mathcal P_n} $ с семейством (по $x \in (0, 1)$) мер $\mu _x $:
$$\mu _x (\lambda )=\frac{x^{\sum\limits_i {\lambda _i } 
}}{F(x)}=\frac{x^{n(\lambda )}}{F(x)}, \text{ где } F(x)=\sum\limits_{n=0}^\infty 
{p(n)x^n},$$
которые в свою очередь индуцируют на пространстве $\mathcal P_n $ исходную равномерную меру  $\mu _n $ (см. утверждение 1 ниже),
во-вторых, в обосновании утверждения указанной теоремы для семейства введенных мер $\mu _x $ при $x$, стремящемся к единице слева, то есть $x \to 1-0$ (что существенно проще, см. далее). Представление меры $\mu_x$ как выпуклой комбинации мер $\mu_n$ (см. утверждение 3 ниже) позволяют делать заключение, что тот же предел (найденный для мер  $\mu _x $ при $x \to 1-0$) будет и для мер $\mu _n $ при $n \to \infty$.

С точки зрения статистической физики $\mathcal P$ 
соответствует (большому) каноническому ансамблю с мерой Гиббса $\mu _x 
(\lambda )=\frac{x^{n(\lambda )}}{F(x)}$, где $F(x)=\sum\limits_\lambda 
{x^{\sum\limits_i {\lambda _i } }} =\sum\limits_{n=0}^\infty 
{p(n)x^{n(\lambda )}} $ -- статистическая сумма. Покажите, что для 
статистической суммы $F(x)$ (или иначе производящей функции последовательности $p(n)$) справедлива формула Эйлера: 
$$F(x)=\sum\limits_{n=0}^\infty {p(n)x^n} =\prod\limits_{k=1}^\infty 
{(1-x^k)^{-1}}.$$

Переход от малого к большому каноническому ансамблю упрощает решение задачи в силу  мультипликативности мер $\mu _x $ (см. утверждение~2 
ниже).

Проверьте справедливость следующих \textbf{утверждений}:

1) $\left. {\mu _x (\lambda )} \right|_{\mathcal P_n } \equiv \frac{\mu _x (\lambda 
)}{\mu _x (\mathcal P_n)}=\mu _n (\lambda )$, то есть мера на $\mathcal P_n $, индуцированная 
мерой $\mu _x $, совпадает с равномерной мерой $\mu _n $.

2) $r_1 (\lambda ),\;r_2 (\lambda ),\ldots$ как случайные величины на $\mathcal P$ 
независимы относительно мер $\mu _x $.

\textbf{Указание}: Проверьте, что $\mu _x \left\{ {\lambda :\;r_k (\lambda 
)=s} \right\}=x^{ks}(1-x^k)$, а совместное распределение является произведением маргинальных.

3) Мера $\mu _x $ является выпуклой комбинацией мер $\mu _n $.

\textbf{Указание}: $\mu _x =\sum\limits_{n=0}^\infty 
{\frac{x^np(n)}{F(x)}\mu _n } $.

Согласно описанной выше схеме (идее) доказательства, проверим выполниммость
$$\mathop {\lim }\limits_{x\to 1-0} \mu_{x} \left\{ {\lambda :\quad \mathop {\sup }\limits_t \left| {\tilde {\phi }_\lambda^{\sqrt {n(\lambda)}} (t)-C(t)} 
\right|<\varepsilon } \right\}=1.$$
Для этого нужно, во-первых,
$$
\mathop {\lim }\limits_{x\to 1-0} \mathbb E_{x}\tilde {\phi }_\lambda^{\sqrt {n(\lambda)}} (t)=C(t),
$$
во-вторых,
$$\mathop {\lim }\limits_{x\to 
1-0} \mathbb E_{x } \left[ {\left( {\tilde {\phi }_\lambda^{\sqrt {n(\lambda)}} (t)-C(t)} 
\right)^2} \right]=0.$$
Заметим, что в отличие от случая, когда усреднение берется по мере $\mu_n$, $n(\lambda)$ -- случайная величина, и пользоваться линейностью математического ожидания, например, при вычислении $$\mathbb E_x \left[ \frac{1}{\sqrt {n(\lambda)}} \sum_{k\ge t\sqrt {n(\lambda)}} r_k(\lambda) \right],$$ нельзя. Преодоление этой трудности аналогично выбору перевального контура в методе Лапласа: вместо предела $x\to 1-0$ выберем последовательность $x_n$ такую, что, с одной стороны, $x_n \to 1-0$ при ${n\to\infty}$, с другой стороны, для каждого $n$ мера $\mu_{x_n}$ концентрируется на  $\mathcal P_n$. Для этого будем оптимизировать (по $x$) вероятность того, что случайное (по мере $\mu_x$) разбиение дает фиксированное исходное значение $n$:
$$ \mu_{x} \left\{ \sum_k kr_k(\lambda) = n \right \} = p(n)\frac{x^n}{F(x)}\to\max\limits_x, $$
так что  $x=x_n$ является корнем уравнения
$x\frac{d}{dx}\left[ {\ln 
F(x)} \right]=n$ (покажите это).
Последнее уравнение можно интерпретировать немного иначе, а именно выбирать значение $x=x_n$ из совпадения математического ожидания (по мере $\mu_x$) случайной величины $n(\lambda)$ с исходным значением $n$:
$$ \mathbb E_x \left[ {\sum\limits_k {kr_k (\lambda )} } 
\right]=n. $$
Обоснуйте (воспользовавшись предствалением статистической суммы $F(x)$ согласно формуле Эйлера и тем, что $r_k(\lambda)$ имеют геометрическое распределение (утверждение 2 выше)), что последнее также эквивалентно $x\frac{d}{dx}\left[ {\ln 
F(x)} \right]=n$. Покажите, что это уравнение имеет единственное решение $x_n \in (0,1)$, подставляя которое обратно в уравнение имеем нормировку:
$$
1 = \sum_{k=1}^n \frac{k}{\sqrt n} \frac {(x_n^{\sqrt n})^{k/{\sqrt n}}}{1-(x_n^{\sqrt n})^{k/{\sqrt n}}}\frac{1}{\sqrt n} \to \int\limits_0^\infty {u\frac{e^{-\sqrt {\varsigma (2)} 
u}}{1-e^{-\sqrt {\varsigma (2)} u}}du}
$$
и при этом $$\mathop {\lim }\limits_{n\to \infty } \mathbb E_{x_n } \left[ {\left( 
{\frac{n(\lambda )}{n}-1} \right)^2} \right]=0.$$ Таким образом, получаем 
$$
\mathbb E_{x_n} \tilde {\phi }_\lambda^{\sqrt n} (t) = \frac{1}{\sqrt n} \sum_{k\ge t\sqrt n} \mathbb E_{x_n} r_k = $$
$$ =\sum_{\frac{k}{\sqrt n}\ge t}\frac {(x_n^{\sqrt n})^{\frac{k}{\sqrt n}}}{1-(x_n^{\sqrt n})^{\frac{k}{\sqrt n}}}\frac{1}{\sqrt n}\to \int\limits_t^\infty {\frac{e^{-\sqrt {\varsigma (2)} 
u}}{1-e^{-\sqrt {\varsigma (2)} u}}du},
$$
что и является нормированной кривой $C(t)$.

Вырожденность предельной формы обосновывается соотношением \[\mathop {\lim }\limits_{n\to 
\infty } E_{x_n } \left[ {\left( {\tilde {\phi }_\lambda^{\sqrt n} (t)-C(t)} 
\right)^2} \right]=0.\]

Таким образом, имеет место
$$\mathop {\lim }\limits_{n\to \infty} \mu _{x_n} \left\{ {\lambda :\quad 
\mathop {\sup }\limits_t \left| {\tilde {\phi }_\lambda^{\sqrt n}  (t)-C(t)} 
\right|<\varepsilon } \right\}=1.$$

Стоит отметить, что ''догадаться''  до ответа о предельной форме диаграмм Юнга можно с помощью 
следующих ''нестрогих'' рассуждений. (См. видеотеку летней школы ''Современная математика'', Дубна, 2013 г., В.А. Клепцын ''Асимптотические задачи комбинаторики'').
Пусть в первом квадранте на целочисленной решетке отмечены пары точек $(x_i 
,y_i )$, так что
\[
0=x_1 <x_2 <\ldots <x_k ,\quad y_1 <y_2 <\ldots <y_k =0.
\]
Число возможных диаграмм Юнга (без фиксации $n$), продолженная граница которых проходит через 
заданные пары точек, приближенно выражается так:
\[
\exp \left\{ {\sum\limits_i {\left( {\Delta _i x+\Delta _i y} \right)H\left( 
{\frac{\Delta _i x}{\Delta _i x+\Delta _i y}} \right)} } \right\},
\]
где $\Delta _i x=x_{i+1} -x_i $, $\Delta _i y=y_{i+1} -y_i $, $H(p)=-p\ln 
p-(1-p)\ln (1-p)$. Проверьте это, воспользовавшись тем, что $C^{[pN]}_N$, где $p\in(0,1)$, без учёта константого множителя приближенно равно $e^{NH(p)}$ для больших $N$. (Считайте, что $\Delta _i x$ и $\Delta _i y$ большие, см. ниже.)

Соедините указанные точки ломаной, продлив ее осями за концевые точки. 
Произведите сжатие по каждой из осей в $\sqrt n$ раз, так чтобы площадь под ломаной (приближённо) равнялась единице. Перейдите в новую систему координат поворотом старой против часовой стрелки 
на $45^{\circ}$, сделав растяжение в $\sqrt{2} $ раза. Ломаная теперь задается как 
график функции $u=f(t)$. Наконец, с ростом числа точек
$$
\exp \left\{ {\sum\limits_i {\left( {\Delta _i x+\Delta _i y} \right)H\left( 
{\frac{\Delta _i x}{\Delta _i x+\Delta _i y}} \right)} } \right\}\to 
$$
$$
\exp 
\left\{\sqrt n {\int {H\left( { \frac{ 1 + {f}'(t)}{2}} \right)dt} } \right\}.
$$
Для нахождения предельной формы нужно найти максимум функционала (с 
энтропийным лагранжианом) при ограничениях $f(t)\ge \vert t\vert $, $\vert 
{f}'\vert \le 1$, $\int {\left( {f(t)-\vert t\vert } \right)dt} =1$. Затем 
вернуться в исходную систему координат.
\end{ordre}
\end{problem}

\begin{problem}\DStar(Статистика выпуклых ломаных \cite{7}.)
\label{convcurve}
\imgh{40mm}{curves1.pdf}{К задаче <<Статистика выпуклых ломаных>> ($c=3$)}
Рассмотрим на плоскости выпуклую ломаную $\Gamma $, выходящую из нуля, у которой вершины являются целыми точками и угол наклона каждого звена неотрицателен и не превосходит ${\pi \mathord{\left/ {\vphantom {\pi  2}} \right. \kern-\nulldelimiterspace} 2} $ ($\Gamma $ -- график кусочно-линейной функции $x_2 = \Gamma(x_1)$). Выпуклость понимается в том смысле, что наклон последовательных звеньев ломаной строго возрастает. Пространство всех таких ломаных обозначим через $\Pi $, а через $\Pi (m_{1} ,m_{2} )$ --- множество ломаных, оканчивающихся в точке $(m_{1} ,m_{2} )$. 

Введем равномерную меру на пространстве $\Pi (m_{1} ,m_{2} )$:
\[P_{m_{1} ,m_{2} } (\Gamma )=\frac{1}{\left|\Pi (m_{1} ,m_{2} )\right|} =\frac{1}{N(m_{1} ,m_{2} )} .\] 

Рассмотрим детерминированную кривую $t_2 = L_{c}(t_1) $, $t_1 \in [0,1]$, заданную уравнением $\left(ct_{1} +t_{2} \right)^{2} =4ct_{2} $ (см. рис. \ref{Fig:curves1.pdf}). Заметим, что точки $(0;0)$ и $(1;c)$ являются соответственно левым и правым концами кривой, также для этой кривой справедливо \[\left. \left({dt_{2} \mathord{\left/ {\vphantom {dt_{2}  dt_{1} }} \right. \kern-\nulldelimiterspace} dt_{1} } \right)\right|_{t_{1} =0} =0, \; \left. \left({dt_{2} \mathord{\left/ {\vphantom {dt_{2}  dt_{1} }} \right. \kern-\nulldelimiterspace} dt_{1} } \right)\right|_{t_{1} =1} =\infty \]
и кривая инвариантна относительно преобразования $$(t_1, t_2) \mapsto (1 - t_2/c, c(1 - t_1)).$$

Покажите справедливость закона больших чисел для выпуклых ломаных $\Gamma\in \Pi(m_{1} ,m_{2} )$:

\textbf{Теорема (Я.Г. Синай).} {\itДля любого $\delta > 0$ справедливо}
$$
\lim_{m_{1},m_2 \to \infty, \frac{m_2}{m_1}\to c} P_{m_{1} ,m_{2} } \left \{ \max_{t_1\in[0;1]}\left | \frac{1}{m_1}\Gamma(m_1 t_1) - L_c(t_1) \right |\le \delta \right\} = 1.
$$
{\it Иными словами, в результате масштабного преобразования (скейлинга) $t_{1} ={x_{1} \mathord{\left/ {\vphantom {x_{1}  m_{1} }} \right. \kern-\nulldelimiterspace} m_{1} } $ и $t_{2} ={x_{2} \mathord{\left/ {\vphantom {x_{2}  m_{1} }} \right. \kern-\nulldelimiterspace} m_{1} } $ форма случайных ломаных становится при $m_{1}, m_2 \to \infty$, $m_2 / m_1 \to c$ детерминированной.}

\end{problem}

\begin{ordre} 

Согласно работе: Синай Я.Г. Вероятностный подход к анализу статистики выпуклых ломаных // Функциональный анализ и его приложения.  1994.  Т.~28, вып.~2.  С.~41--48 (см. также книгу \cite{7}), рассмотрим подход к исследованию асимптотических вероятностных свойств выпуклых ломаных, основанный на понятиях статистической физики -- микроканонического и большого канонического ансамблей (см. предыдущую задачу). С такой точки зрения, распределение вероятностей $P_{m_{1} ,m_{2}}$ на $\Pi(m_{1} ,m_{2})$ вводится как условное распределение, индуцированное подходящей мерой $Q$ (см. ниже), заданной на
пространстве $\Pi$. При этом те или иные свойства ломаных (в частности, закон
больших чисел) устанавливаются сначала по отношению к $Q$, а затем переносятся
на случай $P_{m_{1} ,m_{2}}$  с помощью соответствующей локальной предельной теоремы.

Сначала рассмотрим множество $X$ всех пар взаимно простых целых чисел, и пусть $C_{0} (X)$ -- пространство финитных функций на $X$  с неотрицательными целыми значениями. Покажите, что каждой такой функции естественным образом отвечает некоторая выпуклая (конечнозвенная) ломаная, и наоборот. Введем на $C_{0} (X)$ мультипликативную статистику $Q = Q_{z_1,z_2}$ (параметры $z_1$ и $z_2$, $0<z_1,z_2<1$, определим ниже) -- распределение случайного поля $\nu = \nu(\cdot)$  на $X$ с независимыми значениями при разных $x = (x_1, x_2)\in X$  и распределением:
\[Q_{z_{1} ,z_{2} } (\nu )=\prod _{x=(x_{1} ,x_{2} )\in X}\left[\left(z_{1}^{x_{1} } z_{2}^{x_{2} } \right)^{\nu (x)} \left(1-z_{1}^{x_{1} } z_{2}^{x_{2} } \right)\right] =\] 

\[=\prod _{x=(x_{1} ,x_{2} )\in X}\left(z_{1}^{x_{1} } z_{2}^{x_{2} } \right)^{\nu (x)}  \prod _{x=(x_{1} ,x_{2} )\in X}\left(1-z_{1}^{x_{1} } z_{2}^{x_{2} } \right) .\] 
Это означает, что каждая с.в. $\nu (x)$ имеет геометрическое распределение с параметрами $z_{1}^{x_{1} } z_{2}^{x_{2} } $.

Проверьте, что на пространстве $P_{m_{1} ,m_{2}}$ введенная выше мера $Q_{z_1,z_2}$ индуцирует распределение, не зависящее от параметров $z_1$ и $z_2$, в данном случае равномерное:
\vskip-2mm
\[Q_{z_{1} ,z_{2} } \left(\nu |\Pi (m_{1} ,m_{2} )\right)=\frac{1}{N(m_{1} ,m_{2} )} =P_{m_{1} ,m_{2} } (\nu ).\]

С точки зрения статистической физики, распределение $Q_{z_1,z_2}$ играет роль большого канонического распределения Гиббса, а распределение $P_{m_{1} ,m_{2}}$ -- микроканонического распределения. Если какое-либо событие по отношению к распределению $Q_{z_{1} ,z_{2} } $ имеет малую вероятность, гораздо меньшую, чем вероятность $Q_{z_{1} ,z_{2} } \left(\Pi (m_{1} ,m_{2} )\right)$, то оно имеет и малую вероятность по отношению к распределению $P_{m_{1} ,m_{2} } $. Основная идея заключается в том, чтобы подобрать такие значения параметров $z_{1} ,z_{2} \in \left[0,1\right]$, при которых вероятность $Q_{z_{1} ,z_{2} } \left(\Pi (m_{1} ,m_{2} )\right)$ приняла бы возможно большее значение, то есть чтобы распределение $Q_{z_{1} ,z_{2} }$ концентрировалось на пространстве $\Pi(m_{1} ,m_{2})$:
$$Q_{z_{1} ,z_{2} } \left(\Pi (m_{1} ,m_{2} )\right) = Q_{z_{1} ,z_{2} } \left \{ \sum _{x\in X}\nu (x)x_{1} = m_1, \sum _{x\in X}\nu (x)x_{2} = m_2 \right \} =$$
$$=\frac{N(m_1,m_2)z_1^{m_1}z_2^{m_2}}{\prod_{x=(x_{1} ,x_{2} )\in X}(1-z_{1}^{x_{1} } z_{2}^{x_{2} })^{-1}}\to\max_{z_{1} ,z_{2} \in \left[0,1\right]}.$$
Покажите, что такой выбор согласуется с интуитивным выбором из фиксации ``в среднем''\ по мере $Q_{z_{1} ,z_{2} }$ правого конца случайной ломаной:
\[\begin{array}{l} {\Exp_{z_{1} ,z_{2} } \left(\sum _{x\in X}\nu (x)x_{1}  \right)=\sum _{x\in X}\frac{x_{1}z_1^{x_1}z_2^{x_2}}{1-z_1^{x_1}z_2^{x_2}}=m_{1} ,} \\ {\Exp_{z_{1} ,z_{2} } \left(\sum _{x\in X}\nu (x)x_{2}  \right)=\sum _{x\in X}\frac{x_{2}z_1^{x_1}z_2^{x_2}}{1-z_1^{x_1}z_2^{x_2}}=m_{2}. } \end{array}\] 
\noindent 
Это в свою очередь эквивалентно системе уравнений (см.~предыдущую задачу):
\[\begin{array}{l} 

{z_1 \frac{d}{dz_1} \ln \left(\prod_{x=(x_{1} ,x_{2} )\in X}(1-z_{1}^{x_{1} } z_{2}^{x_{2} })^{-1}\right) =m_{1} ,} \\ 

{z_2 \frac{d}{dz_2} \ln \left( \prod_{x=(x_{1} ,x_{2} )\in X}(1-z_{1}^{x_{1} } z_{2}^{x_{2} })^{-1} \right)=m_{2}. } 

\end{array}\] 

Детали нахождения $z_1$ и $z_2$ см. в указанной выше литературе. Здесь приводится только результат, что $z_{1}$ и $z_{2} $ равны соответственно 
\[1-\left(\frac{\zeta (3)c}{\zeta (2)m_{1} } \right)^{\frac{1}{3}} (1+o(1)), \; \; 1-\left(\frac{\zeta (3)}{\zeta (2)c^{2} m_{1} } \right)^{\frac{1}{3}}  (1+o(1)).\] 
Дзета-функция Римана появляется в связи с тем, что рассматриваются только пары взаимно простых чисел -- пространство $X$ (см. задачу \ref{sec:z_func_riman} из раздела \ref{hard}).

Зафиксируем последовательность $0<\tau _{1} <\tau _{2} <\cdots <\tau _{N} $. Впоследствии $\tau _{1} \to 0$, $\tau _{N} \to \infty $, $\mathop{\max }\limits_{j} (\tau _{j+1} -\tau _{j} )\to 0$. Введем с.в.

\[\zeta _{k}^{(j)} (\nu )=\sum _{x\in X:\; \tau _{j} \le {\frac{x_2}{x_1}} \le \tau _{j+1} }x_{k} \nu (x) ,\quad k=1,2.\] 

\noindent $\zeta _{1}^{(j)} (\nu )$ и $\zeta _{2}^{(j)} (\nu )$ есть приращение по осям $x_{1} $, $x_{2} $ соответственно той части ломаной, где тангенс угла наклона звеньев заключен между $\tau _{j} $ и $\tau _{j+1} $. По отношению к распределению вероятностей $Q_{z_{1} ,z_{2} } $ случайные величины $\zeta _{k}^{(j)} (\nu )$ взаимно независимы при разных $j$. Их математические ожидания:

\[\Exp_{z_{1} ,z_{2} } \zeta _{k}^{(j)} (\nu )=\sum _{\tau _{j} \le {\frac{x_2}{x_1}} \le \tau _{j+1} }x_{k} \frac{z_{1}^{x_{1} } z_{2}^{x_{2} } }{1-z_{1}^{x_{1} } z_{2}^{x_{2} } }  .\] 

С учетом полученных выше параметров $z_{1}, z_{2} $ и заменой переменных $x_1 = m_1^{1/3} t_1$, $x_2 = m_1^{1/3} t_2$, где $t_1$ и $t_2$ меняются с шагом $m_1^{-1/3}$,  из этих соотношений можно получить дифференциальное уравнение второго порядка (см. указанную выше литературу), решением которого (с учётом граничных точек) и является детерминированная кривая (см. рис. \ref{Fig:curves1.pdf}):

\[\left(ct_{1} +t_{2} \right)^{2} =4ct_{2} .\] 

Далее нужно воспользоваться законом больших чисел по отношению к распределениям $Q_{z_{1} ,z_{2} } $ и $P_{m_{1} ,m_{2} } $.

\end{ordre}
\begin{remark}
Помимо уже указанных выше источников см. статью Вершик А.М. Предельная форма выпуклых целочисленных многоугольников и близкие вопросы // Функциональный анализ и его приложения.  1994.  Т. 28, вып.~1.  С. 16--25.
\end{remark}

\section{Предельные теоремы}
\label{zb4}

\begin{problem}
\label{contin}
Пусть для любой $f \in C_b^{\infty}\left(\mathbb{R}\right)$ (класс бесконечно гладких ограниченных функций на $\mathbb{R}$) выполнено соотношение: 
 \[
 \int f(x)\, dF_n(x) \rightarrow \int f(x) \, dF(x), \, n\rightarrow \infty,
 \]
\noindent где $F_n$, $F$ -- функции распределения с.в. $X_n$ и $X$ соответственно. Докажите, что в таком случае имеет место сходимость по распределению  $X_n$ к $X$, то есть выполнено следующее:
\begin{equation*}
F_n(x) \to F(x)\text{ в любой точке непрерывности }F(x) \text{ при } n\to \infty.
\end{equation*}
 
\end{problem}
\begin{remark}
Понятие сходимости по распределению обобщается на случай, когда случайные величины являются вещественнозначными функциями на некотором метрическом пространстве $S$ (т.е. если $X$ - с.в., то $X = X(s),\, s\in S$). В таком случае вполне достаточно рассматривать выполнимость соотношения, приведенного выше, для класса $C_b(S)$ (непрерывных ограниченных функций). 
\par Иногда такую сходимость называютс слабой сходимостью мер. 

См. Биллингсли П. Сходимость вероятностных мер. М.: Наука, 1977.  352 с.
\end{remark}

\begin{problem}
Пусть $\{\xi_n\}_{n=1}^{\infty}$ последовательность с.в и выполнено следующее равенство:
\begin{equation*}
    \mathbb{P} (\xi_n = -n) = \mathbb{P} (\xi_n = n) = \frac{1}{2}. 
\end{equation*}
\par Покажите, что данная последовательность не сходится по распределению.
\par Имеет ли место сходимость соответствующей последовательности функций распределений? Если ответ положительный, то будет ли функцией распределения предельная функция? Будет ли поточечно сходиться соответствующая последовательность характеристических функций? Приведите пример, когда имеет место поточечная сходимость характеристичесих функций, но нет слабой сходимости (сходимости по распределению) соответствующих случайных величин.
\end{problem}

\begin{problem}\Star(Теорема Леви--Крамера.)
\label{levi_kramer}
Пусть $\left\{ {\xi _n } 
\right\}_{n\in { \mathbb{N}}} $ -- последовательность с.в., а $\left\{ {F_n \left( x 
\right)} \right\}_{n\in { \mathbb{N}}} $ и $\left\{ {\phi _n \left( t \right)} 
\right\}_{n\in { \mathbb{N}}} $ -- соответствующие последовательности функций 
распределений и характеристических функций. Покажите, что верны следующие утверждения.
\begin{enumerate}
\item Если существует с.в. $\xi $ с характеристической функцией $\phi 
\left( t \right)$ и такой функцией распределения $F\left( x \right)$, что 
$F_n(x) \to F(x), \, n\to\infty$ в точках непрерывности $F\left( 
x \right)$, то $\phi _n \left( t \right)\to \phi \left( t \right), \, n\to\infty$ 
равномерно на каждом конечном интервале.
\item Если существует такая непрерывная в нуле функция $\phi \left( t 
\right)$, что $\phi _n \left( t \right)\to \phi \left( t \right)$ при $n\to\infty$, то 
существует с.в. $\xi $ с характеристической функцией $\phi \left( t 
\right)$ и функцией распределения $F\left( x \right)$, такой что (при $n\to$\linebreak $\to\infty$)  $F_n \left( x 
\right)\to F\left( x \right)$ в точках непрерывности $F\left( x \right)$, то 
есть $\xi _n \buildrel d \over \longrightarrow \xi $. Более того, 
$F_n \left( x \right) \to F\left( x \right), \, n\to \infty$ равномерно на любом конечном или 
бесконечном множестве точек непрерывности функции $F\left( x \right)$, а 
также $\phi _n \left( t \right)\to \phi \left( t \right)$ равномерно на 
каждом конечном интервале.
\end{enumerate}

\end{problem}
\begin{remark}
См. книгу Сачков В.Н. Вероятностные методы в комбинаторном анализе.  М.: Наука, 1978, а также \cite{Gupta,21,stoianov}.
\end{remark}

\begin{problem}
Пусть с.в. $X_k$ нормально распределены и существует такая с.в. $X$, что $X_n \to X$ в $L_2$, т.е. $\mathbb{E} \left[|X_n-X |^2\right] \to 0$ при $n\to\infty$. Покажите, что с.в. $X$  нормально распределена.
\end{problem}

\begin{problem}(Теорема Линдберга.)
Пусть   
$
\xi_{1}, \ldots , \xi_{n}
$, $n\in\mathbb{N}$, -- 
 независимые с.в., причем  $\Exp(\xi_{i}) = 0$ и $\Var(\xi_{i}) = \sigma^2_{i} < \infty$. Введем обозначения: 
 \[
 S_n = \xi_{1} + \ldots + \xi_{n},
 \]
 \[
 s^2_n = \sigma^2_{1} + \ldots + \sigma^2_{n}.
 \]
  
\noindent Покажите, используя метод телескопических сумм (см. указание), что если $\forall\,\epsilon > 0$ 
   \begin{equation} \label{lindeberg}
  \frac{1}{s_n^2} \underset{k=1}{\overset{n}{\sum}} \underset{x \geq \epsilon s_n}{\int} x^2 dF_{\xi_{k}}(x) \rightarrow 0, \: n \rightarrow \infty,
  \end{equation} 
  
  \noindent то при $n\to\infty$
  $$
  \frac{S_n}{s_n} \overset{d}{\longrightarrow} Z, \, \text{где } Z\in \N(0,1)
  $$
  Покажите также, что из условия \eqref{lindeberg} следует следующее условие на частичную дисперсию (условие пренебрегаемости):
  $$
    \max_{j\leq n}\dfrac{\sigma_j^2}{s_n^2} \rightarrow 0, \, n\rightarrow \infty.
  $$
\end{problem}

\begin{ordre} 
Используя результат задачи \ref{contin}, необходимо доказать, что при $n\to\infty$
\[
\Exp\bigg[f\bigg(\frac{S_n}{s_n}\biggr)\biggr] \rightarrow \Exp(f(N)),
\] 
где $N\sim \mathcal{N}(0,1)$,  $f \in C^{\infty}_b\left(\mathbb{R}\right)$.

Метод телескопических сумм состоит в  последовательной заменe $\xi_{k}$ в сумме $S_n$ на случайные величины $\eta_{k} \in N(0, \sigma_{k}^2)$. В результате чего происходит ступенчатый переход вида
\[
\Exp\bigg[f\bigg(\frac{S_n}{s_n}\biggr)\biggr],\,
\ldots,\,\]
\[\Exp\bigg[f\bigg(\frac{(\xi_{1} + \ldots +\xi_{k-1}) + (\eta_{k} + \ldots + \eta_{n})}{s_n}\biggr)\biggr],\,
\ldots,\,
\Exp[f(N)].
\]
Модуль разности первого и последнего членов данной последовательности можно оценить сверху суммой модулей разности пар последовательно идущих элементов
\begin{align*}
    \sum\limits_{j=1}^n 
    \left|\mathbb{E}\bigg[f\bigg(\dfrac{\ldots+\xi_j + \eta_{j+1}+\ldots}{s_n}\bigg) - 
    f\bigg(\dfrac{\ldots+\xi_{j-1} + \eta_{j}+ \ldots}{s_n}\bigg)
    \bigg] \right|.
\end{align*}

\noindent Используя формулу Тейлора, покажите, что
\begin{equation}\label{tayl_lind}
f(x + h_1) - f(x+h_2)  - \bigg[f'(x)(h_1 - h_2) + \frac{1}{2} f''(x)(h_1^2 - h_2^2)\biggr] \leq g(h_1) + g(h_2),
\end{equation}
\noindent где 
$$g(h) = \underset{x}{\sup}\bigg| f(x+h) - f(x) - f'(x)h - \frac{1}{2} f''(x)h^2\biggr| < O(\min\{h^2, |h|^3\}).$$
Считая, что с.в. $\{\xi_k\}_{k=1}^n$ и $\{\eta_k\}_{k=1}^n$ независимы в совокупности, докажите следующее неравенство:
\[
\biggr|\Exp\bigg[f\bigg(\frac{S_n}{s_n}\biggr)\biggr] - \Exp(f(N))\biggr| \leq \underset{k=1}{\overset{n}{\sum}} \Exp \bigg[g\bigg(\frac{\xi_{k}}{s_n}\biggr)\biggr] + \underset{k=1}{\overset{n}{\sum}} \Exp \bigg[g\bigg(\frac{\eta_{k}}{s_n}\biggr)\biggr].
\] 
Таким образом, задача сводится к оценке остаточных членов формулы Тейлора \eqref{tayl_lind}, где и необходимо воспользоваться условием Линдеберга \eqref{lindeberg}. Разбив каждое математическое ожидание на множествах $(\xi_n <$\linebreak $< \epsilon s_n)$ и $(x_n \geq \epsilon s_n)$, докажите, что при $n\rightarrow \infty$
\begin{align*}
&\underset{k=1}{\overset{n}{\sum}} \Exp \bigg[g\bigg(\frac{\xi_{k}}{s_n}\biggr)\biggr]
    \rightarrow 0,\\
&\underset{k=1}{\overset{n}{\sum}} \Exp \bigg[g\bigg(\frac{\eta_{k}}{s_n}\biggr)\biggr]
\rightarrow 0.
\end{align*}
\end{ordre}

\begin{remark}

Классический вариант ЦПТ может быть получен как следствие теоремы Линдберга~\cite{5}: пусть $\xi_{i}$, $i=1,2,\dots$, -- независимые и одинаково распределенные случайные величины с нулевым математическим ожиданием и конечной дисперсией, тогда при $n\to\infty$
\[
\frac{1}{\sigma \sqrt{n}} \underset{k=1}{\overset{n}{\sum}} \xi_{k} \overset{d}{\longrightarrow}
Z,\, \text{где } Z\in \N(0, 1).
 \]
 Отметим также, что с помощью описанной схемы можно оценить и скорость сходимости в ЦПТ.
 
Можно показать, что если 
\[
\lim_{n\to\infty}\max_{1\le k \le n} \Exp\left[\xi_{k}^2\right]=0\text{ или }\lim_{n\to\infty}\max_{1\le k \le n} \PR\left(|\xi_{k}|\ge\epsilon\right)=0,
\]
то условие \eqref{lindeberg} является не только достаточным, но и необходимым для выполнения ЦПТ. Несложно проверить, что условие \eqref{lindeberg} влечет первое из этих условий, из которого, в свою очередь, следует второе.

См. Биллингсли П. Сходимость вероятностных мер.  М.: Наука, 1977.  352 с.
\end{remark}

\begin{problem}
Пусть $\xi _{1} ,\xi _{2} ,...$ --  независимые одинаково распределенные с.в. с конечной ненулевой дисперсией. \mbox{Обозначим $S_{n} =\sum _{i=1}^{n}\xi _{i}$}. Выяснить, при каких значениях $c$ имеет или не имеет место сходимость:
$$
\PR\left(\frac{S_{n} }{n} \leq c\right)\to I\left(\Exp\xi _{1} \leq c\right),\, 
n\to \infty.
$$

\begin{remark}
См.  Cherny~A.S. The Kolmogorov student's competitions  on probability theory. MSU.
\end{remark}

\end{problem}

\begin{problem}
Пусть имеются независимые одинаково распределенные с.в. $X_i$, $i=1,\dots,n$:

\[X_{i} =\left\{\begin{array}{cc} {1,} & {p} \\ {-1,} & {q=1-p.} \end{array}\right. \] 

Покажите, что верна {\it локальная предельная теорема} (см. \cite{2}) для суммы независимых с.в.  $X_i$: равномерно по всем $x=O\left(\sqrt{n} \right)$ таким, что $(p-q)n+x$ -- целое неотрицательное число, при больших $n$ выполнено
\[\PR\left(\sum _{i=1}^{n}X_{i} =(p-q)n+x \right)\approx \frac{1}{\sqrt{2\pi npq} } \exp \left\{-\frac{x^{2} }{2npq} \right\}.\] 
\end{problem}
\begin{ordre} Покажите, используя \textit{теорему Коши} из курса ТФКП, что
$$\PR\left(S_{n} =\left\lfloor \alpha n\right\rfloor \right)=\frac{1}{2\pi i}     \ointctrclockwise _{|z|=\rho }\frac{(pz+qz^{-1} )^{n} }{z^{\left\lfloor \alpha n\right\rfloor +1} } dz=$$
$$=\frac{1}{2\pi i} \ointctrclockwise  _{|z|=\rho }e^{n\left[\ln (pz+qz^{-1} )-\alpha \ln z\right]} e^{\left[\alpha n-\left\lfloor \alpha n\right\rfloor -1\right]\ln z} dz.$$

Примените \textit{метод перевала} для аппроксимации полученного интеграла. Выберите радиус окружности $\rho $ так, чтобы точка перевала находилась на пересечении этой окружности с положительной вещественной осью: $z=\rho $; для этого найдите максимум функции $\ln (pz+qz^{-1} )-\alpha \ln z$. Замените интеграл по всей окружности на интеграл по $\delta $-дуге, содержащей точку перевала \cite{27}.
\end{ordre}

\begin{problem}\Star(Коралов--Синай \cite{7}.)
Пусть задано гильбертово пространство
$\rm H = L^2\left( {{\rm R},{\rm B},\mu _G } \right)$, где $\mu _G$ --  
гауссовская мера со скалярным произведением 
$$
\langle {f,g} 
\rangle=\int\limits_{-\infty }^{+\infty } 
{f(x)g(x)d\mu_G(x)}.
$$ 

\noindent Пусть $\xi _i $, $i=1,\dots,n$, -- независимые с.в. с нулевым математическим ожиданием и 
плотностью распределения 
$$
p_h =\frac{1}{\sqrt {2\pi } }\left( {1+h(x)} 
\right)e^{-\frac{x^2}{2}}, \, x\in \mathbb{R}^1
$$ где $h$ -- такой элемент $\rm H$, что
\begin{enumerate}
\item величина $\left\| h \right\|$ достаточно мала, 
где 
$$\left\| h \right\|^2=\frac{1}{\sqrt {2\pi } }\int\limits_{-\infty }^{+\infty } 
{h^2(x)e^{-\frac{x^2}{2}}dx},$$

\item $\langle h(x),\mathbbm{1} \rangle=0$,

\item $\langle {h(x),x} \rangle=0$.

\end{enumerate}
Покажите, что последовательность с.в. 
$$
\zeta _n 
=2^{-\frac{n}{2}}\sum\limits_{i=1}^{2^n} {\xi _i } 
$$ 
сходится 
по распределению к нормальной с.в. с нулевым математическим ожиданием и 
дисперсией 
$$
\sigma ^2(p_h )=\frac{1}{\sqrt {2\pi } }\int\limits_{-\infty 
}^{+\infty } {x^2p_h (x)dx}.$$

\end{problem}
\begin{remark} 
Для решения задачи можно воспользоваться \textit{методом ренормгруппы}. Заметим, что
$$\zeta _{n+1} =\frac{1}{\sqrt 2 }\left( {\zeta '_n 
+\zeta ''_n } \right),$$
где 
$$
\zeta '_n 
=2^{-\frac{n}{2}}\sum\limits_{i=1}^{2^n} {\xi _i },\quad 
\zeta ''_n 
=2^{-\frac{n}{2}}\sum\limits_{i=2^n+1}^{2^{n+1}} {\xi _i }
$$
-- независимые 
одинаково распределенные с.в., поэтому 
$$
p_{n+1} (x)=Tp_n (x),
$$ 
где $p_n (x)$~-- плотность распределения $\zeta _n $, а $T$~-- нелинейный 
оператор (действующий в пространстве плотностей), такой что 
$$
Tp(x)=\sqrt 2 
\int\limits_{-\infty }^{+\infty } {p(\sqrt 2 x-u)p(u)du}, \, x\in \mathbb{R}^1
$$ 
В этих терминах 
в задаче нужно показать, что при $n\to \infty$
\[
T^n p_h (x)\to 
\frac{1}{\sqrt {2\pi } \sigma (p_h )}{\rm e}^{-\frac{x^2}{2\sigma ^2(p_h )}}.
\] 

Рассмотрите нелинейный оператор $\tilde{L}$ на пространстве ${\rm H}$, связанный с 
оператором $T$ следующим образом:
$$
T p_h (x)=\frac{1}{\sqrt {2\pi } }\left( 
{1+\tilde{L}\left( {h(x)} \right)} \right){\rm e}^{-\frac{x^2}{2}}.
$$ 
Линеаризуйте 
оператор $\tilde {L}$, показав, что $\tilde {L}h=Lh+O\left( {\left\| h 
\right\|^2} \right)$, где 
\[L(h)(x)=\frac{2}{\sqrt \pi }\int\limits_{-\infty 
}^{+\infty } {e^{-\left( {\frac{x^2}{2}-\sqrt 2 xu+u^2} \right)}h(u)du} .
\]

Далее покажите, что линейный оператор $L$ имеет полное множество собственных 
векторов 
$$h_k (x)=e^{\frac{x^2}{2}}\left( {\frac{d}{dx}} 
\right)^ke^{-\frac{x^2}{2}},\quad k\ge 0,
$$
с собственными значениями $\lambda 
_k =2^{1-\frac{k}{2}}$, $k\ge 0$. 

Пусть ${\rm H}_k $~-- одномерное 
подпространство, натянутое на $h_k (x)$. Тогда $L$~-- сжимающий оператор на 
${\rm H} \setminus \left( {{\rm H}_0 \oplus {\rm H}_1 \oplus {\rm H}_2 } \right)$.

Покажите, что $\tilde {L}$~-- сжимающий оператор на ${\rm H} \setminus \left( {{\rm 
H}_0 \oplus {\rm H}_1 } \right)$, имеющий единственную неподвижную точку: 
\[f_h (x)=\frac{1}{\sigma (p_h )}e^{\frac{x^2}{2}-\frac{x^2}{2\sigma ^2(p_h 
)}}-1.
\]
\end{remark}

\begin{problem}
В игре в рулетку колесо разделено на 38 равных секторов: 18 красных, 18 белых и два сектора (0 и 00) зеленого цвета. Обозначим через $X_{i} $ выигрыш в $i$-й игре. Тогда $X_{1} ,X_{2} ,X_{3} ,...$ -- независимые с.в., имеющие распределение 
\[X_{i} =\left\{\begin{array}{cc} {+1,} & p = 18/38, \\ {-1,} & 1-p = 20/38. \end{array}\right. \] 
Пусть сыграно $n=19^2=361$ партий. С помощью ЦПТ и неравенства Берри--Эссеена оцените вероятность того, что в удастся обыграть казино.
\end{problem}
\begin{remark}
Приведем формулировку теоремы  Берри--Эссеена \cite{19}, \cite{21} Т. 1.
\label{sec:BerryEssen}
 Пусть $\xi_1, \xi_2\dots$ -- независимые одинаково распределенные с.в., причем $\mathbb{E}\xi_i = m$, 
 $\mu_3={\mathbb E}|\xi_i - {\mathbb E}\xi_i|^3<\infty$, $\sigma^2=\mathbb D \xi_i$.
Близость с.в. $\frac{\sum_{i=1}^{n}\xi_i-nm}{\sigma\sqrt{n}}$ к стандартной нормально распределенной с.в. (согласно ЦПТ) в смысле 
близости их функций распределения определяется неравенством Берри--Эссеена: 
$$
\sup\limits_{t\in \mathbb{R}} \left| {\mathbb P}\Bigl( \frac{\sum_{i=1}^{n}\xi_i-nm}{\sigma\sqrt{n}}<t \Bigr) - \Phi(t) 
\right| \le \frac{C_0 \mu_3}{\sigma^3 \sqrt{n}} , 
$$
где $\left(2\pi\right)^{-1/2}\le C_0 < 0.7056$,
$\Phi(t)=\int_{-\infty}^t \frac{e^{-\tau^2/2}}{\sqrt{2\pi}}\, d\tau, \, t\in \mathbb{R}$.

Неравентсво Берри--Эссеена в общем случае неулучшаемо (см., например, Сенатов В.В. Центральная предельная теорема: Точность аппроксимации и асимптотические разложения. М.: Книжный Дом ``ЛИБРОКОМ'', 2009). Действительно, если сохранить вопрос в задаче (об оценке вероятности обыграть казино), но изменить вероятность успеха на $p=0.5$, то легко показать, что неравенство Берри--Эссеена превратится в равенство с $C_0\approx\left(2\pi\right)^{-1/2}$ (см., например, \cite{21} Т. 1).

Однако, в случае когда расматриваются маловероятые события (большие уклонения) неравенство Берри--Эссеена может быть довольно грубым. В частности, в данном примере имеет смысл также посмотреть, что дает использование неравентсв концентрации меры, см., например, Boucheron S., Lugosi G., Massart P. Concentration inequalities: A nonasymptotic theory of independence. – Oxford University Press, 2013. Подробнее об этом будет написано в Части 2. Далее на физическом уровне строгости будет показано, что если не налагать больше никаких условий на моменты с.в., то с точностью до логарифмического множителя порядок оценки в неравенстве Берри--Эссеена оптимальный и с учетом оценок вероятностей больших уклонений, т.е. это неравентсво равномерно оптимально по $t\in \mathbb{R}$.

Пусть независимые одинаково распределенные с.в. $\left\{ {\xi _i } 
\right\}_{i=1}^n $ удовлетворяют условию (для простоты положим $m = 0$)
$$
\PR\left(\xi_i > x \right) = V\left( x \right) = O\left(x^{-\alpha} \right), \quad \alpha > 2.
$$
Тогда (см. задачу \ref{Cramer_Zone} раздела \ref{zb4})
$$
\PR\left( {\sum\limits_{i=1}^n {\xi _i } \ge x} \right)\mathop \simeq 
\limits_{n\gg 1} 1-\Phi \left( {\frac{x}{\sigma}} \right)+ n V\left( x 
\right),
\quad
\Phi \left( x \right)=\frac{1}{\sqrt {2\pi } }\int\limits_{-\infty }^x 
{e^{-{y^2} \mathord{\left/ {\vphantom {{y^2} 2}} \right. 
\kern-\nulldelimiterspace} 2}dy} .
$$
Используя 
$$
0.2e^{-{2x^2}/\pi}\le 1-\Phi \left( x \right)\le e^{-{x^2}/2}, \quad x\gg 1,
$$
отсюда можно получить, что в режиме ЦПТ
$$\PR\left( {\sum\limits_{i=1}^n {\xi _i } \ge x} \right)\mathop \simeq 
\limits_{n\gg 1} 1-\Phi \left( {\frac{x}{\sigma}} \right),\;\;x\le 
\sqrt {\left( {\alpha -2} \right)\sigma^2 n\ln n} ,$$
в режиме тяжелых хвостов
$$\PR\left( {\sum\limits_{i=1}^n {\xi _i } \ge x} \right)\mathop \simeq 
\limits_{n\gg 1} nV\left( x \right),\;\;x>\sqrt {\left( {\alpha -2} 
\right)\sigma^2 n\ln n}.$$
Неравенство Берри--Эссеена может быть объяснено исходя из предположения $\mu_3 < \infty$. Это соотношение означает (точнее говоря, для того, чтобы оно имело место достаточно), что $3 - (\alpha + 1) < -1$, т.е. $\alpha > 3$. При таком $\alpha$ имеем
$$n V\left(\sqrt{\left( {\alpha -2} \right)\sigma^2 n\ln n}\right) = n O\left(n^{-3/2}\right)=O\left(n^{-1/2}\right),$$ что поясняет зависимость $\sim n^{-1/2}$ в правой части неравентсва Берри--Эссена.
\end{remark}

\begin{problem} (Петербургский парадокс \cite{19}, \cite{book12}.)
\label{piter}
Рассмотрим следующую игру (традиционное название игры -- ''орлянка``): 
\par Игрок делает ставку и подкидывает симметричную монету. Если выпадает «орел», то игрок забирает удвоенную ставку, иначе теряет все свои деньги. Броски повторяются независимо в каждой игровой серии. Изначально игрок имеет в банке сумму, условно равную 1 и в каждой серии разыгрывает все деньги, которые у него есть.
\par Пусть $X_{1} ,X_{2} ,X_{3},\dots$  -- независимые с.в., обозначают состояние счета игрока после 1,~2,~3~\dots\ игр соответственно. Нетрудно заметить, что распределение $X_i$ выглядит следующим образом:
\begin{equation*}
\PR\left(X_{i} =2^{k} \right)=2^{-k}, \, k\in \mathbb{N}\cup\{0\}, \, k \leq i.
\end{equation*}
\par Иначе говоря, если в игре в орлянку $k$ раз подряд выпал «орел», то выигрыш будет~$2^{k}$.
\par Справедливой ценой за игру называют математическое ожидание выигрыша при бесконечной игре. Нетрудно проверить, что $\Exp X_{\infty} =+\infty $. Тем не менее, покажите, что выполнено следующее: 
$$
\frac{S_{n} }{n\log _{2} n} \mathop{\to }\limits^{p} 1\quad \text{при~} n\to \infty ,
$$ 
где $S_{n} =\sum _{k=1}^{n}X_{k}$. Проинтерпретируйте этот результат, введя цену за $n$ игр.

\begin{remark} 

Введем с.в. $X_{nk} \, 1\le k\le n$ такие, что для каждого фиксированного $n$ с.в. $X_{nk}$, $1\le k\le n$ -- независимы. Определим также величны $b_n$ так, что $b_{n} >0, \, b_{n} \to \infty$ при $n\to\infty$. Введем дополнительно величины  $\bar{X}_{nk} =X_{nk}  I \left\{X_{nk} \le b_{n} \right\}$ и предположим, что при $n\to \infty$ выполнены следующие условия:
\begin{enumerate}
\item 
$
\sum _{k=1}^{n}\PR\left(\left|X_{nk} \right|>b_{n} \right)\to 0;
$ 
\item 
$
\dfrac{1}{b_{n} ^{2} }\sum _{k=1}^{n}\Var\left( \bar{X}_{nk}\right)   \to 0.
$ 
\end{enumerate}
Тогда верно следующее:
$$\frac{1}{b_{n} }\left(\sum _{k=1}^{n}X_{nk}  -\sum _{k=1}^{n}\Exp \left(\bar{X}_{nk}\right)  \right) \mathop{\to }\limits^{p} 0\quad \text{при~} n\to \infty. 
$$

\noindent Для решения задачи положите $X_{nk} =X_{k} $. В качестве $b_{n} >0$ возьмите $b_{n} =2^{m(n)} $, где $m(n)$ -- целое число, которое можно представить в виде $m(n)=\log _{2} n+K(n)$, $K(n)\to \infty $ при $n\to \infty $. Например, если \mbox{$K(n)\le \log_2 (\log_2 n)$,} то 
\begin{equation*}
\frac{S_{n} }{n\log _{2} n} \mathop{\to }\limits^{p} 1 \text{ при }n\to \infty .
\end{equation*}

\end{remark} 

\end{problem}

\begin{problem}
Пусть в некоторой игре размер выигрыша есть с.в., которая принимает значение $2^{k} -1$ с вероятностью 
$$p_{k} =\frac{1}{2^{k} k(k+1)} \quad \text{для~} k=1,2,3,...$$ 
и значение\textit{ $-1$} с вероятностью 
$$p_{0} =1-\sum _{k=1}^{\infty }p_{k} .$$ 
Проверьте, что математическое ожидание выигрыша равно нулю. Применив теорему из предыдущей задачи, покажите, что при $n\to \infty $  для суммарного размера выигрыша за $n$ партий ($S_{n} $) справедливо
$$\frac{S_{n} }{n \log _{2} n}\mathop{\to }\limits^{p} 1.$$

\begin{remark}  
В замечании к задаче \ref{piter} положите $b_{n} =2^{m(n)} $, где $$m(n)=\min \left\{m\in \mathbb{N}:\; 2^{-m} \frac{1}{\sqrt{m^{3}}} \le n^{-1} \right\}.$$ Следует также обратиться к книге \cite{19}.
\end{remark} 

\end{problem}

\begin{problem}\label{max_stable} (max-устойчивые распределения: Гумбеля, Фреше, Вейбулла.)
Распределение $G\left(x\right),\, x\in \mathbb{R}^1$ называется max-устойчивым, если для любого $n=1,2,...$ существуют $a_{n} >0$ и $b_{n} \in {\mathbb R}$, такие что $G^{n} \left(a_{n} x+b_{n} \right)=G\left(x\right)$.
Пусть есть независимые одинаково распределенные с.в. $X_{1} ,...,X_{n} $ с распределением $F\left(x\right)$. Обозначим\linebreak $X_{\left(n\right)} =\max \left\{X_{1} ,...,X_{n} \right\}$. Легко видеть, что распределение с.в. $X_{(n)}$ есть $F_{X_{\left(n\right)} } \left(x\right)=\left[F\left(x\right)\right]^{n} $.
Обозначим также через $x_{(n)}$ решение уравнения $$\mathbb{P}(X>x_{(n)}) = 1-F(x_{(n)}) = 1/n$$ (так называемое характеристическое наибольшее значение).
\begin{enumerate}
\item Пусть при $x\to\infty$ функция распределения $X_i$ имеет вид\linebreak $F(x) = 1-Ax^{-\alpha}$, $\alpha>0$. Покажите, что функция распределения $Y = X_{(n)}/x_{(n)}$ поточечно стремится к пределу $F_{Y}(y) = {\rm e}^{{-y}^{-\alpha}}$ при $y\geq 0$, $F_{Y}(y) = 0$ при $y<0$ (распределение Фреше).
\item Пусть с.в. $X_i$ ограничены сверху ($X_i\leq 0$) и в окрестности нуля при $x\to 0-$ функция распределения $X_i$ имеет вид $F(x) = 1-A|x|^{\alpha}$, $\alpha>0$, $F(x)=1$ при $x>0$. Покажите, что функция распределения с.в. $Y = X_{(n)}/x_{(n)}$ стремится к пределу $F_{Y}(y) = {\rm e}^{{-|y|}^{\alpha}}$ при $y\leq 0$, $F_{Y}(y) = 1$ при $y>0$ (распределение Вейбулла).
\item Пусть при $x\to\infty$ функция распределения $X_i$ имеет вид $F(x) =$\linebreak $= 1- {\rm e}^{-\lambda x}$. Покажите, что функция распределения с.в. $Y = X_{(n)}-x_{(n)}$ стремится к пределу $F_{Y}(y) = {\rm e}^{-{\rm e}^{-\lambda y}}$ (распределение Гумбеля).

\item Покажите, что распределения Вейбулла, Гумбеля, Фреше, являются max-устойчивыми.

\item Покажите, что если $X_i$, $i=1,2,\dots$ -- независимые одинаково распределенные стандартные
\\

\indent 1)  гауссовские (т.е. нормальные) с.в., то верно
\[
\lim_{n\to\infty}\PR\left(2m_n[\max_{i=1,\dots,n}{X_i} - m_n]\leq z \right)= \exp(-\exp(-z)),
\, z\in \mathbb{R}^1
\]
где 
\[
  m_n = \left[2\log\frac{n+1}{\sqrt{8\pi}\log(n+1)}\right]^{1/2};
\]
\indent 2) лапласовские с.в. $\left(\text{с.в. с плотностью } \exp(-|x|/2), x\in \mathbb{R} \right)$, то верно
\[
\lim_{n\to\infty}\PR\left(\max_{i=1,\dots,n}{X_i} - \log(n/2)\leq z \right)= \exp(-\exp(-z)),
\, z\in \mathbb{R}^1
.
\]
\end{enumerate}

\begin{remark}
Рассмотрим следующий пример (см. Лагутин М.Б. Наглядная математическая статистика.  М.: БИНОМ. Лаборатория знаний, 2007). Пусть $X_1,\dots,X_n$~-- независимые одинаково распределенные с.в. с распределением 
$$
F(x) = \left(1-\frac{1}{\ln x}\right)I\left(x>e\right).  
$$
Обозначим $X_{(n)} = \max\{X_1,\dots,X_n\}$. Оценим $\gamma = P(X_{(n)}>10^7)$ при $n=4$. Из независимости и одинаковой распределенности $X_i$ 
$$P(X_{(n)}\leq x) = [F(x)]^n.$$
Поскольку $\ln 10 \approx 2.3$, а $(1-\epsilon)^n\approx 1-\epsilon n$ при малых $\epsilon$, получаем
$$\gamma = 1 - \left(1-\frac{1}{7\ln 10}\right)^4 \approx 1/4.$$
Таким образом, примерно в каждом четвертом случае значение $X_{(4)}$ будет превышать $10^7$. 
Оказывается, что из-за того, что функция $F(x)$ имеет <<сверхтяжелый>> правый <<хвост>>, распределение с.в. $X_{(n)}$ чрезвычайно быстро с ростом $n$ уходит на бесконечность и никаким линейным преобразованием не удается <<вернуть>> его в конечную область. Точнее, невозможно подобрать такие константы $a_n$ и $b_n>0$, чтобы последовательность $(X_{(n)}-a_{n})/b_n$ сходилась бы по распределению к невырожденному закону.

См.  лекции А.Н.~Соболевского в НМУ  http://www.mccme.ru/ium/ s09/probability.html, а также Лидбеттер~М., Линдгрен~Г., Рот\-сен~X. Экстремумы случайных последовательностей и процессов.  М.: Мир, 1989.

Отметим также, что распределения Гумбеля, Фреше, Вейбулла исчерпывают все возможные типы предельных распределений в классе max-устойчивых распределений. Приложения распределения Гумбеля отражено в задачах \ref{gumbel}, \ref{gibbs} раздела~\ref{hard} и задаче 15 раздела 6.
\end{remark}
\end{problem}

\begin{problem}

Пусть $X_{1} ,...,X_{n}$ -- независимые одинаково распределенные с.в. со стандартным распределением Коши. Распределение Коши задается следующей функцией плотности:
$$
f\left(x\right)=\frac{1}{\pi }\frac{1}{1+x^{2}}, \, x\in \mathbb{R}^1.  
$$ 
Воспользовавшись предыдущей задачей, найдите предельное распределение для должным образом нормированных с.в. $X_{\left(n\right)} = \max\limits_{i}(X_i)$.
\end{problem}
\begin{remark}
На тему свойств распределения Коши рекомендуется также посмотреть задачу \ref{cauchy_gen} раздела \ref{standart} и задачу \ref{cauchy_rad_emitt} раздела \ref{zb4}.
\end{remark}


\begin{problem}
Пусть ${X}_{n} \in {\mathbb R}^{m} $, $n=1,2,...$ -- независимые одинаково распределенные случайные векторы, причем  
$$
\Exp{X}_{n} =0, \quad \Exp({X}_{n} {X}_{n}^{T}) =R
$$ ($R$ -- неотрицательно определенная матрица (по определению, см. задачу \ref{normal} раздела \ref{standart}), однако для наглядности дополнительно будем считать, что $R$ положительно определенная матрица). С помощью аппарата характеристических функций докажите, что для любого борелевского множества $B\subseteq {\mathbb R}^{m} $ верно
\[\mathop{\lim }\limits_{N\to \infty } \PR\left(\frac{1}{\sqrt{N} } \sum _{n=1}^{N}{X}_{n}  \in B\right)=\dfrac{1}{(2\pi)^{m/2}\det(R)^{1/2}} \int _{B}{\rm e}^{-\frac{1}{2} \langle{x},R^{-1}{x}\rangle} d {x} .\] 
\noindent Переформулируйте и решите задачу для случая, когда матрица\;$R$ положительная полуопределенная.
\end{problem}
\begin{remark}
Аппарат характеристических функций -- мощная конструкция для изучения свойств распределений случайных величин, т.к. известно, что есть взаимно однозначное соответствие между вероятностой мерой и её фурье-образом. Читателю можно порекомедовать книгу   Боровкова~А.А. \cite{1b} для более подробного изучения свойств характеристических функций. Также большое количество показательных упражнений на характеристические функции можно найти в книге Ширяева~А.Н.   Т.~1 [\ref{chiraiev}].\\
\indent Развивая историческую часть доказательства ЦПТ, стоит также упомянуть, что классически существует два подхода доказательства указанной выше теоремы. Первый способ -- \textit{метод Линденберга}\linebreak (см.~задачу~5 из этого раздела) и второй более распространенный -- метод характеристических функций (\textit{метод Ляпунова}), с которым можно ознакомиться, например, по книге~\cite{6}. Возникает вопрос: какой из способов доказательства более общий? Оказывается, что первый, т.е. метод Линдеберга.
\end{remark}

\begin{problem}

Пусть $X_1,\ldots,X_n$~-- независимые одинаково распределенные с.в. Пусть также характеристическая функция с.в. $X_k$ представляется 
в окрестности $t=0$ в виде 
$$
\varphi_{X_k}(t)={\mathbb E}(e^{it X_k})=1+imt+o(t). 
$$
 Верно ли, что при $n\to\infty$ 
$$
\frac{1}{n}\sum\limits_{i=1}^{n} X_i \xrightarrow{p} m?
$$
\end{problem}
\begin{remark}
См. Cherny A.S. The Kolmogorov student's competitions on probability theory. MSU.
\end{remark}

\begin{problem}
Пусть при каждом $n\geqslant 1$ независимые с.в. $X_{1n}, X_{2n},\ldots, X_{nn}$ таковы, что $X_{kn}\in \Be(p_{kn})$, где 
$$\lim_{n\to\infty}\max\limits_{1\leqslant k\leqslant n} p_{kn}=0\quad\text{и ~} \lim_{n\to\infty}\sum\limits_{k=1}^{n}p_{kn}=\lambda.$$ 
Тогда при $n\to\infty$
\begin{equation*}
\label{TPois}
{\mathbb P}(S_n=m)\to  e^{-\lambda}\frac{\lambda^m}{m!}, \quad m=0,1,2,\ldots, 
\end{equation*}
 где $S_n=\sum\limits_{k=1}^{n} X_{kn}$. 
\begin{remark}
Если $p=\lambda/n$, то $$C_{n}^{m}p^{m}(1-p)^{n-m} = \frac{\lambda^m}{n!}e^{-\lambda}\left( 1 + O\left(\frac{m^2+\lambda^2}{n}\right)\right).$$
Кроме того, имеет место следующая оценка сходимости по вариации (Прохорова--Ле Кама): 
$$\sum_{m=0}^{\infty}{\left|\mathbb P(S_n=m) - e^{-\lambda}\frac{\lambda^{m}}{m!}\right|} \leq 2 \sum_{k=1}^{n}{p_{kn}^2}.$$
В качестве основных ссылок для этой задачи отметим \cite{19}, Т.~1 [\ref{chiraiev}].
\end{remark}
\end{problem}

\begin{problem}
В течение дня игрок в казино участвует в $N=100$ независимых розыгрышах. В каждом розыгрыше он выигрывает с вероятностью 
$p=0.01$. Оцените вероятность события:\\
\indent а) игрок ни разу не выиграет,\\
\indent б) выиграет ровно один раз, \\
\indent в) выиграет ровно три раза.
\par Предположим, что описанная выше игра из $N$ розыгрышей повторяется в течении $n=100$ дней. 
Оцените вероятность того, что за эти $100$ дней в общей сложности реализуется: не менее $100$ выигрышей; не менее $300$ выигрышей.

\end{problem}

\begin{problem}(Дельта метод.) 
Пусть $T_n$~-- последовательность с.в., такая, что  при $n\to\infty$
$$
\sqrt{n}(T_n-\theta)\xrightarrow{d}Z,\, \text{где }Z\sim \N(0,\sigma^2(\theta)),\quad \sigma(\theta)>0.
$$
Пусть отображение $g:\mathbb{R}\to\mathbb{R}$ дифференцируемо в $\theta$ и $g^{\prime}(\theta)\not=0$. Покажите, что при $n\to\infty$
$$
\sqrt{n}\left(g(T_n)-g(\theta)\right)\xrightarrow{d} Z, \, \text{где } Z\sim\mathcal{N}(0,[g^{\prime}(\theta)]^2\sigma^2(\theta)).
$$
\end{problem}
\begin{remark} 
Пусть  $g^{\prime}(\theta) = 0$, в таком случае распределение $g(T_n)$ определяется третьим членом разложения Тейлора, то есть
$$
g(T_n) = g(\theta) + \frac{(T_n-\theta)^2}{2}g^{\prime\prime}(\theta)+o\left((T_n-\theta)^2\right),
$$ 
где в данной задаче $o(1)$ -- величина, сходящаяся с вероятностью $1$ к~$0$. Поэтому при $n\to\infty$
$$
\sqrt{n}(g(T_n) -g(\theta)) = n\frac{(T_n-\theta)^2}{2}g^{\prime\prime}(\theta)+o(1)\xrightarrow{d} \frac{g^{\prime\prime}(\theta)\sigma^2(\theta)}{2}\chi^2(1),
$$
где $\chi^2(1)$ -- случайная величина с распределением $\chi^2$ с одной степенью свободы.

Для решения задачи рекомендуется ознакомиться с литературой по дельта-методу \cite{Gupta}.
\end{remark}

\begin{problem}
Пусть $\{T_n\}$ -- последовательность $k$-мерных случайных векторов, таких, что при $n\to\infty$
$$
\sqrt{n}(T_n-\theta)\xrightarrow{d}Z, \, \text{где }Z\sim\mathcal{N}(0,\Sigma(\theta)).
$$ 
Пусть $g:\mathbb{R}^{k}\to \mathbb{R}^m$ дифференцируемо в $\theta$ с $\nabla g(\theta)$. 
\begin{enumerate}
\item

Покажите, что при $n\to\infty$
$$
\sqrt{n}\left(g(T_n)-g(\theta)\right)\xrightarrow{d} Z, \, \text{где } Z\sim\mathcal{N}(0,\nabla g(\theta)^{T}\Sigma(\theta)\nabla g(\theta)),
$$
если $\nabla g(\theta)^{T}\Sigma(\theta)\nabla g(\theta)$ положительно определена.
\item
Рассмотрите пример, когда   $X_1,\dots,X_n$ -- независимые одинаково распределенные случайные величины с математическим ожиданием $\mu$ и дисперсией $\sigma^2$, а $T_n = \frac{1}{n}\sum_{i=1}^n X_n$. Покажите, что тогда при $n\to\infty$ 
$$
\sqrt{n}(T_n^2-\mu^2)\xrightarrow{d}Z, \, \text{где }Z\sim\mathcal{N}(0,4\mu^2\sigma^2).
$$
\end{enumerate}
\end{problem}

\begin{remark}
См. также книгу: Боровков А.А. Математическая статистика.  СПб.: Лань, 2010.
\end{remark}

\begin{problem}(ЦПТ для стационарных последовательностей.) 
Пусть $X_i$, $i=$\linebreak $=1,2,\,\dots$ -- стационарная последовательность с $\mathbb{E}(X_i)=\mu$ и $\Var (X_i) =$\linebreak $= \sigma^2<\infty$, обладающая  следующим свойством: для некоторого фиксированного $m$ выполнено, что $(X_1,\dots,X_i)$ и $(X_j,X_{j+1}\dots,)$ независимы при $j-i>m$. Покажите, что $n\to\infty$

$$
\frac{\dfrac{1}{n}\sum\limits_{i=1}^n X_i-\mu}{\sqrt{n}}\overset{d}{\to} \mathcal{N}(0,\tau^2),
$$
где $\tau^2 = \sigma^2+2\sum\limits_{i=2}^{m+1}{\rm cov}(X_1,X_i).$

\end{problem}
\begin{remark}
В главах 9 и 10 книги \cite{Gupta} содержится большое количество вариаций ЦПТ для случайных  последовательностей с зависимыми элементами (например, стационарных последовательностей, марковских последовательностей и др.).
\end{remark}

\begin{problem}
Пусть $X_1$, $X_2$,\dots -- независимые одинаково распределенные с.в. с $\mathbb{E}|X_1|^4 < \infty$. Пусть $\mathbb{E}(X_1)=\mu$ и $\mathbb{D}(X_1)=
\sigma^2$ и $\bar{X} = \frac{1}{n}(X_1 + \dots + X_n)$. Пусть $g$ -- функция с равномерно ограниченной четвертой производной. Покажите, что 
\begin{enumerate}
\item $\mathbb{E}[g(\bar{X})] = g(\mu)+\frac{g^{(2)}(\mu)\sigma^2}{2n}+O(n^{-2})$, 
\item
$\mathbb{D}[g(\bar{X})] = \frac{(g^{\prime}(\mu))^2}{n}+O(n^{-2}).$
\end{enumerate}
\end{problem}
\begin{remark}
Основные ссылки по указанным результатам содержатся в \cite{Gupta}.
\end{remark}

\begin{problem}

Показать, что при бросании симметричной монеты $n$ раз отношение числа выпадений герба к числу выпадений решки почти наверное стремится 
к $1$ при $n\to\infty$, а вероятность того, что число выпадений герба в точности равняется числу выпадений решки, при четном числе 
бросаний стремится к $0$ при $n\to\infty$. 
\end{problem}

\begin{problem}
Пусть при любом $\lambda >0$ с.в. $\xi _{\lambda } $ имеет распределение $\Po(\lambda)$. Докажите, что $(\sqrt{\lambda})^{-1}(\xi _{\lambda } -\lambda)  $ по распределению сходится  к стандартному нормальному распределению при $\lambda \to \infty $.
\end{problem}

\begin{problem}
\label{gamma}
Пусть с.в. $X_n\in \Gamma(\lambda,n)$, т.е. неотрицательная случайная величина  с плотностью распределения при $x\geq 0$
$$
p(x) = \dfrac{\left(x/\lambda\right)^{n-1}{\rm e}^{-x/\lambda}}{\Gamma(n)}.
$$
Напомним, что для целых $n$ гамма-функция принимает значения $\Gamma(n) = (n-1)!$.
Покажите, что из ЦПТ следует 
$$
\frac{X_n-m(\lambda)\cdot n}{\sigma(\lambda)\cdot\sqrt{n}} \xrightarrow{d} Z, \, 
\text{где }Z\sim\mathcal{N}(0,1) \text{ при } n\to\infty . 
$$
Найдите $m(\lambda)$, $\sigma(\lambda)$. 
\end{problem}
\begin{remark}
См. замечание к задаче \ref{laplace} раздела \ref{standart}.
\end{remark}

\begin{problem}
Пусть $X_n$ -- последовательность независимых с.в., сходящаяся по вероятности к с.в. $X$:  $X_n\xrightarrow{p}X$ при $n\to\infty$. Докажите, 
что с.в. $X$ вырождена, т.е. $X\equiv c$, где $c$ -- некоторое число. 
\end{problem}

\begin{remark}

Справедливы следующие утверждения:

\begin{enumerate}
\item
Из любой сходящейся по мере (в частности, по вероятностной) последовательности 
измеримых функций (в частности, с.в.) можно выделить подпоследовательность, сходящуюся почти всюду (п.н.). 

\item
Из закона нуля и единицы Колмогорова следует, что для всякого разбиения прямой ${\mathbb R}$ на борелевские множества $\{ B_m\}_{m\geqslant 1}$ 
ровно для одного $m=m_0:$ $\ {\mathbb P}(A_{B_{m_0}})=1$, для остальных $m:\  {\mathbb P}(A_{B_m})=0$, где 
$$
A_{B_m}=\{ \omega: \, X=\lim\limits_{k\to\infty} X_{n_k}\in B_m \} . 
$$
\end{enumerate}

\noindent См.  Cherny~A.S. The Kolmogorov student's competitions  on probability theory. MSU, 
 а также \cite{22}, \cite{220}.
\end{remark}


\begin{problem}
В некотором городе прошел второй тур выборов. Выбор был между двумя кандидатами $A$ и $B$ (графы <<против всех>> на этих выборах не было). 
Сколько человек надо опросить на выходе с избирательных участков, чтобы исходя из ответов можно было определить долю проголосовавших 
за кандидата $A$ с точностью $3\%$ и с вероятностью, не меньшей $0.99$?
\end{problem}

\begin{problem} (Асимптотическое распределение числа инверсий в случайной перестановке.)
На множестве $n!$ перестановок $n$ различных элементов задано равномерное распределение. Обозначим через $\xi_k$ случайную величину, 
равную числу инверсий, образованных элементом с номером $k$, т.е. равную числу элементов с номерами, меньшими, чем $k$, 
которые стоят в перестановке правее элемента с номером $k$. Покажите, что 
$$
\frac{\sum\limits_{k=1}^{n}\xi_k -\left.n^2\right/4}{\left.n^{3/2}\right/6}\xrightarrow{d} \N(0,1) \quad \text{ при } n\to\infty . 
$$
\end{problem}

\begin{remark} Для решения этой и следующей задачи полезно ознакомиться с книгой Сачков В.Н. Комбинаторные методы дискретной математики.  М.: МЦНМО, 2004.
$ $

\begin{enumerate}

\item
 Введем с.в. 
$$
\xi_{k,i} = I(\text{<<$k$ находится левее числа $i$>>}) , 
$$

И обозначим
$$
\xi_k=\xi_{k,1}+\xi_{k,2}+\ldots +\xi_{k,k-1}.
$$

Покажите, что
$${\mathbb E}\xi_k=\left.(k-1)\right/2,$$

$$
{\mathbb E}\xi_k^2={\mathbb E}\xi_{k,1}^2+\ldots+{\mathbb E}\xi_{k,k-1}^2+2\sum\limits_{i<j<k}{\mathbb E}(\xi_{k,i}\xi_{k,j})=
$$
$$
=\frac{k-1}{2}+2\cdot\frac{(k-1)(k-2)}{2}\cdot\frac{1}{3}=\frac{2k^2-3k+1}{6}.
$$
Более того, с.в. $\xi_k$ и $\xi_m$ некоррелированы, $k\ne m$.

\item 
Для каждого $k$ характеристическая функция с.в. $X_k=\xi_k-{\mathbb E}\xi_k$ имеет вид 
\[
\varphi_{X_k}(t)=1-\frac{t^2 \Var\xi_k}{2}+ o(t^2), \, \text{при малых } t.
\]

\end{enumerate}
\end{remark}

\begin{problem}(Асимптотическое распределение числа циклов в случайной перестановке.)
\label{permutation}
Перестановка $\pi $ на $n$ элементах задается пошагово следующим образом: $\pi (k)$ выбирается случайно равновероятно из $\left\{1,\ldots ,n\right\}\backslash \left\{\pi (1),\ldots ,\pi (k-1)\right\}$. Ясно, что вероятность получения фиксированной перестановки будет $(n!)^{-1} $. Пусть с.в. $\xi _{kn} $, $1\le k\le n$, равна 1, если в качестве $\pi (k)$ выбран элемент, замыкающий какой-либо цикл в перестановке, и равна 0 в противном случае (несложно заметить, что $X_n = \sum\limits_{k=1}^n \xi_{kn}$ равняется количеству циклов в случайно выбранной перестановке).\\
\indent а) Покажите, что величины $\xi_{kn}$ являются попарно независимыми, т.е. если $A_k$ -- множество перестановок, у которых на $k$-м шаге заканчивается цикл (см.~выше), то $A_i$ и $A_j$ -- независимы для $i < j$. Покажите также, что $P\left(A_k\right) = \left(n-k+1\right)^{-1}$.\\
\indent б) Воспользовавшись ЦПТ в форме Ляпунова, покажите, что при $n\to \infty$ число циклов в случайной перестановке $$X_{n} =\sum\limits_{k=1}^{n}\xi _{kn}\sim \mathcal{N}\left(\ln(n)+\gamma, \ln(n)+\gamma - \dfrac{\pi^2}{6}\right),$$ где $\gamma = 0.57721\dots$ -- константа Эйлера.
\end{problem}

\begin{remark}
Напомним, как выглядит ЦПТ в форме Ляпунова. Пусть выполнены условия Линденберга (см. задачу 5) и для с.в. $\{X_i\}_{i=1}^n$ выполнено $\mathbb{E}|X_i|^3 < +\infty$. Тогда определена последовательность
$$
    r_n^3 = \sum\limits_{i=1}^n \mathbb{E}(|X_i - \mathbb{E}X_i|^3).
$$
Если предел $\lim\limits_{n\rightarrow \infty}\dfrac{r_n}{s_n} = 0$ (условие Ляпунова), то
$$
    \dfrac{S_n - m_n}{s_n} \xrightarrow{d} Z, \, \text{где }Z\sim\mathcal{N}(0,1) ,\, n\rightarrow \infty,
$$
где 
\begin{align*}
    &S_n = \sum\limits_{i=1}^n X_i, \, m_n = \sum\limits_{i=1}^n \mathbb{E}X_i ,\\
    &s_n = \sum\limits_{i=1}^n \sigma_i^2, \, \sigma_i^2 = \mathbb{D}\left(X_i\right).
\end{align*}

Интересные результаты об асимптотических свойствах группы перестановок содержатся в работах А.М.~Вершика, А.А.~Шмидта (см.~ссылки далее).
В качестве модели перестановки $g$ длины $n$ можно использовать последовательность чисел $i_1=1,i_2,\dots,i_n$, сгруппированных в циклы, причем каждый цикл начинается с наименьшего\linebreak в нем числа, и концы циклов указаны. Упорядочим циклы по возрастанию наименьших чисел. 
Пример множеств $A_k$: для перестановки $g =$\linebreak $= ((1,5),(2,3,8,6),(4,7))$, $g\in A_2 \bigcap A_6 \bigcap A_8$ (см.~презентацию Вершика А.М. \verb|http://www.mathnet.ru:8080/PresentFiles/231/v231.pdf|).

Определим отображение $g\to (x_1(g),\dots,x_n(g))$, где $x_i(g)$ -- длина $i$-го цикла, нормированная на $n$, которая переводит меру, заданную на перестановках (по условию вероятность получения фиксированной перестановки будет $(n!)^{-1}$),  в дискретную меру на единичном  симплексе. Оказывается, у такой последовательности мер при $n$, стремящемся к бесконечности, существует слабый  предел. Предельная мера порождается случайным рядом $\eta_i$, $i=1,2,\dots,$ с неотрицательными значениями и суммой $1$, причем значения $\eta_k/(\eta_k+\eta_{k+1}+\dots)$ независимы и равномерно распределены  на отрезке $[0,1]$.\\
\indent Иллюстрацией к приведенному результату является задача о ``ломании палочки'' (см. также Часть~2). Отрезок делится с равномерной вероятностью (точка деления, распределенная равномерно на $[0,1]$). Левый отрезок фиксируется и затем ломается правый отрезок. Точки разлома выбираются последовательно с помощью последовательности независимых равномерно распределенных на $[0,1]$ случайных величин 
В работах  А.М.~Вершика показано, что такое случайное разбиение отрезка задает распределение длин отрезков, совпадающее с распределением длин циклов в рассмотренной модели перестановок. 

Доказательство этого факта технически нетривиально, например, используется принцип инвариантности 
Донскера--Прохорова.  
Подробные доказательства имеются в статьях: Вершик~А.М., Шмидт~А.А.  Предельные меры, возникающие в асимптотической теории симметрических групп.~I //   ТВП. 1977.  Т.~22:1. С. 72–-88;  Вершик~А.М., Шмидт~А.А. Предельные меры, возникающие в асимптотической теории симметрических групп.~II //  ТВП. 1978. Т.~23:1.  С.~42–-54.

\end{remark}

\begin{problem}(Предельные меры, А.М. Вершик и др.) 
\label{permutation1}
В качестве множества элементарных исходов рассматривается группа всевозможных 
подстановок (перестановок) $\mathbb{S}_n $ (симметрическая группа), $n\gg 1$. В этой 
группе $n!$ элементов. Припишем каждой подстановке одинаковую вероятность $1 
\mathord{\left/ {\vphantom {1 {n!}}} \right. \kern-\nulldelimiterspace} 
{n!}$.

\begin{enumerate}

\item Покажите, что математическое ожидание числа циклов есть\linebreak $\approx 
\ln n$ (см. также задачу \ref{permutation}).

\item\Star В каком смысле нормированные длины циклов случайной подстановки 
убывают со скоростью геометрической прогрессии со знаменателем $e^{-1}$ ?

\item\Star Положим $$\rho _n \left( a \right)={\left| {\left\{ {g\in \mathbb{S}_n:\;n_{\max } \left( g \right)\le an} \right\}} \right|} \mathord{\left/ 
{\vphantom {{\left| {\left\{ {g\in S_n :\;n_{\max } \left( g \right)\le an} 
\right\}} \right|} {n!}}} \right. \kern-\nulldelimiterspace} {n!},$$ где 
\mbox{$n_{\max } \left( g \right)$~-- длина} максимального цикла в подстановке $g$. 
Покажите, что $\rho _n \left( a \right)$ удовлетворяет \textit{уравнению Дикмана--Гончарова} (40-е годы XX 
века):
\[
\rho _n \left( a \right)=\int\limits_0^a {\rho _n \left( {\frac{a}{1-t}} 
\right)dt}.
\]
\item\Star\,\,\Star Покажите, что начиная с некоторого большого числа $N$ 99{\%} 
натуральных чисел $n$, больших, чем $N$, обладают свойством
\[
n^{0.99}<p_1 \cdot ...\cdot p_{11} ,
\]
где $p_1, \dots, p_{11}$ -- первые наибольшие 11 простых делителей $n$. Наибольшие простые делители
$\{p_i\}_{i=1}^k$ числа $n$ определяются единственным образом через разложение на простые множители 
$$
    n = p_1\cdot p_2\dots p_k, \, \text{где } p_1\geq p_2\dots\geq p_k,
$$
при условии, что $p_1, p_2,\dots, p_k$ -- простые числа (Пример: $36 = 3\cdot 3 \cdot 2 \cdot 2$, тогда по определению $p_1 = 3,\, p_2=3, \, p_3 = 2, \, p_3 = 2$).
\end{enumerate}
\end{problem}
\begin{ordre} Решение задач всех пунктов сводится к задаче о ``ломании палочки'' (см. замечание  к предыдущей задаче). В пункте (г) условие можно переформулировать следующим образом: у (99{\%}) натуральных чисел основная часть 
(99{\%}) числа есть произведение первых 11-ти наибольших простых делителей. Вопрос этого же пункта  можно свести 
к случаю симметрической группы, если заметить, что разложение на набольшие простые делители эквивалентно разложению перестановки на упорядоченное произведение циклов. Чтобы сведение было корректным, надо правильно определить понятие длины ``цикла'' в таком разложении числа.
Это должна быть некоторая аддитивная монотонная, неотрицательная функция циклов $\{p_i\}_{i=1}^k$ (Например: $\log n = \log (p_1) + \dots + \log (p_k)$).
\end{ordre}

\begin{remark} См. 
Вершик А.М., Шмидт А.А.  Предельные меры, возникающие в асимптотической теории симметрических групп.~I //   ТВП. 1977.  Т.~22:1. С. 72–-88;
Вершик А.М., Шмидт А.А.  Предельные меры, возникающие в асимптотической теории симметрических групп.~ II //  ТВП. 1978. Т.~23:1.  С. 42--54; Вершик А.М. Асимптотическое 
распределение разложений натуральных чисел на простые делители // ДАН. 1986. 
Т. 289:2. С. 269--272; Tenenbaum G. Introduction to analytic and 
probabilistic number theory. Cambridge Univ. Press, 1995; 
Arra\-tia~R., 
Ba\-rbour~A.D., Tavar\'{e}~S. Logarithmic combinatorial structures: A~pro\-ba\-bi\-lis\-tic approach. Series ``EMS Monographs in Mathematics''. Z\"{u}rich: Eur. Math. Soc.,  
2003.

Можно более точно сформулировать результаты этой задачи и предыдущей. Для этого определим распределение Пуассона--Дирихле (см. также задачу \ref{sobord}) на множестве числовых рядов (для наглядности можно сичтать, что члены каждого ряда отсоритрованы по убыванию) с неотрицательными элементами, сумма членов которых равна 1. Другми словами рассматривается бесконечномерный единичный симплекс. На симплексе в размерности $n$ можно задать распределение Дирихле $\Dir(\lambda/n,...,\lambda/n)$ (см. задачу \ref{laplace} раздела \ref{standart}). Под \textit{распределением Пуассона--Дирихле} $PD(\lambda)$ с параметром $\lambda$ на бесконечномерном единчном симплексе понимается (слабый) предел распределений $\Dir(\lambda/n,...,\lambda/n)$ при $n\to\infty$. 

Рассмотрим разложение натурального числа $n$ на простые множители
$$
n = p_1\cdot p_2\dots p_{k(n)}, \, \text{где } p_1\geq p_2\dots\geq p_{k(n)}.
$$
Выберем случайно $n$ в соответствии с равномерным распределением на множсетве натуральных чисел от $1$ до $N$ и рассмотрим случайный вектор
$$v_n = \left(\frac{\ln p_1}{\ln n},...,\frac{\ln p_k(n)}{\ln n}\right)^T.$$
Оказывается, этот вектор слабо сходится при $N\to\infty$ к распределению $PD(1)$.

Аналогично, на множестве всех перестановок из $N$ элементов введем равномерную меру, затем выбирем случайную перестановку $n$ ($1\le n \le N!$) и рассмотрим соответсвующий ей случаный вектор (здесь $l_i$ -- длина $i$-го по величине цикла в выбранной перестановке, а $k(n)$ -- число циклов в выбранной перестановке)
$$
v_n = \left(\frac{l_1}{n},...,\frac{l_{k(n)}}{n}\right)^T.
$$
Тогда $v_n$ слабо сходится при $N\to\infty$ к распределению $PD(1)$.

Описанный в предыдущей задаче марковский процесс ``ломания палочки'' имеет своим слабым пределом также $PD(1)$ и позволяет понять свойства распределения $PD(1)$ (см., например, распределение максимального элемента -- Дикмана--Гончарова, убывание для типичных представителе $PD(1)$ элементов ряда со скоростью геометрической прогрессии со знаменталем $e^{-1}$). 

\end{remark}

\begin{problem}(Закон повторного логарифма.)
Пусть $X_n$~-- независимые одинаково распределённые с.в. с нулевым математическим ожиданием и единичной дисперсией. Пусть $S_n = X_1+\ldots+ X_n$. Докажите, что почти наверное
\[\underset{n \rightarrow \infty}{\overline{\lim} } \frac{S_n}{\sqrt{n \log_2 (\log_2 n)}} = \sqrt{2},\]
\[\underset{n \rightarrow \infty}{\underline{\lim} } \frac{S_n}{\sqrt{n \log_2 (\log_2 n)}} = -\sqrt{2}.\]
\end{problem}

\begin{remark} 

Хотя величина $\frac{S_n}{\sqrt{n \log_2 \log_2 n}}$ будет меньше, чем любое заданное  $\varepsilon$  с вероятностью, стремящейся к единице (это следует из~ЦПТ), она будет бесконечное число раз приближаться сколь угодно близко к любой точке отрезка [$-\sqrt{2}, \sqrt{2}$] почти наверное [\ref{chiraiev}] Т.~2.

Приведем также некоторые смежные результаты из \cite{Gupta}. Пусть $X_1$, $X_2$, \dots -- независимые одинаково распределенные с.в. с общей функцией распределения $F$ и пусть $\gamma_n$ -- расходящаяся числовая последовательность (к $+\infty$). Пусть $Z_{n,\gamma}  = \frac{S_n}{\gamma_n}$, $n\geq 1$, и $B(F,\gamma)$ -- множество всех предельных точек $Z_{n,\gamma}$. Тогда существует неслучайное множество $A(F,\gamma)$, такое, что с вероятностью $1$ $B(F,\gamma)$ совпадает c $A(F,\gamma)$.  В~частности:
\begin{itemize}
\item Если $\gamma_n = n^{\alpha}$, $0<\alpha<1/2$, тогда для всех $F$, не вырожденных в $0$, $A(F,\gamma)$ равно всей расширенной вещественной оси, если оно содержит хотя бы одно конечное вещественное число.
\item Если $\gamma_n = n$ и если $A(F,\gamma)$ содержит хотя бы два конечных вещественных числа, то тогда оно должно содержать $\pm\infty$.
\item Если $\gamma_n=1$, тогдa $A(F,\gamma) \in \{\pm\infty\}$ тогда и только тогда, когда для некоторого $a>0$, $\int_{a}^a \text{Re}\{1/(1-\psi(t))\}\,dt\leq \infty$, где $\psi(t)$ -- характеристическая функция $F$.
\item Пусть $\mathbb{E}(X_1)=0$, $\mathbb{D}X_1 = \sigma^2\leq{\infty}$. Тогда при $\gamma_n=$\linebreak $=\sqrt{2n\log_2(\log_2 n)}$, $A(F,\gamma) = [-\sigma,\sigma]$, а если $\mathbb{E}(X_1)=0$, $\mathbb{D}X_1 = \infty$, и выбрать $\gamma_n=\sqrt{2n\log_2(\log_2 n)}$, то $A(F,\gamma)$ содержит по крайней мере одну из $\pm\infty$.
\end{itemize}

Пусть $m_n$ таковы, что $\PR(S_n\leq m_n)\geq 1/2$ и $\PR(S_n\geq m_n)\geq 1/2$. Тогда для того, чтобы существует такая положительная последовательность $\gamma_n$, удовлетворяющая почти наверное соотношениям
$$
-\infty< \underset{n \rightarrow \infty}{\underline{\lim} } \frac{S_n-m_n}{\gamma_n}<
\underset{n \rightarrow \infty}{\overline{\lim} }\frac{S_n-m_n}{\gamma_n}<\infty
$$ необходимо и достаточно, чтобы для всех $c\geq 1$
$$
\underset{x \rightarrow \infty}{\underline{\lim} }\frac{\PR(|X_1|>cx)}{\PR(|X_1|>x)}\leq c^{-2}.
$$
\end{remark}

\begin{problem}
\label{lognormal}
Рассмотрим с.в. с лог-нормальным распределением (см. задачу~\ref{lognorm} из раздела \ref{hard}) \[
X \in \text{Log}\N(m, \sigma^2) \Leftrightarrow \log (X) \in \N(m, \sigma^2). \] 
При каких значениях $\sigma$ лог-нормальное распределение похоже на степенное (при каких значениях  $x$)?
Покажите, что для произведения с.в. характерна сходимость к лог-нормальному распределению. Докажите это же свойство для длин отрезков, получающихся в результате деления единичного отрезка с равномерным многократным выбором точек раздела (см. задачу \ref{permutation} и \ref{permutation1}).     
\end{problem}


\begin{problem}(Пуассоновский процесс \cite{1}.)
\label{sec:poisson}
Пусть необходимо оценить, сколько билетов на метро одного вида $K(T)$ 
продается за одну рабочую смену длительностью $T$. Имеет место формула 
$$
K(T)=\max\Bigl\{ n:\; \sum\limits_{k=1}^{n} X_k<T \Bigr\} , 
$$
где $X_1, X_2, X_3,\ldots$~-- независимые одинаково распределенные по закону ${\rm Exp}(\lambda)$ с.в. ($X_k$ интерпретируется как время между 
$k-1$ и $k$ сделкой (продажей)). Покажите, что 
\begin{enumerate}
\item вероятность ${\mathbb P}(K(T+t)-K(T)=k)$, где $t\geqslant 0$ и $k=0,1,2,\ldots$, не зависит от $T\geqslant 0$; 

\item $\forall\,n\geqslant 1$, $0\leqslant t_1\leqslant t_2\leqslant \ldots\leqslant t_n$ 
 с.в. $\bigl\{ K(t_k)-K(t_{k-1})\bigr\}_{k=1}^{n}$ независимы; 

\item ${\mathbb P}(K(t)>1)=o(t),\quad t>0$; 

\item $\forall\,n\geqslant 1$, $0\leqslant t_1\leqslant t_2\leqslant \ldots\leqslant t_n$, 
$0\leqslant k_1\leqslant k_2\leqslant \ldots\leqslant k_n$, $k_1,\ldots, k_n\in {\mathbb N}\cup \{ 0\}$ 
\begin{multline*}
{\mathbb P}(K(t_1) = k_1,\ldots, K(t_n)=k_n)= \\
=e^{-\lambda t_1} \frac{(\lambda t_1)^{k_1}}{k_1!}\cdot 
e^{-\lambda(t_2-t_1)} \frac{(\lambda(t_2-t_1))^{k_2-k_1}}{(k_2-k_1)!}\cdot \ldots \\
\ldots \cdot e^{-\lambda(t_n-t_{n-1})}\frac{(\lambda(t_n-t_{n-1}))^{k_n-k_{n-1}}}{(k_n-k_{n-1})!} , 
\end{multline*}
в частности, $K(T)\in \Po(\lambda T)$. 

\end{enumerate}
\end{problem}
\begin{remark}
Обратите внимание на задачу \ref{exp_eps} из раздела 2.
\end{remark}

\begin{problem}(Сложный пуассоновский процесс \cite{1}, \cite{51}.)
\label{sec:cpoisson}
В течение рабочего дня фирма осуществляет $K(T)\in \Po(\lambda T)$ сделок ($K(T)$ -- с.в., имеющая распределение Пуассона с параметром 
$\lambda T$, где $\lambda = 100~[\text{сделок/час}]$). Каждая сделка приносит доход $V_n\in R[a,b]$ ($V_n$ -- с.в., имеющая 
равномерное распределение на отрезке $[a,b]=[10, 100]$, $n$ -- номер сделки). Считая, что $K$, $V_1$, $V_2$, $\ldots$ -- 
независимые в совокупности с.в., найдите математическое ожидание и дисперсию выручки за день $Q(T)=\sum\limits_{k=1}^{K(T)} V_k$. Докажите следующее соотношение для характеристической функции $Q(T)$:
$$
\varphi_{Q(T)}(s)=\exp\{ \lambda T(\varphi_{V_k}(s)-1)\} = \exp\left\{\lambda T\left(\int{e^{isx}dF_{V_k}(x)} - 1\right)\right\}. 
$$
Найдите (см. \cite{5}) такие $m(T)$ и $\sigma(T)$, что при $T\to\infty$
$$\frac{Q(T)-m(T)}{\sigma(T)}\xrightarrow{d} Z, \, \text{где }Z\sim\mathcal{N}(0,1).$$
\end{problem}

\begin{ordre}

Примените формулу для условного математического ожидания (см. задачу \ref{cond} раздела \ref{hard})
$$
{\mathbb E}Y= \Exp{( {\mathbb E}(Y|X) )}. 
$$
Установите справедливость следующего соотношения:
$$
\Var Y=\Var({\mathbb E}(Y|X))+{\mathbb E}(\Var(Y|X)) . 
$$

\end{ordre}

\begin{problem}
В течение трех лет фирма из предыдущей задачи работала\linebreak $N=1000$ дней (длина рабочего дня и параметры спроса не менялись). 
Оцените (см. \cite{5}) распределение с.в. 
$$Q^N=\sum\limits_{k=1}^{N} Q_k(T), 
$$
где $Q_k(T)$ -- выручка за $k$-й день. Верно ли, что с.в. $Q^N$ и $Q_k(NT)$ одинаково распределены? 
\end{problem}

\begin{problem}(Пуассоновский поток событий \cite{27}, \cite{202}.)
\label{Poisson}
Рассмотрим интервал 
$\left[ {-N,N} \right]$ и бросим на него независимо и случайно (равномерно) \mbox{$M=\left[ {\rho N} \right]$} точек, где $\rho >0$ -- некоторая 
константа, называемая плотностью. Легко вычислить биномиальную вероятность 
$\PR_{N,M} \left( {k,I} \right)$ того, что в конечный интервал $I\subset 
\left[ {-N,N} \right]$ попадет ровно $k$ точек. Покажите, что для $\PR_{N,M} 
\left( {k,I} \right)$  при $N\to \infty $ верно следующее
\[
\PR\left( {k,I} \right)\mathop =\limits^{def} \mathop {\lim }\limits_{N\to 
\infty } \PR_{N,M\left( N \right)} \left( {k,I} \right)=\frac{\left( {\rho 
\left| I \right|} \right)^k}{k!}e^{-\rho \left| I \right|},
\quad
k=0,1,...
\]
Покажите также, что если $I_1 ,I_2 \subset \left[ {-N,N} \right]$ и $I_1 
\cap I_2 =\emptyset $, то
\[
\PR \left( {k_1 ,I_1 ;k_2 ,I_2 } \right)\mathop =\limits^{def} \mathop {\lim 
}\limits_{N\to \infty } \PR_{N,M\left( N \right)} \left( {k_1 ,I_1 ;k_2 ,I_2 } 
\right)=\PR\left( {k,I_1 } \right)\PR\left( {k,I_2 } \right).
\]
\end{problem}

\begin{problem}(Безгранично делимые с.в.)
\label{sec:infdiv}
Случайная величина $X$ называется \textit{безгранично делимой}, если для любого натурального $n$ найдутся $n$ таких независимых одинаково распределенных с.в. $X_{kn}$ (распределение $X_{kn}$ зависит от $n$), что $X = \sum_{k=1}^{n}{X_{kn}}$, где равенство понимается по распределению. Теорема Колмогорова--Леви--Хинчина утверждает, что $X$ -- безгранично делимая с.в. тогда и только тогда, когда существуют (но не обязательно однозначно определяются \cite{28}) 
$$
 (b,c,\nu(dx)): \; 
c\geqslant 0, \;  \nu(dx)\geqslant 0,\; \nu(\{0\})=0 \text{, такие что} 
$$
\[
\varphi_{X}(s)
=\exp \left\{  is b-\frac{c s^2}{2}+
\int_{-\infty}^{\infty} \biggl( e^{i s x}-1-\frac{i s x}{x^2+1}\biggr)\frac{x^2+1}{x^2}\, \nu(dx) 
\right\}, 
\]
где $\varphi_{X}(s)$~-- характеристическая функция $X$.

Верно ли, что любая безгранично делимая с.в. может быть представлена (имеет такое же распределение) как $X + Q$, где $X\sim \mathcal{N}(m,
\sigma^2)$, а $Q$ -- с.в., имеющая сложное распределение Пуассона?
Определите $(b,c,\nu(dx))$ для $Q(T)$ из задачи 
\ref{sec:cpoisson}.
\end{problem}
\begin{ordre}
См. замечние к задаче \ref{Levi}.
\end{ordre}
\begin{remark} 
В литературе часто можно встретить и другие записи теоремы Колмогорова--Леви--Хинчина, например, 
\[
\varphi_{X}(s)
=\exp \Bigl\{  i s b_0-\frac{c s^2}{2}+
\int_{-\infty}^{\infty} \bigl( e^{i s x}-1-i\frac{s|x|}{1+x^2} \bigr)\, \nu(dx) 
\Bigr\}.
\]
Здесь и далее (в других представлениях) используются разные $\nu(dx)$, однако обозначать их будем одинаковым образом.
Другой пример записи -- представление Леви--Хинчина (см. также задачу \ref{sobord}):
\[
\phi _{X} \left( s \right)= \exp \left\{ {i s b-{c s^2} \mathord{\left/ {\vphantom 
{{cs^2} 2}} \right. \kern-\nulldelimiterspace} 2+\int_{-\infty }^\infty 
{\left( {e^{is x}-1-i s xI\left( {\left| x \right|<1} \right)} \right)\nu 
\left( {dx} \right)} } \right\},
\]
где $I(\text{true}) = 1$, $I(\text{false}) = 0$. Еще один пример будет приведен ниже в этом замечании.

В качестве нетривиальных примеров безгранично делимых распределений упомянем следующие: распределение Стьюдента, Коши, гамма-распределение, лог-нормальное, показательное. Последнее является примером безгранично де\-ли\-мо\-го распределения, которое может быть представлено в виде бесконечной свертки (независимых, одинаково распределенных с.в.) не безгранично делимых с.в. Для нормального  распределения и распределения Пуассона такое представление невозможно. Теоремы Крамера и Райкова говорят соответственно, что если свертка двух распределений нормальна (распределена по закону Пуассона), то компоненты должны также иметь нормальное распределение (распределение Пуассона) с другими параметрами. Подробнее об этом см. \cite{stoianov}. 

Основное свойство безгранично делимых с.в. заключается в том, что такие с.в. (и только такие) могут возникать в качестве пределов по распределению сумм с.в. $\sum_{k=1}^{n}{X_{kn}}$, где для любого $n$ с.в. $X_{kn}$ -- независимы и одинаково распределенны (как именно, вообще говоря, заивисит от $n$). 

Приведем для справки также несколько фактов из \cite{Gupta}.
\begin{enumerate}
\item
Пусть имеется безгранично делимая с.в. со средним $m$ и конечной дисперсией $D$ и  $\phi(s)$ -- характеристической функцией. Тогда 
$$
\ln\phi(s)=ims+\int_{-\infty}^{\infty}({\rm e}^{isx}-1-isx)\frac{\nu(dx)}{x^2},
$$
где $\nu$ -- конечная мера на действительной оси, более того, $\nu(\mathbb{R}) = D$.
\item {\bf Теорема Голди--Стьётела.} {\itПусть положительная с.в. имеет плотность распределения $f(x)$, причем $f$ бесконечно дифференцируема и $(-1)^k f^{(k)}(x)\geq 0$, $\forall x>0$. Тогда $X$ имеет безгранично делимое распределение.} 
\end{enumerate}

Также следует обратиться к книге Гнеденко Б.В., Колмо\-го\-ров~А.Н. Предельные распределения для сумм независимых случайных величин.  М.--Л.: ГТТИ, 1949.

\end{remark}

\begin{problem}{\bf(ЦПТ без требования существования дисперсий.)}
{\itПусть $X_1,X_2, 
\dots$ -- независимые одинаково распределенные случайные величины из распределения $F$ на действительной прямой.  Верна следующая теорема 
\cite{Gupta}: cуществуют такие последовательности $\{a_n\},\{b_n\}$, что при $n\to \infty$ $$\frac{\sum_{i=1}^n(X_n-a_n)}{b_n}\xrightarrow{d} Z, \, \text{где }Z\sim\mathcal{N}(0,1)$$ тогда и только тогда, когда функция $v(x) = \int_{[-x,x]}y^2\,dF(y)$ является медленно меняющейся на бесконечности, то есть для нее верно при любом $t>0$, что} $$\lim_{x\to\infty}\frac{v(tx)}{v(x)}=1.$$
\noindent Покажите, что для  распределения Стьюдента с двумя степенями свободы выполняется приведенная теорема. Найдите нормировочные последовательности $\{a_n\}$ и $\{b_n\}$.
\end{problem}
 \begin{remark}
 Напомним, что по определению распределение Стьюдента c $k$ степенями свободы есть распределение случайной величины $$t = \frac{\xi_0}{\sqrt{\frac{1}{k}(\xi^2_1+\dots+\xi^2_k)}},$$
 где $\xi_j$ независимы, $\xi_j\in \mathcal{N}(0,1)$, $j=0,\dots,k$. Плотность распределения Стьюдента 
 выражается следующей функцией:
 $$
 p(x) = \frac{\Gamma((k+1)/2)}{\sqrt{\pi k}\Gamma(k/2)} \left(1+\frac{x^2}{k}\right)^{-(k+1)/2},
 \, x\in \mathbb{R}^1.
 $$
 Можно также показать, что при $k\rightarrow \infty$ плотность распределения $p(x)$ сходится к плотности стандартного нормального распределения (т.е. нормальное распределение можно рассматривать как распределение Стьюдента с бесконечным числом степеней свободы $k = \infty$).
 \par Распределение Стьюдента крайне важное распределение, которое, в частности, возникает во многих классических задачах статистики,  см. например \cite{52}.
 \end{remark}

\begin{problem}\label{bluzd_ust}
Рассмотрим простую и классическую схему блуждания точки на прямой, соответствующую правилам игры в орлянку:
\[\eta (0)=0,\] 
\[\eta (t+1)=\left\{\begin{array}{cc} {\eta (t)+1,} & {p={1\mathord{\left/ {\vphantom {1 2}} \right. \kern-\nulldelimiterspace} 2} }, \\ {\eta (t)-1,} & {p={1\mathord{\left/ {\vphantom {1 2}} \right. \kern-\nulldelimiterspace} 2} }. \end{array}\right. \] 
Занумеруем в порядке возрастания все моменты времени, когда $\eta (t)=0$. Получим бесконечную последовательность $0=\tau _{0} <\tau _{1} <$\linebreak $<\tau _{2} <...$ Рассмотрим разности $\xi _{i} =\tau _{i} -\tau _{i-1} $ -- последовательность независимых одинаково распределенных с.в.
\begin{enumerate}
\item Покажите, что $\PR\left(\xi _{i} =2m\right)\approx\frac{1}{\sqrt{\pi m (2m - 1)^2}}\sim m^{-3/2}$.
\item Покажите, что математическое ожидание с.в. $\xi _{i} =\tau _{i} -\tau _{i-1} $ равно бесконечности.
\end{enumerate}

\begin{remark}
Этот результат можно проинтерпретировать так: среднее время до первого возвращения блуждания в 0 бесконечно. 
Тем не менее суммы $\tau _{n} =\sum _{i=1}^{n}\xi _{i}  $ при надлежащей нормировке имеют предельное распределение. Имеет место такой результат: 
\[\mathop{\lim }\limits_{n\to \infty } \PR\left\{\frac{2\tau _{n} }{\pi n^{2} } <z\right\}=\left\{\begin{array}{cc} {\frac{1}{\sqrt{2\pi } } \int _{0}^{z}e^{-\frac{1}{2x} } x^{-\frac{3}{2} }  dx,} & {z>0}, \\ {0,} & {z\leq 0}. \end{array}\right. \]  
Случайная величина $\zeta$ имеет устойчивый закон распределения, если для любого натурального $N>1$ найдутся независимые случайные величины $\{Z_n\}_{n=1}^N$, распределение которых совпадает с  распределением $\zeta$  и постоянные $a_N$ и $b_N$, такие что $\zeta = \frac{1}{a_N}(\sum_{n=1}^N Z_n - b_N)$. Имеет место следующее утверждение (см. \cite{21}, Т.1): $\zeta$ может быть пределом по распределению с.в. $\frac{1}{a_N}(\sum_{n=1}^N Z_n - b_N)$, где $Z_n$ -- независимые одинако распределенные с.в. (с распределением не зависящим от $N$), тогда и только тогда, когда $\zeta$ имеет устойчивый закон распределения. Также верна теорема о каноническом представлении для устойчивых законов  Леви--Хинчина \cite{1}: для того чтобы с.в. была устойчивой, необходимо и достаточно, чтобы логарифм ее характеристической функции представлялся формулой
\[
\ln \varphi (s)=i\gamma s - c|s|^{\alpha } \left(1+i\beta \text{sign}(s) \omega (s,\alpha )\right),\] 
где $\gamma\in\mathbb{R}$, $-1\le \beta \le 1$, $0<\alpha \le 2$, $c\ge 0$, $\text{sign}(0) = 0$ и
\[\omega (s,\alpha )=\left\{\begin{array}{cc} {\tg\left(\frac{\pi }{2} \alpha \right),} & {\alpha \ne 1,} \\ {\frac{2}{\pi } \ln |s|,} & {\alpha =1}. \end{array}\right.
\]

Полученное выше распределение соответствует каноническому представлению с $\alpha ={1\mathord{\left/ {\vphantom {1 2}} \right. \kern-\nulldelimiterspace} 2} $, $\beta =1$, $\gamma =0$, $c=1$, и принадлежит семейству кривых Пирсона (показано Н.В. Смирновым).

Особенно просто устроены характеристические функции симметричных устойчивых распределений:
$$
\ln \varphi (s) =  - c|s|^{\alpha }. 
$$

П. Леви также принадлежит немного другой способ классификации устойчивых распределений. А именно, независимые одинаково распределенные с.в. $X_1$, $X_2$, $X$ устойчивы по Леви ($L$-устойчивые), если для любых $s_1\ge 0$, $s_2\ge 0$, существует такое $\alpha \in \left(0,2\right]$, что $s^{\alpha} = s_1^{\alpha} + s_2^{\alpha}$ и $s^{\alpha}X = s_1^{\alpha}X_1 + s_2^{\alpha}X_2$, где последнее равенство понимается как равенство по распределению (с.в. в левой части равенства имеет такое же распределение, как с.в. в правой части равенства). Случай $\alpha = 2$ определяет нормальные (гауссовские) с.в., а случай $\alpha = 1$ -- с.в., распределеные по закону Коши (см., например, задачу \ref{cauchy_rad_emitt} раздела \ref{hard}). Этими двумя случаями исчерпываются ситуации, когда можно выписать явную (аналитическую) формулу для плотности (в симметричном случае -- в несимметричном случае имеется еще одна аналитическая ситуация, которая была разобрана выше, в связи с исследованием времени возвращения при игре в орлянку). Тем не менее, при больших $|x|$ и $\alpha\in(0,2)$ имеет место следующая асимптотика для функции плотности распределения:
$$
p_{\alpha,\beta,\gamma,c}(x)\approx \frac{c^{\alpha}(1 + \beta)\sin(\pi \alpha/2)\Gamma(\alpha + 1)}{\pi |x-\gamma/c|^{\alpha + 1}}.
$$

Более того, имеет место следующий результат: пусть последовательность независимых одинаково распределенных с.в. $\left\{Z_n\right\}_n$ имеет степенные хвосты $\PR(Z_n > x)\approx C_{+}x^{-\alpha}$ и $\PR(Z_n < -x)\approx C_{-}x^{-\alpha}$, при $x\to\infty$, где $\alpha \in (0,2)$, $C_{+}, C_{-} > 0$. Тогда существуют такие последовательности $\left\{a_n\right\}_n > 0$ и $\left\{b_n\right\}_n$ (их определение дано, например, в \cite{202}), что $\frac{1}{a_N}(\sum_{n=1}^N Z_n - b_N)$ слабо сходится при $N\to\infty$ к с.в. $\xi_{\alpha,\beta,0,c}$ c $\beta = \left(C_{+} - C_{-}\right)/\left(C_{+} + C_{-}\right)$, $c = C_{+} + C_{-}$, а $\xi_{\alpha,\beta,0,c}$ имеет устойчивый закон распределения с соответствующими значениями параметров (добавим, что $\gamma = 0$).

Для решения этой задачи и следующих двух задач также полезно познакомиться с учебным пособием \cite{28}. Для более глубокого погружения рекомендуется книга: Пет\-ров~В.В. Предельные теоремы для сумм независимых случайных величин.  М.: Наука, 1987. 
\end{remark}
\end{problem}

\begin{problem}
Устойчивые законы  распределения являются также и  безгранично делимыми. 
 Покажите, что пуассоновское распределение безгранично делимо, но не устойчиво. 
\end{problem}
\begin{remark}
См. предыдущую задачу, а также \cite{stoianov}.
\end{remark}

\begin{problem}(Вывод распределения Хольцмарка \cite{28}, \cite{202}.)
\label{Holtsmark}
Рассмотрим шар \mbox{радиуса $R$} с центром в начале координат и $n$ звезд (точек), расположенных в нем случайно и независимо друг от друга. Пусть каждая звезда имеет единичную массу, и звезды распределены равномерно внутри данного шара.  Гравитационное поле $F_{i}$ (вектор), создаваемое $i$-ой звездой в начале координат, определяется по закону Ньютона--Гука:
\begin{equation*}
    F_i = G\dfrac{r_i}{\|r_i\|_2^3}, \, r_i = (x_i,y_i,z_i)\in \mathbb{R}^3,
\end{equation*}
где $(x_i,y_i,z_i)$-координаты $i$-ой звезды, $G$ -- гравитационная постоянная.
\par Суммарное гравитационное поле от всех звед можно выразить как $S_{n} =F_{1} +...+F_{n} $. Устремим $R$ и $n$ к бесконечности так, чтобы (см. также задачу \ref{Poisson}) $$\frac{4}{3} \pi R^{3} n^{-1} =\const = \rho >0,$$
где $\rho > 0$ -- параметр плотности звёзд (в объеме $V$ содержится в среднем $\rho V$ звёзд). 
\par Покажите, что предельная плотность распределения вектора $S(F)$ зависит только от ``модуля'' поля $\|F\|_2$ и стремится к плотности симметричного устойчивого распределения Хольцмарка с $\alpha = 3/2$.
Также покажите, что данное распределение близко к распределению силы притяжения в нуле, вызванной одной ближайшей звездой.
\end{problem}
\begin{remark}
Можно показать, что задача по существу не изменится, если массу каждой звезды считать с.в. моментом порядка $3/2$ равным $1$ и массы различных звезд предполагать взаимно независимыми с.в., не зависящими также от их расположения (см. \cite{28}).
\par Важно, что при выводе данного распределения пренебрегают размером звёзд (звезды -- материальные точки). Поскольку, основной вклад дает ближайшая звезда, то такое приближение в случае больших плотностей едва ли можно считать хорошим. Таким образом, данный результат является приближением, которое можно уточнить, если считать звезды шарами конечных размеров. В свою очередь, конструктивно описать, что получается в таком случае не представляется возможным; см. также книгу Кендалл М., Моран П. Геометрические вероятности.  М.: Наука, 1972.
\end{remark}

\begin{problem} Пусть $n$ единичных масс равномерно распределены в точках $X_1, \ldots, X_n$ на отрезке $[-n, n]$. На единичную массу в начале координат действует гравитационная сила $$f_n = \sum_{k =1}^n \frac{{\rm sign}(X_k)}{X_k^2}.$$

\noindent Покажите, что $$\mathbb{E}\left[\exp{\left(i s f_n\right)}\right] \rightarrow \exp{\left(-c\sqrt{|s|}\right)}, \quad n\to\infty.$$
\end{problem}
\begin{ordre}
Так как $X_k$ распределено равномерно, то (см. \cite{2013})
\begin{equation*}
\mathbb{E}\left[\exp{\left(is \frac{{\rm sign}(X_k)}{X_k^2}\right) }\right] = \int_{-n}^n \exp{\left(is\frac{{\rm sign}(x)}{x^2}\right) }\frac{dx}{2n} = \frac{1}{n}\int_{0}^n \cos{\left(\frac{s}{x^2}\right)}dx.
\end{equation*}
\end{ordre}

\begin{problem} \label{sobord} Покажите, что если $X\left( t \right)$ -- 
процесс Леви, то (см. замечание к задаче \ref{sec:infdiv})
\[
\exists \;\;\left( {b,\;c,\;\nu \left( {dx} \right)} \right),\,
c\ge 0,\,\nu \left( {dx} \right)\ge 0:
\quad
\int_{-\infty }^\infty {\min \left( {1,\;x^2} \right)\;}
\nu \left( {dx} 
\right)<\infty ,\]
\[
\forall \;\;t\ge 0 \quad 
\phi _{X\left( t \right)} \left( s \right)=\mathbb{E}{\rm e}^{i s X\left( t 
\right)}=
\]
\[=\exp \left\{ {t\left[ {i s b-{c s^2} \mathord{\left/ {\vphantom 
{{cs^2} 2}} \right. \kern-\nulldelimiterspace} 2+\int_{-\infty }^\infty 
{\left( {e^{is x}-1-i s xI\left( {\left| x \right|<1} \right)} \right)\nu 
\left( {dx} \right)} } \right]} \right\}.
\]
\end{problem}
\begin{remark}
 Процессом Леви  $\left\{ {X\left( t 
\right)} \right\}_{t\ge 0} $ называется стохастически непрерывный случайный 
процесс, удовлетворяющий следующим условиям (с небольшими оговорками (cadlag--процесс) -- см.  Applebaum~D. L\'{e}vy processes and stochastic calculus.   2nd ed. Cambridge University Press, 2009.):

\begin{enumerate}
\item $X\left( 0 \right)\mathop =\limits^{\text{п.н.}} 0;$
\item для любых $0\leq \tau < t$ распределение с.в. $X\left( t \right)-X\left( \tau \right)$ зависит только от $t-\tau$ (также говорят, что $\left\{ {X\left( t \right)} \right\}_{t\ge 0} $ имеет стационарные приращения или, что $\left\{ {X\left( t \right)} \right\}_{t\ge 0} $  однородный);
\item для любых $n\in \mathbb{N}$, $0\le t_0 \le t_1 \le ...\le t_n $ выполняется: $X\left( {t_1 } \right)-X\left( {t_0 } \right)$, {\ldots}, $X\left( {t_n } \right)-X\left( {t_{n-1} } \right)$ --- независимые в совокупности с.в. (также говорят, что $\left\{ {X\left( t \right)} \right\}_{t\ge 0} $ имеет независимые приращения).
\end{enumerate}
См.  Т.~1 книги Ширяев А.Н. Основы финансовой 
стохастической математики.  М.: МЦНМО, 2016.
\end{remark}

\begin{remark}\label{moran}

Отметим, что мера Леви $\nu$ может быть сигнулярна  в нуле, и такие процессы могут иметь бесконечно много  скачков на любом непустом временном интервале. 
Примером такого процесса является гамма-процесс (см., например, Cherny A.S. The Kolmogorov students' competition on probability theory. http://www-old.newton.ac.uk/preprints/NI05043.pdf), называемый иногда субординатором Морана \cite{202}.

Процесс Леви, с неубывающими траекториями,  называется субординатором. 
 Простейшими примерами неотрицательных процессов Леви являются пуассоновский процесс и сложный пуассоновский процесс. 
При этом скачки обычно происходят довольно редко, поэтому приращения для приведенных примеров чаще всего равны нулю даже на довольно больших промежутках времени. Мера Леви таких этих процессов конечна. 
 
В качестве примера процесса Леви с бесконечной мерой Леви можно привести субординатор с $\nu(dx) = x^{-1}{\rm e}^{-x}dx$ (субординатор Морана). 
Таким образом, субординатор Морана -- это процесс Леви, отвечающий (безгранично-делимому) гамма-распределению \cite{202}. Для него верно, что $\nu(\mathbb{R})=\infty$.

Этот процесс может быть использован для генерирования последовальностей с.в. с распределением Пуассона--Дирихле (см. задачу \ref{permutation1}). Это распределение упорядоченных по убыванию нормированных скачков субординатора Морана. Распределение Пуассона--Дирихле чрезвычайно часто возникает в разнообразных приложениях, особенно в популяционной генетике и экономике (оно является равновесным распределением для ряда эволюционных моделей). С другой стороны, распределение Пуассона--Дирихле появляется как предел распределения Дирихле на конечномерных симплексах при размерности симплекса, стремящейся к бесконечности, и одинаковыми параметрами, нормированными на размерность симплекса. Детали см. в книге \cite{202}.
\end{remark}

\begin{problem}(Модель Кокса--Росса--Рубинштейна.) 
\label{sec:levilong}
Пусть на ``идеализированном'' фондовом рынке имеется всего две ценные 
бумаги, и торговля осуществляется всего в два момента времени. Пусть цена 
первой бумаги $S$ (будем называть её акцией (stock)) известна в первый 
момент. Цена второй бумаги $C$ (будем называть ее call --- опционом 
европейского типа) не известна в первый момент. Пусть с ненулевой 
вероятностью $p>0$ (мы это $p$ не знаем, но от него ничего зависеть в итоге 
не будет) к моменту времени 2 цена акции вырастет в $u>1$ (up) раз и с 
вероятностью $1-p$ цена акции ``вырастет'' в $d<1$ (down) раз, т.е. упадет. 
Пусть также известны возможные цены опциона во второй момент: $C_u $, если 
акция выросла в цене, и $C_d $,  если акция упала в цене. Для простоты будем 
считать, что банк работает с нулевым процентом, т.е. класть деньги в банк, в 
расчете на проценты, бессмысленно. Говорят, что рынок 
безарбитражный, если не существует таких $k_S $, $k_C $, что
\[
X\left( 1 \right)=k_S S+k_C C=0,
\]
\[
\PR\left( {X\left( 2 \right)\ge 0} \right)=\PR\left( {k_S S\left( 2 \right)+k_C 
C\left( 2 \right)\ge 0} \right)=1,\] причем $\PR\left( {X\left( 2 \right)>0} 
\right)>0.
$

Докажите, что рассматриваемый рынок безарбитражный тогда и только тогда, 
когда

$$C=\tilde {p}C_u +\left( {1-\tilde {p}} \right)C_d ,$$ где $\tilde 
{p}=\frac{1-d}{u-d}$.
\end{problem}
\begin{remark}
 Опцион характеризуется датой исполнения (в нашем 
случае -- момент времени 2) и платежами в момент исполнения ($C_u $\linebreak и $C_d 
)$. Причем эти платежи -- заранее известные функции от цены акции в этот 
момент (введение опционов было мотивировано желанием ``хеджироваться'', 
страховаться от нежелательных изменений цен акций). Основная задача 
заключается в установлении ``справедливой'' цены опциона $C$ в момент 
времени 1 (см.  гл.~7, T.~2, {\S}11 в~\cite{21}; двухтомник Ширяев~А.Н. Основы финансовой 
стохастической математики.  М.: МЦНМО, 2016).

Не имея в начальный момент 1 капитала $X\left( 1 \right)=0$, но, проделав 
некоторую махинацию (продав одних ценных бумаг (в зависимости от специфики 
рынка иногда разрешается\ ``вставать в короткую позицию''\ -- продавать 
ценные бумаги, не имея их в наличии; приобретая при этом в долг) и купив на 
вырученные деньги других бумаг), можно в момент времени 2 гарантированно 
ничего не проиграть, и при этом с ненулевой вероятностью выиграть (не 
уточняя, сколько -- поскольку, ``прокручивая'' по имеющемуся арбитражу 
(пропорционально увеличивая коэффициенты $k_S $, $k_C )$ сколь угодно 
большую сумму, можно получить с ненулевой вероятностью сколь угодно большой 
выигрыш).

Числа $\{\tilde{p}, \tilde{q}\}$, где $\tilde{q} = 1-\tilde{p}$, а $\tilde{p} = \frac{1-d}{u-d}$ задают мартингальную вероятностью (смысл такого определения будет раскрыт в следующих задачах) 
или мартингальную меру. Если существует единственная мартингальная мера, то рынок называется 
полным. На полном рынке неизвестная цена опциона $C$ в начальный 
момент определяется однозначно и может интерпретироваться как ``справедливая 
цена''.

Приведем здесь для справки арбитражную теорему, известную также из курса <<Методов оптимизации>> как теорема Штимке (см.~Жа\-дан~В.Г. Методы оптимизации. Ч.~I.  М.: МФТИ, 2015). Детали также можно посмотреть в книге  Ross~S.M. An Elementary Introduction to Mathematical Finance. Cambridge University Press, 2011. \\ 
\noindent\verb|http://catdir.loc.gov/catdir/samples/cam033/2002073603.pdf|

\textbf{Арбитражная теорема.}  {\itПусть имеется $m$ ценных бумаг и $n$ 
возможных исходов, $r_{ij} $ -- платежи {\rm(}убытки{\rm)} по $i$-й ценной бумаги, 
если исход $j$. Под стратегией инвестора понимается вектор $x\in$\linebreak $\in \mathbb{R}^m$ 
{\rm(}$x_i $ -- количество ценных бумаг типа $i${\rm)}. Имеет место одно и только одно из приводимых ниже утверждений}:

{\it{\rm(}a\,{\rm)} $\exists\,p > 0: p^TR=0$, где $R=\left\| {r_{ij} } 
\right\|_{i,j=1}^{m,n} $ {\rm(}мартингальное соот\-но\-ше\-ние{\rm)}{\rm;}

{\rm(}б\,{\rm)} $\exists\,x\in \mathbb{R}^m: Rx \ge 0,\ Rx \ne 0$ {\rm(}условие существования арбитража{\rm)}.}

\end{remark}
\begin{problem}(Биномиальная $n$-периодная 
модель Кокса--Росса--Рубинштейна.)\label{n-period} Предложите обобщение рынка и 
соответствующих понятий из задачи \label{sec:levilong} на $n$-периодный рынок с возможностью класть 
деньги в банк под процент $r-1$ ($d<r<u)$ -- за один период (под такой же 
процент брать деньги из банка). Опцион исполняется в заключительный $(n+1)$-й 
момент. Платежи по опциону в этот момент известны и описываются известной 
функцией $\bar {C}\left( S \right)$ (например, для указанного в предыдущей задаче 
опциона $\bar {C}\left( S \right)=\max \;\left\{ {0,\;S-X} \right\})$, 
т.е. $C_k \left( {n+1} \right)=\bar {C}\left( {S_k \left( {n+1} \right)} 
\right)$, где $k$ -- состояние, в котором находится рынок в момент времени 
$n+1$. Считайте, что $S_k \left( {n+1} \right)=Su^kd^{n-k}$, т.е. $k$  
характеризует то, сколько раз акция поднималась в цене. \\
\indent Как и в предыдущей задаче, 
требуется определить ``справедливую'' цену опциона.  Обоснуйте формулу Кокса--Росса--Рубинштейна:
\begin{equation}
\label{koks}
C=\frac{1}{r^n}\sum\limits_{k=0}^n {\left( {\begin{array}{l}
 n \\ 
 k \\ 
 \end{array}} \right)} \tilde {p}^k\left( {1-\tilde {p}} \right)^{n-k}\bar 
{C}\left( {Su^kd^{n-k}} \right),
\end{equation}
где $\tilde {p}=\dfrac{r-d}{u-d}.$
\end{problem}

\begin{remark}
Отметим, что в условиях задачи $X$ называется ценой исполнения опциона и считается 
известной. Собственно, вид функции $\bar {C}\left( S \right)=\max \left\{ 
{0,\;S-X} \right\}$ проясняет смысл опциона. Опцион дает право купить (у 
того, кто продал нам опцион) в момент исполнения опциона акцию по цене $X$. 
Если акция стоит дороже в этот момент, то, конечно, мы этим правом 
воспользуемся и получим прибыль (продавец опциона обязан продать нам акцию). 
Если же цена акции меньше цены исполнения опциона, то нам уже не выгодно 
покупать акцию по более дорогой цене, чем рыночная, и мы не исполняем 
опцион, т.е. ничего не делаем (ведь опцион дает нам право, ни к чему не 
обязывая). См. также книги Белопольская~Я.И. Теория арбитража в непрерывном времени.  CПб.: СПбГАСУ, 2006 и Булинс\-кий~А.В. Случайные процессы. Примеры, задачи и упражнения.  М: МФТИ, 2010.
\end{remark}
\begin{problem}(Континуальная биномиальная модель Блэка--Шоулса.) 
\label{bw-cont}
Уместим на отрезке времени $\left[ {0,\;t} \right], \, n+1$  моментов 
(промежутки между которыми одинаковы), в которые осуществляется торговля 
как описано задаче \ref{n-period}. Введем два параметра: $a$ -- снос, $\sigma ^2\ge 0$ -- 
волатильность (дисперсия). Положим,
\[
\mu =a+\frac{\sigma ^2}{2},
\quad
r=\exp \left( {\mu \frac{t}{n}} \right),
\quad
u=\exp \left( {\sigma \sqrt {\frac{t}{n}} } \right),
\]
\begin{equation}
\label{bw}
d=\exp \left( {-\sigma \sqrt {\frac{t}{n}} } \right),
\quad
\tilde {p}=\frac{r-d}{u-d}\approx \frac{1}{2}\left( {1+\frac{a}{\sigma 
}\sqrt {\frac{t}{n}} } \right).
\end{equation}

Переходя к пределу при $n\to \infty $ в формуле \eqref{koks} (с $\bar {C}\left( S 
\right)=$\linebreak $=\max\left\{ {0,\;S-X} \right\})$, и используя \eqref{bw}, получите формулу 
Блэка--Шоулса для справедливой цены опциона в ``континуальной биномиальной 
модели''. Почему вводится именно два параметра (а не один, три и т.д.)? 
Почему
\[
r-1\approx \frac{\mu t}{n},
\quad
u-1\approx \sigma \sqrt {\frac{t}{n}} ,
\quad
d-1\approx -\sigma \sqrt {\frac{t}{n}} ?
\]
Возможны ли какие-нибудь другие осмысленные варианты соотношений типа \eqref{bw}, 
при которых будет существовать предел при $n\to \infty $ в формуле \eqref{koks} (для 
простоты вычислений считайте, что $\bar {C}\left( S \right):=S)$?
\end{problem}
\begin{remark} (Донскер--Прохоров--Скороход--Ле\;Кам--Варадарайн.)\\   
Имеет место слабая сходимость при $n\to \infty $ описанного случайного 
блуждания в дискретном времени к случайному процессу (в непрерывном 
времени), называемому геометрическим броуновским движением. Детали см., 
например, в книге Биллингсли П. Сходимость вероятностных мер.  М.: Наука, 
1977. Также см. гл.~12 в \cite{Gupta}. Также для решения задачи полезно использовать \cite{101}.
\end{remark}
\begin{problem} (Броуновское движение (процесс Башелье) и винеровский процесс.) \label{bachelie}
 Исходя из формулы \eqref{koks}, имеем, что ``рынок'' при определении 
``справедливой'' цены опциона считает, что случайный процесс $S\left( m 
\right)$ (цена акции в момент времени $m)$ эволюционирует согласно 
биномиальной модели с неизменными параметрами $d$, $u$, $\tilde 
{p}$. Построим случайный процесс $S\left( t \right)$ (в непрерывном 
времени), исходя из процесса $S\left( m \right)$, заданного в дискретном 
времени предельным переходом, аналогичным задаче \ref{bw-cont}. Как уже отмечалось,
полученный процесс $S\left( t \right)$ называют геометрическим 
броуновским движением (или случайным процессом Башелье--Самуэльсона) с 
параметрами $a$, $\sigma ^2\ge 0$, а случайный 
процесс $B\left( t \right)=\ln \left(\frac{S\left( t \right)}{S\left( 0 \right)}\right)$ 
-- броуновским движением с параметрами $a$, $\sigma ^2\ge 0$. Если 
$a=0$, $\sigma ^2=1$, то такое броуновское движение имеет специальное 
название -- винеровский процесс $W\left( t \right)$.\\
\indent Покажите, что броуновское движение является процессом Леви (см. задачу \ref{sobord}). Найдите триплет $\left( 
{b,\;c,\;\nu \left( {dx} \right)} \right)$.
\end{problem}
\begin{remark}  Один 
из альтернативных способов введения мартингальных вероятностей $\tilde{p}$ 
основывается на так называемых ``риск нейтральных'' или ``мартингальных'' 
соображениях, заключающихся в том, что $\tilde{p}$ выбирается исходя из 
равенства 
$$
\mathbb{E}_{\tilde{p}} \left[ 
Su^{\sum\limits_{k=1}^n {X_k}}
d^{n-\sum\limits_{k=1}^n {X_k}} \right]=Sr^n,
$$ где независимые одинаково распределенные с.в.  $X_k$ имеют распределение $Be\left( {\tilde {p}} \right)$, а $\mathbb{E}_{\tilde{p}}$ -- матожидание по соответствующей мере, или исходя из того, что процесс приведенной 
(продисконтированной) стоимости акции $\tilde {S}\left( m \right) = S\left( m 
\right)/{r^m}$ должен быть 
мартингалом относительно мартингальной меры $\tilde {p}$ (отсюда и 
название), т.е. $$\mathbb{E}_{\tilde {p}} \left( {\left. {\tilde {S}\left( {m+1} 
\right)} \right|\left( {\tilde {S}\left( 1 \right),...,\tilde {S}\left( m 
\right)} \right)} \right)=\tilde {S}\left( m \right).$$

Параметры геометрического броуновского движения имеют следующий смысл: 
$$a=\frac{1}{t}\mathbb{E}\left[ {\ln \frac{S\left( t \right)}{S\left( 0 \right)}} 
\right],$$ 
$$\sigma ^2=\frac{1}{t}\Var\left[ {\ln \frac{S\left( t 
\right)}{S\left( 0 \right)}} \right].$$
\end{remark}
\begin{problem} \label{geomB}(Геометрическое броуновское движение или процесс 
Башелье--Самуэльсона.) В условиях задачи \ref{bachelie} покажите, что геометрическое 
броуновское движение удовлетворяет следующему стохастическому 
дифференциальному уравнению:
\[
dS\left( t \right)=aS\left( t 
\right)dt+\sigma S\left( t \right)dW\left( t \right),
\]
которое определяет $S\left( t \right)$ как случайный процесс, 
удовлетворяющий соотношению
\[
S\left( t \right)=S\left( 0 \right)+a\int\limits_0^t {S\left( \tau \right)dt} +\sigma \int\limits_0^t 
{S\left( \tau \right)dW\left( \tau \right)},
\]
где второй интеграл понимается в смысле Ито (см.  \cite{101}).
\end{problem}
\begin{remark} Случайные процессы, которые задаются 
стохастическими дифференциальными уравнениями наподобие рассмотренного, 
задают по определению диффузионный процесс Ито (см.\,гл.\,7,\,8 книги Оксендаль~Б.
Стохастические дифференциальные уравнения.  М.: Мир, 2003.).
\end{remark}
\begin{problem} (Формула Ито.) Через конечные разности покажите, что дифференциал процесса $S(t)$ (см. предыдущую задачу) выглядит как 
\[
\Delta S\left( t \right)=\Delta \left( {S\left( 0 \right)\exp \left( 
{B\left( t \right)} \right)} \right)=
\]
\[
=aS\left( t \right)\Delta t+\sigma 
S\left( t \right)\Delta W\left( t \right),
\]
где
\[
\Delta S\left( t \right)=S\left( {t+h} \right)-S\left( t \right),
\quad
\Delta W\left( t \right)=W\left( {t+h} \right)-W\left( t \right),
\]
\[
\Delta t=t+h-t=h,
\quad
h>0.
\]
Предложите общий вид формулы для $\Delta g\left( {t,\;W\left( t \right)} 
\right)$ ($dg\left( {t,\;W\left( t \right)} \right))$.
\end{problem}

\begin{remark}
Воспользуйтесь приближением
\[
\Delta \left( {S\left( 0 \right)\exp 
\left( {at+\sigma W\left( t \right)} \right)} \right)\simeq 
\]
\[
\simeq aS\left( 0 \right)\exp \left( {at+\sigma W\left( t \right)} 
\right)\Delta t+\sigma S\left( 0 \right)\exp \left( {at+\sigma W\left( t 
\right)} \right)\Delta W\left( t \right)+
\]
\[
+\frac{\sigma ^2}{2!}S\left( 0 \right)\exp \left( {at+\sigma W\left( t 
\right)} \right)\left( {\Delta W\left( t \right)} \right)^2\simeq 
\]
\[
\simeq \left( {a+\frac{\sigma ^2}{2}} \right)S\left( 0 \right)\exp \left( 
{at+\sigma W\left( t \right)} \right)\Delta t\,+ 
\]
\[+\,\sigma S\left( 0 \right)\exp 
\left( {at+\sigma W\left( t \right)} \right)\Delta W\left( t \right)
\]

См.  гл.~4 книги  Оксендаль~Б.
Стохастические дифференциальные уравнения.  М.: Мир, 2003, а также книгу
 Белопольская~Я.И. Теория арбитража в непрерывном времени.  CПб.: СПбГАСУ, 2006;
 Булинский~А.В. Случайные процессы. Примеры, задачи и упражнения.  М.: МФТИ, 2010.
\end{remark}

\begin{problem}(Сложный пуассоновский процесс -- процесс Леви.) 
\label{Levi}
Покажите, что 
пуассоновский процесс и сложный пуассоновский процесс (см. задачи \ref{sec:poisson} и \ref{sec:cpoisson}) 
являются процессами Леви (см. задачу \ref{sobord}). Найдите триплеты $\left( {b,\;c,\;\nu \left( {dx} 
\right)} \right)$.
\end{problem}

\begin{remark} Исходя из задач \ref{sec:infdiv}, \ref{sec:levilong}
можно выдвинуть гипотезу, что любой процесс Леви может быть получен как 
сумма броуновского движения и сложного пуассоновского процесса. Такого рода 
утверждение действительно имеет место c некоторыми оговорками (см. замечание к задаче \ref{sobord}), и называется представлением 
Леви--Ито. Однако вместо сложного пуассоновского 
процесса в этом 
представлении в общем случае следует брать процесс, который может быть 
получен как ``предел'' сложных пуассоновских процессов (см. 
Sato K.-I. Levy processes and infinitely divisible distributions. Cambridge, 
1999).\\
\indent Сделанное замечание отчасти поясняет важность трех ключевых 
распределений теории вероятностей и трех типов процессов Леви: вырожденного (когда дисперсия равняется нулю), нормального (броуновское движение), распределения 
Пуассона (пределы сложных пуассоновских процессов). Помимо того, что при 
наиболее естественных предположениях для приложений имеет место сходимость к 
одному из этих трех безгранично делимых законов (например, аналогом теоремы 
Пуассона будет теорема Григелиониса \cite{2}), они в некотором смысле являются 
базисом: любое распределение, которое может возникать в пределе при 
суммировании независимых одинаково распределенных с.в., ``может быть 
получено'' исходя из этих трех базовых распределений. Приведенное замечание частично 
проясняет, в каком смысле любое распределение ``может быть так получено'' (см. также задачу 36).
\end{remark}

\begin{problem} (Сложный процесс восстановления \cite{51}.) Если в 
определении сложного пуассоновского процесса (см. задачу \ref{sec:cpoisson}) заменить пуассоновский процесс  общим процессом восстановления $$\tilde {K}\left( t \right)=\max \left\{ 
{k:\;\;\sum\limits_{i=1}^k {T_i } <t} \right\},$$ где $\left\{ {T_i } 
\right\}$ --- независимые одинаково распределенные с.в., но не обязательно, что $T_i \in \mbox{Exp}\left( \lambda 
\right)$, то получится сложный процесс восстановления $\tilde 
{Q}\left( t \right)$ (также играющий важную роль в разнообразных 
приложениях). 

Считая известными $\mathbb{E}\tilde {K}\left( t \right)$, $\Var\tilde 
{K}\left( t \right)$ и $\mathbb{E}V_i $, $\Var V_i $, определите $\mathbb{E}\tilde {Q}\left( t 
\right)$, $\Var\tilde {Q}\left( t \right)$. Найдите (приближенно) распределение сечения процесса $\tilde 
{Q}\left( t \right)$ при $t\gg 1$.
\begin{ordre}
Используйте ЦПТ в форме А.А. Натана \cite{5}, см. также ориганальные работы  Anscombe F. Large sample theory of sequential estimation // Proc. Cambridge Philos. Soc. 1952. V.~48. P.~600--607 и  Renyi A. On the asymptotic distribution of the sum of a random number of
independent random variables // Acta Math. Hung. 1957. V. 8. P. 193--199.
\end{ordre}
\end{problem}


\begin{problem}(Поток Эрланга $m$-го порядка.)\label{GammaFunc} Будет ли процессом
Леви (см.~задачу \ref{sobord}) $m$ раз ``просеянный'' пуассоновский процесс $E_m \left( t \right)$ с 
параметром $\lambda >0$

$$E_m \left( t \right)=\max \left\{ {k:\;\;\sum\limits_{i=1}^k {T_i } <t} 
\right\},$$ 
где i.i.d. с.в. $T_i \in 
\Gamma 
\left( {\lambda^{-1} ,m} \right)$, $\lambda >0$?
\begin{remark}
Определение гамма-распределения можно найти в условии к задаче \ref{gamma} из данного раздела, а также в замечании к задаче \ref{laplace} раздела \ref{standart}: $$\Gamma \left( {\lambda^{-1} ,m} 
\right) \mathop =\limits^d \underbrace {\mbox{Exp}\left( \lambda \right)+...+\mbox{Exp}\left( 
\lambda \right)}_m,$$ где $\mbox{Exp}\left( \lambda \right)+...+\mbox{Exp}\left( 
\lambda \right)$ -- сумма независимых показательных с.в. 
\end{remark}
\end{problem}

\begin{problem} В модели Блэка--Шоулса--Мертона эволюция цены акции 
описывается геометрическим броуновским движением $$S\left( t \right)=S\left( 
0 \right)\exp \left( {at+\sigma W\left( t \right)} \right),$$ где $W\left( t 
\right)$ -- винеровский процесс ($\sigma >0)$. С помощью эргодической теоремы 
для случайных процессов оцените неизвестный параметр $a$, если известна 
реализация процесса $S\left( t \right)$ на достаточно длинном временном 
отрезке $\left[ {0,\;T} \right]$. Предложите способ оценки неизвестного 
параметра $\sigma $. 
\begin{ordre}
Постройте по $S\left( t \right)$ 
процесс $Y\left( t \right)=f\left( {S\left( t \right)} \right)$, то есть подбирать 
функцию $f(\cdot)$ так, чтобы $Y\left( t \right)$ был 
эргодичен по математическому ожиданию \cite{51} и $\mathbb{E}f\left( {S\left( t \right)} 
\right)= a$ или $\mathbb{E}f\left( {S\left( t \right)} 
\right)=\sigma ^2$. Тогда в качестве оценки можно будет использовать $$\frac{1}{T}\int_{0}^{T}{Y(t)dt}.$$
\end{ordre}
\begin{remark}



Модель Блэка--Шоулса--Мертона широко использовалась на практике (см. задачи~\ref{bachelie},~\ref{geomB}). В 1990\:г. М.\:Шоулс и Р.\:Мер\-тон за свою работу были награждены нобелевской премией по 
экономике. Сейчас популярным классом моделей является $S\left( t 
\right)=$\linebreak $=S\left( 0 \right)\exp \left( {L\left( {\tau \left( t \right)} 
\right)} \right)$, где $L(\cdot)$ -- процесс Леви (см. задачу~39),\linebreak $\tau 
\left( t \right)$ -- случайный процесс, не зависящий от $L(\cdot)$, с возрастающими почти наверное траекториями.
\end{remark}
\end{problem}

\begin{problem}(Статистическая эргодическая теорема фон Неймана для динамических систем Т.~2. \cite{21}) Пусть $(X,\sigma(X))$ --- измеримое пространство, $$T:X\to X:\forall B\in \sigma(X)\Rightarrow \mu(T^{-1}(B)) = \mu(B).$$ Тогда 
$$
\forall f\in L_2(X)\Rightarrow\frac{1}{n}\sum_{k=1}^n f(T^k(x))\overset{L_2}{\longrightarrow} \mathbb{E}(f|\Theta), \, n\rightarrow \infty, 
$$
где $\Theta$ порождено функциями (случайными величиными) $$\left\{f:f(x)\overset{L_2}{=}f(T(x))\right\}.$$
Если $$B\in\sigma(X): T^{-1}(B)\overset{L_2}{=}B\Rightarrow B\overset{L_2}{=}\{\emptyset\}\text{ или } X,$$ то $$\mathbb{E}(f|\Theta) = \mathbb{E} f=\int_{X}f(x)d\mu(x).$$

Не ограничивая общности, будем считать, что $\mu(X)=1$. Построим стационарный в широком смысле случайный процесс $Y(k)=$\linebreak $=f(T^kX)$, где $X$ -- случайная величина, распределенная согласно мере~$\mu$. Верно ли, что 
$$\mathbb{E}(f|\Theta) =\mathbb{E}f\Leftrightarrow$$
$$ \Leftrightarrow\frac{1}{n^2}\sum_{i,j=1}^{n}R_{Y}(i,j) {\rightarrow}0\quad \text{при}\quad n\to \infty,$$
где 
$\dfrac{1}{n^2}\sum_{i,j=1}^{n}R_{Y}(i,j) = \dfrac{1}{n^2}\sum_{i,j=1}^{n} \mathbb{E}[\mathring{Y}(i)\mathring{Y}(j)], \, \mathring{Y}(i) = Y(i) - \mathbb{E}Y(i)$?
\end{problem}

\begin{problem}
Пусть задана последовательность  $x_k = \{\alpha k\}$, где  $\alpha$ -- какое-либо иррациональное число, $\{a\}$ -- дробная часть числа $a$. 

а) С помощью предыдущей задачи объясните, почему для любой $f\in C[0,1]$ выполняется при ${n\to\infty}$:
$$
\frac{1}{n}\sum_{k=1}^n f(x_k){\longrightarrow}\int_{0}^1 f(x)\,dx.
$$

б) Пусть функция, определенная в единичном квадрате $[0,1]^2$ на плоскости, имеет вид
$$
f(x,y) = \begin{cases}
 1, & y\leq x \text{ и } x+y\leq 1,\\
 0 & \textit{иначе}. 
\end{cases}
$$
Верно ли, что при ${n\to\infty}$
$$
\frac{1}{n}\sum_{k=1}^n f(x_{2k-1},x_{2k}){\longrightarrow}\int_0^1\int_0^1 f(x,y)\,dxdy ?
$$

\end{problem}

\begin{problem} (Г. Вейль.) Рассмотрим последовательность 
$\left\{ {a_k } \right\}_{k\in {\mathbb{N}}} $, где $a_k $ -- первая цифра в 
десятичной записи числа $2^k$. Положим $$I_m \left( {a_k } \right)=\left\{ 
{\begin{array}{l}
 1,\quad a_k =m, \\ 
 0,\quad a_k \ne m, \\ 
 \end{array}} \right.\quad m=1,\;2,\;...,9.$$ Существует ли $\mathop {\lim 
}\limits_{n\to \infty } \;\frac{1}{n}\sum\limits_{k=1}^n {I_m \left( {a_k } 
\right)}$? Если существует, то найдите его.
\end{problem}

\begin{ordre} Рассмотрим вероятностное пространство $X=\left( \Omega ,\mathcal{F} ,P\right)$, где $\Omega =\left[ {0,\;1} \right)$, $\mathcal{F} $ 
-- $\sigma $-алгебра борелевских множеств (т.е. $\Xi $ -- минимальная 
$\sigma $-алгебра, содержащая всевозможные открытые множества $\Omega 
=\left[ {0,\;1} \right)$) на $\left[ {0,\;1} \right)$, а $P$ -- равномерная 
мера на $\mathcal{F}$, т.е. $P\left( {\left[ {a,\;b} \right)} \right)=b-a$. 
Рассмотрим с.в.
\[
x\left( \omega \right)=\left\{ {\begin{array}{l}
 1,\quad \omega \in \left[ {\log _{10} m,\;\log _{10} \left( {m+1} \right)} 
\right), \\ 
 0,\quad \omega \in \left[ {0,\;\log _{10} m} \right)\cup \left[ {\log _{10} 
\left( {m+1} \right),\;1} \right). \\ 
 \end{array}} \right.
\]
Рассмотрим случайный процесс (в дискретном времени)
$$X_k \left( \omega \right)=x\left( {T^k\omega } \right),$$ где $T:\left[ 
{0,\;1} \right)\to \left[ {0,\;1} \right)$ определяется по формуле 
 $$T\omega =\left( {\omega +\log _{10} 2} 
\right) \bmod 1,$$ 
где $a \bmod 1$ -- дробная часть числа $a$.
Важно заметить, что преобразование $T$ сохраняет меру, т.е. $$\forall 
\;\;A\in \mathcal{F} \to P\left( {T^{-1}A} \right)=P\left( A \right).$$ Собственно, и 
в более общей ситуации, известная из курса случайных процессов эргодическая 
теорема схожим образом переносится на динамические системы, которые задаются 
фазовым пространством $\Omega $ и динамикой $T:\;\Omega \to \Omega $. 
Согласно теореме Крылова--Боголюбова (см. Синай Я.Г. Введение в 
эргодическую теорию.  М.: ФАЗИС, 1996. Лекция~2), если $\Omega $ --  компакт, 
то всегда найдется как минимум одна инвариантная относительно $T$ мера на 
$\mathcal{F} $. Если построенной по такой динамической системе случайный процесс 
окажется эргодическим (т.е. инвариантная мера единственная), то доля времени пребывания динамической системы в 
заданной области просто равняется мере (той самой инвариантной и 
эргодической) этой области. Ввиду вышесказанного интересно заметить, что 
установление эргодичности является трудной задачей. Например, до сих пор строго не 
обоснована ``эргодическая гипотеза Лоренца'' для идеального газа в сосуде 
(см. Козлов В.В. Тепловое равновесие по Гиббсу и Пуанкаре.  М.--Ижевск: РХД, 2002 и Минлос Р. Введение в математическую статистическую 
физику.  М.: МЦНМО, 2002), см. также задачи 1, 22 раздела 6, задачу~15 раздела~7.

Покажите, что случайный процесс $X_k $ -- стационарный в узком смысле. В 
предположении, что этот процесс эргодичен по математическому 
ожиданию (см. \cite{21} Т. 2, гл. 5) 
 найдите искомый предел.
 
Стоит обратить внимание, что в 
эргодической теореме фигурирует сходимость либо в $L_2$, либо в $L_1 $, 
либо п.н. А в данной задаче требуется (для доказательства существования 
предела и его вычисления) сходимость поточечная. Оказывается, для данной 
задачи из сходимости в $L_2 $ легко следует сходимость п.н., откуда (в свою 
очередь) следует поточечная (подробности см. Синай  Я.Г. Введение в 
эргодическую теорию.  М.: ФАЗИС, 1996, лекция 3 и  Корнфельд~И.П.,  Синай Я.Г., 
 Фомин С.В. Эргодическая теория.  М.: Наука, 1980).
\end{ordre}
\begin{problem} (Гаусс--Гильден--Виман--Кузьмин.) Каждое число из промежутка $\Omega =\left[ {0,\;1} 
\right)$ может быть разложено в цепную дробь (вообще говоря, бесконечную):
\begin{equation*}
    \omega = \dfrac{1}{a_1 + \dfrac{1}{a_2 + \frac{1}{a_3 + \cdots}}}.
\end{equation*}
\par Покажите (см. предыдущую задачу), что для почти всех (в 
равномерной мере) точек $\omega \in \left[ {0,\;1} \right)$
\[
\mathop {\lim }\limits_{n\to \infty } \;\frac{1}{n}\sum\limits_{k=1}^n {I_m 
\left( {a_k \left( \omega \right)} \right)} =\frac{1}{\ln \;2}\ln \left( 
{1+\frac{1}{m\left( {m+2} \right)}} \right).
\]
\end{problem}
\begin{ordre}
Цепные дроби (см. также замечание к задаче 108 раздела 2) играют важную роль, например, в различных вычислениях 
(поскольку позволяют строить в определенном смысле наилучшие приближения 
иррациональных чисел рациональными), в теории динамических систем (КАМ-теории). Для рациональных чисел такие дроби конечны, для квадратичных 
иррациональностей -- периодические (см.~пример ниже, в котором период равен~1):
\[
\frac{\sqrt 5 -1}{2}=\frac{1}{a_1 +\frac{1}{a_2 +\frac{1}{a_3 
+...}}}=\frac{1}{1+\frac{1}{1+\frac{1}{1+...}}}.
\]
Чтобы проверить выписанное соотношение, достаточно заметить, 
что $\frac{\sqrt 5 -1}{2}$  является корнем уравнения $x=\frac{1}{1+x}$ 
(причем из принципа сжимающих отображений следует, что последовательность 
$$x_0 =1,\quad x_{n+1} =\frac{1}{1+x_n }$$ сходится именно к этому корню),
см.  Арнольд В.И. Цепные дроби.  М.: МЦНМО, 
2001 и Хинчин А.Я. Цепные дроби.  Л.: Физматгиз, 1961. Покажите, что преобразование $T:\left[ 
0,1 \right)\to \left[ 0,1\right)$:
\[
T\omega =\left\{ {\begin{array}{l}
 \left\{ {\frac{1}{\omega }} \right\},\quad \omega \in \left( {0,\;1} 
\right), \\ 
 0,\quad \omega =0, \\ 
 \end{array}}\right.
\]
где $\left\{ {5.8} \right\}=0.8$ -- дробная часть числа, сохраняет меру Гаусса
\[
\forall \;\;A\in \mathcal{F} \Rightarrow P\left( A \right)=\frac{1}{\ln \;2}\int\limits_A 
{\frac{dx}{1+x}} ,
\]
где $\mathcal{F}$ -- сигма-алгебра на $\Omega$.

Далее рассуждайте аналогично предыдущей задаче. Эргодичность возникшего 
случайного процесса также можно не доказывать. В данном случае эргодичность означает, что мера Гаусса -- единственная мера (среди мер, нормированных на 1, т.е. вероятностных мер), которая инвариантна относительно введенного преобразования. 
\end{ordre}

\begin{problem}(Максимальный показатель Ляпунова \cite{27}.) Пусть имеется последовательность независимых одинаково распределенных случайных матриц $g_k$ (распределение $g_k$ имеет плотность). Покажите, что существует такое $\lambda\in\mathbb{R}$, что (от выбора нормы число $\lambda$ не зависит) 
\begin{equation*}
\lim_{n\to\infty} \frac{1}{n}\|g_{n}\cdot\ldots\cdot g_{1}\| = \lambda.
\end{equation*}
\end{problem}
\begin{remark}
См. книгу Гренандер У. Вероятности на алгебраических структурах.  М.: Мир, 1965.
\end{remark}

\begin{problem}(Теорема Санова.)
\label{sanov}
Рассматриваются реализации слов длины $n$, состоящих из случайных и независимо распределенных букв из конечного алфавита $A$. Вероятности появления в слове символа алфавита обозначим $p_a$, $a\in A$. Предполагается, что число символов в алфавите $|A|>1$. 
Обозначим за $v_a(s)$ с.в., равной относительной частоте буквы $a$ в строке $s$, т.е. число вхождений буквы $a$ в строку $s$, деленное на $n$, а также определим функцию ${v}(s) = (v_1(s), v_2(s), . . . , v_{|A|}(s))$.

Пусть $\Pi$ -- непустое замкнутое подмножество множества распределений вероятности 
$$
{p}=(\dots,p_a,\dots):\, p_{a}> 0 \text{ для всех } a\in A,\, \sum_{a\in A}p_{a}=1,
$$ совпадающее с замыканием своей внутренности, и $$\mathcal{KL}(\Pi|{p}) = \min_{\mu\in\Pi}\mathcal{KL}(\mu|{p}).$$
Покажите, что при $n\to\infty$
$$
-\frac{1}{n}\log\mathbb{P}({v}(s)\in \Pi)\to\mathcal{KL}(\Pi|{p}).
$$

\end{problem}
\begin{remark}
Доказательство базируется на использовании формулы Стирлинга: $n! = n^n{\rm e}^{-n}\sqrt{2\pi n}[1+O(1/n)]$, откуда вероятность реализации слов $s$ длины $n$ с частотой появления букв $v_a(s)=\mu_a$, $a \in A$, равна
\begin{equation*}
\begin{split}
\mathbb{P}(v(s)=\mu) = \frac{n!}{(n\mu_a)!\dots (n\mu_z)!}p_a^{n\mu_a}\cdot\dots\cdot p_z^{n\mu_z}\approx\\
\approx C \exp(-n\mathcal{KL}(\mu|p)) = \exp(-n\mathcal{KL}(\mu|p) + R),
\end{split}
\end{equation*}
где $|R|<(n+1)(\frac 12 \log{n}+\frac 12 \log(\pi)+\frac{1}{12n})$.

\par Отметим, что данная теорема оценивает скорость сходимости эмпирической меры на алфавите $A$ (построенной по наблюдениям $v(s)$) к истинной мере $p = (p_1,p_2, \dots, p_{|A|})$ в терминах расстояния Кульбака--Лейблера.
\par Более подробно этот результат можно посмотреть в работе Санова И.Н. О вероятности больших отклонений случайных величин // Матем. сб. 1957. Т.~42(84):1. С.~11–44 и в лекциях   А.Н.~Соболевского в НМУ http://www.mccme.ru/ium/s09/probability.html.\\

Также о вероятностях больших отклонений для эмпирических мер стохастических процессов см. следующую литературу: Feng~J., Kurtz~T.G. Large deviations for stochastic processes. Series ``Mathematical sur\-veys and monographs''. V.~131. Providence, RI, USA: American Mathe\-ma\-ti\-cal Society,  2006  и Leonard C. A large deviation approach to op\-ti\-mal trans\-port // arXiv:0710.1461v1 (2007).

\end{remark}

\begin{problem}\Star(Крамеровская зона.)
\label{Cramer_Zone}
Рассмотрим следующую сумму с.в.: $S_n = X_1+\cdots+X_n$, где $X_i$, $i=1,\dots,n, \, \mathbb{E}X_i=0$, $\mathbb{E}X_i^2= d<\infty$.
В силу ЦПТ при $n\to\infty$ выполнено следующее: 
\[
\mathbb{P}(S_n\geq x)\approx 1- \Phi\left(\frac{x}{\sqrt{nd}}\right) 
\]
равномерно по $x$ из интервала $(0,\gamma_{n}\sqrt{nd})$, где $\gamma_n\to \infty$ достаточно медленно (см. ниже). 
\par Пусть $X_1, \dots, X_n$ -- независимымые, одинаково распределенные с.в. с распределением $F$, и пусть правый хвост этого распределения
$$
    F_+(t) = \mathbb{P}(X\geq t), \, t\in \mathbb{R},
$$
есть функция вида:
$$
    F_+(t) \equiv V(t) = t^{-\alpha} L(t), \, \alpha \geq 0,
$$
где для функции $L(t)$ выполнено следующее условие:
$$
    \dfrac{L(vt)}{L(t)}\rightarrow 1 \text{ при } t\rightarrow +\infty
$$
для любого фиксированного $v>0$. Отметим, что в литературе (см. замечание) функция $V(t)$ называется правильно меняющейся с показателем $\alpha \geq 0$, а $L(t)$ называется медленно меняющейся функцией.

Покажите, что логарифмическая функция, её степени $\ln^\gamma t, \, \gamma\in \mathbb{R}$, их линейные комбинации, кратные логарифмы, функции со свойством $L(t)\rightarrow L=\const, \, t\rightarrow +\infty$, являются медленно меняющимися.

Покажите, что если
$$
\mathbb{E}(X_1^2; |X| > t) = o(1/\ln t), \, t\rightarrow +\infty, \, x > \sqrt{n},
$$ 
то
$$
    \mathbb{P}(S_n\geq x)\approx 1- \Phi\left(\frac{x}{\sqrt{nd}}\right)+nV(x), \, n \rightarrow +\infty.
$$
Значение $x$, характеризующее зону уклонений $S_n$, где происходит смена асимптотик $\mathbb{P}(S_n\geq x)$ от <<нормальной>> $1-\Phi(x/\sqrt{nd})$ на асимптотику $nV(x)$, описывающую $\mathbb{P}(S_n\geq x)$ при достаточно больших $x$, следующее:
\[
\sigma(n) = V^{-1}(1/n) = \sqrt{n(\alpha-2)d\ln n}.
\]
При $n=1$ полагаем $\sigma(1)=1$.

\end{problem}
\begin{remark}
См. книгу Боровков А.А., Боровков К.А. Асимптотический анализ случайных блужданий.  Т.~1. Медленно убывающие распределения скачков.  М.: Физматлит, 2008.  652 с.
\end{remark}

\begin{problem}(Маловероятные пути блужданий \cite{27}.) Будем говорить, что случайное  блуждание $\{S_t\}_{t\geq 0}$, $S_t = \sum_{i=1}^t X_i$, где
$X_i$ принимает значения ${-1,1}$ c вероятностями $q$ и $p$ соответственно, 
удовлетворяет принципу больших уклонений с функционалом действия $L_{\tau}(v)$, если
$$
\ln \mathbb{P}(A_{[\tau N],\delta})\sim L_{\tau}(v)N, \, \text{при больших } N,
$$
где 
$$
A_{[\tau N],\delta} = \bigg\{\sup_{t=0,1,\dots,[\tau N]}|S_t-vt|\leq \delta N\biggr\}.
$$

Покажите, что для случайного блуждания с  $p=q=1/2$ (соответствующая мера $\mathbb{P}_{0}$) верен принцип больших уклонений с функционалом действия  $L_{\tau}(v) = \tau (-v\lambda(v)+h(\lambda(v)))$,
где $\lambda=\lambda(v)$ -- решение уравнения 
$$
\frac{{\rm e}^{\lambda}-{\rm e}^{-\lambda}}{{\rm e}^{\lambda}+{\rm e}^{-\lambda}}=v.
$$
\end{problem}
\begin{remark}
Идея доказательства состоит в использовании замены меры, то есть вероятности скачков $1/2$ заменяют на
$$
p_{\lambda} = \frac{{\rm e}^{\lambda}}{{\rm e}^{\lambda}+{\rm e}^{-\lambda}},\,
q_{\lambda} = \frac{{\rm e}^{-\lambda}}{{\rm e}^{\lambda}+{\rm e}^{-\lambda}}
$$
так, чтобы $S_t-vt$ имело нулевые средние. Тогда вероятности по мере $\mathbb{P}_{0}$ представляются через средние по новой бернуллиевской мере $\mathbb{P}_{\lambda}$:
$$
\mathbb{P}(A_{[\tau N],\delta}) = \mathbb{E}_{\lambda}[I(A_{[\tau N,\delta]})\exp(-\lambda S_{[N\tau]}+[N\tau]h(\lambda))],
$$
$$
h(\lambda) = \log\left(\frac{{\rm e}^{\lambda}+{\rm e}^{-\lambda}}{2}\right),
$$
где  $I(A_{[\tau N],\delta})$ -- индикатор события $A_{[\tau N],\delta}$.

Доказательство принципа больших уклонений с указанным в условии задачи функционалом действия следует из оценок 
$$
\exp(N L_{\tau}(v))\mathbb{E}_{\lambda}[I(A_{[\tau N],\delta})]\exp(-|\lambda|\delta N)\leq \mathbb{P}(A_{[\tau N],\delta}), 
$$
$$
\mathbb{P}(A_{[\tau N],\delta}) \leq \mathbb{E}_{\lambda}[I(A_{[\tau N],\delta})]\exp(|\lambda|\delta N)\exp(N L_{\tau}(v)).
$$
См. также  Боровков А.А. Асимптотический анализ случайных блужданий. Быстро убывающие распределения приращений.  М.: Физматлит, 2013.  447 с.
\end{remark}

\begin{problem}(Случайное блуждание в полуплоскости \cite{27}.) Пусть частица находится в начальный момент в одной из точек полуплоскости $\mathbb{Z}\times\mathbb{Z}_{+}$ и совершает на каждом шаге скачок из $(k,l)\in \mathbb{Z}\times\mathbb{Z}_{+}$ в одну из четырех соседних точек решетки $(k+i,l+j)$ с вероятностями $p_{ij}$ (если $l>0$) и в одну из трех соседних точек решетки с вероятностями $q_{ij}$ (если $l=0$). Считая, что $\sum_{j}jp_{ij}<0$, опишите движение частицы \cite{27} (движение к границе и последующее движение вдоль границы).
\end{problem}

\section{Марковские модели макросистем}
\label{macrosystems}

\begin{problem}\Star(Модель П. и Т. Эренфестов в форме Кельберта--Сухова \cite{44}.)
Две собаки сидят бок о бок, страдая от $N\gg 1$ блох. Каждая блоха в промежутке времени $\left[ {t,t+h} \right)$ с 
вероятностью $\lambda h+o\left( h \right)$ ($\lambda =1)$, независимо от остальных, перескакивает на соседнюю собаку. Пусть в начальный момент все блохи находятся на  на собаке с номером~1. Покажите, что для 
всех $t\ge c N$ ($c \sim 10)$
\[
\PR\left( {\frac{\left| {n_1 \left( t \right)-n_2 \left( t \right)} 
\right|}{N}\le \frac{5}{\sqrt N }} \right)\ge 0.99,
\]
где $n_1 \left( t \right)$ -- число блох на первой собаке в момент времени 
$t$, а $n_2 \left( t \right)$~-- на второй (случайные величины). То есть 
относительная разность числа блох на собаках будет иметь порядок малости 
${\it O}\left(N^{-1/2}\right)$ при\linebreak $T\ge 
c N$. Покажите также, что математическое ожидание времени первого 
возвращения в начальное состояние (когда только собака с 
номером 1 страдает от блох) равно $2^N$.

\end{problem}

\begin{ordre}
Введем вектор $p\left( t \right)=\left( {p_0 
\left( t \right),\ldots,p_N \left( t \right)} \right)^T$, где $p_i \left( t 
\right)$ -- вероятность того, что на собаке с номером 1 в 
момент времени $t\ge 0$ ровно $i$ блох. Тогда в уравнении 
Колмогорова--Феллера \cite{1}
\[
\frac{d p^T\left( t \right)}{dt}= p^T\left( t \right)\Lambda 
\]
инфинитезимальная матрица $\Lambda $ имеет следующий вид:
\[
\left[ \Lambda \right]_{ij} =\left\{ {\begin{array}{l}
 0,\quad \left| {i-j} \right|>1, \\ 
 \lambda i,\quad j=i-1, \\ 
 \lambda \left( {N-i} \right),\quad j=i+1, \\ 
 -\lambda N,\quad j=i. \\ 
 \end{array}} \right.
\]
Также в описанной марковской динамике выполняется закон сохранения числа блох: $n_1 \left( 
t \right)+n_2 \left( t \right)\equiv N$, единственный закон 
сохранения. Стационарная (инвариантная) мера для этой системы такова (теорема Санова, см. задачу 58 раздела 5):
\[
\nu \left( {n_1 ,n_2 } \right)=\nu \left( {c_1 N,c_2 N} 
\right)=N!\frac{\left( {1 \mathord{\left/ {\vphantom {1 2}} \right. 
\kern-\nulldelimiterspace} 2} \right)^{n_1 }}{n_1 !}\frac{\left( {1 
\mathord{\left/ {\vphantom {1 2}} \right. \kern-\nulldelimiterspace} 2} 
\right)^{n_2 }}{n_2 !}=C_N^{n_1 } 2^{-N}\simeq 
\]
\[\simeq\frac{2^{-N}}{\sqrt {2\pi c_1 
c_2 } }\exp \left( {-N\cdot H\left( {c_1 ,c_2 } \right)} \right),
\]
где $H\left( {c_1 ,c_2 } \right)=\sum\limits_{i=1}^2 {c_i \ln } \;c_i $. 
Такой вид стационарной меры сопутствует тому, что если 
в начальный момент все блохи находились на одной собаке, то математическое 
ожидание времени первого возвращения макросистемы в исходное состояние будет 
порядка $2^N$. Равновесие же макросистемы естественно трактовать как 
состояние, в малой окрестности которого сконцентрирована стационарная мера 
(принцип максимума энтропии Больцмана--Джейнса):
\[
c^\ast =\left( {\begin{array}{l}
 1 \mathord{\left/ {\vphantom {1 2}} \right. \kern-\nulldelimiterspace} 2 \\ 
 1 \mathord{\left/ {\vphantom {1 2}} \right. \kern-\nulldelimiterspace} 2 \\ 
 \end{array}} \right)=\mathop {\arg \min }\limits_{\begin{scriptsize}\begin{array}{c}
 c_1 +c_2 =1 \\ 
 c\ge  0 \\ 
 \end{array} \end{scriptsize}} H\left( c \right).
\]
Полученный результат справедлив и при другом порядке предельных 
переходов (обратном к рассмотренному выше порядку: $t\to \infty $, $N\to 
\infty )$. А именно, сначала считаем, что при $t=0$ существует предел $c_i 
\left( t \right)\mathop {=\joinrel=}\limits^{\mbox{{\scriptsize п.н.}}} \mathop {\lim }\limits_{N\to 
\infty } {n_i \left( t \right)} \mathord{\left/ {\vphantom {{n_i \left( t 
\right)} N}} \right. \kern-\nulldelimiterspace} N$. Тогда (теорема Т. 
Куртца \cite{101}) этот предел существует при любом $t>0$, причем $c_1 \left( t 
\right)$, $c_2 \left( t \right)$ -- детерминированные (не случайные) 
функции, удовлетворяющие СОДУ:
\[
\frac{dc_1 }{dt}=\lambda \cdot \left( {c_2 -c_1 } \right),
\]
\[
\frac{dc_2 }{dt}=\lambda \cdot \left( {c_1 -c_2 } \right).
\]
Глобально устойчивым положением равновесия этой СОДУ будет $c^\ast $, а 
$H\left( c \right)$ -- функция Ляпунова этой СОДУ (убывает на траекториях СОДУ и имеет минимум в точке $c^\ast )$. 

Схожие рассуждения также применимы к общим 
моделям стохастической химической кинетики (см. замечание к задаче 19).
Оценка скорости сходимости макросистемы к равновесию сводится к оценке 
mixing time соответствующей марковской динамики. Здесь может оказаться  кстати
изопериметрическое неравенство Чигера, см., например, \cite{240}. Другой способ, ярко подчеркивающий связь концентрации стационарной меры и оценку
mixing time, использующий дискретную кривизну Риччи, недавно был предложен в работе Joulin~A., Olli\-vier~Y. Curvature, concentration 
and error estimates for Markov chain Monte Carlo // Ann. Prob. 2010. V. 38, N 6. P. 2418\textbf{--}2442.

\end{ordre}

\begin{remark}
С точки зрения моделей стохастической химической 
кинетики (см. замечание к задаче 19) с реакциями $1\buildrel \lambda \over 
\longrightarrow 2$ и $2\buildrel \lambda \over \longrightarrow 1$, 
равновесная конфигурация выражается с помощью принципа максимума энтропии. В свою очередь, 
этот пример демонстрирует связь энтропии  Больцмана (функция Ляпунова прошкалированной кинетической динамики) и энтропии Санова (функционал действия в неравенствах больших уклонений для инвариантной меры).   

В макросистемах возврат к неравновесным макросостояниям допустим, но случаться этому доводится только через очень большое время (циклы Пуанкаре, парадокс Цермело). Несмотря на обратимость во времени
описанного выше процесса, наблюдается 
необратимая динамика относительной разности числа блох на собаках (парадокс Лошмидта). Наглядной аналогией может стать газ, собранный в начальный 
момент в одной половине сосуда и с течением времени равномерно распределяющийся 
по сосуду. Траектории будут обратимыми, но поскольку изначально система 
выведена из равновесия, то с большой вероятностью (тем большей, чем больше 
$N)$ система будет стремиться вернуться в равновесие. Отсюда и возникает ``стрела времени'', см., например, концовку 
монографии Опойцев В.И. Нелинейная системостатика. М.: Наука, 1986.
\end{remark}

\begin{problem}(Ветвящийся процесс.)
В колонию зайцев внесли нового представителя с необычным геном. Обозначим через $p_k $  вероятность того, что в его потомстве ровно $k$ зайчат унаследуют этот ген ($k=0,1,2,\ldots)$. Точно такое же распределение вероятностей характеризует всех последующих потомков, 
унаследовавших необычный ген. Будем считать, что заячья репродукция происходит
один раз в жизни в возрасте одного года (как раз в этом возрасте находился 
самый первый заяц с необычным геном в момент попадания в колонию).

Обозначим через $G\left( z \right)$  производящую функцию распределения 
$p_k $, $k=0,1,2,...$, т.е. $G\left( z \right)=\sum\limits_{k=0}^\infty {p_k 
z^k} $. Пусть $X_n $ --- количество зайцев в возрасте одного года с необычным 
геном спустя $n$ лет после попадания в колонию первого такого зайца. 
Производящую функцию с.в. $X_n $ обозначим как $\Pi _n \left( z \right)=\Exp\left( 
{z^{X_n }} \right)$, $z \in \mathbb{R}$.

а) Получите уравнение, связывающее $\Pi _{n+1} \left( z \right)$ с 
$\Pi _n \left( z \right)$ посредством $G\left( z \right)$.

\begin{ordre}
Покажите, что $\Exp\left( {\left. {z^{X_{n+1} }} \right|X_n 
} \right)=\left[ {G\left( z \right)} \right]^{X_n }$. Затем возьмите 
математическое ожидание от обеих частей равенства.
\end{ordre}

б) Покажите, что вероятность вырождения гена \[q_n =\PR\left( {X_k 
=0;\;k\ge n} \right)=\Pi _n \left( 0 \right).\] Существует ли предел 
$q=\mathop {\lim }\limits_{n\to \infty } \;q_n $? Если существует, то 
найдите его.

\begin{ordre}
Легко видеть, что функция $G\left( z \right)$  
выпуклая. Уравнение $z=G\left( z \right)$ имеет два корня: один в любом 
случае равен 1, другой $q\le 1$. Если $\nu ={G}'\left( 1 \right)>1$, то 
$q<1$. Если $\nu \le 1$, то $q=1$.
\end{ordre}

\end{problem}

\begin{remark}
См. Севастьянов Б.А. Ветвящиеся процессы (серия ``Теория вероятностей и 
математическая статистика'').  М.: Наука, 1971; Калинкин А.В. Марковские ветвящиеся 
процессы с взаимодействием // УМН. 2002. Т. 57:2(344). С. 23--84. Частным, 
но важным случаем ветвящихся процессов, изученных в этой статье, являются 
модели стохастической химической кинетики, которые нам далее часто будут 
встречаться на протяжении всего раздела.
\end{remark}

\begin{problem}(Модель Шлёгля и теорема Т. Куртца.)
Пусть $Y_n \left( t 
\right)\in {\mathbb Z}_+ $, $t\ge 0$ -- случайный процесс рождения и гибели с 
интенсивностями рождения и гибели ($h\to 0+)$:
\[
\PR\left( {Y_n \left( {t+h} \right)=j+1\left| {Y_n \left( t \right)=j} 
\right.} \right)=n\lambda _n \left( {\frac{j}{n}} \right)h+o\left( h 
\right),
\]
\[
\PR\left( {Y_n \left( {t+h} \right)=j-1\left| {Y_n \left( t \right)=j} 
\right.} \right)=n\mu _n \left( {\frac{j}{n}} \right)h+o\left( h \right),
\]
где
\[
\lambda _n \left( x \right)=1+3x\left( {x-\frac{1}{n}} \right),
\quad
\mu _n \left( x \right)=3x+x\left( {x-\frac{1}{n}} \right)\left( 
{x-\frac{2}{n}} \right).
\]
Этот процесс 
представляет собой модель Шлёгля стохастической химической кинетики: 
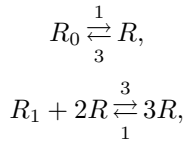
$$R_0 \mathop{\rightleftarrows}\limits_{3}^{1} R,$$ $$R_1 + 2R \mathop{\rightleftarrows}\limits_{1}^{3} 3R,$$ где $Y_n \left( t \right)$ -- количество молекул 
типа $R$ в момент времени $t$. Над (под) стрелками написаны константы 
соответствующих реакций (см.~замечание к~задаче~19).

Введем прошкалированный процесс
\[
X_n \left( t \right)=n^{1 \mathord{\left/ {\vphantom {1 4}} \right. 
\kern-\nulldelimiterspace} 4}\left( {n^{-1}Y_n \left( {n^{1 \mathord{\left/ 
{\vphantom {1 2}} \right. \kern-\nulldelimiterspace} 2}t} \right)-1} 
\right),
\quad
t\ge 0.
\]

а) Положим $E_n =\left\{ {n^{1 \mathord{\left/ {\vphantom {1 4}} 
\right. \kern-\nulldelimiterspace} 4}\left( {n^{-1}y-1} \right):\;y\in {\mathbb 
Z}_+ } \right\}$ и определим (марковскую) полугруппу на банаховом 
пространстве (с нормой $\left\| {\,\cdot \,} \right\|)$ достаточно хороших 
функций $f\left( x \right)$ над $E_n $ ($x\in E_n )$:
\[
T_n \left( t \right)f\left( x \right)\equiv \Exp\left[ {f\left( {X_n \left( t 
\right)} \right)\left| {X_n \left( 0 \right)=x} \right.} \right].
\]
Определим генератор полугруппы (линейный оператор) как (сходимость по норме 
$\left\| {\,\cdot \,} \right\|)$:
\[
G_n f=\mathop {\lim }\limits_{t\to 0+} \frac{1}{t}\left\{ {T_n \left( t 
\right)f-f} \right\}.
\]
Покажите, что
\[
G_n f\left( x \right)=n^{3 \mathord{\left/ {\vphantom {3 2}} \right. 
\kern-\nulldelimiterspace} 2}\lambda _n \left( {1+n^{-1 \mathord{\left/ 
{\vphantom {1 4}} \right. \kern-\nulldelimiterspace} 4}x} \right)\left\{ 
{f\left( {x+n^{-3 \mathord{\left/ {\vphantom {3 4}} \right. 
\kern-\nulldelimiterspace} 4}} \right)-f\left( x \right)} \right\}
+\]
\[+\:n^{3 
\mathord{\left/ {\vphantom {3 2}} \right. \kern-\nulldelimiterspace} 2}\mu 
_n \left( {1+n^{-1 \mathord{\left/ {\vphantom {1 4}} \right. 
\kern-\nulldelimiterspace} 4}x} \right)\left\{ {f\left( {x-n^{-3 
\mathord{\left/ {\vphantom {3 4}} \right. \kern-\nulldelimiterspace} 4}} 
\right)-f\left( x \right)} \right\}.
\]
Положим $Gf\left( x \right)=4{f}''\left( x \right)-x^3f\left( x \right)$. 
Покажите, что этому генератору соответствует диффузионный процесс (отметим, 
что имеет место и единственность): 
\[
dX\left( t \right)=-X^3\left( t \right)dt+2\sqrt 2 dW\left( t \right),
\]
т.е. $[\Exp f(x + dX) - f(x)] / dt \to Gf(x)$ ,
или (если известно $X\left( 0 \right))$
\[
X\left( t \right)=X\left( 0 \right)+2\sqrt 2 W\left( t 
\right)-\int\limits_0^t {X^3\left( s \right)ds} ,
\]
где $W\left( t \right)$ -- винеровский процесс (напомним, в частности, что 
$W\left( t \right)\in$\linebreak $\in \mathcal{N}\left( {0,t} \right))$. Покажите, что для
 $f\in C^2\left( 
{\mathbb R} \right)\cap C_c^1 \left( {\mathbb R} \right)$ (т.е. из пространства дважды гладких функций с финитным носителем производной)
\[
\mathop {\lim }\limits_{n\to \infty } \mathop {\sup }\limits_{x\in E_n } 
\left| {G_n f\left( x \right)-Gf\left( x \right)} \right|=0.
\]
Теорема Т.~Куртца утверждает, что если множество функций $f$, для которых 
это соотношение выполняется, достаточно богато\linebreak  (в нашем случае это так) и 
$X_n \left( 0 \right)$ слабо сходится к $X\left( 0 \right)$, то имеет место 
слабая сходимость соответствующих процессов (см. замечание): $X_n \Rightarrow X$. 

б) Покажите, что имеет место следующее представление в виде 
пуассоновских процессов со случайной заменой времени: 
\[
Y_n \left( t \right)=Y_n \left( 0 \right)+K_+ \left( {n\int\limits_0^t 
{\lambda _n \left( {n^{-1}Y_n \left( s \right)} \right)ds} } \right)- 
\]
\[-K_- 
\left( {n\int\limits_0^t {\mu _n \left( {n^{-1}Y_n \left( s \right)} 
\right)ds} } \right),
\]
$K_+ ,
\quad
K_+ \quad $ -- независимые пуассоновские процессы с одинаковым значением 
параметра, равным 1 (в частности, $K_+ \left( t \right)\in \Po\left( t \right))$. Покажите, что
\[
X_n \left( t \right)=X_n \left( 0 \right)+n^{-3 \mathord{\left/ {\vphantom 
{3 4}} \right. \kern-\nulldelimiterspace} 4}\tilde {K}_+ \left( {n^{3 
\mathord{\left/ {\vphantom {3 2}} \right. \kern-\nulldelimiterspace} 
2}\int\limits_0^t {\lambda _n \left( {1+n^{-1}X_n \left( s \right)} 
\right)ds} } \right)-
\]
\[
-n^{-3 \mathord{\left/ {\vphantom {3 4}} \right. \kern-\nulldelimiterspace} 
4}\tilde {K}_- \left( {n^{3 \mathord{\left/ {\vphantom {3 2}} \right. 
\kern-\nulldelimiterspace} 2}\int\limits_0^t {\mu _n \left( {1+n^{-1}X_n 
\left( s \right)} \right)ds} } \right)+
\]
\[+n^{3 \mathord{\left/ {\vphantom {3 
4}} \right. \kern-\nulldelimiterspace} 4}\int\limits_0^t {\left( {\lambda _n 
\left( {1+n^{-1}X_n \left( s \right)} \right)-\mu _n \left( {1+n^{-1}X_n 
\left( s \right)} \right)} \right)ds} ,
\]
где $\tilde {K}_+ \left( u \right)=K_+ -u$, $\tilde {K}_- \left( u 
\right)=K_- -u$. Используя то, что
\[
\left( {n^{-3 \mathord{\left/ {\vphantom {3 4}} \right. 
\kern-\nulldelimiterspace} 4}\tilde {K}_+ \left( {n^{3 \mathord{\left/ 
{\vphantom {3 2}} \right. \kern-\nulldelimiterspace} 2}\cdot \;} 
\right),n^{-3 \mathord{\left/ {\vphantom {3 4}} \right. 
\kern-\nulldelimiterspace} 4}\tilde {K}_- \left( {n^{3 \mathord{\left/ 
{\vphantom {3 2}} \right. \kern-\nulldelimiterspace} 2}\cdot \;} \right)} 
\right)\Rightarrow \left( {W_+ ,W_- } \right),
\]
где $\left( {W_+ ,W_- } \right)$ -- независимые винеровские процессы, 
покажите, что $X_n \Rightarrow X$, где 
\[
X\left( t \right)=X\left( 0 \right)+W_+ \left( {4t} \right)+W_- \left( {4t} 
\right)-\int\limits_0^t {X^3\left( s \right)ds} .
\]
Покажите, что такое представление предельного процесса совпадает (по 
распределению) с представлением из предыдущего пункта\linebreak (в общем случае это не 
всегда бывает так).

\end{problem}

\begin{remark}
Последовательность случайных процессов $ X_n $ слабо сходится к случайному процессу $X$, если $\mathop {\lim }\limits_{n\to 
\infty } \Exp\left[ {f\left( {X_n } \right)} \right]=\Exp\left[ {f\left( X 
\right)} \right]$ для любого 
ограниченного непрерывного функционала $f$.
Детали имеются в \cite{101}. См. также монографию  Гардинера~К.В. \cite{333}. В последней книге изложены менее строгие, но 
более простые способы получения предельных уравнений для различных 
марковских динамик, которые могут оказаться полезными для решения многих задач 
этого раздела. В рассмотренной задаче мы разобрали простейший пример, 
однако описанная техника также нашла применение в
изучении ветвящихся процессов (глава 9 \cite{101}) и строгом выводе модели 
популяционной генетики Райта--Фишера (глава 10 \cite{101}). 

Популяционная генетика и ее 
окрестности хорошо соответствуют тематике данного раздела, однако ввиду 
достаточно большого количества материала мы лишь ограничимся ссылкой на 
монографию Свирежев Ю.М., Пасеков В.П. Основы математической генетики. М.: Наука, 1982 и популярную 
книгу Горбань А.Н., Хлебопрос Р.Г. Демон Дарвина. Идея оптимальности и естественный отбор. М.: Наука, 
1988 (в электроном виде имеется второе издание: Красноярск, 1998). Также 
можно отметить недавно переведенную книгу Хойл~Ф. Математика 
эволюции. М.--Ижевск: НИЦ ``РХД'', 2012. 
\end{remark}

\begin{problem} (Кинетика социального неравенства и предельные формы.)
\label{soc_ineq}
$\quad$

 а)* В некотором городе проживает  $N\gg 1$ человек (четное число).\linebreak  В начальный 
момент у каждого жителя имеется по $\bar {s}$ монеток. Каждый день жители 
случайно разбиваются на пары. В каждой паре жители скидываются по монетке 
(если один или оба участника банкроты, то банкрот не скидывается, в то 
время как не банкрот в любом случае обязан скинуть монетку). Далее в 
каждой паре случайно разыгрывается победитель, который и забирает ``призовой 
фонд''. Обозначим через $c_s \left( t \right)$  долю жителей города, у 
которых  ровно $s$ ($s=0,...,\bar {s}N)$ монеток на $t$-й день. 
Покажите, что
\[
\exists \;\;a>0:\;\;\forall \;\;\sigma >0,\;t\ge aN\ln N, 
\]
\[
\PR\left( {\left\| 
{c\left( t \right)-c^\ast } \right\|_2 \ge \dfrac{2\sqrt{2} + 4\sqrt{\ln(\sigma^{-1})}}{\sqrt{N}}}  \right)\le \sigma ,
\]
\[
\exists \;\;b,D>0:\;\;\forall \;\;\sigma >0,\;t\ge bN\ln N 
\]
\[\PR\left( {\left\| 
{c\left( t \right)-c^\ast } \right\|_1 \ge D\sqrt {\frac{\ln ^2N+\ln \left( 
{\sigma ^{-1}} \right)}{N}} } \right)\le \sigma ,
\]
где $c_s^\ast \simeq C\exp \left( {-s \mathord{\left/ {\vphantom {s {\bar {s}}}} 
\right. \kern-\nulldelimiterspace} {\bar {s}}} \right)$, а $C\simeq 1 
\mathord{\left/ {\vphantom {1 {\bar {s}}}} \right. 
\kern-\nulldelimiterspace} {\bar {s}}$ находится из условия нормировки 
$\sum\limits_{s=0}^{\bar {s}N} {C\exp \left( {-s \mathord{\left/ {\vphantom 
{s {\bar {s}}}} \right. \kern-\nulldelimiterspace} {\bar {s}}} \right)} =1$. 
Таким образом, кривая (предельная форма) $C\exp \left( {-s \mathord{\left/ 
{\vphantom {s {\bar {s}}}} \right. \kern-\nulldelimiterspace} {\bar {s}}} 
\right)$ характеризует распределение населения по богатству на больших 
временах.

\begin{ordre}
 Для решения этой задачи полезно рассмотреть схожий 
процесс, в котором каждой паре жителей сопоставлен свой (независимый) 
пуассоновский будильник (звонки происходят в случайные моменты времени, 
соответствующие скачкам пуассоновского процесса; параметр интенсивности 
этого пуассоновского процесса называют интенсивностью/параметром 
будильника). У всех будильников одна и 
та же интенсивность $\lambda N^{-1}$. Далее следует погрузить задачу в 
модель стохастической химической кинетики с бинарными реакциями и 
воспользоваться результатом из замечания к задаче 19. Наиболее технически 
сложными моментами в получении указанного в условии задачи результата 
являются оценка mixing time $\sim N\ln N$ и получение поправки под корнем 
$\ln ^2N$. 
\end{ordre}

\begin{remark}
Название модели взято из одноименной статьи К.Ю. Богданова в журнале ``Квант''. В этой 
статье предлагаются и другие правила обмена. Также К.Ю. Богданов является автором одной 
из стохастических динамик, приводящих к модели хищник--жертва и системе 
уравнений Лотки--Вольтера (см. задачу 10).  О возможных обобщениях модели, описанной в предложенной задаче, также можно посмотреть в работах Dragulescu A., Yakovenko V.M. Statistical mechanics of money // 
The European Physical Journal B. 2000. V.~17. P. 723--729 и в других работах 
этих авторов (например, arXiv:0905.1518), а также в работах A.M. Chebotarev'а. Отметим, что в отличие от рассмотренного здесь примера, на практике (и в ряде других моделей) хвосты распределения имеют степенной характер (не экспоненциальный).

Степенные законы в играют важную роль в различных экономических приложениях, не только связанных с описанием расcлоения населения по дохододам. См., например, LeBaron B. Stochastic volatility as a simple generator of apparent financial power laws and long memory // Quantitative Finance. 2001. V. 1:6. P. 621--631; Parisi D.R., Sornette D., Helbing D. Financial price dynamics and pedestrian counter flows: A comparison of statistical stylized facts // Physical Review E. 2013. V. 87. 012804; Challet D., Bochud T. Optimal approximations of power-laws with exponentials: application to volatility models with long memory // Quantitative Finance. 2007. V. 7. no. 6. P. 585--589; http://tuvalu.santafe.edu/~jdf/papers/powerlaw3.pdf. 
\end{remark}

б)* (Задача Булгакова--Маслова о разбрасывании червонцев в варьете.) 
В некотором городе живет $N\gg 1$ жителей (изначально банкротов). Каждый 
день одному из жителей (случайно выбранному в этот день) дается одна 
монетка. В обозначениях п. а) определите предельную форму для $\left\{ {c_s 
\left( {\bar {s}N} \right)} \right\}_{s=0}^{\bar {s}N} $. Сколько надо 
случайно раздать монеток, чтобы с вероятностью $\ge 1-\sigma $ в городе бы 
не было банкротов?

\begin{remark}
Ответ на вопрос задачи помогает 
организовывать (определять, сколько нужно выбрать точек старта) метод 
случайного мультистарта в глобальной оптимизации, см.,~например, Жиглявский~А.А., Жилинскас А.Г. Методы 
поиска глобального экстремума. М.: Наука, 1991.

Покажите, что распределение $\left\{ {c_s \left( {\bar 
{s}N} \right)} \right\}_{s=0}^{\bar {s}N} $  мультиномиальное. Более того, 
оно в точности совпадает со стационарным распределением марковского процесса 
из предыдущего пункта. В работах В.П. Маслова (в основном в журнале ``Математические 
заметки'') собрано большое число тонких результатов о концентрации различных 
(урновых) мер, в том числе возникающих в приложениях. Такие меры, как 
правило, имеют комбинаторную природу и описывают способы распределения 
шаров по урнам при различного рода линейных ограничениях. В нашем случае 
урны отвечают $s=0,...,\bar {s}N$. А ограничения -- закон сохранения числа 
жителей и числа монеток (задают множество (inv) в замечании к задаче 19):
\[
\sum\limits_{s=0}^{\bar {s}N} {n_s \left( t \right)} =N,
\quad
\sum\limits_{s=0}^{\bar {s}N} {sn_s \left( t \right)} =\bar {s}N,
\]
где $n_s \left( t \right)$ -- количество жителей, у которых  ровно\textbf{ 
}$s$ монеток на $t$-й день.
\end{remark}

в)** (Закон Ципфа--Парето и процесс Юла.) 
В городе, в котором проживает большое количество людей, каждый человек изначально не имеет денег, но может участвовать в розыгрыше монеток. База индукции: cначала выбирается один житель, он получает одну монетку. Шаг индукции: в $(k+1) \text{-й}$ день выбирается  новый житель (отличный от $k$ уже выбранных), он получает одну монетку. Еще одна монетка отдается одному из  $k$ предшествующих жителей: с вероятностью  $\alpha < 1$ равновероятно и вероятностью $1 - \alpha$ пропорциональной тому, сколько у жителя уже имеется монеток, т.е. по принципу ``деньги к деньгам'' (в моделях Интернета этот принцип называют ``preferential attachment'').  В соответствии с обозначениями п. а), определите предельную форму для $\left\{ {c_s \left( t 
\right)} \right\}_{s=1}^t $, $t\gg 1$.

\begin{ordre} 
 Изучать стохастическую динамику в этой задаче заметно 
сложнее, чем в двух предыдущих. Поэтому обычно работают в приближении 
среднего поля. В данном контексте это означает, что $n_s \left( t 
\right)\simeq \Exp\left[ {n_s \left( t \right)} \right]$. Далее выписывают на 
$n_s \left( t \right)$ систему зацепляющихся обыкновенных дифференциальных 
уравнений. Автомодельное притягивающее решение этой системы ищут в виде $n_s 
\left( t \right)\sim c_s^\ast t$. После разрешения соответствующих уравнений 
получают степенной закон для зависимости $c_s^\ast \sim s^{-3 + \alpha / (1-\alpha)}$. Такой подход 
применительно к моделям роста Интернета (и изучения степенного закона для 
распределения степеней вершин, см. задачи \ref{prefattach}, \ref{pref_attach} раздела \ref{hard} и задачу \ref{hnmgraph}  раздела 8) 
довольно часто сейчас встречается (в том числе и в учебной литературе). В 
частности, этот подход описан в обзоре Mitzenmacher M. A Brief History of Generative Models 
for Power Law and Lognormal Distributions // Internet Mathematics. 2003. V.~1, N~2. P. 226--251 и рассматривается в книге Blum A., Hopcroft J., Kannan R. Foundations of Data Science // \linebreak e-print. 2016.
\verb|http://www.cs.cornell.edu/jeh/bookJan25_2016.pdf.| По сути, в этом пункте нами был описан 
процесс, возникающий в \linebreak работах 20-х годов XX века по популяционной генетике и 
получивший название {\itпроцесса Юла} (см., например, обзор Newman М.Е.J. Power laws, Pareto 
distributions and Zipf's law // Contemporary Physics. 2005. V. 46, N~5. P.~323--351).
\end{ordre}

\begin{remark}
 Описанные модели восходят к работе В.~Парето конца XIX века, в которой была предпринята попытка объяснить социальное неравенство, 
и к работе Г. Ципфа конца 40-х годов XX века, в которой была отмечена важность 
степенных законов в природе. Эти законы для большей популярности иногда 
преподносят как принцип Парето или принцип 80/20 (80\% результатов 
проистекают всего лишь из двадцати процентов причин). Такие пропорции соответствуют $c_s^\ast \sim s^{-2.1}$ (см. Newman M.E.O.).  Приведем примеры (не 
совсем, правда, точные): 80{\%} научных результатов получили 20{\%} 
ученых, 80{\%} пива выпили 20{\%} людей и т.п. Сейчас много исследований 
во всем мире посвящено изучению возникновения в самых разных приложениях 
степенных законов (распределение городов по населению, коммерческих компаний 
по капитализации, автомобильных пробок по длине). Причиной этому в том числе является возрастание интереса к большим сетям (экономическим, социальным, 
Интернета), см., например, книги и работы M.O. Jackson'а. В~России это 
направление также представлено, например, в Яндексе в отделе А.М.~Райгородского (Гречников--Остроумова--Рябченко--Самосват, 
Леонидов--Мусатов--Савватеев), в ИПМ РАН в группе А.В.~Подлазова 
(собственно, модель п. в) была взята из статьи Подлазова А.В. в сборнике 
``Новое в синергетике. Нелинейность в современном естествознании''. М.: 
ЛИБРОКОМ, 2009. С. 229--256).
\end{remark}

\end{problem}

\begin{problem} (Обезьянка и печатная машинка; закон Ципфа--Мандельброта.)
На печатной машинке $n+1$ символ, один из символов пробел. Обезьянка на 
каждом шаге случайно (независимо и равновероятно) нажимает один из символов. 
Прожив долгую жизнь, обезьянка сгенерировала текст огромной длины. По этому 
тексту составили словарь. Этот словарь упорядочили по частоте встречаемости 
слова (слова -- это всевозможные наборы букв без пробелов, которые хотя бы 
раз встречались в тексте обезьянки между какими-то двумя пробелами). Так, на 
первом месте в словаре поставили самое часто встречаемое слово, на второе 
поставили второе по частоте встречаемости и~т.д. Номер слова в таком словаре 
называется {\itрангом} и обозначается буквой $r$. Покажите, что предельная форма 
кривой, описывающей распределение частот встречаемости слов от рангов, имеет 
вид
\[
\mbox{частота}\left( r \right)\simeq \frac{C}{\left( {r+B} \right)^\alpha 
},
\]
где 
\[
\alpha =\frac{\log \left( {n+1} \right)}{\log \left( n \right)},
\quad
B=\frac{n}{n-1},
\quad
C=\frac{n^{\alpha -1}}{\left( {n-1} \right)^\alpha }.
\]
\end{problem}

\begin{remark}
Такой вывод закона Ципфа (с поправкой) был одним из 
двух, предложенных Б. Мандельбротом. Этот 
вывод можно найти во многих источниках. 

Благодаря Википедии, наиболее 
популярной на эту тему оказалась статья  Li W. Random Text Exibit 
Zipf's--Law--Like Word Frequency Distribution // IEEE Transactions of 
Information Theory. 1992. V. 38(6). P. 1842--1845. В работе (см. также 
Conrad--Mitzenmacher) Бочкарева~В.В., Лернера~Э.Ю. Закон Ципфа для случайных текстов с неравными 
вероятностями букв и пирамида Паскаля // Известия вузов. Математика. 2012. №~12. С. 30--33 приводится обобщение этой модели на случай неравных 
вероятностей букв. В другой работе этих авторов рассматривается обобщение, 
связанное с рассмотрением ситуации, когда последовательность нажимаемых 
символов не i.i.d., а образует цепь Маркова (впрочем, такого рода обобщения 
довольно популярны и рассматривались рядом других авторов). В~работе 
Шрейдер Ю.А. О возможности теоретического вывода статистических закономерностей текста 
(к обоснованию закона Ципфа) // Проблемы передачи информации. 1967. Т. 3, №~1. С. 57--63 была предложена довольно общая схема, приводящая к закону 
Ципфа, в которую можно погрузить и обезьянку с печатной 
машинкой (это полезно показать, а также 
представляется полезным изучение других примеров (в том числе упомянутых) на 
предмет возможности их погружения в эту схему). 
А именно, предположим, что динамика (порождения слов в большом 
тексте) такова, что вероятность того, что в тексте из $N\gg 1$ слов $x_1 $ 
(первое по порядку слово в ранговом словаре) встречалось $N_1 $ раз, $x_2 $ 
(второе по порядку слово в ранговом словаре) встречалось $N_2 $ раза и т.д., 
есть
\[
\sim \frac{N!}{N_1 !\,N_2 !...}\exp \left( {-\eta \sum\limits_{k\in \Nbb} 
{N_k E_k } } \right),
\]
где $E_k $ -- число букв в слове с рангом $k$. Часто считают, что $\eta =0$, 
но зато динамика такова, что число слов и число букв становятся 
асимптотически (по размеру текста) связанными (закон больших чисел). Таким 
образом, к закону сохранения $\sum\limits_{k\in  \Nbb} {N_k =N} $ 
добавляется приближенный закон сохранения $\sum\limits_{k\in  \Nbb} {N_k 
E_k } \simeq \bar {E}N$ ($\bar {E}$ -- среднее число букв в слове). Поиск предельной 
формы приводит к задаче (воспользовались формулой Стирлинга и методом 
множителей Лагранжа):
\[
\sum\limits_{k\in  \Nbb} {\left\{ {N_k \ln N_k +\lambda E_k N_k } 
\right\}} \to \mathop {\min }\limits_{ {\scriptsize \begin{array}{c}
 N_k \ge 0 \\ 
 \sum\limits_{k\in  \Nbb} {N_k =N} \\ 
 \end{array}}} ,
\]
где $\lambda $  либо равняется $\eta $, либо является множителем Лагранжа 
к ограничению $\sum\limits_{k\in \Nbb} {N_k E_k } \simeq \bar {E}N$. 
Стоит отметить, что к аналогичной задаче (с $E_k =k)$ приводит поиск 
предельной формы в модели ``Кинетика социального неравенства'' из п. а) 
задачи 4. Решение нашей задачи дает $N_k = \exp \left( {-\mu -\lambda E_k } 
\right)$, где $\mu $ -- множитель Лагранжа к ограничению $\sum\limits_{k\in 
\Nbb} {N_k } =N$. Далее, считают, что $r\left( E \right)$ --- число различных 
используемых слов с числом букв, не большим $E$, приближенно представимо в 
виде $r\left( E \right)\simeq a^E$. Тогда $N_k \sim k^{-\gamma }$, где 
$\gamma =\lambda \mathord{\left/ {\vphantom {\mu {\ln a}}} \right. 
\kern-\nulldelimiterspace} {\ln a}$. 

Некоторые обобщения закона Ципфа--Мандельброта с обсуждением также и 
лингвистической стороны дела имеются в работе В.П.~Маслова. 
Интересным представляется связь, недавно обнаруженная Ю.И.~Маниным, между 
законом Ципфа и колмогоровской сложностью. Отметим, что ранее на такую связь 
также было указано В.В.~Вьюгиным.

Отметим, что Б. Мандельброт является создателем науки о фракталах. 
Интерес Б. Мандельброта к степенным законом (законам типа Ципфа) был 
закономерен, поскольку степенные законы как раз проявляют свойства 
масштабной инвариантности (самоподобности), свойственной фракталам, и часто 
возникают при описании различных критических явлений (см., например, обзор 
М.Е.J.~Newman'а, 2005). В популярных книгах Б. Мандельброта, большинство из 
которых переведены на русский язык, также периодически встречаются сюжеты на 
эту тему. Вообще, здесь бы хотелось отметить, что очень многие сложные сети, 
возникающие в природе (капиллярная система, ветви деревьев и т.п.), имеют в 
своей основе определенные закономерности, формирующиеся при росте (напомним, 
что социальные сети и Интернет моделируются с помощью модели случайного 
роста ``предпочтительного присоединения''). В результате возникающая 
конфигурация  находится из некоторого 
вариационного принципа, приводящего, так же как и в рассмотренной задаче, к 
степенным законам. Яркий пример имеется, например, в статье West~G.B., Brown~J.H. Enquist~B.J. A General Model 
for the Origin of Allometric Scaling Laws in Biology // SCIENCE. 1997. V.~276. P.~122--126. Представляется интересным искать аналоги ``законов 
физики'' при изучения роста не биологических сетей. То есть небольшой набор 
универсальных правил, которые в различных сочетаниях и в различных 
контекстах дают многообразие возникающих в приложениях сетей больших 
размеров, и отражают их основные (статистические) свойства. Один такой закон 
(предпочтительного присоединения) уже открыт (у этого закона много вариаций, 
за счет которых и удается подогнать ту или иную модель под экспериментальные 
результаты). Работа в этом направлении начата относительно недавно ($\sim 
$15 лет назад), и, кажется, что здесь все самое интересное впереди. 
\end{remark}

\begin{problem}\Star(Heaps' law.) 
Предположим, что сгенерирован достаточно 
большой текст, в котором слова случайно (независимо) выбирались из словаря с 
вероятностями, зависящими от ранга слова по закону Ципфа--Мандельброта (см. 
предыдущую задачу). Покажите, что если сгенерированный текст имеет всего $n$ 
слов, то число различных слов $m$ в нем будет $\sim n^{1 \mathord{\left/ 
{\vphantom {1 \alpha }} \right. \kern-\nulldelimiterspace} \alpha }$. 
Постарайтесь получить более точное описание с.в.~$m$.
\end{problem}
\begin{remark}
 См., например, Grootjen F.A., Leijenhorst D.C., van der Weide Th.P. A formal deviation of Heaps' law // 
Inform. Sciences. 2005. V.~170(24). P.~263--272.
\end{remark}

\begin{problem} (Распределение букв по частоте встречаемости; С.М.~Гусейн-Заде.)
Предположим, что есть некоторая эргодическая динамика 
(детерминированная) на единичном симплексе с равномерной инвариантной мерой, 
описывающая эволюцию частот встречаемости букв в текстах. Эту динамику можно 
понимать как эволюцию состояния некоторого лингвистического мира. Мы 
считаем, что эта динамка развивается в медленном времени (годы). В то же 
время в быстром времени (дни) при заданном ``состоянии мира'' (частотах 
встречаемости букв) случайно (согласно этим частотам) генерируется (с 
постоянной интенсивностью) большое количество текстов. За много лет 
накопилось огромное количество таких текстов. Их объединили все вместе в 
один огромный текст и посчитали частоты встречаемости различных букв в этом 
объединенном тексте. Считая, что всего имеется $n$ различных букв, покажите, 
что $r$-я по величине частота встречаемости приблизительно равна (ранговое 
распределение частот букв)
\[
\frac{1}{n}\sum\limits_{k=r}^n {\frac{1}{k}} \approx \frac{1}{n}\left( {\ln 
\left( {n+1} \right)-\ln r} \right).
\]
\end{problem}

\begin{remark}
 Впрочем, описанный 
закон в определенном смысле ничем не примечателен по сравнению с рядом 
известных конкурирующих законов, см. Гельфанд М.С., Минь Чжао. О ранговых распределениях частот 
букв в~естественных языках // Проблемы передачи информации. 1996. Т. 32, № 2. 
С.~89--95.
\end{remark}

\begin{problem} (Случайный рост диаграмм, предельные формы.)
Под диаграммой будем понимать горку (на плоскости), составленную из одинаковых 
квадратных кирпичей. Слово {\itгорка} означает, что ряд из кирпичей, 
поставленных друг на друга, не ниже ряда из кирпичей, примыкающего справа. 

а)** (Диаграммы Ричардсона.) Пусть горка растет по принципу: новый 
кирпич кладется с равной вероятностью на одну из допустимых позиций, где под 
допустимой позицией понимается ``уголок'' либо две ``крайние'' позиции (над 
самым левым верхнем кирпичом, справа от самого правого нижнего кирпича --- 
далее мы будем называть эти позиции соответственно верхний и нижний 
уголок). Если положить новый кирпич в такой уголок на горке, то он будет 
касаться горки двумя своими гранями (левой и нижней). Пусть случайно, как 
описано выше, положили $n\gg 1$ кирпичей. Перейдя в систему координат с 
прошкалированными фактором $\sqrt n $ осями и разместив левый нижний кирпич 
горки в начале координат, покажите, что в пределе при $n\to \infty $ форма 
горки (предельная форма) будет неслучайной кривой, описываемой уравнением 
$\sqrt x +\sqrt y =\sqrt[4]{6}$, где константа $\sqrt[4]{6}$ была выбрана из 
условий нормировки площади под кривой на 1.

\begin{remark}
См. статью Rost H. Nonequilibrium behavior of a many particle 
process: Density profile and local equilibria // Probability Theory and 
Related Fields. 1981. V. 58, N 1. P. 41--53.
\end{remark}

б)* Как надо правильно прошкалировать оси и какая при этом 
получится предельная форма, если новый кирпич кладется с вероятностью 
$\alpha $ в нижний уголок и с вероятностью $1-\alpha $  равновероятно 
опускается на один из рядов кирпичей, поставленных друг на друга, и далее 
скатывается в уголок, находящийся на этом уровне? Как изменится ответ, если 
новый кирпич равновероятно опускается на один из рядов кирпичей, 
поставленных друг на друга, или нижний уголок, и далее скатывается в уголок, 
находящийся на этом уровне?

\begin{remark}
Следует сравнить с п. в) задачи 4. Предельные формы 
различных диаграмм активно изучались и изучаются в работах А.М.~Вершика и 
его учеников. Часто эти диаграммы имеют вполне содержательную интерпретацию 
(например, диаграммы Юнга). Мера на диаграммах часто задается без описания 
динамики, т.е. не описывается рост диаграмм. Таким примером является 
изучение предельной формы диаграмм Юнга с равномерной мерой (рассматривают и 
другие меры, например, Планшереля), см. задачу \ref{limung} раздела \ref{genF}. 
Сюда же можно отнести и статистику выпуклых ломаных 
(Арнольд--Вершик--Синай), см. задачу \ref{convcurve} раздела \ref{genF}. 
Примечательно, что в обоих случаях предельная форма определяется из решения 
вариационной задачи, которая получается исходя из максимизации (в 
пространстве кривых) должным образом прошкалированного логарифма числа 
диаграмм (выпуклых ломанных), находящихся в малой окрестности заданной 
кривой (аналог действия/энтропии в теоремах о больших уклонениях типа 
Санова). Причем ответ не 
зависит от того, как выбирается эта окрестность (в предположении ее 
малости), важно лишь, чтобы окрестности строились по одним и тем же правилам 
для разных кривых.
Для определения предельной формы диаграмм Юнга мы получаем 
в точности задачу энтропийно-линейного программирования в функциональном 
пространстве, следует сравнить с замечанием к задаче 19.
\end{remark} 

\end{problem}

\begin{problem}\Star(Распределенные модели стохастической химической кинетики.)
На окружности единичного периметра на расстоянии $\varepsilon $ 
друг от друга расположены дома. В  каждом доме живет по одной собаке и 
одной кошке. Для определенности, в начальный момент у каждой собаки $N$ 
блох, а на кошках блох в начальный момент нет. С интенсивностью $\lambda $ 
(см. задачу 1 ``парадокс Эренфестов'') каждая блоха независимо от остальных 
перескакивает с собаки (кошки) на соседа по дому, т.е. кошку (собаку), а с 
вероятностью ${\lambda _{-1} } \mathord{\left/ {\vphantom {{\lambda _{-1} } 
\varepsilon }} \right. \kern-\nulldelimiterspace} \varepsilon +\mu 
\mathord{\left/ {\vphantom {\mu {\varepsilon ^2}}} \right. 
\kern-\nulldelimiterspace} {\varepsilon ^2}$ (${\lambda _{+1} } 
\mathord{\left/ {\vphantom {{\lambda _{+1} } \varepsilon }} \right. 
\kern-\nulldelimiterspace} \varepsilon +\mu \mathord{\left/ {\vphantom {\mu 
{\varepsilon ^2}}} \right. \kern-\nulldelimiterspace} {\varepsilon ^2})$ --\linebreak на 
собаку (кошку) из соседнего дома с номером, на единицу меньшим (большим). 
Покажите (используя теорему Куртца), что при $N\to \infty $, $\varepsilon 
\to 0\,+$ п.н. существует (детерминированный) предел (причем на каждом 
конечном отрезке времени можно говорить о слабой сходимости случайных 
процессов):
\[
c^{dog}\left( {t,x} \right)=\mathop {\lim }\limits_{{\scriptsize \begin{array}{c}
 N\to \infty ,\;\varepsilon \to 0+ \\ 
 x \mathord{\left/ {\vphantom {x \varepsilon }} \right. 
\kern-\nulldelimiterspace} \varepsilon \le k<x \mathord{\left/ {\vphantom {x 
\varepsilon }} \right. \kern-\nulldelimiterspace} \varepsilon +1 \\ 
 \end{array}}} \frac{n_k^{dog} \left( t \right)}{N},
\]
где $n_k^{dog} \left( t \right)  $ -- количество блох у собаки в доме с номером 
$k$ в момент времени $t$, аналогично, через $n_k^{cat} \left( t \right)$, 
определяется $c^{cat}\left( {t,x} \right)$. Покажите, что введенные 
плотности удовлетворяют следующей системе УЧП:
\[
\frac{\partial c^{dog}}{\partial t}=\mu \frac{\partial ^2c^{dog}}{\partial 
x^2}+\left( {\lambda _{-1} -\lambda _{+1} } \right)\frac{\partial 
c^{dog}}{\partial x}+\lambda \cdot \left( {c^{cat}-c^{dog}} \right),
\]
\[
\frac{\partial c^{cat}}{\partial t}=\mu \frac{\partial ^2c^{cat}}{\partial 
x^2}+\left( {\lambda _{-1} -\lambda _{+1} } \right)\frac{\partial 
c^{cat}}{\partial x}+\lambda \cdot \left( {c^{dog}-c^{cat}} \right),
\]
унаследовавшей линейный закон сохранения от стохастической динамики

\[\tag{inv}
\int\limits_0^1 {\left\{ {c^{dog}\left( {t,x} \right)+c^{cat}\left( {t,x} \right)} \right\}dx} \equiv 1.
\] 

Покажите, что стационарное распределение описанного марковского процесса с 
носителем (inv) имеет асимптотическое представление
\[
\sim \exp \left( {-\frac{N}{\varepsilon }\int\limits_0^1 {\left\{ 
{c^{dog}\left( {t,x} \right)\ln \left( {c^{dog}\left( {t,x} \right)} 
\right)+c^{cat}\left( {t,x} \right)\ln \left( {c^{cat}\left( {t,x} \right)} 
\right)} \right\}dx} } \right).
\]
При этом функционал (пространственная энтропия), стоящий при $N 
\mathord{\left/ {\vphantom {N \varepsilon }} \right. 
\kern-\nulldelimiterspace} \varepsilon $, является функционалом Ляпунова 
выписанной ранее системы УЧП. Исследуйте поведение системы на больших 
временах.

\end{problem}

\begin{remark}
Для решения этой задачи (и ряда других задач этого 
раздела) можно рекомендовать монографии \cite{333}  (главы 7, 8), \cite{101}.

В данной задаче рассмотрен, пожалуй, простейший пример пространственно 
распределенной системы. Такие системы часто возникают, например, в 
математической биологии при описании эволюции пространственно распределенных 
взаимодействующих видов. Причем, как правило, за счет более сложного 
взаимодействия и(или) неоднородной (анизотропной) диффузии и(или) сноса (в 
том числе со скоростью/интенсивностью, зависящей от концентраций видов) 
при скейлинге получаются уже нелинейные эволюционные уравнения 
(параболического типа, если есть диффузия). При этом равновесная 
конфигурация, возникающая на больших временах, часто может быть также 
проинтерпретирована в терминах равновесия по Нэшу и Дарвиновского отбора (не 
распределенный вариант, поясняющий отмеченные связи, имеется в задачах 10, 
15, 16). К этому направлению можно отнести и ``принцип эволюционной 
оптимальности'', активно развиваемый в работах В.Н. Разжевайкина.
\end{remark}

\begin{problem}\DStar(Модель хищник--жертва.) 
В сказочном лесу, 
имеющем форму единичного тора (прямое произведение двух окружностей 
единичного периметра), который можно понимать как единичный квадрат с 
отождествленными противоположными сторонами, проведена решетка с квадратными 
ячейками размера $\varepsilon \times \varepsilon $. В узлах этой решетки в 
начальный момент времени случайно (независимо и равновероятно) распределены 
$n_{\text{з}} =\kappa _{\text{з}} \varepsilon ^{-1}$ зайцев и $n_\text{в} =\kappa _\text{в} \varepsilon 
^{-1}$ волков. Волки и зайцы начинают независимо случайно блуждать по 
решетке: с интенсивностью $\lambda _\text{з} \varepsilon ^{-1}$ заяц покидает 
текущий узел решетки, выбирая для перехода равновероятно один из четырех 
соседних узлов, аналогично (с $\lambda _\text{в} \varepsilon ^{-1})$ ведут себя 
волки. Кроме того, с интенсивностью $\mu _\text{з}^+ $ заяц создает в текущем узле 
потомка, который сразу же начинает жить независимой жизнью по общим 
правилам. А волки с интенсивностью $\mu _\text{в}^- $ выбывают из системы. Если в 
результате блуждания волк и заяц окажутся в одном узле, то мгновенно 
происходит реакция: волк съедает зайца и производит потомка, который сразу 
же начинает жить независимой жизнью по общим правилам. В каком смысле можно говорить о том, что в пределе при $\varepsilon \to 0+$ (на конечных отрезках 
времени) динамика изменения числа волков и зайцев соответствует 
модели стохастической химической кинетики (см. замечание к задаче 19):
\[
\left[ \text{З} \right]\buildrel {\mu _\text{з}^+ } \over \longrightarrow 2\left[ \text{З} 
\right],
\quad
\left[ \text{В} \right]+\left[ \text{З} \right]\buildrel K \over \longrightarrow 2\left[ \text{В} 
\right],
\quad
\left[ \text{В} \right]\buildrel {\mu _\text{в}^- } \over \longrightarrow 0,
\]
где константа реакции $K$ определяется по введенным выше параметрам? 
Покажите, что закон действующих масс приводит к системе Гульдберга--Вааге 
(получающейся при каноническом скейлинге, см. замечание к задаче 19):
\[
\frac{dc_\text{з} }{dt}=\mu _\text{з}^+ c_\text{з} -Kc_\text{в} c_\text{з} ,
\]
\[
\frac{dc_\text{в} }{dt}=Kc_\text{в} c_\text{з} -\mu _\text{в}^- c_\text{в} .
\]
Эту систему уравнений принято называть системой Лотки--Вольтера.
\end{problem}

\begin{remark}
 Несложно показать, что в стохастической модели (большое, 
но конечное, число волков и зайцев) с вероятностью 1 на больших временах 
система ``свалится'' в поглощающее состояние (без волков). При этом число 
зайцев в этом состоянии либо будет равняться нулю, либо стремиться к 
бесконечности со временем. При этом система Лотки--Вольтера имеет 
нетривиальное положение равновесия $\left( {c_\text{з}^\ast ,c_\text{в}^\ast } 
\right)=\left( {{\mu _\text{в}^- } \mathord{\left/ {\vphantom {{\mu _\text{в}^- } K}} 
\right. \kern-\nulldelimiterspace} K,{\mu _\text{з}^+ } \mathord{\left/ {\vphantom 
{{\mu _\text{з}^+ } K}} \right. \kern-\nulldelimiterspace} K} \right)$, являющееся 
центром (в окрестности которого система колеблется с частотой $\sqrt {\mu 
_\text{з}^+ \mu _\text{в}^- } )$ причем не только в линейном приближении (локально), но и 
глобально. Это следует из консервативности системы. Первый интеграл имеет 
вид
\[
\mu _\text{в}^- \ln c_\text{з} +\mu _\text{з}^+ \ln c_\text{в} -K\cdot \left( {c_\text{з} +c_\text{в} } \right)\equiv 
\mbox{const}\left( {c_\text{з} \left( 0 \right),c_\text{в} \left( 0 \right)} \right).
\]
Все это означает негрубость системы (чувствительность к возмущениям). Таким 
образом, для данной модели реакций стохастической химической кинетики мы 
имеем неперестановочность порядка взятия предельных переходов $\mathop 
{\lim }\limits_{t\to \infty } \mathop {\lim }\limits_{N\to \infty } \ne 
\mathop {\lim }\limits_{N\to \infty } \mathop {\lim }\limits_{t\to \infty } 
$, где число молекул (агентов) = волков + зайцев стремится к бесконечности. Следует сравнить это 
замечание с замечанием к задаче 19, в котором приводится достаточное условие 
перестановочности. В данной модели можно лишь говорить о том, что система 
(Лотки--)Вольтера аппроксимирует соответствующую модель стохастической 
химической кинетики лишь на конечном промежутке времени. Этот промежуток тем 
больше, чем больше молекул (агентов). Более того, длина промежутка 
неограниченно возрастает с возрастанием числа агентов. Этот конкретный 
пример рассматривался в большом числе книг, переведенных на русский язык 
(Г.~Николис\,--\,И.Р.~Пригожин, К.В. Гардинер, В.Б. Занг).

Если рассматривать более общие биологические модели взаимодействия $n$ 
видов, то в пространственно однородном случае возникает аналогичная СОДУ:
\[
\frac{dc_i }{dt}=c_i f_i \left( {c_1 ,...,c_n } \right),
\quad
i=1,...,n.
\]
Функции $f_i \left( {c_1 ,...,c_n } \right)$  называют {\itмальтузианскими 
функциями}. Если предположить, что имеется устойчивое положение равновесия 
$\bar {c}=$\linebreak $=\underbrace {\left( {c_1 ,...,c_m ,0,...,0} \right)}_n$, 
то необходимо $f_i \left( {\bar {c}} \right)=0$ при $i=1,...,m$ и $f_i 
\left( {\bar {c}} \right)\le 0$ при $i=m+1,...,n$. Это означает (принцип 
эволюционный оптимальности), что
\begin{center}
$из \bar {c}_i >0$ следует $f_i \left( {\bar {c}} \right)=\mathop {\max 
}\limits_{j=1,...,n} \;f_j \left( {\bar {c}} \right)$
\end{center}
или, что то же самое (задача дополнительности, которая во многих важных 
приложениях переписывается как (монотонное) вариационное 
неравенство, которое в случае потенциальности векторного поля $f$ 
переписывается как задача выпуклой оптимизации (см., например, задачу 15)):
\[
\bar {c}_i \cdot \left( {\mathop {\max }\limits_{j=1,...,n} \;f_j \left( 
{\bar {c}} \right)-f_i \left( {\bar {c}} \right)} \right)=0.
\]
Биологический смысл этого условия заключается в том, что выжившие в том или 
ином состоянии виды обязаны иметь максимальное значение мальтузианских 
коэффициентов среди всех потенциально допустимых, вычисленных в этом 
состоянии. Именно эти коэффициенты характеризуют силу видов в ее 
дарвиновском понимании. Подробнее об этом и о возможных обобщениях можно 
посмотреть, например, в учебном пособии  Разжевайкин В.Н. Анализ моделей динамики популяций. 
М.: МФТИ, 2010. 

Отметим, что поскольку матрица Якоби имеет верхнетреугольный блочный 
вид 
$$ J = \left\| {\frac{\partial \left( {c_i f_i \left( c \right)} 
\right)}{\partial c_j }} \right\|_{c=\bar {c}} =\left( 
{{\begin{array}{*{20}c}
 A \hfill & B \hfill \\
 0 \hfill & D \hfill \\
\end{array} }} \right),$$ 
спектр матрицы (по теореме Лапласа) задается 
спектром двух подматриц $A$ и $D=\mbox{diag}\left\{ {f_i \left( {\bar {c}} 
\right)} \right\}_{i=m+1}^n $. Для устойчивости необходима и достаточна 
внутренняя устойчивость (матрицы $A$), означающая устойчивость к изменению 
численностей присутствующих в равновесии видов, и внешняя (матрицы $D)$, 
означающая устойчивость к привнесению в систему извне видов, которые в 
равновесие не представлены. Последнее означает, что $f_i \left( {\bar {c}} 
\right)\le 0$ при $i=m+1,...,n$ (необходимое условие устойчивости матрицы 
$D)$.

В популяционной теории игр рассмотренному выше придается немного другая 
интерпретация. А именно $c_i $ -- концентрация (доля) игроков, использующих 
стратегию $i$, $g_i \left( {c_1 ,...,c_n } \right)$ -- выигрыш от 
использования стратегии $i$ при заданных пропорциях. Каждый игрок может 
выбирать любую из $n$ стратегий. Точнее, каждый игрок (агент), использующий 
стратегию $i$ ($i=1,...,n)$, независимо от остальных с 
интенсивностью $\lambda c_j \left[ {g_j \left( c \right)-g_i \left( 
c \right)} \right]_+ $, меняет свою стратегию на $j\ne i$ (система 
представляет собой модель стохастической химической кинетики с унарными 
реакциями и с непостоянными коэффициентами реакций, см. задачу 15).

Если число игроков стремится к бесконечности так, что корректно определены 
начальные концентрации $c\left( 0 \right)$, то концентрации корректно 
определены и при $t\ge 0$ (теорема Т. Куртца), причем вектор-функция 
$c\left( t \right)$ удовлетворяет выписанной ранее СОДУ с $f_i =g_i -\bar 
{g}$, где $\bar {g}\left( c \right)=\sum\limits_{i=1}^n {c_i g_i \left( c 
\right)} $, которую (для заданного протокола) называют {\itдинамикой 
репликаторов}. В такой интерпретации полученное равновесие приобретает смысл 
равновесия Нэща в популяционной игре. А~принцип эволюционный оптимальности 
можно понимать как условие равновесия по Нэшу. Все необходимые детали 
имеются в книге \cite{222}. 

Рассмотренный в задаче протокол (способ изменения своей стратегии) 
называют {\itпопарным сравнением}. Но можно рассматривать и многие другие 
протоколы, см., например, книгу W. Sandholm'a \cite{222} и задачи~14, 15. Общей 
замечательной особенностью протоколов, отражающих рациональность игроков, 
является (при весьма общих условиях, но не всегда -- см., например, задачу~14) одинаковость аттрактора, что можно грубо сформулировать так: ``все 
дороги ведут в Рим''. То есть при весьма общих условиях, если игроки 
стремятся менять, исходя из текущих наблюдаемых концентраций, свои стратегии 
в сторону увеличения собственного выигрыша, то соответствующая 
стохастическая динамика имеет инвариантную (стационарную) меру, которая, 
часто концентрируется, с ростом числа игроков, около (устойчивых) равновесий 
Нэша, множество которых не зависит от выбранного протокола. В свою очередь, 
множество, на котором концентрируется стационарная мера (в общем случае это 
могут быть и не равновесия Нэша, до которых динамика может ``не дойти'', 
скажем, ``закружившись в цикл'', которому соответствует ``кратер вулкана'' 
-- см., например, задачу~14), можно рассматривать как аттрактор 
соответствующей СОДУ.

Предложенный в задаче скейлинг далеко не единственно возможный (см., 
например, задачу 9). Здесь мы приведем еще один вариант. Если считать, что 
$n_\text{з} =\kappa _\text{з} \varepsilon ^{-2}$, $n_\text{в} =\kappa _\text{в} \varepsilon ^{-2}$, и 
предположить, что находящиеся в одном узле волк и заяц прореагируют не 
мгновенно, а с интенсивностью $K$ (в остальном так же, как и раньше, с 
поправкой, что $\lambda _\text{з} :=\lambda _\text{з} \varepsilon ^{-2}$, $\lambda _\text{в} 
:=\lambda _\text{в} \varepsilon ^{-2})$, то для любого $t\ge 0$ и $x\in \left[ 
{0,1} \right)^2$ существует (при $x_i/\varepsilon \leq  k_i <{x_i }/\varepsilon +1,\;i=1,2$) 
\[
c_\text{з} \left( {t,x} \right)=\mathop {\lim }\limits_{{\scriptsize
 \varepsilon \to 0+ }} \Exp\left[ {n_\text{з} \left( {t;k_1 ,k_2 } \right)} \right],
\]
где $n_\text{з} \left( {t;k_1 ,k_2 } \right)$ -- число зайцев в момент времени $t$ 
в узле решетки с номером $\left( {k_1 ,k_2 } \right)$ (аналогично 
определяется $c_\text{в} \left( {\tau ,x} \right))$. Причем эти функции 
удовлетворяют следующей системе (Лотки--Вольтера с диффузией) уравнений 
реакции-диффузии: 
\[
\frac{\partial c_\text{з} }{\partial t}=\lambda _\text{з} \Delta c_\text{з} +\mu _\text{з}^+ c_\text{з} -Kc_\text{в} 
c_\text{з} ,
\]
\[
\frac{\partial c_\text{в} }{\partial t}=\lambda _\text{в} \Delta c_\text{в} +Kc_\text{в} c_\text{з} -\mu _\text{в}^- 
c_\text{в} .
\]
Детали и обобщения имеются в статье {Малышев В.А., Пирогов С.А. }Обратимость и необратимость в 
стохастической химической кинетике // Успехи мат. наук. 2008. Т. 63, вып. 
1(379). С. 4--36. В частности, если жизнь ``разворачивается'' на $d$-мерном 
торе, то изменения будут только в определении начальных численностей 
$\varepsilon ^{-2}\Rightarrow \varepsilon ^{-d}$.

\end{remark}

\begin{problem}\Star(Введенская--Замятин--Малышев.) 
Ось $x$ разбита на ячейки, расстояние между 
соседними ячейками $\varepsilon $. В ячейках могут находиться частицы ``+'' 
и ``\textbf{--}''. Однако частицы разного знака не могут находиться в одной 
ячейке, поскольку моментально прореагируют и исчезнут. Каждая частица ``+'' 
с интенсивностью $\lambda _+ / \varepsilon$ независимо от остальных перепрыгивает в 
соседнюю справа ячейку (с большей координатой), аналогично, каждая частица 
``\textbf{--}'' с интенсивностью $\lambda_{-} / \varepsilon $ независимо от остальных 
перепрыгивает в соседнюю слева ячейку. Введя плотность числа частиц (с 
учетом знака) $c\left( {t,x} \right)$, считая, что можно корректно  
определить $c\left( {0,x} \right)=  \bigg \{ \begin{array}{l}
 -2 ,\quad  x \leq 0 , \\ 
 \,\,\,\, 1 ,\quad x > 0. \\ 
 \end{array} $, получите в пределе при $\varepsilon \to 0+$ уравнение
\[
\frac{\partial c}{\partial t}+\frac{\partial \phi \left( c \right)}{\partial 
x}=0,
\quad
\phi \left( c \right)=\left\{ {\begin{array}{l}
 c \lambda _+  ,\quad c>0, \\ 
 | c | \lambda _- ,\quad c\le 0. \\ 
 \end{array}} \right.
\]
Как понимать решение (начальной) задачи Коши для этого уравнения? Покажите, 
что обобщенное решение этой задачи Коши будет иметь вид ударной волны. 
Определите скорость этой ударной волны (формула Римана--Рэнкина--Гюгонио).
\end{problem}

\begin{remark}
Подробнее об обобщенных решениях законов сохранения, в 
том числе ударных волнах, написано, например, во второй главе монографии 
Введение в математическое моделирование транспортных потоков. М.: МЦНМО, 
2013.

К сожалению, из-за негладкости в нуле функции $\phi \left( c \right)$ 
напрямую нельзя использовать описанные в этой книге наработки. Тем не менее 
если немного сгладить функцию $\phi \left( c \right)$, например, за счет 
такой зависимости $\lambda \left( c \right)$ при малых $c$, что $\lambda _- 
\left( 0 \right)=\lambda _+ \left( 0 \right)=0$, то обобщенное (по 
Олейник--Кружкову) решение уже будет определено корректно и будет совпадать 
с формально посчитанной (по формуле Римана--Рэнкина--Гюгонио) ударной волной 
в первоначальной постановке.
\end{remark}

\begin{problem} (Игрушечная модель эволюции\ Д. Кьялво.) 
В аудитории собрали большое число студентов. Каждый студент, соответствующий 
определенному биологическому виду, имеет приспособленность (число от нуля до 
единицы). Динамика (эволюция) состоит в том, что на каждом шаге выбирается 
наименее приспособленный вид (студент с наименьшей приспособленностью) и еще 
какие-то два студента, которые выбираются независимо и 
равновероятно. Для каждого из этих трех студентов случайно (независимо 
и с равномерным распределением на $\left[ {0,1} \right])$ выбираются новые 
значения приспособленностей. Определим порог приспособленности $c\in \left[ 
{0,1} \right]$. Будем считать, что число студентов стремится к 
бесконечности, а доля студентов с приспособленностью ниже этого порога 
стремится к нулю. Тогда $c$ определяется из условия, что после одного шага 
эволюции среднее число студентов ниже порога приспособленности не изменится. 
Покажите, что ${c=1} \mathord{\left/ {\vphantom {{c=1} 3}} \right. 
\kern-\nulldelimiterspace} 3$. Предположим, что в начальный момент времени 
все студенты (кроме одного) имеют приспособленности выше пороговой. Лавина начинается тогда, 
когда появляется очередной студент(-ы) с приспособленностью ниже пороговой, а 
заканчивается тогда, когда остается один студент с приспособленностью  выше 
пороговой. Покажите, что распределение лавин по длительности имеет степенной 
закон. Здесь имеется в виду, что если понаблюдать за студентами достаточно 
долго и нарисовать график частота встречаемости лавины / длительность 
лавины, то с хорошей точностью получится степенной закон с показателем 
степени ${\tau =3} \mathord{\left/ {\vphantom {{\tau =3} 2}} \right. 
\kern-\nulldelimiterspace} 2$.

\end{problem}

\begin{ordre}
Отметим, что в исходной постановке Д. Кьялво студенты были 
размещены на окружности и\;\,``еще два студента'' выбирались по принципу 
соседства с наименее приспособленным. Однако это более сложная для 
исследования ситуация. Впрочем, имеют место оценки: $c\simeq 0.667$, $\tau 
\simeq 1.07$, обоснование которых требует более изощренного определения порога 
приспособленности.

Для решения данной задачи воспользуйтесь степенным законом распределением времени 
первого возвращения симметричного случайного блуждания на прямой, см. задачу \ref{bluzd_ust} раздела \ref{zb4}
и обзор Newman M.E.J. Power laws, Pareto distributions and Zipf's law // Contemporary 
Physics. 2005. V. 46, N 5. P. 323--351. Данная задача немного проясняет один 
из аспектов того, почему степенные законы так часто возникают при описании 
зависимости частоты различных катастроф от их масштабов и ряда других 
зависимостей, встречающихся в природе, см. Бак~П. Как работает природа: Теория 
самоорганизованной критичности.  М.: УРСС, 2013. 

\end{ordre}

\begin{problem} (Majority rule model.) 
Популяция состоит из $N\gg 1$ 
агентов, каждый из которых имеет спин (мнение) $\pm 1$. На каждом шаге 
случайно (независимо и равновероятно) выбираются три агента. Если у всех 
агентов одинаковый спин, то ничего не происходит, иначе большинство (два 
агента) ``уговаривает'' (заставляет) оставшегося агента поменять свой спин. 
Опишите, что произойдет с агентами по прошествии достаточно большого 
времени. Постарайтесь получить количественные оценки, в том числе времени 
выхода в равновесие.
\end{problem}

\begin{remark}
Описание более общих моделей можно найти в обзоре 
{Castellano C., Fortunato S., Loreto V. }Statistical physics of social dynamics //   
{arXiv:0710.3256v2} (2009). Так же, как и в предыдущей задаче, модель усложняется, 
если вводить (пространственно) локальные взаимодействия. В частности, многое 
зависит от размерности пространства, в котором расположена решетка со 
спинами. Детали имеются, например, в монографии Лигетт Т.М.  Марковские процессы с 
локальным взаимодействием. М.: Мир, 1989. 
\end{remark}

\begin{problem}\Star(Динамика наилучших ответов для игры 
``Камень--Ножницы--Бумага''.)  Популяция состоит из $N\gg 1$ агентов. 
Стратегией каждого агента является выбор камня, ножниц или бумаги. Выигрыш 
(в единицу времени) агента, выбравшего камень, есть разность [доля агентов, 
выбравших ножницы] -- [доля агентов, выбравших бумагу], выбравшего ножницы есть 
разность [доля агентов, выбравших бумагу] -- [доля агентов, выбравших камень], 
выбравшего бумагу есть разность [доля агентов, выбравших камень] --- [доля 
агентов, выбравших ножницы]. При этом с каждым агентом связан свой 
``пуассоновский будильник''. Такой будильник звонит через 
случайные (одинаково распределенные) промежутки времени, каждый промежуток 
-- независимая показательная случайная величина. Параметр этой случайной 
величины называется параметром будильника. Фактически мы уже сталкивались с 
пуассоновским будильником, например, в задаче ``парадокс Эренфестов'' и в ее 
распределенном варианте. Все будильники имеют одинаковый параметр. Когда 
будильник звонит, агент выбирает исходя из текущей ситуации камень, ножницы 
или бумагу, чтобы получить максимальный выигрыш (best response protocol; 
динамика наилучших ответов). Опишите, что произойдет с агентами по 
прошествии достаточно большого времени. Постарайтесь получить количественные 
оценки.
\end{problem}

\begin{remark}
 Для решения этой и следующей задачи рекомендуется 
посмотреть монографию \cite{222}. Интересно 
отметить, что в этом примере есть чувствительность к выбору ``разумной'' 
динамики. Это в некотором смысле нетипично, если речь идет о сходимости к 
единственному равновесию, как в данной задаче. Скажем, для популярной в 
популяционной теории игр динамики репликаторов в пределе $N\to \infty $ 
(этот предел нужен, чтобы можно было смотреть на стохастическую динамику, 
как на детерминированную) концентрации агентов, использующих разные 
стратегии, будут колебаться вокруг равновесия $\left( {1 \mathord{\left/ 
{\vphantom {1 3}} \right. \kern-\nulldelimiterspace} 3,1 \mathord{\left/ 
{\vphantom {1 3}} \right. \kern-\nulldelimiterspace} 3,1 \mathord{\left/ 
{\vphantom {1 3}} \right. \kern-\nulldelimiterspace} 3} \right)$, и лишь в 
чезаровском смысле будут сходиться к этому равновесию.
\end{remark}

\begin{problem} (Модель равновесного распределения транспортных потоков.)
Ориентированный граф $\Gamma =\left( {V,E} \right)$ -- транспортная сеть 
города\linebreak ($V$~-- узлы сети (вершины), $E\subset V\times V$ -- дуги сети (рёбра 
графа)). Пусть $W=\left\{ {w=\left( {i,j} \right):\;i,j\in V} \right\}$ -- 
множество пар источник--сток; $p=\left\{ {v_1 ,v_2 ,...,v_m } \right\}$ -- 
путь из вершины $v_1 $ в $v_m $, если $\left( {v_k ,v_{k+1} } \right)\in E$, 
$k=1,...,m-1$; $P_w $ -- множество путей, отвечающих корреспонденции $w = (i,j)\in 
W$ (т.е. множество путей, начинающихся в в вершине $i$ и заканчивающихся в вершине $j$); $P=\bigcup\nolimits_{w\in W} {P_w } $ -- совокупность всех путей в сети 
$\Gamma $; $x_p $ -- величина потока по пути $p$, $x=\left\{ {x_p :\;\;p\in 
P} \right\}$; $y_e $ -- величина потока по дуге $e$: $y_e =\sum\limits_{p\in 
P} {\delta _{ep} x_p } $, где $\delta _{ep} =\left\{ {\begin{array}{l}
 1,\mbox{ }e\in p \\ 
 0,\mbox{ }e\notin p \\ 
 \end{array}} \right.$ ($y=\Theta x$, $\Theta =\left\{ {\delta _{ep} } 
\right\}_{e\in E,p\in P} )$; $\tau _e \left( {y_e } \right)$ -- удельные 
затраты на проезд по дуге $e$ (гладкие неубывающие функции); $G_p \left( x 
\right) \quad G_p \left( x \right)=\sum\limits_{e\in E} {\tau _e \left( {y_e } 
\right)\delta _{ep} } $ -- удельные затраты на проезд по пути $p$, $G\left( 
x \right)=\left\{ {G_p \left( x \right):\;\;p\in P} \right\}$ ($G\left( x 
\right)=\Theta ^T\tau \left( y \right))$. Пусть известна матрица 
корреспонденций $\left\{ {d_w } \right\}_{w\in W} $ ($d_{ij}$ -- сколько агентов перемещается из вершины $i$ в вершину $j$). Тогда вектор $x$, 
характеризующий распределение потоков по путям, должен лежать в допустимом множестве: 
$X=\left\{ {x\ge 0:\;\;\sum\limits_{p\in P_w } {x_p } =d_w ,\;w\in W} 
\right\}$. Рассмотрим игру, в которой каждому элементу $w\in W$ 
соответствует свой, достаточно большой ($d_w \gg 1)$, набор однотипных 
агентов (``сидящих на корреспонденции $w$''). Множеством чистых стратегий 
каждого такого агента является $P_w $, а выигрыш (потери со знаком минус) от 
использования стратегии $p\in P_w $ определяется формулой $-G_p \left( x 
\right)$. Переопределим параметры модели
\[
d_w = d_w M \, (M \gg 1), 
\quad
x = x / M, 
\quad
y = y / M, 
\quad
\tau_e(y_e) = \tau_e(y_e/M).   
\]

а)* Покажите, что при $M \to \infty$ динамика наилучших ответов из предыдущей задачи 
``приводит'' к равновесию $x^\ast $ Нэша--Вардропа, которое для данной игры 
загрузок (а стало быть, потенциальной игры) сводится к задаче выпуклой 
оптимизации (потенциала):
\[
\Psi \left( {y\left( x \right)} \right)=\sum\limits_{e\in E} 
{\int\limits_0^{y_e \left( x \right)} {\tau _e \left( z \right)dz} } \to 
\mathop {\min }\limits_{x\in X} .
\]

б) (Энтропийная регуляризация.)** Легко проверить, что если $\tau _e 
^\prime \left( {\,\cdot \,} \right)>0$, $e\in E$, то равновесное 
распределение потоков по ребрам $y^\ast $ единственно (см. п. а)). Но это еще не 
гарантирует единственности $x^\ast $ (напомним связь: $y\left( x 
\right)=\Theta x)$. Когда равновесных распределений потоков по путям может 
быть много, то нужно понять, а какое из этих равновесий все-таки реализуется 
на практике (если реализуется, ведь может быть, например, ни к чему не 
сходящееся блуждание по равновесному множеству). Предположим, что каждый 
агент независимо принимает решения на основе текущей зашумленной информации 
$\arg \mathop {\max }\limits_{p\in P_w } \left\{ {-G_p \left( x \right)+\xi 
_p } \right\}$, где $\xi _p $ -- независимые одинаково распределенные 
случайные величины с нулевым математическим ожиданием и ограниченной 
дисперсией. Будем считать, что $\xi _p $ имеют распределение Гумбеля (см. 
задачу \ref{gumbel} раздела \ref{hard}): $\PR\left( {\xi _p \leq \xi } \right)=\exp 
\left\{ {-e^{-\xi \mathord{\left/ {\vphantom {\xi \sigma }} \right. 
\kern-\nulldelimiterspace} \sigma -E}} \right\}$, $E\approx 0.5772$ -- 
константа Эйлера, а $\Var\xi _p ={\sigma ^2\pi ^2} \mathord{\left/ {\vphantom 
{{\sigma ^2\pi ^2} 6}} \right. \kern-\nulldelimiterspace} 6$. Покажите, что при	 $M \to \infty$
такая динамика (ее называют имитационная логит-динамика) приводит к стохастическому равновесию, которое определяется как 
единственное решение сильно выпуклой (в 1-норме) задачи оптимизации:
\[
\Psi \left( {y\left( x \right)} \right)+\sigma \sum\limits_{w\in W} 
{\sum\limits_{p\in P_w } { {x_p \ln \left( {{x_p } \mathord{\left/ 
{\vphantom {{x_p } {d_w }}} \right. \kern-\nulldelimiterspace} {d_w }} 
\right) } } } \to \mathop {\min }\limits_{x\in X} .
\]
Если $\sigma \to 0+$, то решение этой задачи сходится к 
\[
x^\ast =\arg \mathop {\min }\limits_{x \in X:\;\Theta x = 
y^\ast } \sum\limits_{w\in W} {\sum\limits_{p\in P_w } { {x_p \ln 
\left( {{x_p } \mathord{\left/ {\vphantom {{x_p } {d_w }}} \right. 
\kern-\nulldelimiterspace} {d_w }} \right)} } } ,
\]
где $y^\ast =\arg \mathop {\min }\limits_{y=\Theta x,\;x\in X} \Psi \left( y 
\right)$.

\begin{remark}

 Распределение потоков по путям $x=\left\{ {x_p } 
\right\}\in X$ называется равновесием (Нэша--Вардропа) в популяционной игре 
$\left\langle {\left\{ {x_p } \right\}\in X,\left\{ {G_p \left( x \right)} 
\right\}} \right\rangle $, если из $x_p >0$ ($p\in P_w )$ следует $G_p 
\left( x \right)=$\linebreak $=\mathop {\min }\limits_{q\in P_w } \;G_q \left( x \right)$. 
Или, что то же самое:

\begin{center}
для любых $w\in W,\;\;p\in P_w $ выполняется $x_p \cdot \left( {G_p \left( x 
\right)-\mathop {\min }\limits_{q\in P_w } \;G_q \left( x \right)} 
\right)=0.$
\end{center}

В этой задаче имеется ряд моментов, которые связаны с возможностью 
осуществления предельного перехода по числу агентов, стремящихся к 
бесконечности, и требуют некоторых дополнительных оговорок. Детали имеются, 
например, в уже упоминавшейся книге~\cite{222} (в особенности стоит 
обратить внимание на теорему~11.5.12). Эта задача может быть обобщена на 
многостадийные транспортные модели, которые сейчас популярны на практике:

{Гасников А.В.} Эффективные численные методы поиска равновесий в больших транс-портных сетях. Диссертация на соискание степени д.ф.-м.н. по специальности 05.13.18 – Математическое моделирование, численные методы, комплексы программ. М.: МФТИ, 2016. – 487 с.  http://arxiv.org/find/all/1/all:+gasnikov/0/1/0/all/0/1

Описанная в этой задаче зашумленная динамика наилучших ответов играет важную 
роль при описании поведения ограниченно рациональных агентов: Andersen S.P., de Palma A., Thisse~J.-F. Discrete choice theory of product differentiation. Cambridge: MIT Press,
1992, а~также \cite{222}.

\end{remark}

\end{problem}

\begin{problem}\Star(Эволюция РНК по В.Г. Редько.)
Для простоты рассуждений представим молекулу РНК $S$ в виде 
последовательности битов длины $N \gg 1$. В системе имеется $n$ ($N\ll n\ll 2^N$) молекул РНК. 
В начальный момент все молекулы инициализированы независимо и случайно (из 
равномерного распределения). Будем считать, что самая приспособленная 
молекула $S_0 $ -- это молекула, состоящая из одних нулевых битов. Ее 
приспособленность $f\left( {S_0 } \right)=1$. Приспособленность молекулы $S$ 
есть $f\left( S \right)=\exp \left( {-\rho \left( {S,S_0 } \right)} \right)$, 
где $\rho $ -- расстояние Хэмминга. Динамика системы молекул в дискретном 
времени заключается в следующем. Отбираются $n$ особей в новую популяцию, 
посредством $n$ независимых генераций из дискретного распределения с 
вероятностями исходов $\sim \left\{ {f\left( {S_k } \right)} \right\}_{k=1}^n $ 
(естественный отбор). Далее в получившейся популяции 
каждый бит каждой молекулы независимо подвергается случайной мутации. Вероятность мутации $P\sim 1 
\mathord{\left/ {\vphantom {1 N}} \right. \kern-\nulldelimiterspace} N$. 
Покажите, что после $T={\it O}\left( N \right)$ шагов будет явным 
доминирование наиболее приспособленного вида. Другими словами, с большой 
вероятностью большая часть молекул будет иметь почти все свои биты нулевыми. 
Получите количественные оценки к последнему утверждению.
\end{problem}

\begin{remark}
По поводу постановки задачи и ее окрестностей см. 
{Редько В.Г.} Эволюция, нейронные сети, интеллект: Модели и концепции эволюционной 
кибернетики. М.: УРСС, 2013, а также Эйген М., Шустер П. Гиперцикл. Принципы самоорганизации макромолекул. М.: Мир, 1982; Dyson F. Origins of live. Cambridge University Press, 2004; Резерфорд А. Биография жизни. От первой клетки до генной инженерии. М.: Бином, 2016 и популярные книги и выступления А.В. Маркова (д.б.н., МГУ) на тему происхождения жизни и эволюции. Для оценки скорости сходимости (mixing time) в 
этой и других задачах рекомендуется использовать монографию \cite{240}.

Если функция 
приспособленности многоэкстремальная, а требуется найти глобальный максимум, 
то описанную в задаче динамику можно понимать как один из вариантов 
генетического алгоритма, способного найти глобальный максимум. При этом в 
генетических алгоритмах стараются выбирать $n$ на границе условия $N\ll n$, 
вплоть до $n\simeq N$. Хотя генетические алгоритмы весьма популярны в 
задачах глобальной оптимизации (см., например, Zhigljavsky A., Zilinskas A. Stochastic global 
optimization. Springer Optimization and Its Applications, 2008), однако об 
их сходимости (точнее об оценках скорости сходимости) известно не так много. 
Все же, некоторые результаты есть, см., например, статью Cerf~R. 
Asymptotic convergence of genetic algorithms // Adv. Appl. Prob. 1998. V. 
30, N~2. P. 521--550. 
\end{remark}

\begin{problem}(Эволюция РНК с рекомбинациями, но без отбора.)  
a)*  Для простоты рассуждений рассмотрим молекулу РНК $S$ как цепочку битов длины $N$. В системе имеется $n\gg 2^N$ молекул РНК. В начальный момент все молекулы инициализированы независимо и 
случайно (из равномерного распределения). С интенсивностью $\lambda $ каждый 
$0$ бит каждой молекулы независимо от остальных превращается в $1$ бит, а с 
интенсивностью $\mu $ -- наоборот. Можно  себе представить, что с каждой молекулой может происходить $2N$ различных типов 
реакций, в результате реакции меняется один из битов (соответствующей данной 
реакции). Каждой молекуле и каждому типу реакции соответствует свой 
пуассоновский будильник. Когда будильник звонит, происходит соответствующая 
реакция. Считаем также, что каждой паре молекул $x$ и $y$ соответствует 
$2^N$ однотипных реакций (рекомбинаций). Реакция заключается в том, что 
какой-то фрагмент $x_I $ (цепочка битов, не обязательно идущих подряд) одной 
молекулы копируется и заменяет соответствующий участок $y_I $ в другой 
молекуле. Интенсивность таких реакций положим известной $0\le \phi \left( 
{x_I ,y_I } \right)\le 1$, причем $\phi \left( {x_I ,y_I } \right)=\phi 
\left( {y_I ,x_I } \right)$. Обозначим через $\mu _k \left( t \right)$ -- 
долю молекул с числом единичных битов равным $k$. Опишите поведение системы 
молекул в терминах $\mu _k \left( t \right)$ на больших временах. Отметим, 
что поведение системы не зависит от вида функции $\phi \left( {x_I ,y_I } 
\right)$. Рассмотрите случай, когда $\lambda =1$, $\mu =2$. 
Постарайтесь получить количественные оценки.

б)** Предложите такое обобщение модели из п. а), в 
котором бы учитывался естественный отбор (или такое обобщение модели из 
предыдущей задачи, в котором учитывалась бы рекомбинация) и которое 
позволяет аналитически исследовать поведение на больших временах.
\end{problem}

\begin{remark}
 Для лучшего понимания постановки этой задачи просмотрите цикл недавних работ S.A. Pirogov, A.N. Rybko и др. 
на arxiv.org. Для решения данной задачи полезно выписать с 
помощью теоремы Т. Куртца (впрочем, эта теорема полезна практически во всех 
``эволюционных'' задачах) систему $2^N$ нелинейных дифференциальных 
уравнения на $\mu ^{(x)}\left( t \right)$, где $\mu^{(x)}\left( t \right)$ --- доля 
молекул с набором битов $x$. Эта система получается при предельном переходе 
$n\to \infty $ в предположении, что в начальный момент существуют 
соответствующие предельные пропорции $\mu ^{(x)}\left( 0 \right)$. Детали см. в \cite{101}.
\end{remark}

\begin{problem}\Star(Пуассоновская гипотеза.)  
Имеется $N$ заявок и $M$ 
серверов с процедурой обслуживания FIFO. После обслуживания каждая заявка с 
равной вероятностью отправляется на один из серверов. Времена обслуживания 
заявок --- независимые одинаково распределенные с.в. с ф.р. $F\left( x 
\right)$, $\mu =\int\limits_0^\infty {xdF\left( x \right)} $. Предположим, 
что $N\to \infty $, $M \mathord{\left/ {\vphantom {M N}} \right. 
\kern-\nulldelimiterspace} N\to \lambda $. Исследуйте поведение системы на 
больших временах.
\end{problem}

\begin{remark}
Данная задача является частным случаем более общего 
класса задач, в которых изучаются сети (массового обслуживания) Джексона (и 
их обобщения) при термодинамическом предельном переходе, cм., например, цикл 
работ A. Rybko, S. Shlosman на arxiv.org. Отметим в 
этой связи интересные исследования фазового перехода в таких сетях (см. 
приложение А.А. Замятина, В.А. Малышева в книге Введение в математическое 
моделирование транспортных потоков. М.: МЦНМО, 2013). Однако у приведенной 
задачи есть особенность -- время обслуживания, вообще говоря, не 
предполагается экспоненциальным. Вообще эта задача очень показательна во 
многих отношениях (см. статью А.Н. Рыбко в 4-м номере ``Глобус'').
\end{remark}

\begin{problem} (Теорема Гордона--Ньюэлла и PageRank.)
а) Имеется $N\gg 
1$ пользователей, которые случайно (независимо) блуждают в непрерывном 
времени по ориентированному графу с эргодической инфинитезимальной матрицей 
$\Lambda $. Назовем вектор $p$ (из единичного симплекса) PageRank, если 
$\Lambda p=0$. Обозначим через $n_i \left( t \right)$  число пользователей 
на $i$-й странице в момент времени $t\ge 0$. Покажите, что $n\left( t 
\right)$ асимптотически имеет мультиномиальное распределение с вектором 
параметров PageRank $p$, т.е.
\[
\mathop {\lim }\limits_{t\to \infty } \PR\left( {n\left( t \right)=n} 
\right)\sim \prod\limits_i {\frac{\left( {p_i } \right)^{n_i }}{n_i !}} .
\]
Следовательно (неравенство Хефдинга в гильбертовом пространстве, см. Часть 2),
\[
\mathop {\lim }\limits_{t\to \infty } \PR\left( {\left\| {\frac{n\left( t 
\right)}{N}-p} \right\|_2 \ge \dfrac{2\sqrt{2} +  4\sqrt{\ln(\sigma^{-1})}}{\sqrt{N}} } \right)\le \sigma .
\]

б) Получите тот же результат, что в п. а), рассмотрев 
соответствующую систему унарных химических реакций. Переход одного из 
пользователей из вершины $i$ в вершину $j$  означает превращение одной 
молекулы вещества $i$ в одну молекулу вещества $j$, $n_i \left( t \right)$ 
-- число молекул $i$-го типа в момент времени $t$. Каждое 
ребро графа соответствует определенной реакции (превращению). Интенсивность 
реакций определяется матрицей $\Lambda $ и числом молекул, вступающих в 
реакцию (закон действующих масс). Покажите, что условие $\Lambda p=0$  в 
точности соответствует условию унитарности в стохастической химической 
кинетике (обобщению условия детального баланса, предложенного в 2000~г.
Батищевой--Веденяпиным и независимо Малышевым--Пироговым--Рыбко в 2004~г.). 

\end{problem}

\begin{remark}
 По п. а) полезно посмотреть монографию {Serfozo~R.} Introduction to 
stochastic networks. Springer, 1999, а по п. б) статью {Гасников А.В}., {Гасникова Е.В.} Об 
энтропийно-подобных функционалах {\ldots} // Математические заметки. 2013. 
Т. 94, № 6. С. 816--824 и цитированную там литературу. С примером более 
общих реакций можно познакомиться, например, по задаче ``Кинетика 
социального неравенства'', в которой реакции бинарные. Отметим, что если 
воспользоваться теоремой Санова о больших уклонениях (см. задачу \ref{sanov} 
раздела \ref{zb4}) для мультиномиального распределения, то получим
\[
\frac{N!}{n_1 !\,\cdot ...\cdot n_m !}p_1^{n_1 } \cdot ...\cdot p_m^{n_m } 
=\exp \left( {-N\sum\limits_{i=1}^m {\nu _i \ln \left( {{\nu _i } 
\mathord{\left/ {\vphantom {{\nu _i } {p_i }}} \right. 
\kern-\nulldelimiterspace} {p_i }} \right)} +R} \right),
\]
где $\nu _i ={n_i } \mathord{\left/ {\vphantom {{n_i } N}} \right. 
\kern-\nulldelimiterspace} N$, $\left| R \right|\le \frac{m}{2}\left( {\ln 
N+1} \right)$. Однако последующее применение неравенства Пинскера (см. 
Часть 2) не дает нам равномерной по $m$ оценки в 
1-норме. Как и ожидалось, выписанная в задаче оценка в 2-норме (правая часть 
неравенства под вероятностью) и так полученная оценка в 1-норме, будут 
отличаться по порядку приблизительно в $\sqrt m $ раз, что соответствует 
типичному соотношению между 1- и 2-нормами (см. задачу на теорему Б.С. Кашина 
о поперечниках в Части 2). Слово ``типично'' здесь отвечает, грубо говоря, ситуации, 
когда компоненты вектора одного порядка. Для многих приложений, где 
возникают предельные конфигурации (кривые), описывающиеся вектором с 
огромным числом компонент, имеет место быстрый закон убывания этих 
компонент. Например, для некоторого обобщения (см. Райгородский~А.М. Модели Интернета. 
Долгопрудный: Изд. Дом\;``Интеллект'', 2013) модели Барабаши--Альберт 
случайного роста графа Интернета вектор PageRank с хорошей точностью и с 
высокой вероятностью имеет компоненты, убывающие по степенному закону. Если 
в этих приложениях возникает мультиномиальное распределение с таким 
вектором, то при изучении концентрации этого распределения предпочтительнее 
становится подход с 1-нормой (см., например, задачу 4 ``Кинетика социального 
неравенства'').

К п. б) можно привести следующее обобщение. Предположим, что некоторая 
макросистема может находиться в различных состояниях, характеризуемых 
вектором $n$ с неотрицательными целочисленными компонентами. Будем считать, 
что в системе происходят случайные превращения (химические реакции). 

Пусть $n\to n-\alpha +\beta $, $\left( {\alpha ,\beta } \right)\in J$ -- все 
возможные типы реакций, где $\alpha $ и $\beta $ -- вектора с 
неотрицательными целочисленными компонентами. Введем интенсивность реакции:
\[
\lambda _{\left( {\alpha ,\beta } \right)} \left( n \right)=\lambda _{\left( 
{\alpha ,\beta } \right)} \left( {n\to n-\alpha +\beta } 
\right)=
\]
\[=N^{1-\sum\limits_i {\alpha _i } }K_\beta ^\alpha 
\prod\limits_{i:\;\alpha _i >0} {n_i \cdot ...\cdot \left( {n_i -\alpha _i 
+1} \right)} ,
\]
где $K_\beta ^\alpha \ge 0$ -- константа реакции; при этом 
$\sum\nolimits_{i=1}^m {n_i \left( t \right)} \equiv N\gg 1$. Другими 
словами, $\lambda _{\left( {\alpha ,\beta } \right)} \left( n \right)$ -- 
вероятность осуществления в единицу времени перехода $n\to n-\alpha +\beta 
$. Здесь не предполагается, что число состояний $m=\dim \;n$ и число реакций 
$\left| J \right|$ не зависят от числа агентов $N$. Тем не менее если 
ничего не известно про равновесную конфигурацию $c^\ast $ (типа быстрого 
убывания компонент этого вектора), то дополнительно предполагается, что 
$m\ll N$ -- это нужно для обоснования возможности применения формулы 
Стирлинга при получении вариационного принципа (максимума энтропии) для 
описания равновесия макросистемы $c^\ast $ (в концентрационной форме). 
Однако часто априорно можно предполагать (апостериорно проверив), что 
компоненты вектора $c^\ast $ убывают быстро, тогда это условие можно 
отбросить. Так, например, обстоит дело с ``Кинетикой социального 
неравенства''.

Возникающий марковский процесс считается неразложимым. Далее приводится 
{теорема }(во многом установленная в 1999--2005~гг.\linebreak  в работах Я.Г. 
Батищевой, В.В.~Веденяпина, В.А.~Малышева, С.А.~Пи\-ро\-гова, А.Н.~Рыбко). 

а) $\left\langle {\mu ,n\left( t \right)} \right\rangle \equiv 
\left\langle {\mu ,n\left( 0 \right)} \right\rangle $ (inv) $\Leftrightarrow 
 \quad вектор \mu $  ортогонален каждому вектору семейства $\left\{ {\alpha -\beta } \right\}_{\left( {\alpha ,\beta } 
\right)\in J}$ . 

б) (Т. Куртц.)  Если существует $$\mathop {\lim }\limits_{N\to \infty } {n\left( 0 
\right)} \mathord{\left/ {\vphantom {{n\left( 0 \right)} N}} \right. 
\kern-\nulldelimiterspace} N=c\left( 0 \right),$$ где $K_\beta ^\alpha := K_\beta 
^\alpha \left( {n \mathord{\left/ {\vphantom {n N}} \right. 
\kern-\nulldelimiterspace} N} \right),$ а $m$ и $\left| J \right|$ не зависят от $N$, то для любого $t>0$ с вероятностью 1 существует $$\mathop 
{\lim }\limits_{N\to \infty } {n\left( t \right)} \mathord{\left/ {\vphantom 
{{n\left( t \right)} N}} \right. \kern-\nulldelimiterspace} N=c\left( t 
\right),$$  где $c\left(t\right)$  -- не случайная вектор-функция, удовлетворяющая СОДУ Гульдберга--Вааге:

\[\tag{ГВ}
\frac{dc_i }{dt}=\sum\limits_{\left( {\alpha ,\beta } \right)\in J} {\left( 
{\beta _i -\alpha _i } \right)K_\beta ^\alpha \left( c \right)\prod\limits_j 
{c_j^{\alpha _j } } } .
\]

Гиперплоскость (inv) (с очевидной заменой $n$ на $c$) инвариантна относительно этой динамики.   Более того, случайный процесс ${n\left( t \right)} \mathord{\left/ {\vphantom 
{{n\left( t \right)} N}} \right. \kern-\nulldelimiterspace} N$ слабо сходится при $N\to \infty$ 
 к $c\left( t \right)$  на любом конечном отрезке времени. 

Отметим, что в каком-то смысле жизнь нелинейной динамической 
системы определяется линейными законами сохранения, унаследованными ею при 
скейлинге (каноническом). Этот тезис имеет, по-видимому, более широкое 
применение (см. работы Веденяпин~В.В.).

в) Пусть выполняется условие унитарности (очевидно, что $\xi $, удовлетворяющий условию (U), -- неподвижная точка в (ГВ)):

\[\tag{U}
\exists \;\xi >0:\;\forall \;\beta 
\Rightarrow
 \sum\limits_{\alpha :\;\left( 
{\alpha ,\beta } \right)\in J} {K_\beta ^\alpha \prod\limits_j {\left( {\xi 
_j } \right)^{\alpha _j }} } =\sum\limits_{\alpha :\;\left( {\alpha ,\beta } 
\right)\in J} {K_\alpha ^\beta \prod\limits_j {\left( {\xi _j } 
\right)^{\beta _j }} } .
\]

Тогда неотрицательный ортант ${\mathbb R}_+^m $ расслаивается гиперплоскостями (inv), так что в каждой гиперплоскости (inv) уравнение (U) (положительно) разрешимо притом единственным образом. Стало быть, существует, притом единственная, неподвижная точка $c^\ast \in$ (inv) у системы (ГВ), являющаяся при этом глобальным аттрактором. Система (ГВ) имеет функцию Ляпунова 
$$KL\left( {c,\xi } 
\right)=\sum\limits_{i=1}^m {c_i \ln \left( {{c_i } \mathord{\left/ 
{\vphantom {{c_i } {\xi _i }}} \right. \kern-\nulldelimiterspace} {\xi _i }} 
\right)}. $$

Стационарное (инвариантное) распределение описанного марковского процесса имеет носителем множество (inv) и (с точностью до нормирующего множителя) имеет вид
\[
\frac{N!}{n_1 !\;\cdot \ldots \cdot n_m !}\left( \xi _1 \right)^{n_1 }\cdot \ldots\cdot \left( \xi _m  \right)^{n_m }\sim \exp \left( {-N\cdot KL\left( 
{c,\xi } \right)} \right),
\]
$где \xi $ -- произвольное решение (U), не важно какое именно (от этого, конечно, будет зависеть нормирующий множитель, но это ни на чем не сказывается). При этом условие унитарности (U) является обобщением условия детального равновесия (баланса): 
\[
\exists \;\;\xi >0:\;\;\forall \;\;\left( {\alpha ,\beta } \right)\in J 
\Rightarrow
K_\beta ^\alpha \prod\limits_j {\left( {\xi _j } \right)^{\alpha _j }} 
=K_\alpha ^\beta \prod\limits_j {\left( {\xi _j } \right)^{\beta _j }} ,
\]
принимающего такой вид для мультиномиальной стационарной\linebreak меры.

Используя неравенство Чигера \cite{44} (для оценки mixing time) и неравенство Хефдинга в гильбертовом пространстве (см. Часть 2), отсюда можно получить (зависимость $a\left( {m,c\left( 0 \right)} \right)$ во многих приложениях может быть равномерно ограничена):
\[
\exists \;\;a=a\left( {m,c\left( 0 \right)} \right):\;\;\forall \;\;\sigma 
>0,\;t\ge aN\ln N, 
\]
\[ \PR\left( {\left\| {\frac{n\left( t \right)}{N}-c^\ast } 
\right\|_2 \ge \dfrac{2\sqrt{2} + 4\sqrt{\ln(\sigma^{-1})} }{\sqrt{N}} } 
\right)\le \sigma ,
\]
где (принцип максимума энтропии Больцмана--Джейнса):
\[
c^\ast =\arg \mathop {\max }\limits_{ c \in \left( {\mbox{{\scriptsize inv}}} 
\right)} \left( {-\sum\limits_i {c_i } \ln \left( {{c_i } \mathord{\left/ 
{\vphantom {{c_i } {\xi _i }}} \right. \kern-\nulldelimiterspace} {\xi _i 
}} \right)} \right)=\arg \mathop {\min }\limits_{ c \in \left( 
{\mbox{{\scriptsize inv}}} \right)} KL\left( {c,\xi } \right),
\]
$а \xi $ -- произвольное решение (U), причем $c^\ast $  определяется единственным образом, т.е. не зависит от выбора $\xi $. Геометрически $c^\ast $ -- это\linebreak {KL-проек}\-ция произвольного $\xi $, удовлетворяющего (U), на гиперплоскость (inv), соответствующую начальным данным $c\left( 0 \right)$. Независимость этой проекции от выбора $\xi $ из (U) просто означает, что кривая (U) проходит KL-пер\-пен\-дикулярно через множество (inv).

\textbf{Замечание}. В терминах условия задачи условие унитарности просто 
означает, что в равновесии для любой вершины имеет место баланс числа 
пользователей, входящих в единицу времени в эту вершину с числом 
пользователей, выходящих в единицу времени из этой вершины. В то время как 
условие детального равновесия (см. также задачу~\ref{cript} раздела~\ref{MK}) 
означает, что в равновесии для любой пары вершин число пользователей, 
переходящих в единицу времени из одной вершину в другую, равно числу 
пользователей, переходящих в обратном направлении. Понятно, что второе 
условие является частным случаем первого.

Много интересных 
примеров макросистем, в которых $K_\beta ^\alpha =$\linebreak $=K_\beta ^\alpha \left( {n 
\mathord{\left/ {\vphantom {n N}} \right. \kern-\nulldelimiterspace} N} 
\right)$ и имеет место детальный баланс, собрано в книге  Вайдлих~В. Социодинамика: 
системный подход к математическому моделированию в социальных науках. М.: 
УРСС, 2010. Выше (задача 15) нами был рассмотрен пример на эту тему из книги~\cite{222}. Как правило, такие постановки приводят к функциям 
Ляпунова--Санова вида $\Psi \left( c \right)+\sigma KL\left( {c,\xi } 
\right)$, $\sigma \ge 0$. Причем во всех этих приложениях условие детального 
баланса выполняется точно. Отметим, что с некоторыми дополнительными 
оговорками можно допускать и приближенное выполнение условий детального 
баланса (с точностью ${\it O}\left( {N^{-1}} \right))$.

г) (Е.В. Гасникова.) Верно и обратное утверждение, то есть условие (U) не только достаточное для того, чтобы равновесие находилось из приведенной выше задачи энтропийно-линейного программирования, но и, с небольшой оговоркой (для почти всех $c\left( 0 \right)$), необходимое. Также верно и более общее утверждение, связывающее понимание энтропии в смысле Больцмана (функция Ляпунова прошкалированной кинетической динамики) и Санова (функционал действия в неравенствах больших уклонений для инвариантной меры): если стационарная мера асимптотически представима в виде $\sim \exp \left( {-N\cdot 
V\left( c \right)} \right)$ (вместе с производными до второго порядка включительно), то $V\left( c \right)$  -- функция Ляпунова (ГВ).

К п. г) можно сделать следующее пояснение. Обозначим через $h\left( c 
\right)$ вектор-функцию, стоящую в правой части СОДУ (ГВ). Тогда (см., 
например, книгу К.В. Гардинера \cite{333}) при $N \gg 1$ по теореме Т.~Куртца $n(t)/N$ будет как   
$O\left( {\log N} /{\sqrt N } \right)$-близко (в этом месте 
требуются оговорки, чтобы близость была равномерна по времени) к $x_t 
=x\left( t \right)$ -- решению стохастической системы дифференциальных 
уравнений (с начальным условием $x_0 =c\left( 0 \right))$:
\[
dx_t =h\left( {x_t } \right)dt+\sqrt {\frac{g\left( {x_t } \right)}{N}} dW_t 
,
\]
где функция $g\left( {x_t } \right)>0$ рассчитывается по набору реакций и 
константам реакций (которые могут быть не постоянны и зависеть от 
концентраций). Инвариантная (стационарная) мера $m\left( x \right)$ этого 
однородного марковского процесса удовлетворяет уравнению (см. также задачу \ref{annealing} раздела \ref{MK}):
\[
\frac{1}{2N}\nabla^2{\left(g(x) m(x) \right)}-\mbox{div}\left( {h\left( x \right)m(x)} 
\right)=0,
\]
поскольку плотность распределения $p\left( {t,x} \right)$ процесса $x_t $ 
подчиняется уравнению Колмогорова--Фоккера--Планка:
\[
\frac{\partial p\left( {t,x} \right)}{\partial t}=-\mbox{div}\left( {h\left( 
x \right)p\left( {t,x} \right)} \right)+\frac{1}{2N}\nabla^2{\left( g\left( x  \right) p\left( {t,x} \right) \right)}.
\]
Если известно, что $m\left( x \right)\simeq \mbox{const}\cdot \exp \left( 
{-N\cdot V\left( x \right)} \right)$ (и аналогично для производных до второго порядка включительно), то из уравнения для $m\left( x \right)$ 
имеем при $N \to \infty$
\[
\left\langle {h,\nabla V} \right\rangle =-\frac{1}{2}g\cdot \left( {\nabla 
V} \right)^2+ O\left( \frac{1}{N} \right) \to -\frac{1}{2}g\cdot \left( {\nabla V} 
\right)^2\le 0.
\]
Эта выкладка поясняет (но не доказывает, для доказательства требуются более 
аккуратные рассуждения), почему функция $V\left( c \right)$ может быть 
функцией Ляпунова системы (ГВ) ${dc} \mathord{\left/ {\vphantom {{dc} {dt}}} 
\right. \kern-\nulldelimiterspace} {dt}=h\left( c \right)$. Более того, 
модель стохастической химической кинетики здесь может быть заменена и более 
общими шкалирующимися марковскими моделями.

\end{remark}

\begin{problem}\Star(Теорема Гордона--Ньюэлла; Л.Г. Афанасьева.) 
Рассматривается транспортная сеть, в которой между $N$ 
станциями курсируют $M$ такси. Клиенты пребывают в $i$-й узел в соответствии с 
пуассоновским потоком с параметром $\lambda _i >0$ ($i=1,\;...,\;N)$. Если в 
момент прибытия в $i$-й узел там есть такси, клиент забирает его и с 
вероятностью $p_{ij} \ge 0$ направляется в $j$-й узел, по прибытии в который 
покидает сеть. Такси остается ждать в узле прибытия нового клиента. Времена 
перемещений из узла в узел -- независимые случайные величины, имеющие 
показательное распределение с параметром $\nu _{ij} >0$ для пары узлов 
$\left( {i,j} \right)$. Если в момент прихода клиента в узел там нет такси, 
клиент сразу покидает узел. Считая $p_{ij} =N^{-1}$, $\lambda _i =\lambda $, 
$\nu _{ij} =\nu $, покажите, что вероятность того, что клиент, поступившей в 
узел (в установившемся (стационарном) режиме работы сети), получит отказ, 
равна
\[
p_{\text{отк}} \left( {N,M} \right)=\dfrac{\sum\limits_{k=0}^M 
{\frac{C_{N-2+k}^k \rho ^{M-k}}{\left( {M-k} \right)!}}}
{
\sum\limits_{k=0}^M {\frac{C_{N-1+k}^k \rho ^{M-k}}{\left( 
{M-k} \right)!}
}},
\quad
\rho =N\lambda/ \nu.
\]
Методом перевала покажите справедливость следующей асимптотики при $N\to 
\infty $:
\[
p_{\text{отк}} \left( {N,rN} \right)=1-\frac{2r}{\lambda \mathord{\left/ 
{\vphantom {\lambda \nu }} \right. \kern-\nulldelimiterspace} \nu +r+1+\sqrt 
{\left( {\lambda \mathord{\left/ {\vphantom {\lambda \nu }} \right. 
\kern-\nulldelimiterspace} \nu +r+1} \right)^2-4{\lambda r} \mathord{\left/ 
{\vphantom {{\lambda r} \nu }} \right. \kern-\nulldelimiterspace} \nu } 
}+ O\left( {\frac{1}{N}} \right).
\]

\end{problem}

\begin{remark}
Метод перевала -- очень полезный инструмент исследования 
асимптотик интегралов по параметру. Его подробное изложение имеется, 
например, здесь: Федорюк М.В. Метод перевала. М.: УРСС, 2010.
\end{remark}

\begin{problem}\DStar(Mean field games; Lasry--Lions.)
\label{m_field_games}
Пусть динамика агента 
(игрока) задается стохастическим дифференциальным уравнением
\[
dX_t =-\alpha _t dt+\sigma dW_t ,
\quad
X_0 =x_0, 
\]
где $W_t $ -- винеровский процесс, а $\alpha _t $ -- стратегия (управление), 
$\alpha _t,$\linebreak $ X_t \in {\mathbb R}^d$, Функционал потерь задается формулой
\[
J\left( {x_0 ,\alpha } \right)=\mathop {\underline{\lim }}\limits_{T\to 
\infty } \frac{1}{T}\Exp\left[ {\int\limits_0^T {\left( {L\left( {\alpha _t } 
\right)+F\left( {X_t } \right)} \right)dt} } \right],
\]
где ${L\left( \alpha \right)} \mathord{\left/ {\vphantom {{L\left( \alpha 
\right)} {\left| \alpha \right|}}} \right. \kern-\nulldelimiterspace} 
{\left| \alpha \right|}\to \infty $ при $\left| \alpha \right|\to \infty $. 
Считаем, что функция $F\left( x \right)$ 1-периодическая по каждой 
компоненте вектора $x$. Таким образом, можно ограничиться рассмотрением 
стохастической динамики на $\mbox{T}^d$ -- единичном торе в ${\mathbb R}^d$. 
Введем гамильтониан $H\left( p \right)=$\linebreak $=\mathop {\sup }\limits_{\alpha \in 
{\mathbb R}^d} \left\{ {\left\langle {p,\alpha } \right\rangle -L\left( \alpha 
\right)} \right\}$ и $\nu ={\sigma ^2} \mathord{\left/ {\vphantom {{\sigma 
^2} 2}} \right. \kern-\nulldelimiterspace} 2$.

а) (Fleming--Soner; Bardi--Capuzzo-Dolcetta.) Покажите, что 
если уравнение Гамильтона--Якоби--Беллмана, 
\[
-\nu \Delta V+H\left( {\nabla V} \right)+\lambda =F\left( x \right),
\]
имеет решение (пара $\lambda $, 
$V\left( {\,\cdot \,} \right))$,
то оптимальное значение функционала $\mathop {\inf}\limits_\alpha 
J\left({x_0 ,\alpha}\right)\!\!=\!\!J\left( {x_0 ,\hat {\alpha }} \right)\!\!=\!\!\lambda 
$ и оптимальное управление (в форме синтеза) $\hat {\alpha }\left( x 
\right)=\nabla H\left( {\nabla V\left( x \right)} \right)$. Кроме того, 
марковский (диффузионный) процесс 
\[
dX_t =-\hat {\alpha }_t \left( {X_t } \right)dt+\sigma dW_t 
\] 
имеет инвариантную меру $m\left( x \right)dx$, 
являющуюся решением уравнения Колмогорова--Фоккера--Планка:
\[
\nu \Delta m+\mbox{div}\left( {\nabla H\left( {\nabla V} \right)m} 
\right)=0,
\quad
\int\limits_{\small{\mbox{T}^d}\simeq \left[ {0,1} \right]^d} {m\left( x \right)dx} 
=1,
\quad
m\left( x \right)>0.
\]
Ввиду эргодичности $X_t $, $m\left( x \right)$ -- асимптотическое 
распределение вероятностей положения агента, ведущего себя оптимальным 
образом.

б) (A. Friedman; Bensoussan--Frehse.) Предположим теперь, что 
имеется $N>1$ игроков, и у каждого своя динамика:
\[
dX_t^i =-\alpha _t^i dt+\sigma ^idW_t ,
\quad
X_0^i =x_0^i ,
\]
и своя функция потерь:
\[
J^i\left( {\left\{ {x_0 } \right\},\left\{ \alpha \right\}} \right)=\mathop 
{\underline{\lim }}\limits_{T\to \infty } \frac{1}{T}\Exp \left[ 
{\int\limits_0^T {\left( {L^i\left( {\alpha _t^i } \right)+F^i\left( {X_t^1 
,...,X_t^N } \right)} \right)dt} } \right].
\]
Аналогично видоизменяются определения $H^i\left( p \right)$ и $\nu ^i$. 
Равновесие Нэша $\left\{ {\bar {\alpha }} \right\}$ в такой игре 
определяется из условия (для любого\linebreak $i=1,...,N)$
\[
J^i\left( {\left\{ {x_0 } \right\},\left\{ {\bar {\alpha }} \right\}} 
\right)=\mathop {\min }\limits_{\alpha ^i} J^i\left( {\left\{ {x_0 } 
\right\},\bar {\alpha }^1,...,\bar {\alpha }^{i-1},\alpha ^i,\bar {\alpha 
}^{i+1},...,\bar {\alpha }^N} \right).
\]
Можно показать, что система уравнений ($i=1,...,N)$
\[
-\nu ^i\Delta V_i +H^i\left( {\nabla V_i } \right)+\lambda _i =f^i\left( 
{x;m_1 ,...,m_N } \right),
\]
\[
\nu ^i\Delta m_i +\mbox{div}\left( {\nabla H^i\left( {\nabla V_i } 
\right)m_i } \right)=0,
\]
\[
\int\limits_{\small{\mbox{T}^d}} {m_i \left( x \right)dx} =1,
\quad
m_i \left( x \right)>0,
\quad
\int\limits_{\small{\mbox{T}^d}} {V_i \left( x \right)dx} =0,
\]
где
\[
f^i\left( {x;m_1 ,...,m_N } \right)=\int\limits_{T^{d\left( {N-1} \right)}} 
{F^i\left( {x^1,...,x^{i-1},x,x^{i+1},...,x^N} \right)\prod\limits_{j\ne i} 
{dm_j \left( {x^j} \right)} } ,
\]
имеет решение $\lambda _i $, $V_i \left( {\,\cdot \,} \right)$, $m_i \left( 
{\,\cdot \,} \right)$, $i=1,...,N$. Покажите, что для любого такого решения 
$\hat {\alpha }^i\left( x \right)=\nabla H^i\left( {\nabla V_i \left( x 
\right)} \right)$ -- равновесие Нэша и $\lambda _i =$\linebreak $=J^i\left( {\left\{ {x_0 
} \right\},\left\{ {\hat {\alpha }} \right\}} \right)$ -- соответствующее 
значение игры.

в)  Пусть для всех игроков $\nu ^i=\nu $, $H^i=H$, а $F^i$ зависит 
только от $X^i$ и эмпирической плотности распределения других игроков, т.е.
\[
F^i\left( {x^1,...,x^N} \right)=\Phi \left[ {\frac{1}{N-1}\sum\limits_{j\ne 
i} {\delta _{x_j } } } \right]\left( {x^i} \right),
\]
где
\[
\Phi :\;\left\{ {\mbox{распределения вероятностей на T}^d} \right\}\to 
\]
\[\to
\left\{ {\mbox{липшицевы функции на T}^d} \right\}
\]
-- непрерывный интегральный функционал
\[
\Phi \left[ m \right]\left( x \right)=G\left( {x,\int\limits_{\mbox{T}^d} 
{k\left( {x,y} \right)dm\left( y \right)} } \right).
\]
Можно показать, что последовательность $\lambda _i^{\left( N \right)} $, 
$V_i^{\left( N \right)} \left( {\,\cdot \,} \right), m_i^{\left( N 
\right)} \left( {\,\cdot \,} \right)$, \linebreak $(N\to \infty)$ относительно 
компактна в ${\mathbb R}\times C^2( {\mbox{T}^d} )\times W^{1,p}( 
{\mbox{T}^d} )$ (для любого $1\le p<\infty )$. Следовательно, можно 
выбрать сходящуюся подпоследовательность. Зафиксируем $i$. Покажите, что при 
$N\to \infty $ предел любой такой сходящейся подпоследовательности является 
решением системы ($x\in \mbox{T}^d)$:
\[
-\nu \Delta V+H\left( {\nabla V} \right)+\lambda =\Phi \left[ m 
\right]\left( x \right),
\]
\[
\nu \Delta m+\mbox{div}\left( {\nabla H\left( {\nabla V} \right)m} 
\right)=0,
\]
\[
\int\limits_{\mbox{T}^d} {m\left( x \right)dx} =1,
\quad
m\left( x \right)>0,
\quad
\int\limits_{\mbox{T}^d} {V\left( x \right)dx} =0.
\]
Можно показать, что решение этой системы существует. Если\linebreak $\Phi \left[ m 
\right]\left( x \right)=G\left( {m\left( x \right)} \right)$, $G$ -- 
возрастающая функция, а $H$ -- сильно выпуклая функция, то можно 
гарантировать и единственность.

\end{problem}

\begin{remark}
 Игры среднего поля -- относительно новая и стремительно развивающаяся область исследований (ее возраст около 
10 лет). В этом направлении ведут исследования небольшие группы в 
Екатеринбурге (Ю.В. Авербух) и в Москве (В.Н. Колокольцов), а также французскими и американскими математиками: P.L. Lions, J.M. Lasry, O. Gueant, M. Bardi, A. Bensoussan, J. 
Frehse, D. Lacker и др.
\end{remark}

\begin{problem} (Круговая модель М. Каца.) 
На окружности отмечено $n$ 
равноотстоящих точек, $m$ из которых составляют множество $Q$. В каждой 
точке помещен белый (б) или черный (ч) шар. В единицу времени каждый шар 
переходит против часовой стрелки в соседнюю точку, причем, если шар выходит 
из точки, принадлежащей $Q$, он меняет свой цвет. Пусть $\mu =m 
\mathord{\left/ {\vphantom {m n}} \right. \kern-\nulldelimiterspace} n<1 
\mathord{\left/ {\vphantom {1 2}} \right. \kern-\nulldelimiterspace} 2$. 
Обозначим через $N_\text{б} \left( t \right)$  число белых шаров в момент времени 
$t$, $N_\text{ч} \left( t \right)$ -- число черных шаров в момент времени $t$. 

а) Покажите, что динамика обратима и $2n$-периодична. Если 
предположить, что множество $Q$ выбирается случайно, то обратимость 
теряется. Ограничимся также рассмотрением $t\ll n$. Покажите, что при этих 
предположениях
\[
\Exp_Q \left[ {N_\text{б} \left( t \right)-N_\text{ч} \left( t \right)} \right]=\left( 
{1-2\mu } \right)^t \Exp_Q \left[ {N_\text{б} \left( 0 \right)-N_\text{ч} \left( 0 \right)} 
\right].
\]
б)* Введем в динамику п. а) случайность. Для этого определим 
последовательность с.в., независимую от случайности при выборе 
множества $Q$ ($p<1 \mathord{\left/ {\vphantom {1 2}} \right. 
\kern-\nulldelimiterspace} 2)$:
\[
\chi \left( t \right)=\left\{ {\begin{array}{l}
 1\;\;\;\mbox{с вероятностью }p, \\ 
 -1\;\mbox{с вероятностью }1-p. \\ 
 \end{array}} \right.
\]
Пронумеруем последовательность точек на окружности и сопоставим белому цвету 
число $+1$, черному $-1$. Введем также функцию (случайную, поскольку $Q$ 
случайно)
\[
\delta _k =\left\{ {\begin{array}{l}
 -1\;\;\;k\in Q, \\ 
 +1\;\;\;k\notin Q. \\ 
 \end{array}} \right.
\]
Пусть $C_k \left( t \right)$ обозначает цвет шара в точке $k$ в момент 
времени $t$. В~таких обозначениях новая динамика будет иметь вид
\[
C_k \left( {t+1} \right)=\chi \left( {t+1} \right)\delta _k C_k \left( t 
\right),
\]
т.е. с вероятностью $p$ шар еще может (случайно) изменить цвет (все эти 
обозначения введены для удобства доказательства). Для простоты будем 
считать, что в начальный момент все шары белые. Покажите, что при $t\ll n$
\[
\Exp_{\delta ,\chi } \left[ {\frac{N_\text{б} \left( t \right)-N_\text{ч} \left( t 
\right)}{n}} \right]\sim \left( {1-2\mu } \right)^t\left( {1-2p} \right)^t,
\]
\[
\Var_{\delta ,\chi } \left[ {\frac{N_\text{б} \left( t \right)-N_\text{ч} \left( t 
\right)}{n}} \right]\sim \frac{1}{n}\left( {1-2\mu } \right)^{2t}\left( 
{1-2p} \right)^{2t}.
\]
Покажите также, что если множество $Q$ зафиксировано (не случайно), то 
сомножитель со степенью $\left( {1-2\mu } \right)$ исчезает вместе с 
ограничением $t\ll n$.

\end{problem}

\begin{remark}
 Более подробно об этой модели написано в старых, но до 
сих пор актуальных и популярных книгах М.~Каца, переведенных на русский 
язык (см.,~нап\-ри\-мер,~\cite{20}), а также в монографии {Опой\-цев~В.И.} Нелинейная системостатика. М.: Наука, 1986. 
Отметим, что если  при задании множества $Q$ перейти  по принципу случайного 
выбора $m$ позиций к случайному разыгрыванию каждой позиции: с~вероятностью 
$\mu $ каждая точка (независимо от остальных) принадлежит~$Q$, а~с~вероятностью $1-\mu $ не принадлежит, то получение фактора $\left( {1-2\mu } 
\right)^t$ существенно упрощается. В статистической физике такой переход и 
его различные обобщения называют {\itметодом большого канонического ансамбля}. 
Хотя сама по себе модель Каца в каком-то смысле является ``карикатурной'' 
моделью статистической физики, она 
по-прежнему  вызывает интерес  ведущих специалистов в этой области и смежных 
областях (см., например, работы В.В.~Козлова, С.З.~Ад\-жиева--В.В.~Ве\-де\-ня\-пина). Именно так, по-видимому, и задумывалось  
 М.~Кацем, который пытался обобщить модель из парадокса Эренфестов, 
см.~за\-да\-чу~1. Тем не менее, эта модель оказалась очень важной стартовой 
площадкой для многих исследований ввиду своей простоты с одной стороны и 
возможностью продемонстрировать разнообразные аспекты статистической 
физики с другой: парадокс обратимости, парадокс времени возвращения.

\end{remark}

\section{Метод Монте-Карло}
\label{MK}

\begin{problem}
На плоскости дано ограниченное измеримое по Лебегу множество~$S$. Требуется найти площадь (меру Лебега) этого множества \linebreakс заданной точностью~$\varepsilon$. 

Поскольку по условию множество ограничено, то вокруг него можно описать квадрат со стороной $a$. Выберем декартову систему координат 
в одной из вершин квадрата с осями, параллельными сторонам квадрата. Рассмотрим $n$  независимых с.в. $\{ X_k\}_{k=1}^{n}$,  имеющих 
одинаковое равномерное распределение в этом квадрате, т.е. $X_k\in R([0,a]^2)$. Введем с.в. 
$$
Y_k=I(X_k\in S)=\begin{cases}
1,\quad X_k\in S,\\
0, \quad X_k\notin S.
\end{cases} 
$$
Тогда $\{ Y_k\}_{k=1}^{n}$ -- независимые одинаково распределенные с.в.. Ясно, что $Y_k\in\Be(p(S))$. Следовательно, по УЗБЧ 
$$
\frac{Y_1+\ldots+Y_n}{n} \xrightarrow{\text{ п.н. }} {\mathbb E}(Y_1)=p(S)=\frac{\mu(S)}{a^2} \quad \text{ при  } n\to\infty . 
$$
Оцените сверху следующую вероятность:
$$
{\mathbb P}\Bigl( \Bigl| \frac{Y_1+\ldots+Y_n}{n}-\frac{\mu(S)}{a^2}\Bigr|>\delta \Bigr) . 
$$
\end{problem}

\begin{ordre}
См. раздел 5, а также Часть 2.
\end{ordre}

\begin{problem}(Вычисление значения интеграла.) а) Требуется вычислить с заданной точностью $\varepsilon $ и с заданной доверительной вероятностью $\gamma $ абсолютно сходящийся интеграл 
\[
J=\int _{\left[0,\; 1\right]^{m} }f\left({x}\right)d{x}.
\] Считайте, что $\forall \; \; {x}\in \left[0,\; 1\right]^{m}$  выполнено $\left|f\left({x}\right)\right|\le 1$.

\begin{remark}
Введем случайный $m$-вектор ${X}\in R\left(\left[0,\; 1\right]^{m} \right)$ и с.в.\linebreak $\xi =f\left({X}\right)$. Тогда $\Exp\xi =\int _{\left[0,\; 1\right]^{m} }f\left({x}\right)d{x} =J$. Поэтому получаем оценку интеграла $\bar{J}_{n} =\frac{1}{n} \sum _{k=1}^{n}f\left(\tilde{X}_{k} \right) $, где $\tilde{X}_{k} $, $k=1,...,n$ -- повторная выборка значений случайного вектора ${X}$ (т.е. все $\tilde{x}_{k} $, $k=\overline{1,n}$,  независимы и одинаково распределены: так же, как и вектор $X$). В~задаче требуется оценить сверху число $n$ ($n\gg m$), начиная с которого $\mathbb{P}\left(\left|J-\bar{J}_{n} \right|\le \varepsilon \right)\ge \gamma $.
\end{remark}

\indent б) Решите задачу из п. а) при дополнительном предположении липшицевости функции $f\left({x}\right)$, разбив единичный куб на $n=N^{m} $ одинаковых кубиков со стороной ${1\mathord{\left/ {\vphantom {1 N}} \right. \kern-\nulldelimiterspace} N} $ и используя оценку $\bar{J}_{n} =$\linebreak $=\frac{1}{n} \sum _{k=1}^{n}f\left(\tilde{X}_{k} \right) $, где $\tilde{X}_{k} $ -- имеет равномерное распределение в \textit{k}-м кубике.

\indent в) (Метод включения особенности в плотность.) Решите задачу \linebreak п. а), не предполагая, что $f\left({x}\right)$ -- ограниченная функция на единичном кубе. Предложите способы уменьшения дисперсии полученной оценки интеграла. Как можно использовать информацию об особенностях функции $f\left({x}\right)$?
\end{problem}
\begin{ordre}
См. книгу Соболь И.М. Численный метод Монте-Карло.  М.: Наука, 1973, а также \cite{24}.
\end{ordre}

\begin{problem}Стандартный способ моделирования с.в. -- {\it метод обратной функции}. Покажите, что если с.в. $\eta $ равномерно распределена на отрезке $\left[0,1\right]$, то с.в. $\xi =F^{-1} \left(\eta \right)$ имеет функцию распределения $F\left(x\right)$. Предполагается, что $F\left(x\right)$ непрерывна и строго монотонна. Как выглядит формула для моделирования с.в. из показательного распределения с функцией распределения $F\left(x\right)=\left(1-e^{-\lambda x} \right)I\left(x\geq 0\right)$?
\begin{remark}
См. также задачу \ref{inverse} раздела \ref{standart}.
\end{remark}
\end{problem}

\begin{problem}
Пусть с.в. $\eta _{1} $, $\eta _{2}$ независимы и имеют равномерное распределение на отрезке $\left[0,1\right]$. Докажите, что с.в. $X$ и $Y$, где $X=\sqrt{-2\ln \eta _{1} } \cos \left(2\pi \eta _{2} \right)$, $Y=\sqrt{-2\ln \eta _{1} } \sin \left(2\pi \eta _{2} \right)$ -- независимые и одинаково распределенные: стандартно нормально ${\rm {\mathcal N}}\left(0,1\right)$.

\begin{ordre}
Покажите, что
\[f_{XY} (x,y)=\frac{1}{\sqrt{2\pi } } e^{-\frac{x^{2} }{2} } \frac{1}{\sqrt{2\pi } } e^{-\frac{y^{2} }{2} } =\frac{1}{2\pi } e^{-\frac{x^{2} +y^{2} }{2} } .\] 
Перейдите к полярным координатам, не забыв о якобиане замены переменных.
\end{ordre}

\end{problem}

\begin{problem}

Если $X$ -- с.в., имеющая стандартное нормальное распределение, то $X^{-2} $ имеет устойчивую плотность (см. замечание к задаче~\ref{bluzd_ust} раздела~5):
\[\frac{1}{\sqrt{2\pi } } e^{-\frac{1}{2x} } x^{-\frac{3}{2} }, \quad x>0.\] 
Используя это, покажите, что если $X$ и $Y$ -- независимые нормально распределенные с.в. с нулевым математическим ожиданием и дисперсиями $\sigma _{1}^{2} $ и $\sigma _{2}^{2} $, то величина $Z=\frac{XY}{\sqrt{X^{2} +Y^{2} } } $ нормально распределена с~дисперсией $\sigma _{3}^{2} $, такой, что $\frac{1}{\sigma _{3}^{2} } =\frac{1}{\sigma _{1}^{2} } +\frac{1}{\sigma _{2}^{2} } $.

\end{problem}

\begin{problem}(Теорема Бернштейна \cite{28}, \cite{21} Т.~2.) а) С помощью неравенства Чебышёва установите следующий результат из анализа: 

\[
\forall \; \; f\in C\left[0,1\right]\Rightarrow \left\| f_{n} -f\right\| _{C\left[0,1\right]}
\to 0, \, n\to \infty
\] 

\[
f_{n} \left(x\right)=\sum_{k=0}^{n}f\left(\frac{k}{n} \right) C_{n}^{k} x^{k} \left(1-x\right)^{n-k}. 
\]

\indent б) Исходя из  задачи~9 и п.~а) предложите способ генерирования распределения с.в. $\xi $, имеющей плотность $f_{\xi } \left(x\right)$ с финитным носителем. Пусть, для определенности,  носителем будет отрезок $\left[0,1\right]$.

\end{problem}

\begin{problem}(Метод фон Неймана.) 
Пусть с.в. $\xi $ распределена на отрезке $\left[a,b\right]$, причем ее плотность распределения ограничена: $\mathop{\max }\limits_{x\in \left[a,b\right]} f_{\xi } (x) < C$.\linebreak Пусть с.в. $\eta _{1} $, $\eta _{2} $, \dots \ независимы и равномерно распределены на $\left[0,1\right]$, $X_{i} =a+\left(b-a\right)\eta _{2i-1} $, $Y_{i} =C\eta _{2i} $, $i=1,2,...$, т.е. пары $\left(X_{i} ,Y_{i} \right)$ независимы и равномерно распределены в прямоугольнике $\left[a,b\right]\times \left[0,C\right]$. Обозначим через $\nu $ номер первой точки с координатами $\left(X_{i} ,Y_{i} \right)$, попавшей под график плотности $f_{\xi } (x)$, т.е.\linebreak $\nu =\min \left\{i:\: Y_{i} \le f_{\xi } (X_{i} )\right\}$. Положим $X_{\nu } =\sum _{n=1}^{\infty }X_{n} I\left(\nu =n\right) $.

\begin{enumerate}
\item Покажите, что с.в. $X_{\nu } $ распределена так же, как $\xi $.

\item Сколько в среднем точек $\left(X_{i} ,Y_{i} \right)$ потребуется «вбросить» в прямоугольник $\left[a,b\right]\times \left[0,C\right]$ для получения одного значения $\xi $?

\item Предложите модификацию рассмотренного метода для генерации дискретной случайной величины, принимающей значения $\lbrace 1, 2, ... , k \rbrace$ с одинаковой вероятностью, имея в распоряжении монету (генератор бинарной случайной величины).   
\end{enumerate}
\end{problem}

\begin{problem}
Как с помощью с.в. $\xi $, равномерно распределенной на отрезке $\left[0,1\right]$ ($\xi \in R\left[0,1\right]$), и симметричной монетки построить с.в. $X$, имеющую плотность распределения $f_{X} (x)=\frac{1}{4} \left(\frac{1}{\sqrt{x} } +\frac{1}{\sqrt{1-x} } \right)$, $x\in \left(0,1\right)$?
\end{problem}

\begin{problem}

Пусть $\xi $ распределена на $\left[0,1\right]$ с плотностью $f_{\xi } (x)$, представимой в виде степенного ряда $\sum _{k=0}^{\infty }a_{k} x^{k}  $ с $a_{k} \ge 0$. Положим $p_{k} ={a_{k} \mathord{\left/ {\vphantom {a_{k}  (k+1)}} \right. \kern-\nulldelimiterspace} (k+1)} $. Тогда $f_{\xi } (x)=\sum _{k=0}^{\infty }p_{k} \cdot (k+1)x^{k}  $. Примените \textit{метод суперпозиции} для моделирования с.в. $\xi $.

\begin{ordre}
\textit{Метод суперпозиции}:

\indent а) Разыгрывается значение дискретной с.в., принимающей значения $k=0,1,2,...$ с вероятностями $p_{k}, \, k = 0,1,2\dots$.

\indent б) Для каждого $k$ моделируется с.в. с функцией распределения $F_{k} (x) = x^{k+1}$, $x\in [0,1]$ (например, методом обратной функции).

\end{ordre}

\end{problem}

\begin{problem}
В случае общего положения невозможно построить генератор дискретной с.в., принимающей $n$ значений, быстрее чем за $O(n)$. Однако, существуют процедуры построения генератора такие, что можно, начиная со второго раза, генерировать распределение этой с.в. за $O(\log n)$. Более того, если потребуется сгенерировать распределение немного отличной с.в. (например, в которой поменялось несколько вероятностей исходов и, как следствие, нормировочный множитель всего распределения), то это также можно сделать за $O(\log n)$. Покажите, что если случайная величина принимает разные значения с одинаковыми вероятностями, то  построить генератор такой с.в. можно и за $O(\log n)$.
\end{problem}

\begin{problem}(Алгоритм Кнута–Яо.)
С помощью бросаний симметричной монетки требуется сгенерировать распределение заданной дискретной с.в., принимающей конечное число значений. Обобщите описанную ниже схему на общий случай. Предположим, что нам нужно сгенерировать распределение с.в., принимающей три значения 1, 2, 3 с равными вероятностями 1/3. Действуем таким образом. Два раза кидаем монетку: если выпало 00, то считаем, что выпало значение 1, если 01, то 2, если 11, то 3. Если 10, то еще два раза кидаем монетку и повторяем рассуждения. Покажите, что в среднем с помощью не более чем $\log_2 (n - 1) + 2$ подбрасываний симметричной монетки можно сгенерировать распределение дискретной с.в., принимающей, вообще говоря, с разными вероятностями $n$ значений.
\end{problem}
\begin{remark}
См. Ермаков С.М. Метод Монте-Карло в вычислительной математике.  М.: Бином, 2009.
\end{remark}

\begin{problem}(Метод Уокера.)
Пусть с.в. $\xi$ принимает значения $1, \ldots, n$ с вероятностями $p_1, \ldots, p_n$. Докажите, что 
с.в. $\xi$ можно сгенерировать при помощи смеси дискретных распределений с двумя исходами и с.в. $\xi'$, принимающей значения $1, \ldots, n$ с одинаковыми вероятностями. 
\end{problem}

\begin{ordre}

\begin{enumerate}
\item Покажите, что $\exists i, j \neq i: \;  p_i \leq 1/n, \; p_i + p_j > 1/n$.
\item Зафиксируем пару $i, j$ с указанным свойством и определим с.в. $\eta^{(1)}$ с двумя исходами:
\[
\mathbb{P}(\eta^{(1)}=i)=q_i^{(1)} = n p_i, \; \mathbb{P}(\eta^{(1)}=j)=q_j^{(1)} = 1 - n p_i, \; q_k^{(1)} = 0.   
\]
Покажите, что $\xi$ может быть представлена в виде смеси $\eta^{(1)}$ и случайной величины $\zeta^{(1)}$, принимающей значения из множества $\{1, \ldots, n\} \backslash i$. 
\item Повторите рассуждения пункта б) для с.в. $\zeta^{(1)}$ и таким образом получите, что 
\[
p_l = \frac{1}{n} \sum \limits_{j=1}^n q_l^{(j)}.
\]
\end{enumerate}

См. также  Ермаков С.М. Метод Монте-Карло в вычислительной математике.  М.: Бином, 2009.
\end{ordre}

\begin{problem}(Теорема Пайка.) \label{paika} 
Пусть ${X}_{k} ,\; k=1,...,n$ -- независимые равномерно распределенные на отрезке $\left[0,1\right]$ с.в. Упорядочим эти с.в., введя обозначения
\[\mathop{\min }\limits_{k=1,...,n} {X}_{\left(k\right)} ={X}_{\left(1\right)} \le ...\le {X}_{\left(n\right)} =\mathop{\max }\limits_{k=1,...,n} {X}_{\left(k\right)} .\] 
Пусть ${\eta}_{k} ,\; k=1,...,n+1$ -- независимые показательно распределенные с.в.: 
\[\mathbb{P}\left({\eta}_{k} >t\right)=e^{-t} ,\; t\ge 0.\] 
Покажите, что
\[\left({X}_{\left(1\right)} {\rm ,}\; {X}_{\left(2\right)} {\rm ,...,}\; {X}_{\left(n\right)} \right) \mathop{=}\limits^{d} \left(\frac{{\eta}_{1} }{\sum _{k=1}^{n+1}{\eta}_{k}  } {\rm ,}\; \frac{{\eta}_{1} +{\eta}_{2} }{\sum _{k=1}^{n+1}{\eta}_{k}  } {\rm ...,}\; \frac{\sum _{k=1}^{n}{\eta}_{k}  }{\sum _{k=1}^{n+1}{\eta}_{k}  } \right). \] 
\end{problem}

\begin{remark}
См. Лагутин М.Б. Наглядная математическая статистика.  М.: Бином, 2009 и Кендалл М., Моран П. Геометрические вероятности.  М.: Наука, 1972.
\end{remark}

\begin{problem}(Генерация точек из равномерного распределения.) а) Пусть ${X}$~--  $n$-мерный вектор с независимыми одинаково распределенными компонентами с распределением $\mathcal{N}(0,1)$. Покажите, что ${\xi} = \frac{{X}}{\|{X}\|_2}$ имеет  равномерное распределение на единичной сфере в $\mathbb{R}^n$. 

\indent б) Пусть ${X}$ -- вектор с независимыми компонентами, каждая с распределением Лапласа (т.е. с плотностью $p(x)=\frac{1}{2}\exp(-|x|)$). Какое распределение имеет вектор  ${\xi} = \frac{{X}}{\|{X}\|_1}$?
\end{problem}

\begin{problem}\Star (Сдвиг Бернулли.)
Рассмотрим динамическую систему (ДС):\linebreak $x\to\{2x\}$ ($\{\cdot \}$ -- дробная часть числа), преобразования отрезка $X =$\linebreak $= [0,1]$ в себя. Покажите, что для почти всех (по равномерной мере Лебега на отрезке $X$) точек старта полученная c помощью ДС последовательность точек будет ``квазислучайной'' (схожей с последовательностью независимых одинаково распределенных на отрезке $X$ с.в.), то есть для неё, например, справедливы ЗБЧ и ЦПТ.
\begin{ordre}
 Pokhodzei~B.B. Niederreiter H. Random number generation and quasi-Monte Carlo methods // CBMS-NSF regional conference series in applied mathematics, V. 63. Philadelphia, Penn.: SIAM, 1992.
\end{ordre}
\begin{remark}
В контексте этой задачи рекомендуется ознакомиться также с понятиями непредсказуемой последовательности, типичной последовательности, случайной (сложной) по Колмогорову последовательности, например, по книге Верещагин Н.К., Успенский В.А., Шень А. Колмогоровская сложность и алгоритмическая случайность.  М.: МЦНМО, 2013. Задача отражает то обстоятельство, что случайность ``переносится'' из начальных данных (даже простым в смысле Колмогорова алгоритмом). Рассматриваемая динамическая система последовательно считывает числа после запятой в двоичной записи начальной точки. Почти все числа из отрезка $[0,1]$ сложны (несжимаемы) в смысле Колмогорова, правда для подавляющего большинства конкретных чисел это утверждение не может быть доказано. Студенты ФУПМ имеют возможность познакомиться с колмогоровской сложностью и алгоритмическими вопросами теории вероятностей по курсу (книге) Вьюгин В.В. Колмогоровская сложность и алгоритмическая случайность.  М.: МФТИ, 2012. В частности, рекомендуется ознакомиться со сложностным доказательством закона повторного логарифма и эргодической теоремы Биркгофа\linebreak(--Хинчина), а также с  вопросами эффективности эргодической теоремы. 

Другая причина появления случайности -- это поведение динамических систем (например, рассмотренной) в условиях небольших внешних возмущений. По данной тематике стоит отметить книгу Опойцева~В.И. Нелинейная системостатика.  М.: Наука, 1986 и в цикле недавних работ В.А. Малышева и А.А.~Лыкова. Вообще между динамическими системами и случайными (марковскими) процессами имеется глубокая взаимообогащающая связь. См., например, конструкцию Улама в книге Бланк~М.Л. Устойчивость и локализация в хаотической динамике.  М.: \mbox{МЦНМО}, 2001 и Синай~Я.Г. Как математики изучают хаос // Математическое просвещение. 2001. Вып.~5. С.~32--46. Интересные взгляды на эту науку также имеются в книгах Вент\-цель~А.Д., Фрейд\-лин~М.И. Флуктуации в динамических системах под воздействием малых случайных возмущений.  М.: Наука, 1979, \cite{333}, \cite{101}.

Ну а самое главное, что в приложениях нам нужна не ``природная случайность''. Нам нужно лишь выполнение для полученной квазислучайной последовательности неких тестов типа ЗБЧ, ЦПТ с такими же оценками скорости сходимости (или не сильно худшими). Однако оказывается, что псевдослучайная последовательность может даже увеличить скорость сходимости по сравнению с настоящими случайными данными. Обнаружено это было более 40 лет назад (сюда можно отнести метод выбора узлов Холтона--Соболя, обобщающий идею задачи 2 п. б) настоящего раздела), но по-прежнему активно используется на практике (см. Соболь И.М. Численные методы Монте-Карло.   М.: Наука, 1973) и позволяет вместо $~1/\sqrt{n}$ получать точность $~n^{\varepsilon - 1}$ cо сколь угодно малым $\varepsilon>0$, см. также
blogs.princeton.edu/imabandit/2014/12/22/guest-post-by-\-sasho-\-nikolov-\-beating-\-monte-\-carlo.

В заключение хочется отметить, что в последние десятилетия очень бурно развивается область Theoretical Computer Sсience, связанная с изучением генераторов псевдослучайных чисел (в том числе в связи с проблемой $P \ne NP$), см.  
Разборов А.А. Theoretical Computer Sсience: взгляд математика. 2013. \\ 
\noindent\small{\verb|http://www.people.cs.uchicago.edu/razborov/files/computerra.pdf|}

\end{remark}
\end{problem}

\begin{problem}(PageRank.)
Ориентированный граф $G=\left\langle {V,E} \right\rangle$  сети Интернет представляется в виде набора web-страниц $V$ и ссылок между ними: запись $\left( {i,j} \right)\in E$ означает, что на $i$-й странице имеется ссылка на $j$-ю страницу. 

\begin{enumerate}

\item По web-графу случайно блуждает пользователь. За один такт 
времени пользователь, находящийся на web-странице с номером $i$, с вероятностью $p_{ij} $ переходит по ссылке на web-страницу $j$. 
Пусть из любой web-страницы можно по ссылкам перейти на любую другую 
web-страницу (условие неразложимости) графа $G$. Проверьте, что при бесконечно долгом блуждании доля времени, которую пользователь проведет на web-странице с номером $k$, есть $p_k $, где $ 
{p} =\left( {p_1 ,...,p_n } \right)^T $, ${p}^T = {p}^T P$, 
$P=\left\| {p_{ij} } \right\|_{i,j=1}^{n,n} $ -- стохастическая матрица, $\sum_k p_k = 1$ (решение единственно  ввиду неразложимости $P$). Обратим внимание, что ответ не зависит от того, с какой вершины стартует пользователь.

\item В условиях предыдущего пункта пустим независимо блуждать по web-графу $N$ 
пользователей ($N\gg \left| V \right| \gg 1)$. Пусть $n_i \left( t \right)$ -- число посетителей web-страницы $i$ в момент времени $t$. Считая стохастическую матрицу $P$ неразложимой и апериодической, покажите, 
что
\[
\exists \;\;\lambda _{0.99} >0,\;T_G>0:\;\;\;\forall \;\;t\ge 
T_G,
\]
\[
\PR\left( {\left. {\left| {\frac{n_k \left( t \right)}{N}-p_k } \right|\le 
\frac{\lambda _{0.99} }{\sqrt N }} \right)\ge 0.99} \right.,
\]
где $ {p}^T= {p}^T P$ (решение единственно). Как и в п. а), и здесь ответ не зависит от того, с каких вершин стартуют пользователи.

\end{enumerate}

\end{problem}

\begin{remark}
Если организовать случайные блуждая, следуя п.~б) 
(отметим, что эти блуждания хорошо распараллеливаются по числу 
блуждающих пользователей), то при определенных условиях можно получить решение задачи  ${p}^T = {p}^TP$ значительно быстрее, чем, скажем, ${\it O}\left( {n^2} \right)$. Такой способ численного поиска вектора $ {p}$ основан на методе Markov 
chain Monte Carlo. 
Детали и ссылки см., например, в работе Гасников А.В., Дмитриев Д.Ю. Об эффективных рандомизированных алгоритмах поиска вектора PageRank // ЖВМиМФ. 2015. Т.~55:3. С.~355--371.  arXiv:1410.3120 и ее развитие arXiv:1701.02595. Упомянем также недавнюю работу Belloni~A., Chernozhukov V. On the Computational Complexity of MCMC-based Estimators in Large Samples //  arXiv:0704.2167 (2012),  содержащую строгие результаты об эффективной вычислимости байесовских оценок. Если 
использовать неравенства концентрации меры и иметь оценки на спектральную 
щель матрицы $P$ (см. \cite{44,240}), то приведенный в п.~б) результат можно также сделать  более строгим, а именно точнее 
оценивать скорость сходимости и плотность концентрации (см.  задачу~19 раздела~\ref{macrosystems}).
\end{remark}

\begin{problem}\Star\,\,\Star (Markov Chain Monte Carlo Revolution и состоятельность
оценок максимального правдоподобия; P. Diaconis.)
\label{cript}
В руки опытных криптографов попалось закодированное письмо (10~000 символов). Чтобы это 
письмо прочитать, нужно его декодировать. Для этого берется стохастическая 
матрица переходных вероятностей $P=\left\| {p_{ij} } \right\|$ (линейный 
размер которой определяется числом возможных символов (букв, знаков 
препинания и т.п.) в языке, на котором до шифрования было написано письмо, -- 
этот язык известен и далее будет называться базовым), в которой $p_{ij} $ -- это вероятность появления символа с номером~$j$ сразу после символа 
под номером~$i$. Такая матрица может быть идентифицирована с помощью 
статистического анализ какого-нибудь большого текста, скажем, ``Войны и 
мира'' Л.Н.~Толстого.

Пускай способ (де)шифрования (подстановочный шифр) определяется некоторой 
неизвестной дешифрующей функцией $\bar {f}$ -- преобразование 
(перестановка) множества кодовых букв во множество символов базового языка.

В качестве ``начального приближения'' выбирается какая-то функция $f$, 
например, полученная исходя из легко осуществимого частотного анализа. Далее 
рассчитывается вероятность выпадения полученного закодированного текста 
$ {x}$, сгенерированного при заданной функции $f$ (функция 
правдоподобия):

\[
L\left( {{x};f} \right)=\prod\limits_k {p_{f\left( {x_k } 
\right),f\left( {x_{k+1} } \right)} } . 
\]

Случайно выбираются два аргумента у функции $f$ и значения функции при этих 
аргументах меняются местами. Если в результате получилась такая $f^\ast $, 
что $L\left( { {x};f^\ast } \right)\ge L\left( { {x};f} \right)$, то 
$f:=f^\ast $, иначе независимо бросается монетка с вероятностью выпадения 
орла\linebreak $p={L\left( { {x};f^\ast } \right)} \mathord{\left/ {\vphantom 
{{L\left( { {x};f^\ast } \right)} {L\left( { {x};f} \right)}}} 
\right. \kern-\nulldelimiterspace} {L\left( { {x};f} \right)}$, и если 
выпадает орёл, то $f:=f^\ast $, иначе $f:=f$. Далее процедура повторяется (в 
качестве $f$ выбирается функция, полученная на предыдущем шаге).

Объясните, почему предложенный алгоритм после некоторого числа итераций с 
большой вероятностью и с хорошей точностью восстанавливает дешифрующую 
функцию $\bar {f}$? Почему сходимость оказывается такой быстрой (0.01 с на современном PC)?

\end{problem}

\begin{ordre}
Описанный в задаче пример взят из обзора Diaconis P. The Markov 
chain Monte Carlo revolution // Bulletin (New Series) of the AMS. 2009. V.~49:2. P.~179--205. Детали того, что будет написано далее, можно найти в работах Jerrum~M., Sinclair~A. The Markov chain Monte Carlo method: an approach to 
approximate counting and integration // Approximation Algorithms for NP-hard 
Problems / D.S.~Hochbaum ed. Boston: PWS Publishing, 1996. P. 482--520; Lezaud~P. Chernoff-type bound for finite Markov chain // Annals of Applied Probability. 1998. V. 8. P. 849--867; Joulin~A., Ollivier~Y. Curvature, concentration and 
error estimates for Markov chain Monte Carlo // Ann. Prob. 2010. V. 38:6. 
P. 2418--2442; Paulin~D. Concentration inequalities for Markov chains by 
Marton couplings. 2013. arXiv:1212.2015v2. См. также \cite{44}, \cite{240}.

Для того чтобы построить однородный дискретный марковский процесс с конечным 
числом состояний, имеющий наперед заданную инвариантную (стационарную) меру 
$\pi $, переходные вероятности ищутся в следующем виде: $p_{ij} =p_{ij}^0 
b_{ij} $, $i\ne j$; $p_{ii} =1-\sum\limits_{j:\;\;j\ne i} {p_{ij} } $, где 
$p_{ij}^0 $ -- некоторая ``затравочная'' матрица, которую будем далее 
предполагать симметричной. Легко проверить, что матрица $p_{ij} $ имеет 
инвариантную (стационарную) меру $\pi $, если при $p_{ij}^0 >0$
\[
\frac{b_{ij} }{b_{ji} }=\frac{\pi _j p_{ji}^0 }{\pi _i p_{ij}^0 }=\frac{\pi 
_j }{\pi _i }.
\]
Чтобы найти $b_{ij} $, достаточно найти функцию $F:\;{\mathbb{R}}_+ \to \left[ 
{0,1} \right]$ такую, что

$$\frac{F\left( z \right)}{F\left( {1 \mathord{\left/ {\vphantom {1 z}} 
\right. \kern-\nulldelimiterspace} z} \right)}=z\quad \text{и}\quad b_{ij} =F\left( 
{\frac{\pi _j p_{ji}^0 }{\pi _i p_{ij}^0 }} \right)=F\left( {\frac{\pi _j 
}{\pi _i }} \right).$$

Пожалуй, самый известный пример (именно он и использовался в задаче) такой 
функции $\tilde {F}\left( z \right)=\min \left\{ {z,1} \right\}$ -- алгоритм 
Хастингса--Метрополиса. Заметим, что для любой такой функции $F\left( z 
\right)$ имеем $F\left( z \right)\le \tilde {F}\left( z \right)$. Другой 
пример дает функция $F\left( z \right)=z \mathord{\left/ {\vphantom {z 
{\left( {1+z} \right)}}} \right. \kern-\nulldelimiterspace} {\left( {1+z} 
\right)}$. Заметим также, что $p_{ij}^0 $ обычно выбирается равным $p_{ij}^0 
=1 \mathord{\left/ {\vphantom {1 M}} \right. \kern-\nulldelimiterspace} M_i 
$, где $M_i $ -- число ``соседних'' состояний у $i$, или
\[
p_{ij}^0 =1 \mathord{\left/ {\vphantom {1 {\left( {2M} \right)}}} \right. 
\kern-\nulldelimiterspace} {\left( {2M} \right)},
\quad
i\ne j;
\quad
p_{ii}^0 =1 \mathord{\left/ {\vphantom {1 2}} \right. 
\kern-\nulldelimiterspace} 2,
\quad
i\ne j.
\]
При больших значениях времени $t$, согласно эргодической теореме, имеем, что 
распределение вероятностей близко к стационарному $\pi $. Действительно, при 
описанных выше условиях имеет место условие детального баланса (марковские 
цепи, для которых это условие выполняется, иногда называют обратимыми):
\[
\pi _i p_{ij} =\pi _j p_{ji} ,\;i,j=1,...,n,
\]
из которого сразу следует инвариантность меры $\pi $, т.е.
\[
\sum\limits_i {\pi _i p_{ij} } =\pi _j \sum\limits_i {p_{ji} } =\pi _j 
,\;j=1,...,n.
\]
Основное применение замеченного факта состоит в наблюдении, что время выхода 
марковского процесса на стационарную меру (mixing time) во многих случаях 
оказывается удивительно малым. Притом что 
выполнение одного шага по времени случайного блуждания по графу, отвечающему 
рассматриваемой марковской цепи, как следует из алгоритма Кнута--Яо (см.~задачу~11), также 
может быть быстро сделано. Таким образом, довольно часто можно получать 
эффективный способ генерирования распределения дискретной случайной величины 
с распределением вероятностей $\pi $ за время, полиномиальное от логарифма 
числа компонент вектора $\pi $.

Для лучшего понимания происходящего в условиях задачи отметим, что одним из 
самых универсальных способов получения асимптотически наилучших оценок 
неизвестных параметров по выборке является метод наибольшего правдоподобия 
(И.А. Ибрагимов -- Р.З. Хасьминский, В.Г. Спокойный). Напомним вкратце, в чем он 
заключается. Пусть имеется выборка из распределения, зависящего от 
неизвестного параметра -- в нашем случае выборкой $ {x}$ из 10~000 
элементов будет письмо, а неизвестным\; ``параметром''\; будет функция~$f$. 
Далее считается вероятность (или плотность вероятности в случае непрерывных 
распределений) $L\left( { {x};f} \right)$ того, что выпадет данный $
{x}$ при условии, что значение параметра $f$. Если посмотреть на $L\left( { {x};f} \right)$ 
как на распределение в 
пространстве параметров ($ {x}$  зафиксирован), то при большом объеме 
выборки (размерности~${x}$) при естественных условиях это распределение 
концентрируется в малой окрестности наиболее вероятного значения
\[
f\left( {{x}} \right)=\arg \mathop {\max }\limits_f L\left( { 
{x};f} \right),
\]
которое ``асимптотически'' совпадает с неизвестным истинным значением $\bar 
{f}$.

\end{ordre}

\begin{remark}
Для оценки mixing time нужно оценить спектральную щель 
стохастической матрицы переходных вероятностей, задающей исследуемую 
марковскую динамику, то есть нужно оценить расстояние от максимального 
собственного значения этой матрицы, равного единице (теорема 
Фробениуса--Перрона), до следующего по величине модуля. Именно это число 
определяет основание геометрической прогрессии, мажорирующей исследуемую 
последовательность норм разностей расстояний между 
распределением в данный момент времени и стационарным (финальным) распределением.
Для оценки спектральной щели разработано довольно много 
методов, из которых мы упомянем лишь некоторые: неравенство Пуанкаре 
(канонический путь), изопериметрическое неравенство Чигера (проводимость), с 
помощью техники каплинга (получаются простые, но, как правило, довольно 
грубые оценки), с помощью каплинга Мертона, с помощью дискретной кривизны 
Риччи и теорем о концентрации меры (Мильмана--Громова). Приведем некоторые 
примеры применения MCMC: тасование $n$ карт разбиением приблизительно на 
две равные кучи и перемешиванием этих куч (mixing time$\sim \log _2 
n)$;\footnote{ Здесь контраст проявляется, пожалуй, наиболее ярко. Скажем, 
для колоды из 52 карт пространство состояний марковской цепи будет иметь 
мощность 52! (если сложить времена жизней в наносекундах каждого человека, 
когда либо жившего на Земле, то это число на много порядков меньше 52!). В 
то время как такое тасование: взять сверху колоды карту и случайно 
поместить ее во внутрь колоды, отвечающее определенному случайному 
блужданию, с очень хорошей точностью выйдет на равномерную меру, отвечающую 
перемешанной колоде, через каких-то 200--300 шагов. Если брать тасование 
разбиением на кучки, то и того меньше -- за 8--10 шагов. Это обнаружил P. Diaconis.} Hit and Run 
(mixing time $\sim n^3)$; модель Изинга -- $n$ спинов на отрезке, 
стационарное распределение = распределеню Гиббса, глауберова динамика 
(mixing time $\sim n^{{2\log _2 e} \mathord{\left/ {\vphantom {{2\log _2 e} 
T}} \right. \kern-\nulldelimiterspace} T}$, $0<T\ll 1$ (см. задачу 18)); проблема поиска 
кратчайших гамильтоновых путей; имитация отжига (см. задачу 20) для решения 
задач комбинаторной оптимизации, MCMC для решения задач перечислительной 
комбинаторики. Но, пожалуй, самым известным примером (Dyer--Frieze--Kannan) 
является полиномиальный вероятностный алгоритм (работающий быстрее известных 
``экспоненциальных'' детерминированных) приближенного поиска центра тяжести 
выпуклого множества и вычисления его объема. Одна из работ в этом 
направлении была удостоена премии Фалкерсона -- аналога Нобелевской премии в 
области Computer Science. Близкие идеи используются и при применении 
экспандеров в Computer Science. В 2010 году премия Неванлинны 
была вручена Д. Спилману, в частности, за сублинейное (по числу элементов 
матрицы, отличных от нуля) решение системы линейных уравнений и эффективное использование экспандеров (см.~Часть~2). 

\end{remark}

\begin{problem} (Одномерная модель Изинга \cite{44}, \cite{240}.) 
Рассмотрим конечный отрезок одномерной целочисленной решетки $\{0,1,\dots,n\}$, в каждой вершине $k$ которой находится спин, принимающий два значения $\sigma(k)=\pm1$. При этом считаем, что $\sigma(0)=\sigma(n)=1$. Определим гамильтониан системы $H(\sigma) = \sum_{k=0}^{n-2}(1-\sigma(k)\sigma(k+1))/2$. Определим расстояние Больцмана--Гиббса по формуле $\pi(\sigma) = Z^{-1}\exp(-\beta H(\sigma))$, $\beta = T^{-1}>0$~-- величина, обратная <<температуре>>, а $Z$~-- нормирующий множитель (статсумма). Одним из способов построения однородной дискретной марковской цепи с заданным стационарным распределением $\pi(\sigma)$ является использование распределение Глаубера:
\begin{enumerate}
\item Выбираем $k\in\{1,\dots,n-1\}$ согласно равномерному распределению.
\item С вероятностью $$p=\exp(\beta H(\sigma_{k,+1}))/\left(\exp(\beta H(\sigma_{k,+1}))+\exp(\beta H(\sigma_{k,-1}))\right)$$ новым состоянием будет $\sigma_{k,+1}$, а с вероятностью $1-p$ состояние  $\sigma_{k,-1}$,
где $$\sigma_{k,+1} = \sigma_{k}(i),\,i\not=k,\,\sigma_{k,+1} = 1;$$
$$\sigma_{k,-1} = \sigma_{k}(i),\,i\not=k,\,\sigma_{k,-1} = -1.$$
Покажите, что характерное время выхода на стационарное распределение (mixing time) 
этой марковской цепи оценивается сверху как $n^{2\log_2(\beta\exp(1))}$ в предположении $\beta \gg 1$.
\end{enumerate}
\end{problem}

\begin{problem}\Star (Hit and Run.)
Ряд задач, в которых используется метод Монте-Карло, предполагает возможность случайно равномерно набрасывать точки в некоторое наперед заданное множество (не обязательно выпуклое и связное). Например, при вычислении интеграла с помощью метода
Монте-Карло или при численном решении задач оптимизации, в которых нужно уметь
приближенно находить центр тяжести множества. Исходя из MCMC подхода (см. \cite{240}) обоснуйте,
аккуратно оговорив детали, следующий способ генерации точек.

Берем любую точку внутри множества и проводим случайно направление через эту точку, далее с помощью
граничного оракула случайно генерируем (с помощью равномерного распределения) на
этом направлении внутреннюю точку рассматриваемого множества. Через эту точку снова
проводим случайное направление и т.д. Хорошо ли будет работать Hit and Run для вытянутых множеств или для множеств, имеющих достаточно острые углы? Предложите модификацию (например, с помощью эллипсоидов Дикина) алгоритма Hit and Run для таких “плохих” множеств. Предложите другие способы случайно равномерно набрасывать точки в некоторое наперед заданное множество (например, Shake and Bake). Хорошо ли будет работать метод Shake and Bake для множеств в пространствах большой размерности? Как следует действовать, если рассматриваемое множество имеет простую структуру: $n$-мерный куб, $n$-мерный шар, $n$-мерный симплекс, многогранник?
\end{problem}

\begin{remark}
См. видео выступления Б.Т. Поляка “В окрестностях Монте-Карло” на mathnet.ru и оригинальную статью Lovasz L. and  Vempala S. Hit-and-run from a corner  // SIAM Journal on Computing. 2006. V.~35, N~4. P.~985--1005.
\end{remark}

\begin{problem}\DStar(Глобальная оптимизация и монотонный симметричный 
марковский поиск; Некруткин--Тихомиров.) 
Рассматривается задача глобальной оптимизации $f\left( x \right)\to \mathop {\min }\limits_{x\in 
\mathbb{R}^n} $. Считаем, что глобальный минимум достигается в единственной 
точке $x^\ast $ (причем для любых $\varepsilon >0$ выполняется 
условие\footnote{$B_\varepsilon ^c \left( {x } \right)$ -- дополнение 
шара $B_\varepsilon { \left( x \right)}$ радиуса $\varepsilon $ с 
центром в точке $x $.} $\inf \left\{ {f\left( x \right):\;\;x\in B_\varepsilon^c \left( {x^* } \right)} \right\}>f\left( {x^* } \right))$, где 
$f\left( x \right)$ -- непрерывная функция, дважды гладка в точке 
$x^\ast $, причем матрица Гессе $G$ функции $f\left( x \right)$ в этой точке 
положительно определена. Опишем алгоритм (с точностью до выбора функции 
плотности распределения $g\left( r \right)$, $r\in \left[ {0,\infty } 
\right))$:

\begin{enumerate}
\item (начальный шаг) выбираем точку старта $x_0 =x$;
\item (шаг $k<N$) независимо генерируем с.в. $\xi _k $ из центрально симметричного распределения с заданной плотностью $g\left( r \right)$. Если
	\begin{itemize}
	\item $f\left( {x_k +\xi _k } \right)\le f\left( {x_k } \right),$ то $x_{k+1} =x_k +\xi _k $,
	\item иначе $x_{k+1} =x_k $;
	\end{itemize}

\end{enumerate}
Введем обозначения
\[
M_r = \left\{ {x\in B_r{\left( {x^* } \right)}:\;\;f\left( x 
\right)<f\left( y \right)\;\mbox{для всех }y\in B_r^{c} \left( {x^* } 
\right)} \right\},
\]
\[
\tau _\varepsilon =\min \left\{ {n\in {\rm N}:\quad x_n \in M_\varepsilon } 
\right\},
\quad
\delta \left( x \right)=\inf \left\{ {r\ge 0:\;\;x\in M_r } \right\},
\quad
\]
\[
\Gamma =\prod\limits_{i=1}^n {\left( {\frac{\lambda _i }{\lambda _{\min } }} 
\right)^{1 \mathord{\left/ {\vphantom {1 2}} \right. 
\kern-\nulldelimiterspace} 2}} ,
\]
где $\lambda _i $ -- собственные числа $G$. Покажите, что, как бы 
мы не выбирали функцию плотности $g\left( r \right)$, всегда при 
$\varepsilon <\rho \left( {x,x^\ast } \right)$ имеет место следующая оценка 
снизу:
\[
\Exp\left[ {\left. {\tau _\varepsilon } \right|x_0 =x} \right]\ge \ln \left( 
{{\rho \left( {x,x^\ast } \right)} \mathord{\left/ {\vphantom {{\rho \left( 
{x,x^\ast } \right)} \varepsilon }} \right. \kern-\nulldelimiterspace} 
\varepsilon } \right)+2.
\]
Покажите, что метод с плотностью (есть много других вариантов)
\[
g\left( r \right)=\nu \left( r \right)r^{-d},
\quad
\nu \left( r \right)=\frac{c}{\left( {e+n\left| {\ln r} \right|} \right)\ln 
^2\left( {e+n\left| {\ln r} \right|} \right)},
\]
где $c$  находится из условия нормировки, дает оценку ($\varepsilon 
<{\delta \left( x \right)} \mathord{\left/ {\vphantom {{\delta \left( x 
\right)} 2}} \right. \kern-\nulldelimiterspace} 2)$:

\[\tag{1}
\Exp\left[ {\left. {\tau _\varepsilon } \right|x_0 =x} \right]\le b^n\Gamma \ln 
^2\left( \varepsilon \right)\ln ^2\left( {\ln \left( \varepsilon \right)} 
\right)\left| {\ln \left( {\delta \left( x \right)} \right)} \right|, 
\]
где $b\in \left( {2,3} \right)$ (для простоты восприятия мы привели здесь 
огрубленный вариант). 

\end{problem}

\begin{remark}
Если отказаться от гладкости и(или) положительной 
определенности матрицы $G$, то вместо $\Gamma $, которое в типичных 
ситуациях растет с размерностью пространства экспоненциально быстро, в (*) 
можно использовать $F_{\varepsilon ,x}^{-1} $, где
\[
F_{\varepsilon ,x}  = \mathop {\inf }\limits_{\varepsilon \le r<\delta 
\left( x \right)} \left\{ {{\mbox{vol}\left( {M_r } \right)} \mathord{\left/ 
{\vphantom {{\mbox{vol}\left( {M_r } \right)} {\mbox{vol}\left( {B_r \left( 
{x^\ast } \right)} \right)}}} \right. \kern-\nulldelimiterspace} 
{\mbox{vol}\left( {B_r \left( {x^\ast } \right)} \right)}} \right\}.
\]
Отметим, что все приведенные результаты сохраняются с небольшими поправками 
и для оценок вероятностей больших уклонений, т.е. для
\[
n\left( {x,\varepsilon ,\gamma } \right)=\min \left\{ {n:\;\;\PR\left( {x_n 
\in M_\varepsilon \left| {x_0 =x} \right.} \right)\ge \gamma } \right\}=
\]
\[=\min 
\left\{ {n:\;\;\PR\left( {\tau _\varepsilon \le n\left| {x_0 =x} \right.} 
\right)\ge \gamma } \right\}.
\]
Детали имеются в работах А.С. Тихомирова, опубликованных за последние 20 лет 
в ЖВМ и МФ. 

Изложенные в этой задаче результаты могут вызвать на первых порах удивление. 
И, действительно, как такое возможно, чтобы в задачах глобальной оптимизации 
зависимость числа итераций от точности была логарифмическая, в то время как 
известны нижние оценки, в которых эта точность входит в степени размерности 
пространства (в случае равномерной гладкости высокого порядка степень можно 
понижать) в знаменателе, см., например, Zhigljavsky~A.,  Zilinskas~A.  Stochastic global optimization //
 Optimization and Its Applications. Springer, 2008.  V.~9. P.~262. Тут стоит отметить, что, 
во-первых, нижние оценки получаются для детерминированных методов, ну и 
самое главное, что ``проклятие размерности'' здесь также никуда не делось (см., например, недавнюю работу Б.Т. Поляка и П.С. Щербакова на эту тему arXiv:1603.00311). 
Даже при самом благоприятном раскладе в оценку (1) входит фактор $2^n$, 
экспоненциально растущий с ростом размерности пространства. В отличие от 
глобальной оптимизации, в выпуклой оптимизации такие проблемы можно решать (см.~часть~2).

\end{remark}

 \begin{problem}(Глобальная оптимизация и simulated annealing.)
 \label{annealing}
 Пожалуй, самым популярным сейчас методом глобальной оптимизации (правда, с очень 
плохими на данный момент теоретическими оценками скорости сходимости) 
является simulated annealing (имитация затвердевания или отжига), 
представляющий собой дискретное приближение решения стохастического 
дифференциального уравнения\footnote{ Детально изученного в статье S.~German., C.P.~Hwang. 
Diffusions for global optimization // SIAM J. Control and Optimization. 
1986. V.~24:5. P.~1031--1043.}
\[
dx_t =-\nabla f\left( {x_t } \right)dt+\sqrt {2T\left( t \right)} dw_t ,
\]
где $w_t $ -- винеровский процесс. Покажите, что при неограничительных 
условиях и $T\left( t \right)\equiv T$ траектория $x_t $ имеет при $t\to 
\infty $ стационарное распределение с плотностью Гиббса:
\[
\frac{\exp \left( {-{f\left( x \right)} \mathord{\left/ {\vphantom {{f\left( 
x \right)} T}} \right. \kern-\nulldelimiterspace} T} \right)}{\int {\exp 
\left( {-{f\left( z \right)} \mathord{\left/ {\vphantom {{f\left( z \right)} 
T}} \right. \kern-\nulldelimiterspace} T} \right)dz} },
\]
экспоненциально концентрирующееся в окрестности единственной точки 
глобального минимума $x^\ast $ дважды гладкой функции $f\left( x \right)$ 
при $T\to 0+$. Однако при $T\to 0+$ и время выхода на это стационарное 
распределение неограниченно возрастает, что создает проблемы для 
практического применения. Более правильно брать $T\left( t \right)=$\linebreak $=c 
\mathord{\left/ {\vphantom {c {\ln \left( {2+t} \right)}}} \right. 
\kern-\nulldelimiterspace} {\ln \left( {2+t} \right)}$, где $c$ -- 
достаточно большое число. Покажите, что тогда для любой начальной точки $x_0 
$ траектория процесса $x_t $ имеет в пределе $t\to \infty $ (который 
фактически с хорошей точностью проявляется уже на конечных временах) 
распределение, сосредоточенное в точке~$x^\ast$. 

\end{problem}

\begin{remark} 
Детали и способы дискретизации можно почерпнуть, например, из 
работы Kushner~H.  Asymptotic global behavior for stochastic approximation and 
diffusion with slowly decreasing noise effects: global minimization via 
Monte Carlo // SIAM J. Appl. Math. 1987. V. 47:1. P.~169--183. В СССР в этом направлении работал А.Л. Фрадков.
\end{remark}

\begin{problem} (Multilevel Monte\;Carlo; M.\;Giles.) 
Некоторый диффузионный процесс описывается стохастическим дифференциальным урав\-нением
\[
dS\left( t \right)=a\left( {S(t),t} \right)dt+b\left( {S(t),t} \right)dW\left( t 
\right),
\quad
0\le t\le T,
\quad
S\left( 0 \right)=S_0 ,
\]
где $W\left( t \right)$ -- винеровский процесс. Задана липшицева функция 
$f\left( S \right)$. Требуется предложить численный способ оценивания 
\[
Y=\Exp\left[ {f\left( {S\left( T \right)} \right)} \right].
\]

а)* Дискретизируем задачу по схеме Эйлера:
\[
\hat {S}_{n+1} =\hat {S}_n +a\left( {\hat {S}_n ,t_n } \right) \Delta t +b\left( 
{\hat {S}_n ,t_n } \right)\Delta W_n,
\quad \Delta t = h.
\]
возьмем $N$ независимых реализаций $\left\{ {\hat {S}_n^{\left( i \right)} } 
\right\}$ и положим
\[
\overline {Y}=\frac{1}{N}\sum\limits_{i=1}^N {f\left( {\hat {S}_{T 
\mathord{\left/ {\vphantom {T h}} \right. \kern-\nulldelimiterspace} 
h}^{\left( i \right)} } \right)} .
\]
Покажите, что найдутся такие $C_1 ,C_2 >0$, что
\[
\mbox{MSE}=\mathbb{E}\left[ {\left( {\overline {Y}-Y} \right)^2} \right]\approx C_1 
N^{-1}+C_2 h^2.
\]
Покажите, что если от оценки требуется точность $\varepsilon $ ($\sqrt 
{\mbox{MSE}} ={\it O}\left( \varepsilon \right))$, то оптимально (с точки 
зрения $\mbox{Total}\left( \varepsilon \right)$ -- общего числа 
арифметических операций / генерирования нормальных с.в. $\Delta W_n )$ 
выбирать
\[
h={\it O}\left( \varepsilon \right),
\quad
N={\it O}\left( {\varepsilon ^{-2}} \right),
\quad
\mbox{Total}\left( \varepsilon \right)={\it O}\left( {\varepsilon ^{-3}} 
\right).
\]

б)** Предложим другой (более эффективный) способ оценивания. Для 
этого введем константу $M>1$ и положим
\[
h_l =M^{-l}T,
\quad
\overline {Y}_l =\frac{1}{N_l }\sum\limits_{i=1}^{N_l } {\left( {f\left( {\hat 
{S}_{T \mathord{\left/ {\vphantom {T {h_l }}} \right. 
\kern-\nulldelimiterspace} {h_l }}^{\left( i \right)} } \right)-f\left( 
{\hat {S}_{T \mathord{\left/ {\vphantom {T {h_{l-1} }}} \right. 
\kern-\nulldelimiterspace} {h_{l-1} }}^{\left( i \right)} } \right)} 
\right)} ,
\]
\[
\overline {Y}_0 =\frac{1}{N_0 }\sum\limits_{i=1}^{N_0 } {f\left( {\hat {S}_{T 
\mathord{\left/ {\vphantom {T {h_0 }}} \right. \kern-\nulldelimiterspace} 
{h_0 }}^{\left( i \right)} } \right)} ,
\quad
\overline {Y}=\sum\limits_{l=0}^L {\overline {Y}_l } .
\]
Покажите, что
\[
\mbox{Bias}= \left | \mathbb{E}\left[ {\overline {Y}-Y} \right] \right | ={\it O}\left( {h_L } \right),
\]
\[
V_l =D\left[ {f\left( {\hat {S}_{T \mathord{\left/ {\vphantom {T {h_l }}} 
\right. \kern-\nulldelimiterspace} {h_l }}^{\left( i \right)} } 
\right)-f\left( {\hat {S}_{T \mathord{\left/ {\vphantom {T {h_{l-1} }}} 
\right. \kern-\nulldelimiterspace} {h_{l-1} }}^{\left( i \right)} } \right)} 
\right]={\it O}\left( {h_l } \right),
\]
Обратите внимание, что $V_l $ уменьшается с ростом уровня $l$, что в свою очередь сокращает требуемое число реализаций схемы с каждым новым уровнем.  Мы хотим, чтобы
\[
\sqrt {\mbox{MSE}} =\sqrt {\mbox{Bias}^2+\Var\left[ {\overline {Y}} \right]} \le 
\mbox{Bias}+\sqrt {\Var\left[ {\overline {Y}} \right]} \sim \varepsilon ,
\]
что достигается, если положить
\[
L={\log \left( {\varepsilon ^{-1}} \right)} \mathord{\left/ {\vphantom 
{{\log \left( {\varepsilon ^{-1}} \right)} {\log M}}} \right. 
\kern-\nulldelimiterspace} {\log M}+{\it O}\left( 1 \right),
\quad
\Var\left[ {\overline {Y}} \right]=\sum\limits_{l=0}^L {N_l^{-1} V_l } \sim 
\sum\limits_{l=0}^L {N_l^{-1} h_l } \sim \varepsilon ^2.
\]
Учитывая это, покажите, что решение задачи $$\mbox{Total}\left( \varepsilon 
\right)=\sum\limits_{l=0}^L {N_l h^{-1}_l } \to \mathop {\min }\limits_{\left\{ 
{N_l } \right\}\ge 0}$$ 
при ограничении $\sum\limits_{l=0}^L {N_l^{-1} h_l } 
={\it O}\left( {\varepsilon ^2} \right)$ имеет вид $N_l ={\it O}\left( 
{\varepsilon ^{-2}Lh_l } \right)$. Таким образом, $\mbox{Total}\left( 
\varepsilon \right)={\it O}\left( {\varepsilon ^{-2}\left( {\log \varepsilon 
} \right)^2} \right)$.

\end{problem}

\begin{remark}
Описанный в п. б) метод был предложен относительно 
недавно в контексте разработки эффективных численных методов оценки 
финансовых инструментов на рынке (Giles~M. Multilevel Monte Carlo // Acta Numerica. 2015. V. 24. P. 259--328), поэтому он попал далеко не во все 
классические монографии на эту тему: Glasserman~P. Monte Carlo methods in financial 
engineering. Springer, 2005; Graham~C., Talay~D. Stochastic simulation and Monte Carlo methods: mathematical foundation of stochastic 
simulation. Series ``Stochastic modelling and applied probability''. 2013. V.~68. 
 Тем не менее, мы рекомендуем эти книги для погружения в~область 
численных методов финансовой математики.

\end{remark}

\begin{problem} (Задача об обнаружении сигнала на фоне не случайных помех, 
Фишер--Граничин--Поляк) 
Имеется известный (наблюдаемый) центрированный 
случайный процесс $\phi _n $ (с конечными моментами до четвертого порядка 
включительно). Мы знаем, что реализация $\phi _1 \ne 0$. Имеется неизвестное 
значение некоторого параметра $\theta \in \left[ {0,1} \right]$. Мы 
наблюдаем значения случайного процесса $y_n =\phi _n \theta +\varepsilon _n 
$, где $\left| {\varepsilon _n } \right|\le C$, но ``природа'' $\varepsilon 
_n $ нам неизвестна. Известно лишь, что с.в. $\left\{ {\varepsilon _n } \right\}$ 
и $\left\{ {\phi _n } \right\}$ -- независимы в совокупности. Покажите, что оценка 
\[
\hat {\theta }_n =\frac{\sum\limits_{k=1}^n {\phi _k y_k } 
}{\sum\limits_{k=1}^n {\phi _k^2 } }
\]
состоятельна п.н., т.е. $\hat {\theta }_n 
\overset{\text{п.н.}}{\longrightarrow} 
\theta $.
\end{problem}

\begin{remark}
Приведите $y_n =\phi _n \theta +\varepsilon _n $, 
$n=1,2,....$ к виду
\[
\frac{1}{n}\sum\limits_{k=1}^n {\phi _k y_k } 
=\frac{1}{n}\sum\limits_{k=1}^n {\phi _k^2 \theta } 
+\frac{1}{n}\sum\limits_{k=1}^n {\phi _k \varepsilon _k } ,
\quad
n=1,2,.....
\]
Эта задача является пожалуй одной из самых простых по теме оценивания 
неизвестных параметров при почти произвольных помехах. Идея (восходящая к Р. 
Фишеру) введения дополнительной рандомизации, которая бы устраняла за счет 
усреднения независимые от этой рандомизации помехи не случайной природы, 
весьма плодотворна во многих областях, например, в безградиентной 
оптимизации (см., например, работы arXiv:1412.3890, arXiv:1509.01679, arXiv:1701.03821). Подробнее об этом написано в монографии 
\textit{Граничин О.Н., Поляк Б.Т.} Рандомизированные алгоритмы оценивания и оптимизации при почти произвольных 
помехах. М.: Наука, 2003.
\end{remark}

\section{Вероятностный метод в комбинаторике}
\label{combinatorics}

\begin{problem}
Поверхность некоторой шарообразной планеты состоит из океана и суши (множество мелких островков). Суша занимает больше половины 
площади планеты. Также известно, что суша есть множество, принадлежащее борелевской  $\sigma$-алгебре на сфере. На планету хочет 
совершить посадку космический корабль, сконструированный так, что концы всех шести его ножек лежат на поверхности планеты. 
Посадка окажется успешной, если не меньше четырех ножек из шести окажутся на суши. Возможна ли успешная посадка корабля на планету?
\end{problem}

\begin{ordre}
Введем индикаторную функцию для одной посадки
\[ \xi_i = 
\begin{cases}
1, & i\text{-я ножка оказалась на суше,}\\
0 & \text{иначе.}
\end{cases}
\]
Тогда число ножек, оказавшихся на суше, есть $\xi = \sum \limits_{i=1}^6 \xi_i$.
Покажите, что   $ \Exp \xi > 3$ (усреднение берется по всем возможным посадкам). Значит, существует посадка, для которой   $\xi  > 3$, то есть успешная.
\end{ordre}

\begin{problem} 
Пусть $n\ge 2k$ и семейство $F$ является пересекающимся семейством $k$-элементных подмножеств множества $\left\{0,\ldots ,n-1\right\}$, то есть для любых двух множеств $A,B\in F$ выполняется условие $A\cap B\ne \emptyset $. Найдите с помощью вероятностного метода верхнюю оценку на размер семейства $F$ (а именно, покажите, что $|F|\le C_{n-1}^{k-1} $). Покажите, что эта оценка не улучшаемая.
\end{problem}

\begin{ordre} (Теорема Эрдеша--Ко--Радо.)
Семейство $F$ может содержать не более $k$ множеств вида $A_{s} =\left\{s,s+1,\ldots ,s+k-1\right\}$ (сумма берется по модулю $n$), $0\le s\le n-1$.

Пусть $\sigma $ -- случайная перестановка на множестве $\left\{0,\ldots , n-1\right\}$ и $i$ -- случайное число из множества $\left\{0,\ldots , n-1\right\}$. Пусть $A=$\linebreak $=\left\{\sigma (i),\sigma (i+1),\ldots ,\sigma (i+k-1)\right\}$ (сумма берется по модулю~$n$). Покажите, что, с одной стороны (согласно доказанному выше утверждению), $\PR\left[A\in F\right]\le k/n $, с другой стороны (с учетом равновероятности выбора $A$ из всех $k$-множеств), $\PR\left[A\in F\right]=|F|/C_{n}^{k} $.

\end{ordre} 
\begin{remark}
К этой и последующим задачам данного раздела можно рекомендовать книгу~\cite{15}.
\end{remark}

\begin{problem}(Задача Рамсея.)  Докажите, что для произвольного графа $G=$\linebreak $=(V,E)$ всегда можно 
раскрасить вершины в два цвета таким образом, чтобы не менее половины рёбер 
были ``разноцветными'', то есть соединяли вершины разного цвета.
\end{problem}
\begin{ordre}
Вычислите математическое ожидание числа ``разноцветных'' ребер для случайной раскраски вершин графа в два цвета.
\end{ordre}
\begin{remark}
О применении вероятностного подхода к комбинаторике и теории графов рекомендуется посмотреть также следующие  книги:

Айгнер М., Циглер Г. Доказательства из Книги. Лучшие доказательства со времен Евклида до наших дней.  М.: Мир, 2006.  256~с.

Эссе  Gowers W.T. The Two Cultures of Mathematics: сборник статей ``Математика: границы и перспективы'' / пер. с англ. под ред. Д.В.~Аносова и А.Н.~Паршина.  М.: Фазис, 2005.  606~с.
\end{remark}

\begin{problem}
На турнир приехали $n$ игроков. Каждая пара игроков, согласно регламенту турнира, должна провести одну встречу (ничьих быть не может). Пусть 
$$
C_n^k\cdot (1-2^{-k})^{n-k}<1 . 
$$
Докажите, что тогда игроки могли сыграть так, что для каждого множества из $k$ игроков найдется игрок, который побеждает их всех. 

\end{problem}

\begin{ordre}
Введем $A_K$ --- событие, состоящее в том, что не существует игрока, побеждающего всех игроков из множества $K$. 
Докажите, что 
$$
{\mathbb P}\bigl(\bigcup\limits_{K\subset\{1,..,n\},|K|=k} A_K \bigr)\leqslant C_n^k\cdot (1-2^{-k})^{n-k} . 
$$

\end{ordre}

\begin{problem}
Рассмотрим матрицу $n\times n$, составленную из лампочек, каждая из которых либо включена $(a_{ij}=1)$, либо выключена $(a_{ij}=-1)$. 
Предположим, что для каждой строки и каждого столбца имеется переключатель, поворот которого ($x_i=-1$ для строки $i$ и 
$y_j=-1$ для столбца $j$) переключает все лампочки в соответствующей линии: с <<вкл.>> на <<выкл.>> и с <<выкл.>> на <<вкл.>>. 
Тогда для любой начальной конфигурации лампочек можно установить такое положение переключателей, что разность между числом включенных и 
выключенных лампочек будет не меньше $(\sqrt{2/\pi}+o(1))n^{3/2}$. 
\end{problem}

\begin{ordre}
Рассмотрите  переключатель по столбцам как с.в., принимающую с равной вероятностью значения $1$, $-1$. Каждому переключателю по столбцам необходимо подобрать переключатель по строкам, максимизирующий разность включенных и 
выключенных лампочек. Распределение данной разности можно оценить при помощи ЦПТ.       
\end{ordre}

\begin{problem}
Назовем \textit{турниром} ориентированный граф $T=(V,E)$ такой, что $(x,x)\notin E$ для любой вершины $x\in V$, а для любых двух различных вершин $x\ne y$, $x,y\in V$ либо $(x,y)\in E$, либо $(y,x)\in E$. Множество вершин назовем игроками, каждая пара игроков ровно один раз встречается на матче, если игрок $x$ выигрывает у игрока $y$, то $(x,y)\in E$. Гамильтоновым путем графа назовем перестановку вершин $(x_{1} ,x_{2} ,\ldots ,x_{n} )$, что для всех $i$ игрок $x_{i} $ выигрывает у $x_{i+1} $. Несложно показать, что любой турнир содержит гамильтонов путь. Покажите, что найдется такой турнир на $n$ вершинах, для которого число гамильтоновых путей не меньше, чем $n!/2^{n-1}$.
\end{problem}

\begin{ordre}

Рассмотрите случайный турнир (направление каждого ребра выбирается независимо от других с вероятностью $1/2$). Пусть $X$ -- число гамильтоновых путей в случайном турнире. Для каждой перестановки $\pi $ обозначим через $X_{\pi } $ индикаторную с.в. события, что гамильтонов путь, соответствующий этой перестановке, содержится в случайном турнире. Представьте $X$ в виде суммы таких индикаторных с.в. и, воспользовавшись линейностью математического ожидания, получите, что $\Exp X=n!/2^{n-1}$.
\end{ordre}

\begin{problem}
Дано $k$ перестановок натуральных чисел от 1 до $n$, $n>100$. Оказалось, что этот набор перестановок -- минимальный (по количеству перестановок), обладающий следующим свойством: для любых десяти чисел от 1 до $n$ любую их перестановку можно получить вычеркиванием оcтальных чисел из одной из данных. Докажите, что $\ln n \leq k \leq 10^{100} \ln n$.
\end{problem}

\begin{remark}
Эта задача, а также последующие девять задач нам предоставил Федор Петров (ПОМИ РАН).
\end{remark}

\begin{problem}
Докажите, что числа от 1 до $2^n$ можно покрасить в два цвета так, чтобы не было арифметической прогрессии длины $2n$ одного цвета.
\end{problem}

\begin{problem}
На столе лежат $n$ монет орлами вверх. Каждую минуту Вася равновероятно выбирает одну из монет  и переворачивает ее. Докажите, что вероятность того, что через $k$ минут все монеты
будут лежать решками вверх, не превосходит $2^{1-n}$.
\end{problem}

\begin{problem} а) В алфавите племени УАУАУА только две буквы,
причем никакое слово их языка не является началом другого слова.
Докажите, что $\sum N_i2^{-i}\leq 1$,
где $N_i$ --- количество слов длины $i$ в этом языке.
\indent б) В алфавите племени ОЕЕ только две буквы. Люди
этого племени записывают
предложения без пробелов и это никогда
не приводит к двусмысленности (то есть для
любой конечной последовательности букв есть не более одного способа
разбить их на слова). Докажите, что
$\sum N_i2^{-i}\leq 1$,
где $N_i$ -- количество слов длины $i$ в этом языке.
\end{problem}

\begin{problem}
В таблице $n\times n$ расставлены различные действительные числа.
Докажите, что можно так переставить ее строки, что
ни в одном столбце не будет возрастающей (сверху вниз) последовательности длины $\geq 100 \sqrt{n}$.
\end{problem}

\begin{problem}
В однокруговом волейбольном турнире участвовало тысяча команд.
Всегда ли можно выбрать 21 команду и пронумеровать их так,
чтобы в любой паре из этих команд победила та, номер которой больше?
\end{problem}

\begin{problem}
В двудольном графе меньше, чем $2^n$ вершин, и в каждой
имеется список из $n$ цветов. Докажите, что можно
покрасить каждую вершину в один из цветов ее списка
так, чтобы смежные вершины были разных цветов.
\end{problem}

\begin{problem}

\begin{enumerate}
\item В компании из $n$ человек некоторые пары
дружат, а некоторые другие враждуют, при этом у каждого не более
пяти врагов. Известно, что в любом множестве людей, среди которых нет пар врагов, имеется не более чем
$k$ пар друзей. Докажите, что общее количество пар друзей не превосходит $2^{2011}k$.
\item То же, если вражда (в отличие от дружбы) -- не обязательно симметричное отношение: каждый человек неприязненно относится не более чем к пятерым, и в любом множестве людей, среди которых никто ни к кому не относится неприязненно, не  более чем $k$ пар друзей.
\end{enumerate}
\end{problem}

\begin{problem}
Несколько мальчиков ``раскидывают на морского'', кому водить в игре. Для этого каждый из них одновременно с другими
``выбрасывает'' на пальцах число от 0 до 5. Числа складываются и сумма отсчитывается по кругу начиная с заранее
выбранного мальчика (ему соответствует
ноль). Водить будет тот, на ком остановится счет.
При каком числе мальчиков этот метод является справедливым, то есть
вероятность водить одинакова у всех мальчиков?
\end{problem}

\begin{problem}
\label{triangles}
Рассматривается случайный граф $G(n,p)$ (модель Эрдеша--Реньи см. задачу \ref{sec:erdRenyi} раздела \ref{hard}). С.в. $X$ равна числу треугольников в графе. Покажите, что

\begin{enumerate}

\item если $p(n) \ll  n^{-1} $, то граф $G$ почти всегда свободен от треугольников, то есть $\mathop{\lim }\limits_{n\to \infty } \PR (
X > 0)= 0$;

\item если $p(n)\gg n^{-1} $, то граф $G$ почти всегда содержит треугольник, то есть $\mathop{\lim }\limits_{n\to \infty } \PR (X = 0)=0$.
\end{enumerate}

Говорят, что пороговая функция свойства ``граф  свободен от треугольников'' графа $G(n,p)$ равна $n^{-1} $.

\end{problem}

\begin{ordre}

Для случайной величины $X \geq 0$ справедливы неравенства: 

\begin{enumerate}

\item $\PR( X>0) \le \Exp X$, 

\item $\PR( X=0) \le \PR( |X- \Exp X|\ge \Exp X)$. 

\end{enumerate}

\noindent Введем событие $B_{S}$ -- ``$S$ является треугольником''. Тогда 
\[
X=\sum _{|S|=3} \I [B_{S}] = \sum _{|S|=3} X_S,
\]
где $ X_S = \I [B_{S}]$.
Дисперсию зависимых индикаторов предлагается оценивать неравенством  
\[
\Var X\leq  \Exp X+\sum \cov\left( X_{S_1} , X_{S_2} \right). 
\] 
Здесь суммирование ведется по всем упорядоченным зависимым парам различных трехэлементных множеств:

\[
\sum \cov\left( X_{S_1} , X_{S_2} \right) \leq \sum \PR \left( B_{S_1} , B_{S_2} \right)  =
\]\[
=\sum _{S_1} \PR ( B_{S_1} )  \sum  \PR (B_{S_2} | B_{S_1} )  = \Exp X\sum _{}  \PR (B_{S_2} | B_{S_1} ).
\]

\end{ordre}

\begin{problem} (Парадигма Пуассона.) 
Пусть $G\left(n,\frac{c}{n}\right)$, где $c$ -- некоторая константа, -- случайный граф, построенный по модель Эрдеша--Реньи (см. задачу \ref{sec:erdRenyi} раздела \ref{hard}). С~помощью неравенства Янсона (см. замечание) покажите, что с.в. $X = X(n)$, равная числу треугольников в графе, имеет почти пуассоновское распределение с параметром $$\mu =\lim_{n\to\infty}\Exp X = c^3/6,$$ в частности, 
\[ 
\lim_{n\to\infty}\PR (X=0) = e^{ - \mu} .
\]

\end{problem}

\begin{remark}
Согласно предыдущей задаче, пороговая функция свойства ``граф свободен от треугольников''\  равна $n^{-1}$. С учетом обозначений, введенных в предыдущей задаче, \textit{неравенство Янсона} имеет вид
\[
\prod _{|S|=3} \PR \left( \overline{B}_{S} \right)  \le \PR \left(\mathop{\wedge }\limits_{|S|=3} \overline{B}_{S} \right) \leq e^{-\mu +\frac{\Delta }{2} },
\] 
где $\mu =\sum _{|S|=3} \PR( B_{S} )  $, $\Delta =\sum _{|S\cap T|=2} \PR( B_{S} B_{T} )$.

Заметим, что левое неравенство переходит в равенство для взаимно независимых событий $B_S$. При $n\to\infty$ 
$$\prod _{|S|=3} \PR \left( \overline{B}_{S} \right) = \left (1 - \left( \frac{c}{n} \right)^3  \right )^{C_n^3}\to e^{-\frac{c^3}{6}}.$$

На самом деле события $B_S$, $B_T$ зависимы, если $|S\cap T|=2$. Неравенство Янсона дает поправку для ``почти независимых''\ событий через $\Delta=C_n^4C_4^2 \left( \frac {c}{n}\right)^5 = o(1)$. Таким образом, с.в. $X$ --- сумма большого числа индикаторов ``почти независимы''\ событий, имеет почти пуассоновское распределение. 

Более детальный подход к парадигме Пуассона дает метод ``решета Бруна''\ (см. \cite{15}):
$$
\lim_{n\to\infty}\PR (X=k) = \frac{\mu^k}{k!}e^{ -\mu}.
$$
\end{remark}

\begin{problem} 
Покажите, что пороговая функция события: размер максимальной клики $\omega (G)$ в случайном графе $G(n,p)$ (см. задачу \ref{sec:erdRenyi} раздела \ref{hard}) не меньше 4 -- равна $n^{-\frac{2}{3} } $.
\end{problem}

\begin{ordre}
См. задачу \ref{triangles}. 
\end{ordre}

\begin{problem}
Покажите, что для каждого целого числа $n$ найдется раскраска ребер полного графа $K_{n} $ в два цвета (синий, красный), при которой число одноцветных подграфов $K_{4} $ не превосходит $C_{n}^{4} 2^{-5} $. Предложите детерминированный алгоритм построения такой раскраски за полиномиальное от $n$ время.
\end{problem}

\begin{ordre}

Зададим весовую функцию $W(K_{n} )$ частично раскрашенного графа $K_{n} $ как $W(K_{n} )=\sum _{}w(K) $, где суммирование ведется по всем копиям $K$ графа $K_{4} $ в $K_{n} $ и вес~$w(K)$ подграфа~$K$ равен вероятности того, что копия $K$ окажется одноцветной в случае, когда все бесцветные ребра графа $K_{n} $ будут случайно и независимо раскрашены в два цвета. 

Произвольным образом упорядочим все $C_{n}^{2} $ ребер графа $K_{n}$ и создадим цикл их перебора.   Цвет очередного ребра  выбирается так, чтобы минимизировать получающийся вес, то есть при $W_{red} \leq W_{blue} $ ребро раскашивается в красный цвет, в противном случае --- в синий. Покажите, что в такой процедуре вес графа $K_{n} $ с течением времени не возрастает. 

\end{ordre}

\begin{problem}
Покажите, что можно так раскрасить в два цвета ребра полного графа с $n$ вершинами (т.е. графа (без петель), в котором любые две 
различные вершины соединены одним ребром), что любой его полный подграф с $m$ вершинами, где 
$2C_n^m (\left.1\right/2)^{C_m^2}<1$, имеет ребра разного цвета. 
\end{problem}

\begin{problem}(Концентрация хроматического числа.) 
\label{azuma}
Для произвольных $n$ и $p\in (0, 1)$ покажите, что распределение с.в., равной $\chi (G)$ -- хроматическому числу графа $G(n,p)$ (случайный 
граф в модели Эрдеша--Реньи), является плотно сконцентрированным около 
среднего значения:
\[
\forall \lambda >0  \quad \PR\left\{ {\vert \chi (G)-\Exp\chi (G)\vert >\lambda 
\sqrt {n-1} } \right\}\le 2e^{-\frac{\lambda ^2}{2}}.
\]
\end{problem}
\begin{remark}
    Хроматическое число $\chi(G)$ графа $G = (V,E)$ определяется следующим образом:
    \begin{align*}
        \chi(G) := \min \{
            &N : V = \bigcup\limits_{i=1}^N C_i, \, C_n \cap C_m = \emptyset, m \neq n;\\ 
            &\forall e = (e_1, e_2) \in E \rightarrow \exists C_i, C_j, \, i \neq j : e_1 \in C_i,\, e_2 \in C_j
        \},
    \end{align*}
    т.е. это минимальное число $N$, такое, что множество вершин $V$ графа $G=(V,E)$ можно разбить на непересекающиеся классы 
    $C_1, \dots, C_N$:
    $$
        V = \bigcup\limits_{i=1}^N C_i, \, C_i \cap C_j = \emptyset, \, i\neq j,
    $$
    таких, что любое ребро из $E$ соединяет только вершины из разных классов.
\end{remark}

\begin{ordre}

Воспользуйтесь неравенством Азумы: для мартингальной последовательности $X_0 =X_1,\ldots ,X_m $, удовлетворяющей условию $\vert X_{i+1} -X_i \vert \le 1$ для всех $0\le i< m$, справедливо 
\[
\PR\left( {\vert X_m -c\vert >\lambda \sqrt m } \right)\le 
2e^{-\frac{\lambda ^2}{2}}.
\]
В качестве такого мартингала можно взять мартингал проявления 
вершин, заданный следующим образом: $X_i = \Exp\left[ {\chi (G)} \vert  G_{1:i}\right]$, $X_i$ -- условное математическое ожидание значения $\chi (G)$ при зафиксированном подграфе с вершинами $1,\ldots,i$.

Доказательство не дает рецепта вычисления самого среднего. Оказывается (B. Bollobas, 1988), если $p=1/2$, то почти всегда 
$$
\chi (G) \sim \frac{n}{2\log_{2}n}.
$$
Описанная техника популярна в вероятностном анализе рандомизированных алгоритмов, см. Часть 2 и Dubhashi D.P., Panconesi A. Concentration of Measure for the Analysis of Randomized Algorithms. Cambridge, 2012.
\end{ordre}

\begin{problem}(Максимальный размер клики.)
Для произвольного $n$ покажите, что распределение с.в., равной $\omega (G)$ -- максимальному размеру клики графа $G\left( {n,1/2} 
\right)$ (случайный граф в модели Эрдеша--Реньи см. задачу \ref{sec:erdRenyi} раздела \ref{hard}), удовлетворяет неравенством
\[
\PR\left ( {\omega (G)<k} \right) < e^{-(c+o(1))\frac{n^2}{k^8}},
\]
где $c$ -- некоторая положительная константа.

\end{problem}

\begin{ordre}
 Воспользуйтесь неравенством Азумы (см.~указание к задаче \ref{azuma}) для мартингала проявления ребер $X_0 ,\ldots ,X_m $ (здесь $m=C_n^2 )$, заданного следующим образом: $X_0 = \Exp \left[ {Y(G)} \right]$, $X_i $ -- условное математическое ожидание значения $Y(G)$, при 
условии, что первые $i$ ребер/пропусков фиксированы, $Y(G)$ -- максимальный размер семейства непересекающихся по ребрам $k$-клик в графе. 
\[
X_0 = \Exp\left[ {Y(G)} \right]\ge \left( {9+o(1)} 
\right)\frac{n^2}{2k^4},
\]
\[
\left\{ \omega (G)<k \right\} 
\Leftrightarrow 
\left\{ Y(G)=0 \right\}
\Leftrightarrow
\left\{ X_m =0 \right\}.
\]
Далее осталось применить неравенство Азумы для 
\[
\PR\left( {X_m =0} \right)\le \PR\left( {X_m -X_0 \le -X_0 } \right).
\]

\end{ordre}

\begin{problem}
Пусть для модели Эрдеша--Реньи случайного графа $G(n,p),\linebreak p=n^{-\alpha }$, где $\alpha$ -- фиксированное, $\alpha > 5/6$. Тогда существует\linebreak $u=u(n,p)$ такое, что почти всегда 
\[
u\le \chi (G)\le u+3.
\]
\end{problem}

\begin{ordre}
Восплользуйтесь следующей технической леммой.

\begin{lemma}
\textit{Пусть $\alpha$, $c$~-- фиксированные числа, $\alpha > 5/6$. Пусть 
$p=$\linebreak $=n^{-\alpha }$. Тогда почти наверное при $n\to\infty$ каждые $c \sqrt n $ вершин графа 
$G(n,p)$ могут быть правильно раскрашены в три цвета.}
\end{lemma}

Для доказательства леммы предположим противное. Возьмем случайный граф
$G(n,p)$, пусть $T$ -- подмножество (вершин исходного графа) минимального размера, которое нельзя правильно раскрасить в три цвета. Поскольку для всякого $x\in T$ подграф, порожденный множеством $T\backslash \{x\}$, является 3-раскрашиваемым, а подграф, порожденный $T$, не является таковым, $x$ имеет по меньшей мере трех соседей в подграфе, порожденном $T$. То есть если $\vert T\vert =t$, то подграф, порожденный множеством $T$, имеет по меньшей мере $3t/2$ ребер. Вероятность того, что существует такое $T$ с 
$t\le c\sqrt n $, не превосходит $\sum\limits_{t=4}^{c\sqrt n } {C_n^t
C_{C_t^2 }^{\frac{3t}{2}} p^{\frac{3t}{2}}} $. 
Поскольку $C_n^t \le \left( 
{\frac{ne}{t}} \right)^t$ и $C_{C_t^2 }^{\frac{3t}{2}} \le \left( 
{\frac{te}{3}} \right)^{\frac{3t}{2}}$, то 
\[
C_n^t C_{C_t^2 }^{\frac{3t}{2}} 
p^{\frac{3t}{2}}\le \left[ {\frac{ne}{t}\left( {\frac{te}{3}} 
\right)^{\frac{3}{2}}n^{-\frac{3\alpha }{2}}} \right]^t\le \left[ {c_1
n^{1-\frac{3\alpha }{2}}t^{\frac{1}{2}}} \right]^t\le \left[ {c_2 n^{-\left( 
{\frac{3\alpha }{2}-\frac{5}{4}} \right)}} \right]^t
\]
и вероятность заданного события есть $o(1)$ (что и доказывает справедливость леммы).

Далее для произвольного малого $\varepsilon >0$ выберем $u=u(n,p,\varepsilon)$ -- наименьшее целое число, удовлетворяющее неравенству
\[
\PR\left( {\chi (G) \leq u} \right) > \varepsilon.
\]
Далее покажите, что с вероятностью не меньше $1-\varepsilon$ существует $u$-рас\-краска всех, кроме не более чем $c\sqrt n $ вершин. Для этого воспользуйтесь неравенством Азумы для мартингала проявления вершин с 
теоретико-графовой функцией $Y(G)$, равной минимальному размеру множества 
вершин $S$, для которого граф, индуцированный исходным графом $G$, но без вершин $S$, может быть правильно раскрашен в $u$ цветов:
\[
\begin{array}{l}
 \PR\left( {Y\le \Exp Y-\lambda \sqrt {n-1} } \right)<e^{-\frac{\lambda ^2}{2}}, 
\\ 
 \PR\left( {Y\ge \Exp Y+\lambda \sqrt {n-1} } \right)<e^{-\frac{\lambda ^2}{2}}, 
\\ 
 \end{array}
\]
где $\lambda $ удовлетворяет соотношению $e^{-\frac{\lambda ^2}{2}} < \varepsilon$.
Из определения $u$ имеем $\PR\left(Y=0\right) > \varepsilon$. Значит, $\Exp Y\le \lambda \sqrt {n-1} $ и $\PR\left( {Y\ge 2\lambda \sqrt {n-1} } \right)<\varepsilon$.

Интересные результаты (закон 0 и 1, см. задачу \ref{1and0law} раздела \ref{hard}) имеют место для случайных графов, построенных по модели Эрдеша--Реньи, с $p \ll n^{-\alpha}$ или $p \gg  n^{-\alpha}$ для рациональных $\alpha$: многие свойства случайного графа, выразимые на языке логики первого порядка, в таких случаях имеют место либо с вероятностью стремящейся к 1, либо к 0. 
Spencer J.H. The Strange Logic of Random Graphs, Algorithms Combin., V. 22, Springer-Verlag,
Berlin, 2001.

\end{ordre}

\begin{problem}(Балансировка векторов.) 
Пусть ${\rm B}$ -- произвольное нормированное пространство, $v_1 ,\ldots ,v_n $ -- элементы ${\rm B}$, причем $\left\| {v_i } \right\|_2\le 1$. 
Пусть $\varepsilon _1 ,\ldots ,\varepsilon _n $ -- радемахеровские с.в., то есть 
независимые с.в. с распределением $\PR\left( {\varepsilon _i =+1} 
\right)=\PR\left( {\varepsilon _i =-1} \right)=1/2$. Положим\linebreak 
$Y=\left\| {\varepsilon _1 v_1 +\ldots +\varepsilon _n v_n } \right\|_2$.
Покажите справедливость неравенства для произвольного $\lambda >0$:
\[
\PR\left( {\vert Y-\Exp Y\vert >\lambda \sqrt n } \right)\le 2e^{-\frac{\lambda 
^2}{2}}.
\]

\end{problem}

\begin{ordre}  
Воспользуйтесь неравенством Азумы для мартингала, 
полученного последовательным проявлением $\varepsilon_i$.
\end{ordre}

\begin{problem} (Эрдеш--Кац)
Пусть $\omega (n)\to \infty $ произвольно медленно. Покажите, что число тех $x\in \left\{1,\ldots ,n\right\}$, для которых

\[\left|\nu (x)-\ln (\ln n)\right|>\omega (n)\sqrt{\ln (\ln n)} ,\] 
есть $o(n)$. Здесь $\nu (x)$ -- количество простых чисел $p$, делящих $x$ (без учета кратности).
\end{problem}

\begin{remark} 
Грубо это утверждение говорит, что ``почти все'' $n$ имеют число простых делителей (без учета кратности) ``в некотором смысле близкое'' к $\ln (\ln n)$.
(См. также задачу 54 из раздела 3).
\end{remark} 

\begin{ordre} 
Пусть $x$ случайно выбирается из множества $\left\{1,\ldots ,n\right\}$. Для простого $p$ положим: 

\[
X_{p} =\left\{\begin{array}{cc} {1,} & {x \; \mbox{ делится на }   \; p ,} \\ {0,} & { x \;  \mbox{ не делится на }\; p .} \end{array}\right. 
\] 
$X=\sum X_{p}  $, где сумма ведется по всем простым $p\le M\equiv n^{0.1} $. Так как никакое $x\le n$ не может иметь более 10 простых делителей, больших $M$, то $\nu (x)-10\le X(x)\le \nu (x)$ (то есть границы больших уклонений для $X$ переходят в асимптотически равные им границы для $\nu $).

Покажите, что математическое ожидание и дисперсия с.в. $X$ равны $\ln (\ln n)+O(n^{-0.9})$, учтя соотношение $\sum _{p\le x}1/p =\ln\ln x$, где сумма берется по всем простым $p\le x$.

\end{ordre} 

\begin{remark} 
Справедливо также соотношение
\[
\mathop {\lim }\limits_{n\to \infty } \frac{1}{n}\left| {\left\{ {k\leq n:\;\nu \left( k \right)\ge \ln (\ln n)+\lambda \sqrt {\ln (\ln n)} } 
\right\}} \right|=\frac{1}{\sqrt {2\pi } }\int\limits_\lambda ^\infty 
{e^{-{t^2} \mathord{\left/ {\vphantom {{t^2} 2}} \right. 
\kern-\nulldelimiterspace} 2}dt} .
\]
\end{remark}

\begin{problem}
Пусть ${\cal M}=\left\{ {M_1 ,\ldots ,M_s } \right\}$ -- совокупность, 
состоящая из различных $k$-сочетаний элементов множества $\left\{ {1,\ldots 
,n} \right\}$. Назовем\linebreak  $S\subset \left\{ {1,\ldots ,n} \right\}$ \textit{системой общих представителей }(с.о.п.) 
для ${\cal M}$, если $S\cap M_i \ne \emptyset $ для всех $i=1,\ldots ,s$. 
Интерес представляет минимальная (по мощности) с.о.п., т.е. та с.о.п., на 
которой достигается минимум: 
\[
\tau \left( {\cal M} \right)=\min \left\{ {\left| S 
\right|:\;S-\mbox{с.о.п.}} \right\}.
\]
Ясно, что минимальная с.о.п. может быть не единственной (т.е. минимум в 
предыдущем выражении достигается не на единственной $S)$. Зафиксируем 
параметры $n,s,k$ и введем искусственно равномерную дискретную вероятностную 
меру на множестве всех совокупностей ${\cal M}$ (в силу того, что при 
фиксированных $n,s,k$ число таких совокупностей ${\cal M}$ конечно).

а) Выберем согласно введенной вероятностной мере случайную совокупность 
${\cal M}$, найдем все возможные минимальные с.о.п. для нее, пусть их 
количество равно $N({\cal M})$ (это с.в.). Найдите математическое 
ожидание $N({\cal M})$ при условии, что $\tau ({\cal M})=l$.

б) Далее будем интересоваться величиной 
\[
\zeta (n,s,k)=\mathop {\max }\limits_{\cal M} \tau ({\cal M}),
\]
где максимум берется по совокупностям ${\cal M}$ с фиксированными 
параметрами $n,s,k$.

Для получения нижней границы на значения величины $\zeta (n,s,k)$ можно 
воспользоваться вероятностным методом. Согласно введенному выше 
вероятностному пространству на множестве всех совокупностей ${\cal M}$ 
рассмотрим случайное событие $A=\left\{ {{\cal M}:\;\tau ({\cal M})\le l} 
\right\}$. Покажите, что 
\[
\PR (A)\le G(n,s,k,l) = \frac{C_n^l C_{C_n^k -C_{n-l}^k 
}^s }{C_{C_n^k }^s }.
\]
Если параметры $n,s,k,l$, таковы, что $G(n,s,k,l)<1$, то вероятность 
отрицания события $A$ положительна, т.е. существует такая совокупность 
${\cal M}$, для которой $\tau ({\cal M})>l$, а значит, и $\zeta (n,s,k)>l$. 

\textbf{Замечание. }Если $l=l(n,s,k)\approx \frac{n}{k}\ln \frac{sk}{n}$ (в 
предположении, что $sk>n)$, то можно показать, что $G(n,s,k,l)\mathop \to 
\limits 0$, $n\to \infty$, а значит, ``почти всякая'' совокупность обладает 
огромной по размеру минимальной с.о.п. (т.е. с вероятностью, стремящейся к 
единице, случайная совокупность ${\cal M}$ имеет $\tau ({\cal M})>l\approx 
\frac{n}{k}\ln \frac{sk}{n})$. На самом деле полученная выше нижняя оценка 
на значения величины $\zeta (n,s,k)$ асимптотически точна. Иными словами, 
можно доказать следующую теорему (см. Райгородский А.М. Системы общих 
представителей в комбинаторике и их приложения в геометрии.  М.: МЦНМО, 
2009.  132 с.): для любых $n,s,k$ справедливо неравенство
\[
\zeta (n,s,k)\le \max \left\{ {\frac{n}{k},\frac{n}{k}\ln \frac{sk}{n}} 
\right\}+\frac{n}{k}+1.
\]
\end{problem}

\begin{problem}(Локальная лемма Ловаса (ЛЛЛ).)
Орграф зависимостей $(V,D)$ для набора событий $A_{1} ,\ldots A_{t} $ определяется следующим образом:
$V=$\linebreak $=\{ 1,\ldots ,t\};$ $D$ определяется согласно правилу: $(i_{1} ,k),\ldots ,(i_{s} ,k)\notin D$ эквивалентно тому, что $A_{k} $ не зависит от группы событий $A_{i_{1} } ,\ldots, A_{i_{s} } $.

Пусть $D$ -- множество дуг орграфа зависимостей набора событий $A_{1} ,\ldots, A_{n} $. Пусть нашлись такие $x_{1} ,\ldots ,x_{n} \in (0,1)$, что $\forall\,k\in$\linebreak$\in \left\{1,\ldots ,n\right\}$ выполнено неравенство
\[\PR(A_{k} )\le x_{k}\cdot \prod _{i:\; (i,k)\notin D}(1-x_{i} ) .\] 
Тогда
\[\PR\left(\bar{A}_{1} \bigcap \ldots \bigcap \bar{A}_{n} \right)\ge (1-x_{1} )\cdot\ldots\cdot (1-x_{n} )>0.\] 
Примените ЛЛЛ для оценки диагональных чисел Рамсея (см. задачу 3 этого раздела) $R(s,s)>n$, где $n=\left\lfloor 2^{0.5s} \right\rfloor $, то есть покажите, что существует граф на $n$ вершинах, у которого нет ни клик, ни независимых множеств (н.м.) размера $s$. 
\end{problem}

\begin{ordre}
Рассмотрите случайный граф на $n$ вершинах с вероятностью проведения ребра $1/2 $. Для фиксированного подмножества $s$-вершин $U$ определите событие $A_{U} $ --- множество вершин $U$ образует либо клику, либо н.м. Для выбранной модели случайного графа $\PR(A_{U} )=2^{1-C_{s}^{2} } $. Заметьте, что на событие $A_{U} $ влияют только те $A_{U'} $, для которых $\left|U\bigcap U'\right|\ge 2$, то есть на фиксированное событие $A_{U} $ влияют менее $C_{s}^{2} C_{n-2}^{s-2} $ событий.
\end{ordre}

\begin{problem}
Оцените с помощью неравенства Талаграна (см. замечание, указание и \cite{15}) вероятность того, что случайный граф 
$G\left( {n,1/2} \right)$ не имеет клик размера~$k$ и сравните результат с задачей 22 этого же раздела.

\begin{remark}
Пусть множество элементарных исходов $\Omega$ является произведением множеств $\Omega = \prod_{i=1}^n \Omega_i$ (вероятностная мера является соотвественно произведением мер). Для любого множества $A\subseteq\Omega$ и для любого $x\in\Omega$ определим расстояние Талаграна как ``хитро'' взвешенную метрику Хэмминга:
$$
\rho(x, A) = \sup_{\alpha \in \mathbb R^n_+, \|\alpha\|_2 = 1} \inf_{y\in A}\sum_{i=1}^n \alpha_i \mathbb I [x_i  \neq y_i]
$$
(стоит отметить, что веса $\alpha$ могут зависеть от $x$). Введем $t$-окрест\-ность  множества $A\subseteq\Omega$:
$$
A_t = \{ x\in \Omega: \rho(x,A) \le t \}.
$$
Неравенство Талаграна (1995):
$$
\PR(A)(1 - \PR(A_t))\le e^{-t^2 / 4}.
$$
В случае, когда $\Omega = \{0,1\}^n$ расстояние Талаграна
$\rho(x,A)$ совпадает с евклидовым расстоянием от точки $x$ до выпуклой оболочки множества $A$.

Неравенство Талаграна можно примениться к теоретико-графовым функциям $X(x)$ $x\in\Omega$, удовлетворяющим следующим двум условиям:

1) для любых $x$ и $y$, отличающихся не более чем в одной координате, выполнено неравенство $|X(x) - X(y)| \le 1$;

2) $X(x)$ является проверяемой со сложностью $f: \mathbb N \to \mathbb N$, то есть если $X(x) \ge s$, то найдётся подмножество индексов $I \subseteq \{1, \ldots, n\}$, мощности не более, чем $f(s)$, что для всех $y\in\Omega$, совпадающих с $x$ в координатах из $I$, выполнено неравенство $X(y) \ge s$.
\end{remark}

Тогда для любых $b \in \mathbb N$ и $t\in \mathbb R_+$ справедливо неравенство:
$$
\PR(X \le b - t\sqrt{f(b)}) \PR(X\ge b)\le e^{-t^2 / 4}.
$$
\begin{ordre} 
В качестве  $X(x)$  возьмите максимальное число непересекающихся по ребрам $k$-клик. Проверьте, что эта теоретико-графовая 
функция $X(x)$ удовлетворяет необходимым условиям, причем  $f(s)= s C_k^2 $. В качестве $b$ возьмите медиану $X$, а $t$ подберите так, чтобы из неравенства Талаграна можно было получить оценку для
\[
\PR\left( {\omega (G)<k} \right)=\PR\left( {X\le 0} \right),
\]
где $\omega (G)$ -- размер максимальной клики в графе $G$.
\end{ordre}
\end{problem}

\begin{problem}(Теорема Шеннона.) 
Используя вероятностный метод, докажите теорему 
Шеннона. Пусть $\Sigma =\{0,1\}$ -- двухбуквенный алфавит, 
$p<1/2$, $\varepsilon >0$ сколь угодно мало. Существует схема 
кодирования со скоростью передачи данных, превосходящей $1-H(p)-\varepsilon 
$, и вероятностью ошибки при передаче меньшей, чем $\varepsilon $. Здесь $H(p)$ -- энтропия: $$H(p) = -p \log_2 p - (1-p) \log_2 (1-p).$$
Под схемой кодирования подразумеваются функции кодирования $f$ и 
декодирования $g$: 
\[f:\quad \{0,1\}^k\to \{0,1\}^n; \quad g:\quad \{0,1\}^n\to 
\{0,1\}^k.\] Скорость передачи данных в такой схеме определяется отношением 
$k/n$. Пусть $e=\left( {e_1 \ldots e_n } \right)$ --- случайный шумовой
вектор, компоненты которого независимы и имеют распределение Бернулли с 
параметром $p$. В предположении, что случайное сообщение $x$ имеет 
равномерное на $\{0,1\}^k$ распределение (наименее благоприятное априорное распределение, см. Часть 2), определим вероятность правильной передачи данных как $\PR\left( {g\left( {f(x)+e} \right)=x} \right)$ (здесь сложение берется по $\mod {2})$.
\end{problem}

\begin{ordre}
Рекомендуется ознакомиться с книгами \cite{444}, \cite{15} и Рома\-щенко~А., Румянцев~А., Шень~А.  Заметки по теории кодирования.  М.: МЦНМО, 2011.  80~с.

Для больших значений $n$ выберем значение $k=\left\lceil {n\left( {1-H(p)-\varepsilon } \right)} \right\rceil $, обеспечивающее 
нужную скорость передачи данных. Зададим $\delta >0$ такое, что $p+\delta <1/2$ и $H(p+\delta )<H(p)+\varepsilon/2$. Тогда утверждение теоремы есть применение теоремы Варшамова--Гилберта (см. Часть 2): в пространстве $\{0,1\}^n$ с метрикой Хэмминга существует $n(p+\delta )$-сеть мощности не менее чем $2^k$. Основываясь на вероятностном методе докажите 
последнее. 

Функцию кодирования $f: \{0,1\}^k\to 
\{0,1\}^n$ зададим случайно, выбирая для каждого $x$ значение $f(x)$ 
согласно равномерному распределению на $\{0,1\}^n$. Функцию декодирования определим следующим образом: $g(y)=x$, если $x\in \{0,1\}^k$ --- это единственный вектор такой, что $f(x)$ находится на расстоянии Хэмминга не более чем $n(p+\delta )$ от вектора $y\in \{0,1\}^n$. Если указанных векторов $x\in \{0,1\}^k$ нет или напротив больше одного, то декодирование назовем некорректным. Покажите, что вероятность некорректного декодирования 
экспоненциально мала.

Воспользуйтесь неравенством, оценивающим объем шара радиуса $an$, где 
$a<1/2$, в пространстве $\{0,1\}^n$ с метрикой Хэмминга: 
$\sum\limits_{i=0}^{na} {C_n^i } \le 2^{n(H(a)+o(1))}$. Это следствие 
применения неравенства Чернова (см. Часть 2 и п. в) задачи \ref{mark_cheb} раздела \ref{standart}) к вероятности того, что биномиальная 
с.в. с параметрами~$n$, $1/2$ принимает значения, превосходящие $n-na$.
\end{ordre}

\begin{problem}(Игра лжецa.) 
Пусть Пол пытается угадать число $x\in \{1,\ldots ,n\}$ 
у лживой Кэрол, сопротивляющейся этому. Пол может задать $q$ вопросов вида ``Верно ли, что $x\in S?$'', где $S$ --- произвольное подмножество $\{1,\ldots ,n\}$. Вопросы задаются последовательно, $i$-й вопрос Пола может зависеть от предыдущих ответов. Кэрол может лгать, но она не может солгать больше $k$ раз. Покажите с помощью вероятностного подхода, что при условии
\[
n>\frac{2^q}{\sum\limits_{i=0}^k {C_q^i } }
\]
существует выигрышная стратегия у Кэрол. С помощью метода дерандомизации (см. \cite{15}) опишите явную стратегию Кэрол.
\end{problem}

\begin{ordre}
Заметьте, что в связи с тем, что Кэрол сопротивляется 
правильному отгадыванию числа Полом, здесь можно считать, что ее выигрышная стратегия такова, что она в действительности не загадывает число заранее, а выбирает ответы на каждый вопрос Пола так, что все ее ответы согласуются с более чем одним $x$\linebreak (с учетом информации о том, что она может за всю игру солгать не более чем $k$ раз). Тогда зафиксируем стратегию Пола. Пусть Кэрол играет случайно, то есть с вероятностью $1/2$ отвечает либо ``да'', либо ``нет'' на каждый вопрос Пола. Пусть $I_x $ -- индикатор того, что число $x$ согласуется с ответами Кэрол. Покажите, что $\Exp\left( {I_x } \right)=2^{-q}\sum_{i=0}^k {C_q^i } $. Таким образом, в силу 
линейности математического ожидания, среднее значение числа таких $x$, что согласуются с ответами Кэрол, есть $n2^{-q}\sum_{i=0}^k {C_q^i } $, что по условию больше $1$.

Переход от вероятностного доказательства существования (выигрышной стратигии Кэрол) к явному построению называется процедурой дерандомизации (см. Часть 2).
Введем вес игровой ситуации $W$, равный математическому ожиданию числа таких $x$, что согласуются с ответами Кэрол, при условии, что она играет случайно. Тогда стратегия Кэрол: на каждом шаге максимизировать вес игровой ситуации. Если для произвольной игровой ситуации с весом  $W$ и некоторым ходом Пола обозначить за $W^y$ и $W^n$ веса игровой ситуации после ответа Кэрол ``да'', ``нет'' соответственно, то выполняется соотношение $W=(W^y+W^n)/2$. Итак, стратегия Кэрол не дает весу уменьшиться. Если в начале игры вес был 
больше единицы, значит, он будет больше единицы и в конце игры, что 
соответствует выигрышу Кэрол.

\end{ordre}

\begin{problem} (Модель Бакли--Остгуса.)
\label{hnmgraph}
Рассмотрим случайную последовательность графов $\{G_n\}$ (модель роста Интернета), полученную следующим образом. Зафиксируем параметр $a>0$. Пусть $G_1$ -- граф с одной вершиной $1$ и одной петлей $(1,1)$. Далее предположим, что граф $G_{n}$ уже построен. Граф $G_{n+1}$ получается путем добавления к графу $G_{n}$ одной вершины $n+1$ и одного ребра. С вероятностью $$\frac{a}{\left(a+1\right)n+1}$$ это ребро будет направлено из $n+1$ в $n+1$, с вероятностью $$\frac{deg_{n}(i)+a-1}{\left(a+1\right)n+1}$$ будет добавлено ребро $(n+1,i), i=1,...,n$. Здесь $deg_{n}(i)$ -- степень вершины $i$ в графе $G_{n}$. Нетрудно видеть, что выбор того, куда проводить новое ребро, зависит от того, как были проведены все предыдущие ребера. Покажите, что случайный граф $G_n$ можно задать с помощью $n$ независимых с.в.
\end{problem}

\begin{ordre}
Рассмотрите последовательность независимых с.в. $\{\xi_i\}$:
$$\PR\left(\xi_i = 2j - 1 \right) = \frac{a}{\left(a + 1 \right)i - 1},\qquad j = 1,...,i,$$
$$\PR\left(\xi_i = 2j \right) = \frac{1}{\left(a + 1 \right)i - 1},\qquad j = 1,...,i-1.$$
Интерпретация следующая: $\xi_i$ -- второй конец ребра, выходящего из вершины $i$. Если $\xi_i = 2j - 1$, то ребро идет просто в вершину $j$. Если же $\xi_i = 2j$, то ребро из вершины $i$ идет в ту же вершину, что и ребро из вершины $j$. Значение $\xi_j$ само может быть четным ($j = 2l$). Тогда снова перенаправляем конец ребра, выходящего из $i$ согласно $\xi_l$.

Представление $G_n$  с помощью  так введенных независимых с.в. позволяет воспользоваться неравенством Талаграна (см. \cite{15}) для подсчета различных статистик графа. В частности, имеет место следующее неравенство для количества вершин $X_n(d)$ с размером второй окрестности (вторая окрестность заданной вершины формируется из всех таких вершин, отличных от выбранной, в которые можно попасть из выбранной, пройдя путь, состоящий ровно из двух ребер), равным $d = O(n^{1/(4+a) - \delta})$:
\[
\PR(|X_n(d) - \Exp X_n(d)| > (\Exp X_n(d))^{1-\varepsilon} ) \to 0, \quad n \to \infty.
\]  
Пропущенные здесь деталм см. в книге книге  Райгородский А.М. Модели Интернета.  Долгопрудный: ИД <<Интеллект>>, 2013.  64 c.
\end{ordre}

\begin{remark}
Частным случаем приведенной модели графа является модель Bollobas--Riordan при $a = 1$, свойства которой хорошо изучены. К такого типа графам для уменьшения разреженности применяется постобработка: граф $G_{nm}$ преобразуется к графу $G_n^{m}$ путем ``склейки'' (объединения) равных групп из $m$ подряд идущих вершин в одну.

Считается, что модель графа должна удовлетворять следующим эмпирически выявленным характеристикам (см., например, Newman М.Е.J. Networks: An Introduction. Oxford University Press, 2010):
\begin{enumerate}
\item число ребер пропорционально числу вершин, в то время как число треугольников на порядок больше числа ребер;
\item одна большая компонента связности небольшого ($\sim 6$) диаметра;
\item граф устойчив к случайному удалению вершин, в то время как удаление вершин максимальной степени приводит к разбиению графа на компоненты;
\item степени, вторые степени (размеры второй окрестности) и\linebreak PageRank вершин подчиняются степенному распределению;
\item кластерный коэффициент -- вероятность того, что соседи случайной вершины  сами соединены между собой -- имеет неизменное с ростом числа вершин значение;
\item средняя степень соседей случайной вершины степени $d$ имеет распределение $d^{\delta}$, где $\delta < 0$ характерно для веб-графов,  в то время как $\delta > 0$ характерно для социальных сетей;
\item количество ребер между вершинами заданных степеней имеет специфический вид распределения;
\item  веб-графу свойственно наличие выраженных двудольных подграфов 
(любители--любимые сайты, покупатели--продавцы ссылок).  
\end{enumerate}

Одним из основных недостатков модели $G_n^{m}$ является низкий кластерный коэффициент (убывает с ростом $n$) и, как следствие, недостаточное число треугольников. Недавно сотрудниками Яндекс была предложена модель, лишенная последнего недостатка (Ryabchenko--Samosvat--Ostroumova, arXiv:1205.3015). 

Данную задачу рекомендуется сопоставить с задачей \ref{pref_attach} из раздела \ref{hard} и задачей \ref{soc_ineq} из раздела \ref{macrosystems}. 

Имеются различные обобщения моделей preferential attachment (предпочтительного присоединения). Некоторые из этих моделей имеют яркую геометрическую структуру; см.
https://blogs.princeton.edu/ \\ imabandit/2015/01/30/some-pictures-in-geometric-probability/ 

\end{remark}

\atoc{$ $}









\thispagestyle{plain}
\newpage

\addtocontents{toc}{\vspace*{3mm}\noindent {\bf Литература} \hfill \bf 210 \par}   

\renewcommand\refname{Литература}
\makeatletter
\renewcommand{\@biblabel}[1]{#1.}
\makeatother

\pagenumbering{gobble}

\clearpage

\thispagestyle{empty}

{ \center \large  Учебное издание

\vspace{21mm}

{\small
{\bfseries Бузун}~\,Назар\;Олегович\\
{\bfseries Гасников}~\,Александр\;Владимирович\\
{\bfseries Гончаров}~\,Фёдор\;Олегович\\
{\bfseries Горбачев}~\,Олег\;Геннадьевич\\
{\bfseries Гуз}~\,Сергей\;Анатольевич\\
{\bfseries Крымова}~\,Екатерина\;Александровна\\
{\bfseries Натан}~\,Андрей\;Александрович\\
{\bfseries Черноусова}~\,Елена\;Олеговна\\

}
\vspace{13mm}

{\large СТОХАСТИЧЕСКИЙ АНАЛИЗ  \\[1pt] В ЗАДАЧАХ\\
[18pt]  Часть I 

}

}
\vspace{27,5mm}

{\parindent=0mm \small

\small {Редактор \emph{И.\,А.~Волкова}. Корректор \emph{Н.\,Е.~Кобзева}}

\small {Компьютерная верстка \emph{Н.\,Е.~Кобзева}}

\vspace*{3mm}

Подписано в печать 11.04.2016. Формат 60$\times$84\,\!$^{1}\!/\!_{16}$. 

Усл. печ. л.  13{,}25. Уч.-изд. л. 12{,}75. Тираж 400~экз. Заказ №\,140.

\vspace*{3mm}

Федеральное государственное автономное образовательное учреждение

высшего профессионального образования

<<Московский физико-технический институт>>

141700, Московская обл., г. Долгопрудный, Институтский пер., 9

Тел. (495) 408-58-22, e-mail: rio@mipt.ru

\rule[2pt]{\textwidth}{0.2pt}

Отдел оперативной полиграфии <<Физтех-полиграф>>

141700, Московская обл., г. Долгопрудный, Институтский пер., 9

Тел. (495) 408-84-30, e-mail: polygraph@mipt.ru



}



\end{document}